\documentclass[11pt]{amsart}
\usepackage[margin=1in]{geometry}               

\usepackage{cancel}
\usepackage{graphicx}
\usepackage{amssymb}
\usepackage{amsmath, amsthm}
\usepackage{epstopdf}
\usepackage[hidelinks]{hyperref}
\usepackage{lscape}
\usepackage{kbordermatrix}

\usepackage{wrapfig}
\usepackage{enumitem}

\usepackage[font=footnotesize,labelfont=bf, footnotesize]{caption}

\usepackage{tikz}
\usetikzlibrary{shapes}
\usetikzlibrary{positioning}

\usetikzlibrary{shapes.misc}
\tikzset{cross/.style={cross out, draw=black, minimum size=2*(#1-\pgflinewidth), inner sep=0pt, outer sep=0pt},
cross/.default={1pt}}

\usepackage{calc}

\usetikzlibrary{decorations.markings}
\tikzstyle{vertex}=[circle, draw, inner sep=0pt, minimum size=6pt]
\newcommand{\vertex}{\node[vertex]}

\usetikzlibrary{decorations.pathreplacing}
\usepgflibrary{decorations.markings}

\newcommand\cons{\text{c}}

\newcommand{\sign}{\text{sign}}

\newtheorem{thm}{Theorem}[subsection]
\newtheorem{lemma}[thm]{Lemma}
\newtheorem{cor}[thm]{Corollary}
\newtheorem{conj}[thm]{Conjecture}

\theoremstyle{definition}
\newtheorem{defn}[thm]{Definition}
\newtheorem{trule}[thm]{Rule}
\newtheorem{procedure}[thm]{Procedure}
\newtheorem{rem}[thm]{Remark}
\newtheorem{example}[thm]{Example}

\newtheorem{innercustomthm}{Theorem}
\newenvironment{customthm}[1]
  {\renewcommand\theinnercustomthm{#1}\innercustomthm}
  {\endinnercustomthm}

\numberwithin{equation}{subsection}

\title{Combinatorics of the double-dimer model}
\author{Helen Jenne}
\begin{document}

\begin{abstract}
We prove that the partition function for tripartite double-dimer configurations of a planar bipartite graph satisfies a recurrence related to the Desnanot-Jacobi identity from linear algebra. A similar identity for the dimer partition function was established nearly 20 years ago by Kuo and has applications to random tiling theory and the theory of cluster algebras. This work was motivated in part by the potential for applications in these areas. Additionally, we discuss an application to a problem in Donaldson-Thomas and Pandharipande-Thomas theory.
%which will be the subject of a forthcoming paper. 
The proof of our recurrence requires generalizing work of Kenyon and Wilson; specifically, lifting their assumption that the nodes of the graph are black and odd or white and even. 
\end{abstract}

\maketitle

\section{Introduction}

\begin{wrapfigure}{r}{0.35\textwidth}
%%Double-dimer configuration
\centering

\begin{tikzpicture}[scale=.5]
\def\maxX{7}
%\foreach \x [count=\n] in {0,...,\maxX}{
 %   \foreach \y in {0,...,7}{
 %       \draw[line width=.5pt] (\x,0) -- (\x,7);
 %       \draw[line width=.5pt] (0,\y) -- (\maxX,\y);
 %   }
%}

\node at  (-0.5,0) {$1$};
\node at  (3,-0.5) {$2$};
\node at  (6,-0.5) {$3$};
\node at  (7.5,2) {$4$};
\node at (7, 7.5) {$5$};
\node at (4, 7.5) {$6$};
\node at (1, 7.5) {$7$};

\node at  (-.5,5) {$8$};

%path from node 1 to node 8
   \draw[line width = .25mm]  (0, 0) -- (1, 0);
   \draw[line width = .25mm] (2, 0) -- (2, 1); 
   \draw[line width = .25mm]  (1, 1) -- (0, 1); 
   \draw[line width = .25mm] (1, 3) -- (1, 4); 
   \draw[line width = .25mm] (0, 2)-- (1, 2);
   \draw[line width = .25mm] (0, 5) -- (1, 5);
   
   \draw[line width = .25mm] (-.1, 3) -- (-.1, 4);
      \draw[line width = .25mm] (0.1, 3) -- (0.1, 4);
   
   \draw[line width = .25mm] (-.1, 6) -- (-.1, 7);
      \draw[line width = .25mm] (.1, 6) -- (.1, 7);
   
   \draw[line width = .25mm]  (2.1, 2) -- (2.1, 3);
      \draw[line width = .25mm]  (1.9, 2) -- (1.9, 3);
   
    \draw[line width = .25mm] (1.9, 4) -- (1.9, 5);
       \draw[line width = .25mm] (2.1, 4) -- (2.1, 5);
    
   \draw[line width = .25mm]  (3.1, 3) -- (3.1, 4);
    \draw[line width = .25mm]  (2.9, 3) -- (2.9, 4);
   
%path from node 1 to node 8
   \draw[line width = .25mm] (1, 0) -- (2, 0);
 \draw[line width = .25mm]  (2, 1) -- (1, 1);
 \draw[line width = .25mm] (0, 1) -- (0, 2);
  \draw[line width = .25mm]  (1, 2) -- (1, 3);
   \draw[line width = .25mm] (1, 4) -- (1, 5);

    %node 7 to node 6 
   \draw[line width = .25mm]  (1, 7) -- (2, 7); 
       \draw[line width = .25mm] (3, 7) -- (4, 7); 
    \draw[line width = .25mm] (3, 7) -- (4, 7); 

   \draw[line width = .25mm]   (3, 6) -- (4, 6);
     \draw[line width = .25mm]  (5, 6) -- (6, 6);
    \draw[line width = .25mm] (3, 5) -- (4, 5);
    \draw[line width = .25mm]   (5, 5) -- (6, 5);
    
       \draw[line width = .25mm] (1, 5.9) -- (2, 5.9); 
          \draw[line width = .25mm] (1, 6.1) -- (2, 6.1); 
       
     \draw[line width = .25mm]   (5, 7.1) -- (6, 7.1);
     \draw[line width = .25mm]   (5, 6.9) -- (6, 6.9);

  \draw[line width = .25mm] (2, 7) -- (3, 7); 
\draw[line width = .25mm]  (4, 6) -- (5, 6);
   \draw[line width = .25mm]  (4, 5) -- (5, 5);
\draw[line width = .25mm]  (6, 6) -- (6, 5);
 \draw[line width = .25mm]  (3, 5) -- (3, 6);

     %node 2 to node 5
   \draw[line width = .25mm]  (3, 0) -- (3, 1);
   \draw[line width = .25mm]   (3, 2) -- (4, 2); 
   \draw[line width = .25mm]   (4,3)-- (4, 4);
   \draw[line width = .25mm]  (5, 4) -- (6, 4);
   \draw[line width = .25mm] (7, 4) -- (7,5);
   \draw[line width = .25mm]   (7, 6) -- (7,7);
      %node 2 to node 5
\draw[line width = .25mm] (3, 1) -- (3, 2);
\draw[line width = .25mm]  (4, 2) -- (4,3);
\draw[line width = .25mm]  (4, 4) -- (5, 4);
\draw[line width = .25mm]  (6, 4) -- (7, 4);
\draw[line width = .25mm]  (7, 5) -- (7,6);

   %node 3 to node 4
      \draw[line width = .25mm]  (6, 0) -- (7, 0);
   \draw[line width = .25mm]   (7, 1) -- (6, 1);
   \draw[line width = .25mm]   (6, 2) -- (6, 3);
   \draw[line width = .25mm]  (7, 3) -- (7, 2);
   
   \draw[line width = .25mm] (4.9, 2) -- (4.9, 3);
      \draw[line width = .25mm] (5.1, 2) -- (5.1, 3);
      
   \draw[line width = .25mm]  (5, 0) -- (5, 1);
   \draw[line width = .25mm] (4, 1) -- (4, 0);
   
   \draw[line width = .25mm]  (7, 0) -- (7, 1);
\draw[line width = .25mm] (6, 1) -- (6, 2);
\draw[line width = .25mm]  (6, 3)-- (7, 3);
\draw[line width = .25mm] (4, 0) -- (5, 0);
\draw[line width = .25mm]   (5, 1) -- (4, 1);

\foreach \x [count = \n] in {0, 2, 4, 6}{
\foreach \y in {0,  2, 4, 6}{
       \filldraw[fill=black, draw=black] (\x,\y) circle (0.15cm); 
         \filldraw[fill=white, draw=black] (\x+1,\y) circle (0.15cm); 
        }
        }
\foreach \x [count = \n] in {1, 3, 5, 7}{
\foreach \y in {1, 3, 5, 7}{
        \filldraw[fill=black, draw=black] (\x,\y) circle (0.15cm); 
          \filldraw[fill=white, draw=black] (\x-1,\y) circle (0.15cm); 
        }
        }
\end{tikzpicture}

\caption{A double-dimer configuration on a grid graph with $8$ nodes. In this configuration, the pairing of the nodes is 
$((1, 8), (3, 4), (5, 2), (7, 6))$.}
\label{fig:DDconfig}
\end{wrapfigure}
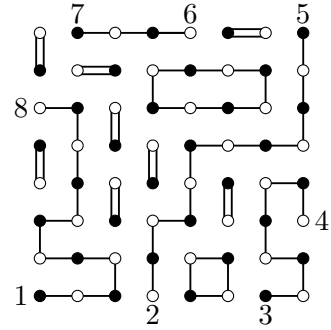

Let $G = (V_1, V_2, E)$ be a finite edge-weighted bipartite planar graph embedded in the plane with $|V_1| = |V_2|$. Let ${\bf N}$ denote a set of special vertices called {\em nodes} on the outer face of $G$ numbered consecutively in counterclockwise order. A {\em double-dimer configuration} on $(G, {\bf N})$ is a multiset of the edges of $G$ with the property that each internal vertex is the endpoint of exactly two edges, and each vertex in ${\bf N}$ is the endpoint of exactly one edge. In other words, it is a configuration of paths connecting the nodes in pairs, doubled edges, and
disjoint cycles of length greater than two (called {\em loops}).
% except the subgraph of degree 2 at the internal vertices and degree 1 at the nodes, except for possibly having some doubled edges. 
Define a probability measure Pr where the probability of a configuration is proportional to the product of its edge weights 
%(doubled edges are counted twice) plan
times $2^{\ell}$, where $\ell$ is the number of loops in the configuration.
%This is the {\em double-dimer model.} 
Kenyon and Wilson initiated the study of the double-dimer model in \cite{KW2006}, when they showed how to compute the probability that a random double-dimer configuration has a particular node pairing.

Before going into the details of Kenyon and Wilson's work, we will describe Kuo's recurrence for {\em dimer configurations}, which is the motivation for this paper,  and state one of our main results. 
% by writing the probability as a polynomials in variables that can be thought of as boundary measurements.
%If $G$ is a finite edge-weighted planar bipartite graph, 
A dimer configuration (or perfect matching) of $G$ is a collection of the edges that covers all of the vertices exactly once. The weight of a dimer configuration is the product of its edge weights. Let $Z^{D}(G)$ denote the sum of the weights of all possible dimer configurations on $G$.
In \cite{Kuo}, Kuo proved that $Z^{D}(G)$ satisfies an elegant recurrence. 
%about the dimer model using a technique called {\em graphical condensation}. 
\begin{thm}\cite[Theorem 5.1]{Kuo}
\label{thm:kuo}
Let $G = (V_1, V_2, E)$ be a planar bipartite graph with a given planar embedding in which $|V_1| = |V_2|$. Let vertices $a, b, c,$ and $d$ appear in a cyclic order on a face of $G$. If $a, c \in V_1$ and $b, d \in V_2$, then
\begin{equation}
\label{eqn:kuo}
Z^{D}(G)Z^{D}(G -\{a, b, c, d\}) =Z^{D}(G  - \{a, b\})Z^{D}(G - \{c, d\})   +  Z^{D}(G  - \{a, d\})Z^{D}(G - \{b, c\}).
\end{equation}
\end{thm}
%In fact, when $G$ is an edge-weighted graph, the same recurrence holds for $W(G)$, defined to be the sum of the weights of all possible perfect matchings on $G$ \cite[Theorem 5.1]{Kuo}\footnotemark. 
%Kuo \cite{Kuo} proved that if $G = (V_1, V_2, E)$ is a plane bipartite graph in which $|V_1| = |V_2|$, and the vertices $a, b, c,$ and $d$ appear in a cyclic order on a face of $G$ such that $a, c \in V_1$ and $b, d \in V_2$, then
%\begin{equation}
%\label{eqn:kuo}
%Z^{D}(G)Z^{D}(G -\{a, b, c, d\}) =Z^{D}(G  - \{a, b\})Z^{D}(G - \{c, d\})   +  Z^{D}(G  - \{a, d\})Z^{D}(G - \{b, c\}), 
%\end{equation}
%where $Z^{D}(G)$ denotes the number of dimer covers (or perfect matchings) of $G$. 
His proof uses a technique called {\em graphical condensation}, which is named for its resemblance to {\em Dodgson condensation}, a method for computing the determinants of square matrices.

%\begin{figure}[htb]
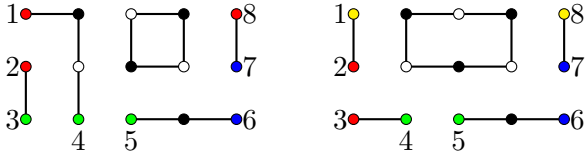
\begin{wrapfigure}{l}{0.5\textwidth}
\centering
\begin{tikzpicture}[scale = .7]
\filldraw[fill=green, draw=black] (0,0) circle (0.1cm); %green node
\filldraw[fill=green, draw=black] (1,0) circle (0.1cm); %green
\filldraw[fill=green, draw=black] (2,0) circle (0.1cm); %green 
\filldraw (3,0) circle (0.1cm);
\filldraw[fill=blue, draw=black] (4,0) circle (0.1cm); %blue 
\filldraw[fill=red, draw=black] (0,1) circle (0.1cm); %red
\filldraw[fill=white, draw=black] (1,1) circle (0.1cm); %white
\filldraw (2,1) circle (0.1cm);
\filldraw[fill=white, draw=black] (3,1) circle (0.1cm); %white
\filldraw[fill=blue, draw=black] (4,1) circle (0.1cm); %blue
\filldraw[fill=red, draw=black]  (0,2) circle (0.1cm); %red
\filldraw (1,2) circle (0.1cm);
\filldraw[fill=white, draw=black] (2,2) circle (0.1cm);  %white
\filldraw (3,2) circle (0.1cm);
\filldraw[fill=red, draw=black]  (4,2) circle (0.1cm); %red

\node at  (-.25,2) {$1$};
\node at  (-.25,1) {$2$};
\node at  (-.25,0) {$3$};

\node at  (4.25,2) {$8$};
\node at  (4.25,1) {$7$};

\node at  (1,-.4) {$4$};
\node at  (2,-.4) {$5$};
\node at  (4.25,0) {$6$};

\draw[thick] (2.1,0) -- (3.9,0); %5 to 6
\draw[thick]  (1, 0.1) -- (1, .9);%4 to 1
\draw[thick] (1, 1.1) -- (1, 2) -- (0.1, 2);%4 to 1
 \draw[thick]  (0, 0.1) -- (0, .9); %2 to 3
 \draw[thick] (4, 1.1) -- (4, 1.9); %7 to 8
\draw[thick]  (2,1) -- (2.9,1);
\draw[thick]  (3,1.1)-- (3,2) --(2.1, 2);
\draw[thick]  (2,1.9) -- (2, 1); %cycle
\end{tikzpicture} \hspace{.5cm}
\begin{tikzpicture}[scale = .7]
\filldraw[fill=red, draw=black] (0,0) circle (0.1cm); %red node
\filldraw[fill=green, draw=black] (1,0) circle (0.1cm); %green
\filldraw[fill=green, draw=black] (2,0) circle (0.1cm); %green 
\filldraw (3,0) circle (0.1cm);
\filldraw[fill=blue, draw=black] (4,0) circle (0.1cm); %blue 
\filldraw[fill=red, draw=black] (0,1) circle (0.1cm); %red
\filldraw[fill=white, draw=black] (1,1) circle (0.1cm); %white
\filldraw (2,1) circle (0.1cm);
\filldraw[fill=white, draw=black] (3,1) circle (0.1cm); %white
\filldraw[fill=blue, draw=black] (4,1) circle (0.1cm); %blue
\filldraw[fill=yellow, draw=black]  (0,2) circle (0.1cm); %yellow
\filldraw (1,2) circle (0.1cm);
\filldraw[fill=white, draw=black] (2,2) circle (0.1cm);  %white
\filldraw (3,2) circle (0.1cm);
\filldraw[fill=yellow, draw=black]  (4,2) circle (0.1cm); %red

\node at  (-.25,2) {$1$};
\node at  (-.25,1) {$2$};
\node at  (-.25,0) {$3$};

\node at  (4.25,2) {$8$};
\node at  (4.25,1) {$7$};

\node at  (1,-.4) {$4$};
\node at  (2,-.4) {$5$};
\node at  (4.25,0) {$6$};

\draw[thick] (2.1,0) -- (3.9,0); %5 to 6
\draw[thick]  (0, 1.9) -- (0, 1.1); %1 to 2
 \draw[thick]  (0.1, 0) -- (.9, 0); %3 to 4
 \draw[thick] (4, 1.1) -- (4, 1.9); %7 to 8
\draw[thick]  (1.1,1) -- (2.9,1);
\draw[thick]  (3,1.1)-- (3,2) --(2.1, 2);
\draw[thick]  (1.9, 2) -- (1,2) -- (1, 1.1); %cycle
\end{tikzpicture} \hspace{.1cm}

\caption{ \small{Two double-dimer configurations on a grid graph.  The pairing of the nodes on the left is a tripartite pairing because the nodes can be colored contiguously using three colors so that no pair contains nodes of the same RGB color. The pairing on the right is not a tripartite pairing because four colors are required.}}
\label{fig:tripartite}
\end{wrapfigure}
%\end{figure}

In this paper, we will show that when $\sigma$ is a {\em tripartite pairing}, a similar identity to (\ref{eqn:kuo}) holds for $Z^{DD}_{\sigma}(G, {\bf N})$, the weighted sum of all double-dimer configurations on $(G, {\bf N})$ with pairing $\sigma$. % when the pairing is a {\em tripartite pairing.} 

A planar pairing $\sigma$ is a tripartite pairing if the nodes can be divided into three circularly contiguous sets $R, G$, and $B$ so that no node is paired with a node in the same set (see Figure~\ref{fig:tripartite}). We often color the nodes in the sets red, green, and blue, in which case $\sigma$ is the unique planar pairing in which like colors are not paired. 
%We will call a pairing {\em strictly tripartite} if three colors are required
%An example of the type of recurrence we prove is the following. 

The following double-dimer version of equation (\ref{eqn:kuo}) is a corollary to Theorem~\ref{thm:cond} in Section \ref{sec:organization}.

\begin{thm}
\label{cor:cond}
Let $G = (V_1, V_2, E)$ be a finite edge-weighted planar bipartite graph with a set of nodes {\bf N}. Divide the nodes into three circularly contiguous sets $R$, $G$, and $B$ such that $|R|$, $|G|$ and $|B|$ satisfy the triangle inequality and let $\sigma$ be the corresponding tripartite pairing\footnotemark.\footnotetext{If $|R|, |G|$, and $|B|$ do not satisfy the triangle inequality, there is no corresponding tripartite pairing $\sigma$.}
Let $x, y, w, v$ be nodes appearing in a cyclic order such that the set
$\{x,y,w,v\}$ contains at least one node of each RGB color\footnotemark. 
%at least one of these four nodes is red, at least one is green, and at least one is blue.
% such that $\{x, y\}$ and $\{v, w\}$ are pairs of $\sigma$.
%If $x$ and $w$ are both black and $y$ and $v$ are both white, then
 If $x, w \in V_1$ and $y, v \in V_2$ then
%Contiguously color the nodes red, green and blue using the minimum number of colors so that no pair contains nodes of the same color, and assume that the minimum number of colors required is 3. 
%Assume that the nodes are contiguously colored red, green, and blue (a color may occur zero times). 
%Assume the nodes are are numbered consecutively in counterclockwise order starting with the red nodes so that node 1 is the first red node.
%Let the nodes $x,y, w, v$ appear in a cyclic order. If $x, w$ are black and $y, v$ are white, and if in addition $\sigma(x) = y$ and $\sigma(w) = v$ then
%\small
\begin{eqnarray*}
Z^{DD}_{\sigma}(G, {\bf N}) Z^{DD}_{\sigma_{xywv}}(G, {\bf N} - \{x, y, w, v\}) \hspace{-.2cm} &=& \hspace{-.2cm}
Z^{DD}_{\sigma_{xy}}(G, {\bf N} - \{x, y\})  Z^{DD}_{\sigma_{wv}}(G, {\bf N} - \{w, v\})\\ 
&  &\hspace{-.1cm}  + 
 Z^{DD}_{\sigma_{xv}}(G, {\bf N} - \{x, v\})  Z^{DD}_{\sigma_{wy}}(G, {\bf N} - \{w, y\}), 
 \end{eqnarray*}
 \normalsize
%Let $x, y, w, v$ be nodes appearing in a cyclic order such that $\{x, y\}$ and $\{v, w\}$ are pairs of $\sigma$. If $x, w \in V_1$ and $y, v \in V_2$ then
%Contiguously color the nodes red, green and blue using the minimum number of colors so that no pair contains nodes of the same color, and assume that the minimum number of colors required is 3. 
%Assume that the nodes are contiguously colored red, green, and blue (a color may occur zero times). 
%Assume the nodes are are numbered consecutively in counterclockwise order starting with the red nodes so that node 1 is the first red node.
%Let the nodes $x,y, w, v$ appear in a cyclic order. If $x, w$ are black and $y, v$ are white, and if in addition $\sigma(x) = y$ and $\sigma(w) = v$ then
%\begin{eqnarray*}
%Z^{DD}_{\sigma}(G, {\bf N}) Z^{DD}_{\sigma_5}(G, {\bf N} - \{x, y, w, v\}) 
% &= & Z^{DD}_{\sigma_1}(G, {\bf N} - \{x, y\})  Z^{DD}_{\sigma_2}(G, {\bf N} - \{w, v\})  \\
% & &  +
% Z^{DD}_{\sigma_3}(G, {\bf N} - \{x, v\})  Z^{DD}_{\sigma_4}(G, {\bf N} - \{w, y\}) 
% \end{eqnarray*}
% where $\sigma_i$ is the unique planar pairing on the corresponding node set in which like colors are not paired together. 
where for $i, j \in \{x, y, w, v\}$, $\sigma_{ij}$ is the unique planar pairing on ${\bf N} - \{i, j\}$ %corresponding node set 
 in which like colors are not paired together.
%where $\sigma_{xy}$ is the unique planar pairing on ${\bf N} - \{x, y\}$ %corresponding node set 
% in which like colors are not paired together, and the pairings $\sigma_{wv}, \sigma_{xv}$, etc. are defined similarly. 
 \end{thm}

 \footnotetext{The nodes of $G$ have two colors: the black-white coloring from the bipartite assumption, and the RGB coloring. The coloring we are referring to is often clear from context, but to avoid ambiguity we will often write RGB color to emphasize that we are referring to the red, green, blue coloring of the nodes rather than the black-white coloring.}

We illustrate Theorem~\ref{cor:cond} with an example. 
\begin{example}
\label{ex:thmillustration}
If $G$ is a graph with eight nodes colored red, green, and blue as shown below, then $\sigma = ((1, 8), (3, 4), (5, 2), (7, 6))$. 
If $x = 8, y = 1, w= 2, v= 5$, then by Theorem~\ref{cor:cond}, 
\begin{center}
\small
\hspace{.2cm} 
 $ Z^{DD}_{\sigma}({\bf N}) Z^{DD}_{\sigma_{1258}}({\bf N} - \{1, 2, 5, 8\})$\hspace{.025cm}
$=$ 
\hspace{.025cm}
 $Z^{DD}_{\sigma_{18}}({\bf N} -\{ 1, 8\})  Z^{DD}_{\sigma_{25}}({\bf N} - \{2, 5\})$
%\hspace{.15cm}  
$+$ 
%\hspace{.15cm}
 $Z^{DD}_{\sigma_{12}}({\bf N} - \{1, 2\})  Z^{DD}_{\sigma_{58}}({\bf N} - \{5, 8\}) $
 \end{center}
\normalsize

\hspace{-.6cm}
\begin{minipage}{.175\textwidth}
%%DOUBLE DIMER CONFIGURATION WITH ALL NODES
\begin{tikzpicture}[scale=.3]
\def\maxX{7}

\node at  (-0.5,0) {\color{red}\footnotesize{$1$}};
\node at  (3,-0.5) {\color{red}\footnotesize{$2$}};
\node at  (6,-0.5) {\color{red}\footnotesize{$3$}};
\node at  (7.5,2) {\color{green}\footnotesize{$4$}};
\node at (7, 7.5) {\color{green}\footnotesize{$5$}};
\node at (4, 7.5) {\color{green}\footnotesize{$6$}};
\node at (1, 7.5) {\color{blue}\footnotesize{$7$}};
\node at  (-.5,5) {\color{blue}\footnotesize{$8$}};
\normalsize

%path from node 1 to node 8
   \draw (0, 0) -- (2, 0) -- (2, 1) -- (0, 1) -- (0, 2) -- (1, 2) -- (1, 5) -- (0, 5);
   \draw (0, 3) -- (0, 4);
   \draw (0, 6) -- (0, 7);
   \draw (2, 2) -- (2, 3);
    \draw (2, 4) -- (2, 5);
        \draw (3, 3) -- (3, 4);
  
  %node 7 to node 6 
  \draw (1, 7) -- (4, 7); 
   \draw (1, 6) -- (2, 6); 
   \draw (3, 6) -- (6, 6) -- (6, 5) -- (3, 5) -- (3, 6);
   \draw (5, 7) -- (6, 7);
   
   %node 2 to node 5
   \draw (3, 0) -- (3, 2) -- (4, 2) -- (4, 4) -- (7, 4) -- (7,7);
   
   %node 3 to node 4
   \draw (6, 0) -- (7, 0) -- (7, 1) -- (6, 1) -- (6, 3)-- (7, 3) -- (7, 2);
   \draw (5, 2) -- (5, 3);
   \draw (4, 0) -- (5, 0) -- (5, 1) -- (4, 1) -- (4, 0);
   	
\foreach \x [count = \n] in {0, 2, 4, 6}{
\foreach \y in {0,  2, 4, 6}{
       \filldraw[fill=black, draw=black] (\x,\y) circle (0.15cm); 
         \filldraw[fill=white, draw=black] (\x+1,\y) circle (0.15cm); 
        }
        }
\foreach \x [count = \n] in {1, 3, 5, 7}{
\foreach \y in {1, 3, 5, 7}{
        \filldraw[fill=black, draw=black] (\x,\y) circle (0.15cm); 
          \filldraw[fill=white, draw=black] (\x-1,\y) circle (0.15cm); 
        }
        }
        
 %        \filldraw[fill=red, draw=black] (0, 0) circle (0.15cm); %1    
%\filldraw[fill=red, draw=black] (3, 0) circle (0.15cm);  %2    
 % \filldraw[fill=red, draw=black] (6, 0) circle (0.15cm); %3 
% \filldraw[fill=green, draw=black] (7, 2) circle (0.15cm);  %4
 %\filldraw[fill=green, draw=black] (7, 7) circle (0.15cm);  %5  
% \filldraw[fill=green, draw=black] (4, 7) circle (0.15cm);  %6  
%\filldraw[fill=blue, draw=black] (1, 7) circle (0.15cm);   %7
% \filldraw[fill=blue, draw=black] (0, 5) circle (0.15cm);   %8

\end{tikzpicture}
%%DOUBLE DIMER CONFIGURATION WITH NODES 1, 2, 5 AND 8 DELETED
\end{minipage}
\begin{minipage}{.175\textwidth}
\begin{tikzpicture}[scale=.3]
\def\maxX{7}

%\node at  (-0.5,0) {\color{red}\footnotesize{$1$}};
%\node at  (3,-0.5) {\color{red}\footnotesize{$2$}};
\node at  (6,-0.5) {\color{red}\footnotesize{$3$}};
\node at  (7.5,2) {\color{green}\footnotesize{$4$}};
%\node at (7, 7.5) {\color{green}\footnotesize{$5$}};
\node at (4, 7.5) {\color{green}\footnotesize{$6$}};
\node at (1, 7.5) {\color{blue}\footnotesize{$7$}};
%\node at  (-.5,5) {\color{blue}\footnotesize{$8$}};
\normalsize

%path from node 1 to node 8
   \draw (0, 0) -- (3, 0) -- (3, 1) -- (0, 1) -- (0, 0); 
   \draw (0, 6) -- (0, 7);
   
          \draw (0, 2) -- (1, 2);
       \draw (2, 2) -- (3, 2);
    \draw (4, 2) -- (5, 2);
          \draw (0, 3) -- (1, 3);
       \draw (2, 3) -- (3, 3);
    \draw (4, 3) -- (5, 3);

        \draw (1, 5) -- (2, 5);
  
  %node 7 to node 6 
  \draw (1, 7) -- (4, 7); 
   \draw (1, 6) -- (2, 6); 
   \draw (3, 6) -- (6, 6) -- (6, 5) -- (3, 5) -- (3, 6);
   \draw (5, 7) -- (6, 7);
   
   %node 5 to node 8
     \draw (0, 5) -- (0, 4);
     \draw (1, 4) -- (2, 4); 
     \draw (3, 4) -- (4, 4); 
     \draw (5, 4) -- (6, 4);
     \draw (7, 4) -- (7,5);
       \draw (7, 6) -- (7,7);
   
   %node 3 to node 4
   \draw (6, 0) -- (7, 0) -- (7, 1) -- (6, 1) -- (6, 3)-- (7, 3) -- (7, 2);
 %  \draw (5, 2) -- (5, 3);
   \draw (4, 0) -- (5, 0) -- (5, 1) -- (4, 1) -- (4, 0);
   	
\foreach \x [count = \n] in {0, 2, 4, 6}{
\foreach \y in {0,  2, 4, 6}{
       \filldraw[fill=black, draw=black] (\x,\y) circle (0.15cm); 
         \filldraw[fill=white, draw=black] (\x+1,\y) circle (0.15cm); 
        }
        }
\foreach \x [count = \n] in {1, 3, 5, 7}{
\foreach \y in {1, 3, 5, 7}{
        \filldraw[fill=black, draw=black] (\x,\y) circle (0.15cm); 
          \filldraw[fill=white, draw=black] (\x-1,\y) circle (0.15cm); 
        }
        }
        
 %        \filldraw[fill=red, draw=black] (0, 0) circle (0.15cm); %1    
%\filldraw[fill=red, draw=black] (3, 0) circle (0.15cm);  %2    
 % \filldraw[fill=red, draw=black] (6, 0) circle (0.15cm); %3 
% \filldraw[fill=green, draw=black] (7, 2) circle (0.15cm);  %4
 %\filldraw[fill=green, draw=black] (7, 7) circle (0.15cm);  %5  
% \filldraw[fill=green, draw=black] (4, 7) circle (0.15cm);  %6  
%\filldraw[fill=blue, draw=black] (1, 7) circle (0.15cm);   %7
% \filldraw[fill=blue, draw=black] (0, 5) circle (0.15cm);   %8
      
\end{tikzpicture}
\end{minipage}
%%DOUBLE DIMER CONFIGURATION WITH NODES 1 AND 8 DELETED
\begin{minipage}{.175\textwidth}
\begin{tikzpicture}[scale=.3]
\def\maxX{7}
%\foreach \x [count=\n] in {0,...,\maxX}{
 %   \foreach \y in {0,...,7}{
 %       \draw[line width=.5pt] (\x,0) -- (\x,7);
 %       \draw[line width=.5pt] (0,\y) -- (\maxX,\y);
 %   }
%}

%\node at  (-0.5,0) {\color{red}\footnotesize{$1$}};
\node at  (3,-0.5) {\color{red}\footnotesize{$2$}};
\node at  (6,-0.5) {\color{red}\footnotesize{$3$}};
\node at  (7.5,2) {\color{green}\footnotesize{$4$}};
\node at (7, 7.5) {\color{green}\footnotesize{$5$}};
\node at (4, 7.5) {\color{green}\footnotesize{$6$}};
\node at (1, 7.5) {\color{blue}\footnotesize{$7$}};
%\node at  (-.5,5) {\color{blue}\footnotesize{$8$}};
\normalsize

%path from node 1 to node 8
   \draw (0, 0) -- (2, 0) -- (2, 1) -- (0, 1)  -- (0, 0);
   \draw (0, 2) -- (1, 2) -- (1, 5) -- (0, 5) -- (0, 2);
   \draw (0, 3) -- (0, 4);
   \draw (0, 6) -- (0, 7);
   \draw (2, 2) -- (2, 3);
    \draw (2, 4) -- (2, 5);
        \draw (3, 3) -- (3, 4);
  
  %node 7 to node 6 
  \draw (1, 7) -- (4, 7); 
   \draw (1, 6) -- (2, 6); 
   \draw (3, 6) -- (6, 6) -- (6, 5) -- (3, 5) -- (3, 6);
   \draw (5, 7) -- (6, 7);
   
   %node 2 to node 5
   \draw (3, 0) -- (3, 2) -- (4, 2) -- (4, 4) -- (7, 4) -- (7,7);
   
   %node 3 to node 4
   \draw (6, 0) -- (7, 0) -- (7, 1) -- (6, 1) -- (6, 3)-- (7, 3) -- (7, 2);
   \draw (5, 2) -- (5, 3);
   \draw (4, 0) -- (5, 0) -- (5, 1) -- (4, 1) -- (4, 0);
   	
\foreach \x [count = \n] in {0, 2, 4, 6}{
\foreach \y in {0,  2, 4, 6}{
       \filldraw[fill=black, draw=black] (\x,\y) circle (0.15cm); 
         \filldraw[fill=white, draw=black] (\x+1,\y) circle (0.15cm); 
        }
        }
\foreach \x [count = \n] in {1, 3, 5, 7}{
\foreach \y in {1, 3, 5, 7}{
        \filldraw[fill=black, draw=black] (\x,\y) circle (0.15cm); 
          \filldraw[fill=white, draw=black] (\x-1,\y) circle (0.15cm); 
        }
        }
        
 %        \filldraw[fill=red, draw=black] (0, 0) circle (0.15cm); %1    
%\filldraw[fill=red, draw=black] (3, 0) circle (0.15cm);  %2    
 % \filldraw[fill=red, draw=black] (6, 0) circle (0.15cm); %3 
% \filldraw[fill=green, draw=black] (7, 2) circle (0.15cm);  %4
 %\filldraw[fill=green, draw=black] (7, 7) circle (0.15cm);  %5  
% \filldraw[fill=green, draw=black] (4, 7) circle (0.15cm);  %6  
%\filldraw[fill=blue, draw=black] (1, 7) circle (0.15cm);   %7
% \filldraw[fill=blue, draw=black] (0, 5) circle (0.15cm);   %8

\end{tikzpicture}
\end{minipage}
\hspace{-.75cm}
%%DOUBLE DIMER CONFIGURATION WITH NODES 2 AND 5 DELETED
\begin{minipage}{.175\textwidth}
\begin{tikzpicture}[scale=.3]
\def\maxX{7}
%\foreach \x [count=\n] in {0,...,\maxX}{
 %   \foreach \y in {0,...,7}{
 %       \draw[line width=.5pt] (\x,0) -- (\x,7);
 %       \draw[line width=.5pt] (0,\y) -- (\maxX,\y);
 %   }
%}

\node at  (-0.5,0) {\color{red}\footnotesize{$1$}};
%\node at  (3,-0.5) {\color{red}\footnotesize{$2$}};
\node at  (6,-0.5) {\color{red}\footnotesize{$3$}};
\node at  (7.5,2) {\color{green}\footnotesize{$4$}};
%\node at (7, 7.5) {\color{green}\footnotesize{$5$}};
\node at (4, 7.5) {\color{green}\footnotesize{$6$}};
\node at (1, 7.5) {\color{blue}\footnotesize{$7$}};
\node at  (-.5,5) {\color{blue}\footnotesize{$8$}};
\normalsize

%path from node 1 to node 8
   \draw (0, 0) -- (2, 0) -- (2, 1) -- (0, 1) -- (0, 2) -- (1, 2) -- (1, 5) -- (0, 5);
   \draw (0, 3) -- (0, 4);
   \draw (0, 6) -- (0, 7);
   \draw (2, 2) -- (2, 3);
    \draw (2, 4) -- (2, 5);
      %  \draw (3, 3) -- (3, 4);
  
  %node 7 to node 6 
  \draw (1, 7) -- (4, 7); 
   \draw (1, 6) -- (2, 6); 
   \draw (3, 6) -- (6, 6) -- (6, 5) -- (3, 5) -- (3, 6);
   \draw (5, 7) -- (6, 7);
   
   %node 2 to node 5
   \draw (7, 7) -- (7, 6); 
   \draw (7, 5) -- (7, 4);
   \draw (6, 4) -- (5, 4); 
   \draw (3, 2) -- (4, 2) -- (4,4) -- (3, 4) -- (3, 2);
   \draw (3, 0) -- (3, 1);
   %\draw (3, 0) -- (3, 2) -- (4, 2) -- (4, 4) -- (7, 4) -- (7,7);
   
   %node 3 to node 4
   \draw (6, 0) -- (7, 0) -- (7, 1) -- (6, 1) -- (6, 3)-- (7, 3) -- (7, 2);
   \draw (5, 2) -- (5, 3);
   \draw (4, 0) -- (5, 0) -- (5, 1) -- (4, 1) -- (4, 0);
   	
\foreach \x [count = \n] in {0, 2, 4, 6}{
\foreach \y in {0,  2, 4, 6}{
       \filldraw[fill=black, draw=black] (\x,\y) circle (0.15cm); 
         \filldraw[fill=white, draw=black] (\x+1,\y) circle (0.15cm); 
        }
        }
\foreach \x [count = \n] in {1, 3, 5, 7}{
\foreach \y in {1, 3, 5, 7}{
        \filldraw[fill=black, draw=black] (\x,\y) circle (0.15cm); 
          \filldraw[fill=white, draw=black] (\x-1,\y) circle (0.15cm); 
        }
        }
        
 %        \filldraw[fill=red, draw=black] (0, 0) circle (0.15cm); %1    
%\filldraw[fill=red, draw=black] (3, 0) circle (0.15cm);  %2    
 % \filldraw[fill=red, draw=black] (6, 0) circle (0.15cm); %3 
% \filldraw[fill=green, draw=black] (7, 2) circle (0.15cm);  %4
 %\filldraw[fill=green, draw=black] (7, 7) circle (0.15cm);  %5  
% \filldraw[fill=green, draw=black] (4, 7) circle (0.15cm);  %6  
%\filldraw[fill=blue, draw=black] (1, 7) circle (0.15cm);   %7
% \filldraw[fill=blue, draw=black] (0, 5) circle (0.15cm);   %8

\end{tikzpicture}
\end{minipage}
%%DOUBLE DIMER CONFIGURATION WITH NODES 1 AND 2 DELETED
\begin{minipage}{.175\textwidth}
\begin{tikzpicture}[scale=.3]
\def\maxX{7}
%\foreach \x [count=\n] in {0,...,\maxX}{
 %   \foreach \y in {0,...,7}{
 %       \draw[line width=.5pt] (\x,0) -- (\x,7);
 %       \draw[line width=.5pt] (0,\y) -- (\maxX,\y);
 %   }
%}

%\node at  (-0.5,0) {\color{red}\footnotesize{$1$}};
%\node at  (3,-0.5) {\color{red}\footnotesize{$2$}};
\node at  (6,-0.5) {\color{red}\footnotesize{$3$}};
\node at  (7.5,2) {\color{green}\footnotesize{$4$}};
\node at (7, 7.5) {\color{green}\footnotesize{$5$}};
\node at (4, 7.5) {\color{green}\footnotesize{$6$}};
\node at (1, 7.5) {\color{blue}\footnotesize{$7$}};
\node at  (-.5,5) {\color{blue}\footnotesize{$8$}};
\normalsize

%path from node 1 to node 8
   \draw (0, 0) -- (3, 0) -- (3, 1) -- (0, 1) -- (0, 0); 
   \draw (0, 6) -- (0, 7);
    \draw (4, 2) -- (4, 3);
      \draw (3, 2) -- (3, 3);
   \draw (2, 2) -- (2, 3);
      \draw (1, 2) -- (1, 3);
        \draw (0, 2) -- (0, 3);
        \draw (1, 5) -- (2, 5);
  
  %node 7 to node 6 
  \draw (1, 7) -- (4, 7); 
   \draw (1, 6) -- (2, 6); 
   \draw (3, 6) -- (6, 6) -- (6, 5) -- (3, 5) -- (3, 6);
   \draw (5, 7) -- (6, 7);
   
   %node 5 to node 8
     \draw (0, 5) -- (0, 4) -- (7, 4) -- (7,7);
   
   %node 3 to node 4
   \draw (6, 0) -- (7, 0) -- (7, 1) -- (6, 1) -- (6, 3)-- (7, 3) -- (7, 2);
   \draw (5, 2) -- (5, 3);
   \draw (4, 0) -- (5, 0) -- (5, 1) -- (4, 1) -- (4, 0);
   	
\foreach \x [count = \n] in {0, 2, 4, 6}{
\foreach \y in {0,  2, 4, 6}{
       \filldraw[fill=black, draw=black] (\x,\y) circle (0.15cm); 
         \filldraw[fill=white, draw=black] (\x+1,\y) circle (0.15cm); 
        }
        }
\foreach \x [count = \n] in {1, 3, 5, 7}{
\foreach \y in {1, 3, 5, 7}{
        \filldraw[fill=black, draw=black] (\x,\y) circle (0.15cm); 
          \filldraw[fill=white, draw=black] (\x-1,\y) circle (0.15cm); 
        }
        }
        
 %        \filldraw[fill=red, draw=black] (0, 0) circle (0.15cm); %1    
%\filldraw[fill=red, draw=black] (3, 0) circle (0.15cm);  %2    
 % \filldraw[fill=red, draw=black] (6, 0) circle (0.15cm); %3 
% \filldraw[fill=green, draw=black] (7, 2) circle (0.15cm);  %4
 %\filldraw[fill=green, draw=black] (7, 7) circle (0.15cm);  %5  
% \filldraw[fill=green, draw=black] (4, 7) circle (0.15cm);  %6  
%\filldraw[fill=blue, draw=black] (1, 7) circle (0.15cm);   %7
% \filldraw[fill=blue, draw=black] (0, 5) circle (0.15cm);   %8

\end{tikzpicture}
\end{minipage}
\hspace{-.5cm}
\begin{minipage}{.175\textwidth}
\begin{tikzpicture}[scale=.3]
\def\maxX{7}
%\foreach \x [count=\n] in {0,...,\maxX}{
 %   \foreach \y in {0,...,7}{
 %       \draw[line width=.5pt] (\x,0) -- (\x,7);
 %       \draw[line width=.5pt] (0,\y) -- (\maxX,\y);
 %   }
%}

\node at  (-0.5,0) {\color{red}\footnotesize{$1$}};
\node at  (3,-0.5) {\color{red}\footnotesize{$2$}};
\node at  (6,-0.5) {\color{red}\footnotesize{$3$}};
\node at  (7.5,2) {\color{green}\footnotesize{$4$}};
%\node at (7, 7.5) {\color{green}\footnotesize{$5$}};
\node at (4, 7.5) {\color{green}\footnotesize{$6$}};
\node at (1, 7.5) {\color{blue}\footnotesize{$7$}};
%\node at  (-.5,5) {\color{blue}\footnotesize{$8$}};
\normalsize

%path from node 1 to node 8
   \draw (0, 0) -- (2, 0) -- (2, 1) -- (0, 1) -- (0, 2) -- (1, 2) -- (1, 5) -- (0, 5);
   \draw (0, 3) -- (0, 4);
   \draw (0, 6) -- (0, 7);
   \draw (2, 2) -- (2, 3);
    \draw (2, 4) -- (2, 5);
        \draw (3, 3) -- (3, 4);
  
  %node 7 to node 6 
  \draw (0, 5) -- (0, 6);
  \draw (0, 7) -- (1, 7);
  \draw (2, 7) -- (3, 7); 
  \draw (4, 7) -- (5, 7);
  \draw (6, 7) -- (7, 7);
   \draw (1, 6) -- (2, 6); 
   \draw (3, 6) -- (6, 6) -- (6, 5) -- (3, 5) -- (3, 6);
   \draw (5, 7) -- (6, 7);
   
   %node 2 to node 5
   \draw (3, 0) -- (3, 2) -- (4, 2) -- (4, 4) -- (7, 4) -- (7,7);
   
   %node 3 to node 4
   \draw (6, 0) -- (7, 0) -- (7, 1) -- (6, 1) -- (6, 3)-- (7, 3) -- (7, 2);
   \draw (5, 2) -- (5, 3);
   \draw (4, 0) -- (5, 0) -- (5, 1) -- (4, 1) -- (4, 0);
   	
\foreach \x [count = \n] in {0, 2, 4, 6}{
\foreach \y in {0,  2, 4, 6}{
       \filldraw[fill=black, draw=black] (\x,\y) circle (0.15cm); 
         \filldraw[fill=white, draw=black] (\x+1,\y) circle (0.15cm); 
        }
        }
\foreach \x [count = \n] in {1, 3, 5, 7}{
\foreach \y in {1, 3, 5, 7}{
        \filldraw[fill=black, draw=black] (\x,\y) circle (0.15cm); 
          \filldraw[fill=white, draw=black] (\x-1,\y) circle (0.15cm); 
        }
        }
        
   %        \filldraw[fill=red, draw=black] (0, 0) circle (0.15cm); %1    
%\filldraw[fill=red, draw=black] (3, 0) circle (0.15cm);  %2    
 % \filldraw[fill=red, draw=black] (6, 0) circle (0.15cm); %3 
% \filldraw[fill=green, draw=black] (7, 2) circle (0.15cm);  %4
 %\filldraw[fill=green, draw=black] (7, 7) circle (0.15cm);  %5  
% \filldraw[fill=green, draw=black] (4, 7) circle (0.15cm);  %6  
%\filldraw[fill=blue, draw=black] (1, 7) circle (0.15cm);   %7
% \filldraw[fill=blue, draw=black] (0, 5) circle (0.15cm);   %8

\end{tikzpicture}
\end{minipage}

%\vspace{.25cm}

%\noindent $Z_{\sigma}({\bf N}) Z_{\sigma_1}({\bf N} - 1, 2, 5, 8)$ means the number of double-dimer configurations
%that pair all eight nodes multiplied by the number of double-dimer configurations that pair only the nodes $3, 4, 6, 7$ (with the pairings shown above). 
%Similarly, $Z_{\sigma_2}({\bf N} - 1,8) Z_{\sigma_1}({\bf N} - 2, 5)$ means the number of double-dimer configurations that pair all nodes except 1 and 8 multiplied by the number of double-dimer configurations that pair all nodes except 2 and 5 (with the pairings shown). \\

\end{example}

We were motivated to find an analogue of Theorem~\ref{thm:kuo} by its potential applications, which we discuss in the next section.

\subsection{Applications}

Kuo's work has a variety of applications. 
%Some are computational: 
For example, Kuo uses graphical condensation to give a new proof that the number of tilings of an order-$n$ Aztec diamond is $2^{n(n+1)/2}$ \cite[Theorem 3.2]{Kuo} and a new proof for MacMahon's generating function for plane partitions that are subsets of a box \cite[Theorem 6.1]{Kuo}. His results also have applications to random tiling theory (see \cite[Section 4.1]{Kuo}) and the theory of cluster algebras.

Cluster algebras are a class of commutative rings introduced by Fomin and Zelevinsky \cite{FZ} to study total positivity and dual canonical bases in Lie theory. The theory of cluster algebras has since been connected to many areas of math, including quiver representations, Teichm{\"u}ller theory, Poisson geometry, and integrable systems \cite{Williams}. 
In \cite{Lai1, Lai2}, Tri Lai and Gregg Musiker study toric cluster variables for the quiver associated to the cone over the del Pezzo surface $dP_3$, giving algebraic formulas for these cluster variables as Laurent polynomials. Using identities similar to Kuo's Theorem~\ref{thm:kuo}, they give combinatorial interpretations of most of these formulas
\cite{Lai1}.

We expect Theorem~\ref{cor:cond} to have similar applications. In addition, by using both Theorem~\ref{thm:kuo} and Theorem~\ref{cor:cond} we can give a direct proof of a problem in Donaldson-Thomas and Pandharipande-Thomas theory.

\subsubsection{Application to Donaldson-Thomas and Pandharipande-Thomas theory.}
\label{sec:DTapp}

Donaldson-Thomas (DT) theory, Pandharipande-Thomas (PT), and Gromov-Witten (GW) theory are branches of enumerative geometry closely related to mirror symmetry and string theory.
The DT and GW theories give frameworks for counting curves\footnotemark~on a threefold $X$. One of the conjectures in \cite{MNOP1, MNOP} gives a correspondence between the DT and GW frameworks, which has been proven in special cases, such as when $X$ is toric \cite{MNOP2}.

PT theory gives a third framework for counting curves when $X$ is a nonsingular projective threefold that is Calabi-Yau. The correspondence between the DT and PT frameworks was first conjectured in \cite{PT16} and was proven in \cite{Bridgeland}, which is closely related to the work in \cite{Toda}. 
 Specifically, let $X$ be a toric Calabi-Yau 3-fold. Define $Z_{DT}(q) = \sum\limits_{n} I_{n} q^{n}$, where $I_{n}$ counts length $n$ subschemes of $X$, and $Z_{PT}(q) = \sum\limits_{n} P_{n} q^n$, where $P_{n}$ counts stable pairs on $X$ (see \cite{PT16}). 
Bridgeland proved that these generating functions coincide up to a factor of $M(q)  = \prod\limits_{n = 1}^{\infty} \dfrac{1}{(1 - q^{n})^{n} }$, which is the total $q$-weight of all plane partitions \cite{MacMahon}. 

 \footnotetext{The frameworks differ in what is meant by a curve on $X$.}
 
 \begin{thm}\cite[Theorem 1.1]{Bridgeland}
\label{conj32}
$Z_{DT}(q) = Z_{PT}(q) M(q) \footnotemark$.
\end{thm}

\footnotetext{In \cite{Bridgeland, PT2009, MNOP1, MNOP2} and elsewhere in the geometry literature, the formulas have $-q$ rather than $q$. The sign is there for geometric reasons which are immaterial to us.}
%\begin{wrapfigure}{r}{0.3\textwidth}
%\begin{figure}[h]
%\centering
%\vspace{-1.5cm}
%\includegraphics[width=1.75in]{boundarypartitionsnomargins}
%\caption{A plane partition asymptotic to $\lambda=(3, 2)$, $\mu=(3, 1)$, $\nu=(3, 1, 1)$. Image credit: Okounkov, Reshetikhin, Vafa \cite{ORV}.}
%\vspace{-1cm}
%\label{fig:boundarypartitions}
%\end{wrapfigure}

The application of Theorem~\ref{cor:cond} that we describe relates to Theorem~\ref{conj32} at the level of the topological vertex. 
Define $V_{\lambda, \mu, \nu} =q^{c(\lambda, \mu, \nu)}  \sum\limits_{\pi} q^{|\pi|}$, where the sum is taken over all plane partitions $\pi$ {\em asymptotic to} $(\lambda, \mu, \nu)$.
%(see Figure~\ref{fig:boundarypartitions}). 
Maulik, Nekrasov, Okounkov, and Pandharipande \cite{MNOP1, MNOP} proved that $Z_{DT}(q) = V_{\lambda, \mu, \nu}$ and thus $V_{\lambda, \mu, \nu}$ is called the DT topological vertex. 
Let $W_{\lambda, \mu, \nu}=  q^{c(\lambda, \mu, \nu)} \sum\limits_{i} d_{i} q^{i}$ where $d_{i}$ is a certain weighted enumeration of {\em labelled box configurations} of length $i$ \cite{PT2009}. In \cite[Theorem/Conjecture 2]{PT2009} Pandharipande and Thomas conjecture that $W_{\lambda, \mu, \nu}$ is the stable pairs vertex, i.e. that $Z_{PT}(q) = W_{\lambda, \mu, \nu}$. 

%The Calabi-Yau cause of Conjecture 4 \in \cite{PT2009} is $V_{\lambda, \mu, \nu} = W_{\lambda, \mu, \nu} M(-q)$.

%The following is the Calabi-Yau case of Conjecture 4 in \cite{PT2009}.

In a forthcoming paper with Gautam Webb and Ben Young (for an extended abstract, see \cite{JWY}), we prove that
\begin{conj}\cite[Calabi-Yau case of Conjecture 4]{PT2009}
\label{conj:vwm}
$V_{\lambda, \mu, \nu} = W_{\lambda, \mu, \nu} M(q)$.
\end{conj}

Pandharipande and Thomas remark that a straightforward (but long) approach to this conjecture using DT theory exists \cite{PT2009}. 
Our proof interprets $V_{\lambda, \mu, \nu}$ using the dimer model and $W_{\lambda, \mu, \nu}$ using the double-dimer model,
%The latter requires writing down a weight-preserving correspondence between double-dimer configurations and labeled box configurations. Once this is done, 
and then uses Theorems~\ref{thm:kuo} and~\ref{cor:cond} to show that both $V_{\lambda, \mu, \nu}/M(q)$ and $W_{\lambda, \mu, \nu}$ satisfy the same recurrence. Conjecture~\ref{conj:vwm}, taken together with a substantial body of geometric work, proves the aforementioned Theorem/Conjecture 2 of \cite{PT2009}. For further details, see \cite{JWY}.
% with the same initial conditions.

\subsection{Proof of Theorem~\ref{cor:cond}}

\label{sec:proofsketch}
%Having motivated Theorem~\ref{cor:cond}, we now discuss the main ideas of the proof. 
Presently, we discuss the main ideas behind the proof of Theorem~\ref{cor:cond}. We start by giving an overview of the results from \cite{KW2006, KW2009} that are needed for our work.

\subsubsection{Background}
\label{sec:KWwork}

%First, we give an overview of the results from \cite{KW2006, KW2009} that are needed for our work.
Kenyon and Wilson gave explicit formulas for the probability that a random double-dimer configuration has a particular node pairing $\sigma$. When $\sigma$ is a tripartite pairing, this probability is proportional to the determinant of a matrix.

To be more precise, we need to introduce some notation and definitions. 
Since $G$ is bipartite, we can color its vertices black and white so that each edge connects a black vertex to a white vertex. 
Let $G^{BW}$ be the subgraph of $G$ formed by deleting the nodes except for the ones that are black and odd or white and even. 
Define $G^{WB}$ analogously, but with the roles of black and white reversed.
Let $G^{BW}_{i, j}$ be the graph $G^{BW}$ with nodes $i$ and $j$ included if and only if they were not included in $G^{BW}$. 
For convenience, Kenyon and Wilson assume the nodes alternate in color, so all nodes are black and odd or white and even. (If a graph $G$ does not have this property, we can add edges of weight 1 to each node that has the wrong color to obtain a graph whose double-dimer configurations are in a one-to-one weight-preserving correspondence with double-dimer configurations of $G$.)

%Let $Z^D(G^{BW})$ and $Z^D(G^{BW}_{i, j})$ be the weighted sum of dimer configurations on $G^{BW}$ and $G^{BW}_{i, j}$, respectively.

For each planar pairing $\sigma$, Kenyon and Wilson showed the normalized probability
 $$\widehat{ \Pr }(\sigma) := \Pr(\sigma) \dfrac{ Z^D(G^{WB})}{Z^D(G^{BW})} = \dfrac{Z^{DD}_{\sigma}(G, {\bf N}) }{(Z^D(G^{BW}))^2}$$
 that a random double-dimer configuration has pairing $\sigma$ is an integer-coefficient homogeneous polynomial in the quantities $X_{i, j} := \dfrac{Z^D(G^{BW}_{i, j})}{Z^D(G^{BW})}$  \cite[Theorem 1.3]{KW2006}.

For example, the normalized probability $\widehat{\Pr}$ that a random double-dimer configuration on eight nodes has the pairing $((1, 8), (3, 4), (5, 2), (7, 6))$ (see Figure~\ref{fig:DDconfig}) is
\begin{eqnarray*}
%\begin{split}
\setlength{\arraycolsep}{2.5pt}
\widehat{\Pr}
\footnotesize
\left(\hspace{-.1cm}
\begin{array}{ c| c | c | c}
1 & 3 & 5 & 7\\
8 & 4 & 2 & 6
\end{array}
\hspace{-.1cm}\right)
\normalsize
 &= & X_{1, 8} X_{3, 4} X_{5,2} X_{7, 6} - X_{1, 4} X_{3, 8} X_{5, 2} X_{7, 6} + X_{1, 6}X_{3, 4}X_{5, 8}X_{7, 2} - X_{1, 8}X_{3, 6}X_{5, 2}X_{7, 4} \\
 & &- X_{1, 4}X_{3, 6}X_{5, 8}X_{7, 2} + X_{1, 6}X_{3, 8}X_{5, 2}X_{7, 4}.
 %\end{split}
\end{eqnarray*}

%The computation of the polynomials requires the following definition: 
%for any odd-even pairing $\tau$, let
%$$X'_{\tau} = (-1)^{\text{\# crosses in }\tau} \prod\limits_{i \text{ odd} } X_{i,\tau(i)}.$$

Kenyon and Wilson gave an explicit method for computing these polynomials: they defined a matrix $\mathcal{P}^{(DD)}$ with rows indexed by planar pairings and columns indexed by odd-even pairings. They showed how to calculate the columns of the matrix completely combinatorially and proved that for any planar pairing $\sigma$, 
 \begin{equation}
 \label{eqn:kwthm14}
\widehat{\Pr}(\sigma) = \sum\limits_{\text{ odd-even pairings } \tau} \mathcal{P}^{(DD)}_{\sigma, \tau} X'_{\tau}. 
\end{equation}
where
$X'_{\tau} = (-1)^{\text{\# crosses of }\tau} \prod\limits_{i \text{ odd} } X_{i,\tau(i)}$ \cite[Theorem 1.4]{KW2006}.

% that projects a vector space with basis vectors indexed by odd-even pairings to a vector space whose basis vectors are indexed by planar pairings.

%To prove Theorem \ref{KWthm1.3}, Kenyon and Wilson construct a projection matrix $\mathcal{P}^{(DD)}$ that projects a vector space with basis vectors indexed by odd-even pairings to a vector space whose basis vectors are indexed by planar pairings.

%For odd-even pairings $\tau$ they defined $X'_{\tau} = (-1)^{\text{\# crosses of }\tau} \prod\limits_{i \text{ odd} } X_{i,\tau(i)}$ and they
 %proved that for any planar pairing $\sigma$, 

%The assertion that the polynomials have integer coefficients follows from their proof that the columns of the matrix $\mathcal{P}^{(DD)}$ can be computed combinatorially \cite[Theorem 1.4]{KW2006}. 

%This statement of their theorem has been abbreviated. For the full statement, see \cite{KW2006}. I
In the case where $\sigma$ is a tripartite pairing, $\widehat{ \Pr }(\sigma)$ is a determinant of a matrix whose entries are $X_{i, j}$ or $0$. 
% \cite{KW2009}. 

\begin{thm}\cite[Theorem 6.1]{KW2009}
\label{thm:kw61}
Suppose that the nodes are contiguously colored red, green, and blue (a color may occur zero times), and that $\sigma$ is the (unique) planar pairing in which like colors are not paired together. Let $\sigma(i)$ denote the item that $\sigma$ pairs with item $i$. We have
$$\widehat{ \Pr }(\sigma) = \det [1_{i, j \text{ RGB-colored differently } } X_{i, j} ]^{i = 1, 3, \ldots, 2n-1}_{j = \sigma(1), \sigma(3), \ldots, \sigma(2n-1) }.$$
\end{thm}

Initially, it seems that Theorem~\ref{cor:cond} will follow immediately from combining
Theorem~\ref{thm:kw61} with the Desnanot-Jacobi identity.

\begin{thm}[Desnanot-Jacobi identity]
Let $M = (m_{i, j})_{i, j=1}^{n}$ be a square matrix, and for each $1 \leq i, j \leq n$, let $M_{i}^{j}$ be the matrix that results from $M$ by deleting the $i$th row and the $j$th column. Then
$$\det(M) \det(M_{i, j}^{i, j}) = \det(M_{i}^{i}) \det(M_{j}^{j}) - \det(M_{i}^{j}) \det(M_{j}^{i})$$
\end{thm}

However, we run into some technical obstacles, which we illustrate with an example. 

\subsubsection{Example}

Suppose we wish to prove the equation from Example~\ref{ex:thmillustration}:
 $$ Z^{DD}_{\sigma}({\bf N}) Z^{DD}_{\sigma_{1258}}({\bf N} - \{1, 2, 5, 8\})\hspace{-.1cm}
=\hspace{-.1cm}
 Z^{DD}_{\sigma_{18}}({\bf N} -\{ 1, 8\})  Z^{DD}_{\sigma_{25}}({\bf N} - \{2, 5\})
 +
 Z^{DD}_{\sigma_{12}}({\bf N} - \{1, 2\})  Z^{DD}_{\sigma_{58}}({\bf N} - \{5, 8\}) $$
 where recall that $\sigma = ((1,8), (3,4), (5, 2), (7, 6))$. 
Then the matrix $M$ from Theorem~\ref{thm:kw61} is 
$$M =
\begin{pmatrix}
X_{1,8} &X_{1, 4} & 0 & X_{1, 6} \\
X_{3, 8} & X_{3, 4} & 0 & X_{3, 6}  \\
X_{5, 8} & 0 & X_{5, 2} & 0  \\
0 & X_{7, 4} & X_{7, 2} & X_{7, 6}   \\
\end{pmatrix}. $$

\begin{figure}[htb]
%\begin{minipage}{.25\textwidth}
\centering
\begin{tikzpicture}[scale=.5]
\def\maxX{7}
%\foreach \x [count=\n] in {0,...,\maxX}{
 %   \foreach \y in {0,...,7}{
 %       \draw[line width=.5pt] (\x,0) -- (\x,7);
 %       \draw[line width=.5pt] (0,\y) -- (\maxX,\y);
 %   }
%}

\node at  (-0.5,0) {{$1$}};
%\node at  (3,-0.5) {\footnotesize{$2$}};
\node at  (6,-0.5) {{$3$}};
\node at  (7.5,2) {{$4$}};
%\node at (7, 7.5) {\footnotesize{$5$}};
\node at (4, 7.5) {{$6$}};
\node at (1, 7.5) {{$7$}};
\node at  (-.5,5) {{$8$}};

\node at (0, -1.2) {};

%path from node 1 to node 8
   \draw (0, 0) -- (2, 0) -- (2, 1) -- (0, 1) -- (0, 2) -- (1, 2) -- (1, 5) -- (0, 5);
   \draw (0, 3) -- (0, 4);
   \draw (0, 6) -- (0, 7);
   \draw (2, 2) -- (2, 3);
    \draw (2, 4) -- (2, 5);
      %  \draw (3, 3) -- (3, 4);
  
  %node 7 to node 6 
  \draw (1, 7) -- (4, 7); 
   \draw (1, 6) -- (2, 6); 
   \draw (3, 6) -- (6, 6) -- (6, 5) -- (3, 5) -- (3, 6);
   \draw (5, 7) -- (6, 7);
   
   %node 2 to node 5
   \draw (7, 7) -- (7, 6); 
   \draw (7, 5) -- (7, 4);
   \draw (6, 4) -- (5, 4); 
   \draw (3, 2) -- (4, 2) -- (4,4) -- (3, 4) -- (3, 2);
   \draw (3, 0) -- (3, 1);
   %\draw (3, 0) -- (3, 2) -- (4, 2) -- (4, 4) -- (7, 4) -- (7,7);
   
   %node 3 to node 4
   \draw (6, 0) -- (7, 0) -- (7, 1) -- (6, 1) -- (6, 3)-- (7, 3) -- (7, 2);
   \draw (5, 2) -- (5, 3);
   \draw (4, 0) -- (5, 0) -- (5, 1) -- (4, 1) -- (4, 0);
   	
\foreach \x [count = \n] in {0, 2, 4, 6}{
\foreach \y in {0,  2, 4, 6}{
       \filldraw[fill=black, draw=black] (\x,\y) circle (0.15cm); 
         \filldraw[fill=white, draw=black] (\x+1,\y) circle (0.15cm); 
        }
        }
\foreach \x [count = \n] in {1, 3, 5, 7}{
\foreach \y in {1, 3, 5, 7}{
        \filldraw[fill=black, draw=black] (\x,\y) circle (0.15cm); 
          \filldraw[fill=white, draw=black] (\x-1,\y) circle (0.15cm); 
        }
        }
        
      %   \filldraw[fill=red, draw=black] (0, 0) circle (0.15cm); %1    
%\filldraw[fill=red, draw=black] (3, 0) circle (0.15cm);  %2    
%  \filldraw[fill=red, draw=black] (6, 0) circle (0.15cm); %3 
% \filldraw[fill=green, draw=black] (7, 2) circle (0.15cm);  %4
% \filldraw[fill=green, draw=black] (7, 7) circle (0.15cm);  %5  
% \filldraw[fill=green, draw=black] (4, 7) circle (0.15cm);  %6  
%\filldraw[fill=blue, draw=black] (1, 7) circle (0.15cm);   %7
% \filldraw[fill=blue, draw=black] (0, 5) circle (0.15cm);   %8
      
\end{tikzpicture}\hspace{1cm}
%\end{minipage}
%\begin{minipage}{.25\textwidth}
%\begin{center}
\begin{tikzpicture}[scale=.5]
\def\maxX{7}
%\foreach \x [count=\n] in {0,...,\maxX}{
 %   \foreach \y in {0,...,7}{
 %       \draw[line width=.5pt] (\x,0) -- (\x,7);
 %       \draw[line width=.5pt] (0,\y) -- (\maxX,\y);
 %   }
%}

\node at  (-0.5,0) {$1$};
%\node at  (3,-0.5) {$2$};
\node at  (6,-0.75) {\cancel{$3$}};
\node at  (6.75,-1) {\color{red}$2$};
\node at  (7.5,2) {\cancel{$4$}};
\node at  (8.25,1.75) {\color{red}$3$};
%\node at (7, 7.5) {$5$};
\node at (4, 7.75) {\cancel{$6$}};
\node at (4.75, 7.75) {\color{red}$4$};
\node at (1, 7.75) {\cancel{$7$}};
\node at (1.75, 7.75) {\color{red}$5$};
\node at  (-.5,5) {\cancel{$8$}};
\node at  (-1.25,5) {\color{red}$6$};

%path from node 1 to node 8
   \draw (0, 0) -- (2, 0) -- (2, 1) -- (0, 1) -- (0, 2) -- (1, 2) -- (1, 5) -- (0, 5);
   \draw (0, 3) -- (0, 4);
   \draw (0, 6) -- (0, 7);
   \draw (2, 2) -- (2, 3);
    \draw (2, 4) -- (2, 5);
      %  \draw (3, 3) -- (3, 4);
  
  %node 7 to node 6 
  \draw (1, 7) -- (4, 7); 
   \draw (1, 6) -- (2, 6); 
   \draw (3, 6) -- (6, 6) -- (6, 5) -- (3, 5) -- (3, 6);
   \draw (5, 7) -- (6, 7);
   
   %node 2 to node 5
   \draw (7, 7) -- (7, 6); 
   \draw (7, 5) -- (7, 4);
   \draw (6, 4) -- (5, 4); 
   \draw (3, 2) -- (4, 2) -- (4,4) -- (3, 4) -- (3, 2);
   \draw (3, 0) -- (3, 1);
   %\draw (3, 0) -- (3, 2) -- (4, 2) -- (4, 4) -- (7, 4) -- (7,7);
   
   %node 3 to node 4
   \draw (6, 0) -- (7, 0) -- (7, 1) -- (6, 1) -- (6, 3)-- (7, 3) -- (7, 2);
   \draw (5, 2) -- (5, 3);
   \draw (4, 0) -- (5, 0) -- (5, 1) -- (4, 1) -- (4, 0);
   	
\foreach \x [count = \n] in {0, 2, 4, 6}{
\foreach \y in {0,  2, 4, 6}{
       \filldraw[fill=black, draw=black] (\x,\y) circle (0.15cm); 
         \filldraw[fill=white, draw=black] (\x+1,\y) circle (0.15cm); 
        }
        }
\foreach \x [count = \n] in {1, 3, 5, 7}{
\foreach \y in {1, 3, 5, 7}{
        \filldraw[fill=black, draw=black] (\x,\y) circle (0.15cm); 
          \filldraw[fill=white, draw=black] (\x-1,\y) circle (0.15cm); 
        }
        }
        
      %   \filldraw[fill=red, draw=black] (0, 0) circle (0.15cm); %1    
%\filldraw[fill=red, draw=black] (3, 0) circle (0.15cm);  %2    
%  \filldraw[fill=red, draw=black] (6, 0) circle (0.15cm); %3 
% \filldraw[fill=green, draw=black] (7, 2) circle (0.15cm);  %4
% \filldraw[fill=green, draw=black] (7, 7) circle (0.15cm);  %5  
% \filldraw[fill=green, draw=black] (4, 7) circle (0.15cm);  %6  
%\filldraw[fill=blue, draw=black] (1, 7) circle (0.15cm);   %7
% \filldraw[fill=blue, draw=black] (0, 5) circle (0.15cm);   %8
      
\end{tikzpicture}
%\end{minipage}

\caption{Left: A double-dimer configuration on a grid graph with node set ${\bf N} - \{2, 5\}$ and pairing $((1, 8), (3, 4), (7, 6))$. Right: The same double-dimer configuration after relabeling the nodes.}
\label{fig:examplefromtalk}
\end{figure}
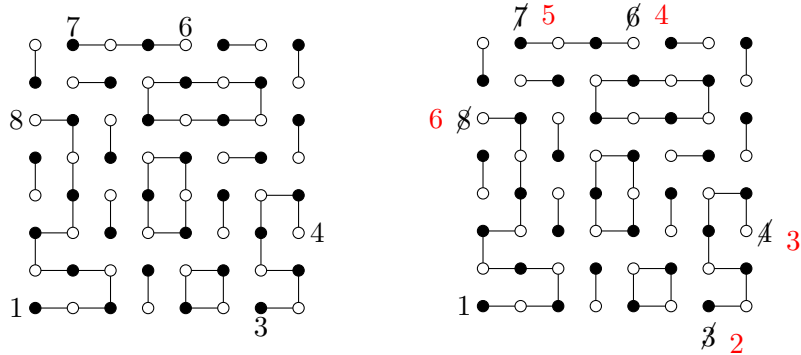

Since the first row and column of $M$ correspond to nodes 1 and 8, respectively, and the third row and column correspond to nodes 5 and 2, 
%Since in Example~\ref{ex:thmillustration} we chose the pairs $(1, 8)$ and $(2, 5)$ of $\sigma$, 
we apply the Desnanot-Jacobi identity with $i = 1$ and $j = 3$:
%and 3 correspond to nodes 1 and 5, and by Theorem~\ref{thm:kw61}, columns 1 and 3 correpond to nodes $\sigma(1) = 8$ and $\sigma(5) = 2$. 
%We have
$$\det(M) \det(M_{1, 3}^{1, 3}) = \det(M_{1}^{1}) \det(M_{3}^{3}) - \det(M_{1}^{3}) \det(M_{3}^{1}).$$
By Theorem~\ref{thm:kw61}, 
$$\det(M) = \dfrac{Z^{DD}_{\sigma}(G, {\bf N}) }{(Z^{D}(G^{BW}))^2 }.$$
We also need to prove, for example, that
\begin{equation}
\label{eqn:exampleminor}
\det(M_{3}^{3} )  =  \dfrac{Z^{DD}_{\sigma_{25}}(G, {\bf N} - \{2, 5\}) }{(Z^{D}(G^{BW}))^2 },
\end{equation}
where
$$M_{3}^{3}  =
\begin{pmatrix}
X_{1,8} &X_{1, 4} & X_{1, 6} \\
X_{3, 8} & X_{3, 4} &  X_{3, 6}  \\
0 & X_{7, 4} & X_{7, 6} 
\end{pmatrix}.$$
 An example of a
%and $Z^{DD}_{\sigma_2}(G, {\bf N} - \{2, 5\}) $ is the generating function for 
double-dimer configuration counted by $Z^{DD}_{\sigma_{25}}(G, {\bf N} - \{2, 5\}) $
% on a grid graph $G$ with node set ${\bf N} - \{2, 5\}$ and pairing $\sigma_2 = ((1, 8), (3, 4), (7, 6))$ 
 is shown in Figure~\ref{fig:examplefromtalk}. 

%An example of a double-dimer configuration on $(G, {\bf N} - \{2, 5\})$ 
%where $G$ is the $8 \times 8$ grid graph with nodes placed as in %Example~\ref{ex:thmillustration}

% One such configuration is shown in Figure~\ref{fig:examplefromtalk}. 

%We cannot apply Theorem~\ref{thm:kw61} to $M_{3}^{3}$ because the nodes are not numbered consecutively. 
We cannot apply Theorem~\ref{thm:kw61} to prove equation~(\ref{eqn:exampleminor})
%$M_{3}^{3}$ 
because the nodes are not numbered consecutively. We might hope to resolve this by relabeling the nodes, as shown in Figure 3.
But since Kenyon and Wilson assume that all nodes are black and odd or white and even, in order to satisfy the assumptions of Kenyon and Wilson's theorem, we need to add edges of weight 1 to nodes $2$ and $3$. Call the resulting graph $\widetilde{G}$ and let $\widetilde{X}_{i,j} = \dfrac{Z^{D}(\widetilde{G}^{BW}_{i, j}) }{Z^{D}(\widetilde{G}^{BW})}$. The matrix from Theorem~\ref{thm:kw61} is
%$$\tilde{M} 
%=  \begin{pmatrix}
%\dfrac{Z^D(\tilde{G}^{BW}_{1,6})}{Z^D(\tilde{G}) } & \dfrac{ Z^D(\tilde{G}%^{BW}_{1, 3} } & \dfrac{ Z^D(\tilde{G}^{BW}_{1, 4} \\
%\dfrac{ Z^D(\tilde{G}^{BW}_{2, 6} & Z^D(\tilde{G}^{BW}_{2, 3} &  Z^D(\tilde{G}^{BW}_{2, 4}  \\
%0 & Z^D(\tilde{G}_{5, 3} & Z^D(\tilde{G}_{5, 4}  
%\end{pmatrix}$$
$$\widetilde{M} =
\begin{pmatrix}
\widetilde{X}_{1,6} &0& \widetilde{X}_{1, 4} \\
\widetilde{X}_{3, 6} & \widetilde{X}_{3, 2} &  0 \\
0 & \widetilde{X}_{5, 2} & \widetilde{X}_{5, 4}  
\end{pmatrix}.
$$

%Recalling that when $i$ and $j$ have the same parity, $X_{i,j} =  \dfrac{Z^{D}(G^{BW}_{i, j}) }{Z^{D}(G^{BW})} = 0$, we see that this approach fails. 
To prove equation (\ref{eqn:exampleminor}) it suffices to show
\begin{equation}
\label{eqn:messy}
(Z^{D}(\widetilde{G}^{BW}))^2 \det(\widetilde{M}) = (Z^{D}(G^{BW}))^2 \det(M_{3}^{3}),
\end{equation}
%which proves (\ref{eqn:exampleminor}), 
since $\det(\widetilde{M}) = 
  \dfrac{Z^{DD}_{\sigma_2}(\widetilde{G}, {\bf N} - \{2, 5\}) }{(Z^{D}(\widetilde{G}^{BW}))^2 }$ by Theorem~\ref{thm:kw61}. 
  
 Verifying equation (\ref{eqn:messy}) is a straightforward computation, but as we consider graphs with more nodes, the computations quickly become more involved.
%Although it is possible to verify similar computations
%Computations similar to the one required to show equation~\ref{eqn:messy} quickly become messy as we consider graphs with more nodes.
To be able to
 %The observation is that in order to 
 interpret the minors of Kenyon and Wilson's matrix outside of small examples, we need to lift their assumption that the nodes of the graph are black and odd or white and even. 

Notice that under the assumption that the nodes of the graph are black and odd or white and even, $X_{i, j} = \dfrac{Z^{D}(G^{BW}_{i, j})}{Z^{D}(G^{BW})}
=  \dfrac{Z^{D}(G_{i, j})}{Z^{D}(G)}$. This suggests that the correct generalization of Kenyon and Wilson's matrix will have entries $\dfrac{Z^{D}(G_{i, j})}{Z^{D}(G)}$. 

%The key observation is that we would not have encountered this issue if
%the matrix entries in Theorem~\ref{thm:kw61} were $\dfrac{Z^{D}(G_{i, j})}%{Z^{D}(G)}$ rather than $\dfrac{Z^{D}(G^{BW}_{i, j})}{Z^{D}(G^{BW})}$. 

%We might hope to apply Theorem~\ref{thm:kw61} after
%We might hope to achieve this by 
%relabeling the nodes and corresponding matrix entries. If we do this, we have the labeling shown in Figure~\ref{fig:examplefromtalk} and the matrix
%$$\begin{pmatrix}
%X_{1,6} &X_{1, 3} & X_{1, 4} \\
%X_{2, 6} & X_{2, 3} &  X_{2, 4}  \\
%0 & X_{5, 3} & X_{5, 4}   \\
%\end{pmatrix}. $$
%Recalling that when $i$ and $j$ have the same parity, $X_{i,j} =  \dfrac{Z^{D}(G^{BW}_{i, j}) }{Z^{D}(G^{BW})} = 0$, we see that this approach fails. 

%The key observation is that we would not have encountered this issue if the matrix entries in Theorem~\ref{thm:kw61} were $\dfrac{Z^{D}(G_{i, j})}{Z^{D}(G)}$ rather than $\dfrac{Z^{D}(G^{BW}_{i, j})}{Z^{D}(G^{BW})}$. 

%Our approach
% to proving equations like (\ref{eqn:theissue}) in general
% is to prove 

\subsubsection{Our approach}
\label{sec:ourapproach}
The previous remark motivates our approach, 
which is to define
$Y_{i, j} :=\dfrac{Z^D(G_{i, j})}{Z^D(G)}$
and $\widetilde{\Pr}(\sigma) = \dfrac{ Z^{DD}_{\sigma}(G, {\bf N}) }{ (Z^{D}(G))^{2} }$. 
When $G$ is a graph with nodes that are either black and odd or white and even, $Z^D(G) = Z^D(G^{BW})$, so $Y_{i, j} = X_{i,j}$ and
$\widetilde{\Pr}(\sigma) = \widehat{\Pr}(\sigma)$.

In this paper, we will prove analogues of many of Kenyon and Wilson's results from \cite{KW2006, KW2009} in the variables $Y_{i, j}$. 
  Once we have established our generalization of Theorem~\ref{thm:kw61}, we will be able to apply the Desnanot-Jacobi identity to prove Theorem~\ref{cor:cond}.

\subsection{Organization of paper}
\label{sec:organization}

This paper is structured as follows.

In Section 2, we generalize some of Kenyon and Wilson's results from \cite{KW2006}.  
%namely their theorem that shows the probability that a double-dimer configuration has a particular pairing is a polynomial function of certain boundary measurements. 
The main result of Section 2 is an analogue of \cite[Theorem 1.4]{KW2006}: we show that we can write $\widetilde{\Pr}(\sigma)$ as an integer-coefficient homogeneous polynomial in the quantities $Y_{i, j}$. 
 %$Y_{i, j} :=\dfrac{Z^D(G_{i, j})}{Z^D(G)}.$
 % When $G$ is a graph with nodes that are either black and odd or white and even, $Y_{i, j} = X_{i,j}$ and $(Z^D(G))^2 = (Z^D(G^{BW}))^2$, so this result generalizes Kenyon and Wilson's Theorem 1.4 \cite{KW2006}. 
  To this end, we define
%Following Kenyon and Wilson's work, for any black-white pairing $\rho$, we define
$$Y'_{\rho} = (-1)^{\text{\# crosses of }\rho} \prod\limits_{i \text{ black} } Y_{i,\rho(i)}$$
for any black-white pairing $\rho$. 
Note that we work with black-white pairings rather than odd-even pairings since we are not requiring that the nodes are either black and odd or white and even. In \cite{KW2006, KW2009}, black-white pairings and odd-even pairings coincide,
%. In our general setting, $i$ and $j$ may have the same parity, but if $i$ and $j$ are the same color then there are no dimer configurations of $G_{i, j}$, so $Y_{i,j} = 0$. 
so $X_{i, j} = 0$ when $i$ and $j$ have the same parity, which occurs exactly when they have the same color\footnotemark. In our general setting, $Y_{i, j}$ may be nonzero when $i$ and $j$ have the same parity, but if $i$ and $j$ are the same color then there are no dimer configurations of $G_{i, j}$, so $Y_{i,j} = 0$. 

\footnotetext{Here, and elsewhere in Section~\ref{sec2},``same color'' refers to the black-white coloring from the bipartite assumption.}

Our analogue of Kenyon and Wilson's matrix $\mathcal{P}^{(DD)}$ (see equation (\ref{eqn:kwthm14})) is $\mathcal{Q}^{(DD)}$. 
%To define $\mathcal{Q}^{(DD)}$, we use Kenyon and Wilson's work as a road map, proving analogues of Lemmas $3.1-3.5$ and Theorem 3.6 from \cite{KW2006}.
 The rows of $\mathcal{Q}^{(DD)}$ are indexed by planar pairings and columns are indexed by black-white pairings.
 To prove that $\mathcal{Q}^{(DD)}$ is integer-valued, we show that the columns of this matrix can be computed combinatorially,
%It has rows indexed by planar pairings and columns indexed by black-white pairings. 
%To show that it is integer-valued, we show that the columns of this matrix can be computed combinatorially 
and in Section~\ref{sec:firstmajorproof} we prove the following theorem:
\begin{thm}
\label{thm:thm1}
Let $G$ be a finite edge-weighted planar bipartite graph with a set of nodes. 
For any planar pairing $\sigma$,
%\dfrac{Z^{DD}_{\sigma}(G, {\bf N})}{(Z^D(G))^2}
$$
\widetilde{\Pr}(\sigma)
= \sum_{\text{black-white pairings } \rho} \mathcal{Q}^{(DD)}_{\sigma, \rho} Y'_{\rho},$$
where the coefficients $\mathcal{Q}^{(DD)}_{\sigma, \rho}$ are all integers. 
\end{thm}
  
  %This work shows that using the variables $X_{i, j}$, while simplying the arguments drastically, was unnecessary. 

To prove Theorem~\ref{thm:thm1}, we use Kenyon and Wilson \cite{KW2006} as a road map, 
proving analogues of Lemmas $3.1-3.5$ and Theorem 3.6 from \cite{KW2006}.
%of their lemmas. 
Because we follow their work so closely, before presenting each of our lemmas we state the corresponding lemma from \cite{KW2006}. In some cases the proofs are very similar. In others, substantially more work is required.

In Section 3, we use our results from Section 2 to generalize Kenyon and Wilson's determinant formula from Theorem~\ref{thm:kw61}. Before stating our version of their formula, we observe that 
\small
$$\det [1_{i, j \text{ RGB-colored differently } } X_{i, j} ]^{i = 1, 3, \ldots, 2n-1}_{j = \sigma(1), \sigma(3), \ldots, \sigma(2n-1) }
= \sign_{OE}(\sigma)\det [1_{i, j \text{ RGB-colored differently } } X_{i, j} ]^{i = 1, 3, \ldots, 2n-1}_{j = 2, 4, \ldots, 2n },$$
\normalsize
where
$\sign_{OE}(\sigma)$
 is the 
parity of the permutation 
$\begin{pmatrix}
\frac{\sigma(1)}{2} & \frac{\sigma(3)}{2} & \cdots & \frac{\sigma(2n-1)}{2}
\end{pmatrix}$
written in one-line notation. 

We prove that
\begin{thm}
\label{thm61}
Let $G$ be a finite edge-weighted planar bipartite graph with a set of nodes.
Suppose that the nodes are contiguously colored red, green, and blue (a color may occur zero times), and that $\sigma$ is the (unique) planar pairing in which like colors are not paired together.  We have
%$$\dfrac{Z^{DD}_{\sigma}(G, {\bf N})}{(Z^D(G))^2}
$$\widetilde{\Pr}(\sigma)= \sign_{OE}(\sigma) \det [1_{i, j \text{ RGB-colored differently } } Y_{i, j} ]^{i = b_1, b_2, \ldots, b_{n}}_{j = w_1, w_2, \ldots, w_{n} },$$
where $b_1 < b_2 < \cdots < b_n$ are the black nodes and $w_1 < w_2 < \cdots < w_n$ are the white nodes.
\end{thm}

%In Section 3 we also prove our main result:
By combining Theorem~\ref{thm61} with the Desnanot-Jacobi identity, we prove our main result:

\begin{thm}
\label{thm:cond}
Let $G= (V_1, V_2, E)$ be a finite edge-weighted planar bipartite graph with a set of nodes {\bf N}.
Divide the nodes into three circularly contiguous sets $R$, $G$, and $B$ such that $|R|, |G|,$ and $|B|$ satisfy the triangle inequality and let $\sigma$ be the corresponding tripartite pairing.
If $x, w \in V_1$ and $y, v \in V_2$ then
\begin{eqnarray*}
& & 
 \sign_{OE}(\sigma) \sign_{OE}(\sigma'_{xywv})Z^{DD}_{\sigma}(G, {\bf N}) Z^{DD}_{\sigma_{xywv}}(G, {\bf N} - \{x, y, w, v\}) \hspace{.4cm}\\ &=& 
\sign_{OE}(\sigma'_{xy}) \sign_{OE}(\sigma'_{wv})
Z^{DD}_{\sigma_{xy}}(G, {\bf N} - \{x, y\})  Z^{DD}_{\sigma_{wv}}(G, {\bf N} - \{w, v\}) \\ 
&& -  \sign_{OE}(\sigma'_{xv}) \sign_{OE}(\sigma'_{wy})
 Z^{DD}_{\sigma_{xv}}(G, {\bf N} - \{x, v\})  Z^{DD}_{\sigma_{wy}}(G, {\bf N} - \{w, y\}),
 \end{eqnarray*}
  where for $i, j \in \{x, y, w, v\}$, $\sigma_{ij}$ is the unique planar pairing on ${\bf N} - \{i, j\}$ %corresponding node set 
 in which like colors are not paired together, and $\sigma_{ij}'$ is the pairing after the the node set ${\bf N} - \{i, j\}$ has been relabeled so that the nodes are numbered consecutively. 
\end{thm}

Theorem~\ref{cor:cond} follows as a corollary; 
the additional assumptions in
Theorem~\ref{cor:cond} lead to a nice simplification of the signs in Theorem~\ref{thm:cond}.

As discussed, Theorems~\ref{thm:thm1} and \ref{thm61} generalize the combinatorial results of \cite{KW2006, KW2009, KW11}. 
 The main questions of interest in these bodies of work involve asymptotic and probabilistic properties of the double-dimer model, which were further studied in \cite{K14, Dubedat, GR}. In \cite{KP}, Kenyon and Pemantle give a connection between the double-dimer model and cluster algebras.
  %appears in \cite{KP}. 
None of these results required taking minors of the matrices from Theorem~\ref{thm:kw61}, so the assumption that the nodes of $G$ are black and odd or white and even was convenient and suitable for their purposes.

\section{Proof of Theorem \ref{thm:thm1}}
%\section{Lemmas from Section 3 of Kenyon and Wilson}

\label{sec2}

In this paper, $G$ always denotes a finite edge-weighted bipartite planar graph embedded in the plane with a set of $2n$ nodes ${\bf N}$ on the outer face of $G$ numbered consecutively in counterclockwise order. Kenyon and Wilson \cite{KW2006, KW2009} assume that the nodes alternate in color so that the black nodes are odd and the white nodes are even. We allow the nodes to have any coloring, as long as ${\bf N}$ has an equal number of black and white nodes. 

To prove Theorem \ref{thm:thm1}, we need to prove analogues of Lemmas $3.1 -3.5$ and Theorem 3.6 from Kenyon and Wilson \cite{KW2006} in this more general setting. 
 For ease of exposition, we prove our lemmas in a slightly different order. 
 
\subsection{Lemma 3.4 from Kenyon and Wilson}

\label{sec:lem34}

The purpose of this section is to prove an analogue of the following lemma from Kenyon and Wilson \cite{KW2006} for black-white pairings. 
%graphs that do not necessarily have the property that all odd nodes are black and all even nodes are white. 

\begin{lemma}\cite[Lemma 3.4]{KW2006}
\label{lem:kw34}
For odd-even pairings $\rho$, 
$$\sign_{OE}(\rho) \prod\limits_{(i, j) \in \rho} (-1)^{(|i-j|-1)/2} = (-1)^{\# \text{ crosses of } \rho}.$$
\end{lemma}

A {\em cross} of a pairing $\rho$ is a set of two pairs $(a, c)$ and $(b, d)$ of $\rho$ such that $a < b < c < d$. Recall from Section~\ref{sec:organization}
 that the sign of an odd-even pairing
$\rho = ((1, \rho(1)), (3, \rho(3)), \ldots, (2n-1, \rho(2n-1)))$
 is the 
parity of the permutation 
$\begin{pmatrix}
\frac{\rho(1)}{2} & \frac{\rho(3)}{2} & \cdots & \frac{\rho(2n-1)}{2}
\end{pmatrix}$
written in one-line notation.

For our version of this lemma,
% will hold for black-white pairings that may not be odd-even. Therefore, 
we need to define the sign of a black-white pairing $\rho$, which we will denote $\sign_{BW}(\rho)$. 

\begin{defn}
\label{def:signbw}
 If $\rho$ is a black-white pairing, then we can write $\rho = ((b_1, w_1), (b_2, w_2), \ldots, (b_n, w_n))$, where $b_1 < b_2 < \cdots < b_n$. 
 Let $r: \{w_1, \ldots, w_n\} \to \{1, \ldots, n\}$ be the map 
 defined by \\
\mbox{$r(k) = \#\{i:w_i \leq w_k\}.$}
 %that relabels $w_1, \ldots, w_n$ with the labels $1, \ldots, n$ so that $\min\limits_{1 \leq k \leq n} \{w_{k} \}$ is relabeled $1,\ldots,$ $\max\limits_{1 \leq k \leq n} \{w_{k} \}$ is relabeled $n$. 
 Then the sign of $\rho$, denoted $\sign_{BW}(\rho)$, is the parity of the permutation
$$\sigma_{\rho} =  \begin{pmatrix}
 r(w_1) & r(w_2) & \cdots & r(w_n)
 \end{pmatrix}$$
 written in one-line notation. 
% The sign of $\rho$, denoted $\sign_{BW}(\rho)$, is the parity of the permutation
% $w_1 w_2 \cdots w_n$, when $w_1, \ldots, w_n$ are relabeled with the labels $1, \ldots, n$ so that $\min\limits_{1 \leq k \leq n} \{w_{k} \}$ is relabeled $1,\ldots,$ $\max\limits_{1 \leq k \leq n} \{w_{k} \}$ is relabeled $n$.
\end{defn}

When $\rho$ is a pairing that is both black-white and odd-even, these signs agree.

\begin{lemma}
\label{lem:OEandBWsigns}
If $\rho$ is a black-white pairing that is also odd-even, then $\sign_{OE}(\rho) = \sign_{BW}(\rho)$.
\end{lemma}

The proof of Lemma~\ref{lem:OEandBWsigns} is straightforward, but it is postponed to Section~\ref{sec:OEandBW} for clarity of exposition.

In Lemma~\ref{lem:kw34}, the sign of a pair $(i, j)$ of $\rho$ is $(-1)^{(|i-j|-1)/2}$.  
If $\rho$ is a black-white pairing that is not odd-even and $(b, w)$ is a pair in $\rho$, it is not necessarily the case that $\frac{ |b-w| - 1}{2}$ is an integer. Therefore we need a different way to define the sign of a pair. 

To motivate this definition, notice that if two nodes of the opposite color $b$ and $w$ have the same parity, it cannot be the case that the nodes between $b$ and $w$ alternate black and white. Therefore we must keep track of the number of consecutive nodes of the same color between $b$ and $w$. Consecutive nodes of the same color appear in pairs. For example, if we have a graph with eight nodes so that nodes $1, 3, 4,$ and $6$ are black and nodes $2, 5, 7, 8$ are white, there are two pairs of consecutive nodes of the same color: $(3, 4)$ and $(7, 8)$. Since we frequently use the term pair when describing pairings of the nodes, we will refer to pairs of consecutive nodes as {\em couples of consecutive nodes} instead.

\begin{defn}
If $(b, w)$ is a pair in a black-white pairing,
let $a_{b, w}$ be the number of couples  of consecutive nodes of the same color in the interval $[\min\{b, w\}, \min\{b, w\} + 1, \ldots, \max\{b, w\}]$. 
\end{defn}

We note that
a triple of consecutive nodes that are all the same color contributes 2 to $a_{b,w}$.

\begin{rem}
\label{rem:signiswelldefined}
If $(b,w)$ is a pair in a black-white pairing, then
 $\frac{|b-w| +a_{b, w} -1}{2}$
 is an integer. 
 \end{rem}

\begin{proof}
Let $(n_1, n_1 + 1), (n_2, n_2 + 1), \ldots, (n_{2k}, n_{2k} + 1)$ 
be a complete list of couples of consecutive nodes of the same color in ${\bf N}$ so that $n_1 < n_2 < \cdots < n_{2k}$, where it is possible that $n_{i+1} = n_{i} + 1$.
Every time we reach a couple of consecutive nodes, the black nodes and white nodes switch parity.
That is, if the black nodes in the interval $[n_{\ell} + 1,n_{\ell} + 2, \ldots, n_{\ell+1}]$ are odd, then the black nodes
in the interval $[n_{\ell+1} + 1,n_{\ell+1} + 2, \ldots, n_{\ell+2}]$ are even. (Note that these intervals could be length 1). 
It follows that if $b$ and $w$ are the same parity, then there are an odd number of couples of consecutive nodes in the interval $[\min\{b, w\}, \min\{b, w\} + 1, \ldots, \max\{b, w\}]$.
So in this case $\frac{|b-w| +a_{b, w} -1}{2}$ is an integer. 
If $b$ and $w$ are opposite parity, then there are an even number of couples of consecutive nodes in the interval $[\min\{b, w\}, \min\{b, w\} + 1, \ldots, \max\{b, w\}]$.
So $\frac{|b-w| +a_{b, w} -1}{2}$ is an integer in this case as well. 
\end{proof}

\begin{defn}
\label{def:signpair}
If $(b, w)$ is a pair in a black-white pairing, let
\[
\text{sign}(b, w) = 
(-1)^{(|b-w|+ a_{b, w}-1)/2} .
\]
\end{defn}

We observe that
when the nodes of $G$ alternate black and white, $a_{b, w} = 0$ for all pairs $(b, w)$, so this definition of the sign of a pair agrees with Kenyon and Wilson's definition.

\begin{rem}
\label{rem:notation}
For the remainder of the paper, we use the following notation. We let %from the proof of Remark \ref{rem:signiswelldefined}. 
%That is, we let
\begin{itemize}
\item $(n_1, n_1 + 1), (n_2, n_2 + 1), \ldots, (n_{2k}, n_{2k} + 1)$ 
be a complete list of couples of consecutive nodes of the same color so that $n_1  < \cdots < n_{2k}$,
\item $(s_1, s_1 + 1), (s_2, s_2 + 1), \ldots, (s_k, s_k+1)$ be a complete list of couples of consecutive black nodes so that $s_1 < \cdots < s_{k}$, and
\item $(u_1, u_1 + 1), (u_2, u_2 + 1), \ldots, (u_k, u_k+1)$ be a complete list of couples of consecutive white nodes so that $u_1 < \cdots < u_k$.
\end{itemize}
Note that we could have $n_{i+1} = n_{i} + 1$, $s_{i+1} = s_i + 1$, or $u_{i+1} = u_i + 1$. 
\end{rem}

Since we are allowing arbitrary node colorings, many of our results
contain a global sign that depends on the order in which the couples of consecutive nodes appear.
% we need to define a global sign that depends on the node coloring. This sign  
 For example, suppose a node set {\bf N} has two couples of consecutive nodes: a couple of consecutive black nodes $(s, s+1)$ and a couple of consecutive white nodes $(u, u+1)$. Then the global sign will be $1$ if $u < s$ and $-1$ if $s < u$. To emphasize that this sign only depends on the relative ordering of the couples of consecutive nodes of the same color, we use the notation $\sign_{\cons}({\bf N})$. 

\begin{defn} 
\label{def:nodesign}

Using the notation from Remark~\ref{rem:notation}, 
if node $1$ is black, 
define the map $\varphi: \{n_1, n_2, \ldots, n_{2k} \} \to \{1,2, \ldots, 2k\}$ by
$$\varphi(n_j) = 
\begin{cases}
2i-1 & \mbox{if } n_j = u_i \\
2i & \mbox{if } n_j = s_i
\end{cases}.
$$
Then the image of $\{n_1, n_2, \ldots, n_{2k}\}$ under the map $\varphi$ can be considered as a permutation in one-line notation:
$$\sigma_{\bf N} = \begin{pmatrix}
\varphi(n_1) & \varphi(n_2) & \cdots & \varphi(n_{2k})
\end{pmatrix}.$$
 Define $\sign_{\cons}({\bf N})$ to be the sign of this permutation. Note that if $u_1 < s_1 < u_2 < s_2< \cdots < u_k < s_k$ then $\sigma_{\bf N} = 
\begin{pmatrix}
1 & 2 & \cdots & 2k
\end{pmatrix}$,
so $\sign_{\cons}({\bf N})= 1$.

If node $1$ is white,
define the map $\varphi: \{n_1, n_2, \ldots, n_{2k} \} \to \{1,2, \ldots, 2k\}$ by
$$\varphi(n_j) = 
\begin{cases}
2i-1 & \mbox{if } n_j = s_i \\
2i & \mbox{if } n_j = u_i
\end{cases}.
$$
As above, the image of $\{n_1, n_2, \ldots, n_{2k}\}$ under the map $\varphi$ can be considered as a permutation in one-line notation and we define $\sign_{\cons}({\bf N})$ to be the sign of this permutation. Note that if $s_1 < u_1 < s_2 < u_2< \cdots < s_k < u_k$, $\sign_{\cons}({\bf N})= 1$. 
%If there is some other ordering of $u_1, \ldots, u_k, s_1, \ldots, s_k$, write $u_1, \ldots, u_k, s_1, \ldots, s_k$, down in order from smallest to largest. Replace $s_1$ with $1$, $u_1$ with $2$, $s_2$ with $3$, $u_2$ with $4$, $\ldots$, $s_k$ with $2k-1$ and $u_k$ with $2k$. This can be considered as a permutation in one-line notation. Define $\sign_{\cons}({\bf N})$ to be the sign of this permutation.

In the case where there are no consecutive nodes of the same color, we define $\sign_{c}({\bf N}) =1$. 
\end{defn}

%\begin{rem*}
In Definition \ref{def:nodesign}, if node $1$ is black, it is possible that $s_k = 2n$.
Similarly, if node $1$ is white, it is possible that $u_k = 2n$. 
%\end{rem*}

\begin{defn}
\label{defn:inversioninnodecolors}
 Since the image of $\{n_1, n_2, \ldots, n_{2k} \}$ under the map $\varphi$ can be considered a permutation in one-line notation, we say that a pair $(u_\ell, s_m)$ is an inversion with respect to the node coloring of {\bf N} if $(\varphi(u_{\ell}), \varphi(s_m))$ is an inversion of $\sigma_{{\bf N}}$.
\end{defn}

\begin{example} Let ${\bf N}$ be a set of nodes where node 1 is black.
\begin{itemize}
\item If {\bf N} has four couples of consecutive nodes of the same color with $u_1 < s_1 < s_2 < u_2$, then $\sigma_{\bf N} = \begin{pmatrix} 1& 2& 4 & 3 \end{pmatrix}$, so $\sign_{\cons}({\bf N}) = -1$. The pair $(s_2, u_2)$ is an inversion with respect to the node coloring.
% because getting from the ordering $u_1 < s_1 < u_2 < s_2$ to $u_1 < s_1 < s_2 < u_2$ requires swapping the locations of $s_2$ and $u_2$ (1 transposition). 
\item If instead $s_1 < u_1 < s_2 < u_2$, then $\sigma_{\bf N} = \begin{pmatrix} 2& 1& 4 & 3 \end{pmatrix}$, so
$\sign_{\cons}({\bf N}) = 1$. The pairs $(s_1, u_1)$ and $(s_2, u_2)$ are inversions. % getting from the ordering $u_1 < s_1 < u_2 < s_2$ to $u_1 < s_1 < s_2 < u_2$ requires swapping the locations of $s_1$ and $u_1$ and swapping the locations of $s_2$ and $u_2$ (2 transpositions).
\end{itemize}
\end{example}

\begin{example} Let ${\bf N}$ be a set of nodes where node 1 is white.
If {\bf N} has six couples of consecutive nodes of the same color with $s_1 < s_2 < u_1 < u_2 < u_3 < s_3$, then  $\sigma_{\bf N} = \begin{pmatrix} 1&3& 2 & 4 & 6 & 5 \end{pmatrix}$. The pairs $(s_2, u_1)$ and $(u_3, s_3)$ are inversions. 
\end{example}

\begin{rem}
\label{rem:inversioninnodecolors}
 If node 1 is black, $(u_\ell, s_m)$ is an inversion with respect to the node coloring when $u_\ell < s_m$ and $\ell > m$. The pair $(s_m, u_\ell)$ is an inversion when $s_m < u_\ell$ and $m \geq \ell$. 
If node 1 is white, $(u_\ell, s_m)$ is an inversion with respect to the node coloring when $u_\ell < s_m$ and $\ell \geq m$. The pair $(s_m, u_\ell)$ is an inversion when $s_m < u_\ell$ and $m > \ell$. 
\end{rem}

We have now established the definitions needed for
%Now we are ready to state 
our version of Kenyon and Wilson's lemma. 

\begin{lemma}[analogue of Lemma 3.4 from \cite{KW2006}]
\label{lemma34}
If $\rho$ is a black-white pairing,
% on a graph $G$ with node set ${\bf N}$, 
\begin{equation*}
\sign_{\cons}({\bf N}) \sign_{BW}(\rho) \prod\limits_{(b, w) \in \rho} \sign(b, w) =  (-1)^{\# \text{ crosses of } \rho}.
\end{equation*}
%where $\text{sign}(b, w)$ is defined in Definition \ref{def:signpair}. 
%Specifically, if the nodes of $G$ meet the assumptions in Assumptions \ref{assumption2}, 6(b), (c), (d), or (f), we have
%\begin{equation*}
%\label{eqn1:lemma34}
%-\sign_{BW}(\rho) \prod\limits_{(b, w) \in \rho} \sign(b, w) =  (-1)^{\# \text{ crosses of } \rho}
%\end{equation*}
%If the nodes meet the assumptions in Assumptions \ref{assumption3}, \ref{assumption4}, \ref{assumption5}, or 6(a) or (e),
%\begin{equation*}
%\label{eqn2:lemma34}
%\sign_{BW}(\rho) \prod\limits_{(b, w) \in \rho} \sign(b, w) =  (-1)^{\# \text{ crosses of } \rho}.
%\end{equation*}
\end{lemma}

We remark that in Kenyon and Wilson's case where all black nodes are odd and all white nodes are even, there are no consecutive nodes of the same color, so for all $(b, w) \in \rho$, $a_{b, w} = 0$ and thus $\sign(b, w)= (-1)^{(|b-w|+ a_{b, w}-1)/2} = (-1)^{(|b-w| -1)/2}$. If all black nodes are odd and all white nodes are even, a black-white pairing is also an odd-even pairing, and by Lemma~\ref{lem:OEandBWsigns}, $\sign_{BW}(\rho) = \sign_{OE}(\rho)$. Finally, by Definition~\ref{def:nodesign}, $\sign_{\cons}({\bf N}) = 1$. So in this case, Lemma~\ref{lemma34} agrees exactly with Lemma~\ref{lem:kw34}. 

Before proving Lemma~\ref{lemma34}, we will prove the following:

\begin{lemma}
\label{firstlemma34}
%Let $G$ be a graph with node set {\bf N}. 
There exists a planar black-white pairing $\rho$ such that
\begin{equation*}
 \sign_{BW}(\rho) \prod\limits_{(b, w) \in \rho} \sign(b, w) =  \sign_{\cons}({\bf N}) .
\end{equation*}
%where $\text{sign}(b, w)$ is defined in Definition \ref{def:signpair}. 
\end{lemma}

\subsubsection{Proof of Lemma \ref{firstlemma34}}

We will prove Lemma~\ref{firstlemma34} by induction on $k$, where {\bf N} is a set of $2n$ nodes with $2k$ couples of consecutive nodes of the same color. The following lemma is the base case $k=1$. 

\begin{lemma}[Base case of Lemma~\ref{firstlemma34}] 
\label{lem:planarBWrho}
For any node coloring such that there are exactly two couples of consecutive nodes of the same color, there is a planar black-white pairing $\rho$ such that 
$$\sign_{BW}(\rho) \prod\limits_{(b, w) \in \rho} \sign(b, w) = \sign_{\cons}({\bf N}).$$
\end{lemma}

\begin{proof}
%There are several cases to consider.
%As usual, 
Let $(n_1, n_1 + 1), (n_2, n_2 + 1)$ be the list of the couples of consecutive nodes of the same color so that $n_1 < n_2$. 
%Let $(s, s+1)$ denote the couple of consecutive black nodes and let $(u, u+1)$ denote the couple of consecutive white nodes. 
There are two cases to consider: Either $n_1$ and $1$ are opposite colors, or $n_1$ and $1$ are the same color.

If $n_1$ and $1$ are opposite colors, 
the pairing  $\rho = ((1,2), (3, 4), \ldots, (2n-1, 2n))$ is black-white. 
To see this, note that since $n_1$ and $1$ are opposite colors, $n_1$ is even, so the only pairs of adjacent nodes that are both the same color are of the form 
$(x, x+1)$, where $x$ is even, or $(2n, 1)$.
Since all pairs of $\rho$ are of the form $(i, i +1)$ where $i$ is odd and $i+1$ is even,
$\rho$ is a black-white pairing. By the previous definitions we deduce that $\sign_{BW}(\rho) = 1$ and  $\prod\limits_{(b, w) \in \rho} \sign(b, w) = 1$. Since we assumed $n_1$ and $1$ are opposite colors, either $u_1 < s_1$ (if $1$ is black) or $s_1 < u_1$ (if $1$ is white) and thus $\sign_{\cons}({\bf{N}}) =1$ by Definition~\ref{def:nodesign}. Therefore the claim holds.

If $n_1$ and $1$ are the same color,
the pairing  $\rho = ((2n, 1), (2, 3), \ldots, (2n-2, 2n-1))$ is black-white. 
The reasoning is analogous to the previous case: 
%note that since $n_1$ and $1$ are the same color, 
$n_1$ is odd, so the only pairs of adjacent nodes that are both the same color are of the form 
$(x, x+1)$, where $x$ is odd. 
In this case, 
$\sign_{BW}(\rho) = (-1)^{n-1}$ and $\prod\limits_{(b, w) \in \rho} \sign(b, w) = (-1)^{(2n-1 + a_{2n, 1} -1)/2} = (-1)^{n-1}(-1)^{a_{2n, 1}/2}$, so 
$$\sign_{BW}(\rho) \prod\limits_{(b, w) \in \rho} \sign(b, w) = (-1)^{a_{2n, 1}/2} =(-1)^{2/2} = -1.$$ Similar to the previous case, since we assumed $n_1$ and $1$ are the same color, $\sign_{\cons}({\bf{N}}) = -1$, and thus the claim holds.
\end{proof}

\begin{defn}
\label{def:invofrho}
Suppose $\rho$ is a black-white pairing. 
Then recall 
that
we can write \\ $\rho = ((b_1, w_1), (b_2, w_2), \ldots, (b_n, w_n))$, where $b_1 < b_2 < \cdots < b_n$. 
 We say that $(w_i, w_j)$ is an inversion of $\rho$ if $i < j$ and $w_i > w_j$. 
 Note that $(w_i, w_j)$ is an inversion of $\rho$ if and only if $(r(w_i), r(w_j) )$ is an inversion of $\sigma_{\rho}$ (see Definition \ref{def:signbw}). 
 \end{defn}

 \begin{defn}
 \label{defn:invofpi}
 Similarly, if $\pi$ is an odd-even pairing,  then we can write\\
$\pi = ((1, \pi(1)), (3, \pi(3)), \ldots, (2n-1, \pi(2n-1)))$
and we will say $(\pi(i), \pi(j))$ is an inversion of $\pi$ if $i < j$ and $\pi(i) > \pi(j)$. This is equivalent to defining an inversion of $\pi$ to be an inversion of the permutation
$\begin{pmatrix}
\frac{\pi(1)}{2} & \frac{\pi(3)}{2} & \cdots & \frac{\pi(2n-1)}{2}
\end{pmatrix}$.
 \end{defn}%}

 \begin{figure}[h!]
\centering
\begin{tikzpicture}[scale = .5]
	\vertex[fill] (n1) at (0, 0) [label=below:$1$] {};
	\vertex[fill]  (n2) at (2,0) [label=below:$2$] {};
	\vertex (n3) at (4,0) [label=below:$3$] {};
	\vertex (n4) at (6, 0) [label=below:$4$] {};
	
	%arcs of pi
	\draw  (n2) arc (0:180:1cm);
	\draw  (n4) arc (0:180:1cm);
	
			\vertex[fill=white] (n3) at (4,0) {};
	\vertex[fill=white] (n4) at (6, 0) {};
	
	%arcs of rho

\end{tikzpicture} \hspace{0.5cm}
\begin{tikzpicture}[scale = .5]
	\vertex[fill] (n1) at (0, 0) [label=below:$1$] {};
	\vertex[fill]  (n2) at (2,0) [label=below:$2$] {};
	\vertex (n3) at (4,0) [label=below:$3$] {};
	\vertex (n4) at (6, 0) [label=below:$4$] {};
	
	%arcs of pi
	\draw  (n3) arc (0:180:1cm);
	\draw  (n4) arc (0:180:3cm);
			\vertex[fill=white] (n3) at (4,0) {};
	\vertex[fill=white] (n4) at (6, 0) {};
	%arcs of rho

\end{tikzpicture} \hspace{0.5cm}
\begin{tikzpicture}[scale = .5]
	\vertex[fill] (n3) at (4, 0) [label=below:$1$]  {};
	\vertex (n4) at (6,0)[label=below:$2$]  {};
	\vertex[fill] (n5) at (8,0) [label=below:$3$]  {};
	\vertex[fill] (n6) at (10, 0) [label=below:$4$]  {};
	\vertex (n7) at (12,0) [label=below:$5$] {};
	\vertex[fill] (n8) at (14,0) [label=below:$6$] {};
	\vertex (n11) at (16,0) [label=below:$7$] {};
	\vertex (n12) at (18, 0) [label=below:$8$] {};
	
	%arcs of pi
	\draw  (n12) arc (0:180:7cm);
	\draw  (n11) arc (0:180:1cm);
	\draw  (n7) arc (0:180:3cm);
	\draw  (n6) arc (0:180:1cm);
	
	%arcs of rho
	%\draw  (n3) arc (180:360:4cm);
	%\draw  (n4) arc (180:360:1cm);
%	\draw  (n6) arc (180:360:4cm);
%	\draw  (n8) arc (180:360:1cm);
	
		\vertex[fill=white] (n4) at (6,0) {};
		\vertex[fill=white]  (n7) at (12,0) {};
			\vertex[fill=white] (n11) at (16,0)  {};
	\vertex[fill=white] (n12) at (18, 0) {};

\end{tikzpicture}

\caption{An inversion of a planar pairing $\pi$ corresponds to a nesting. Left: the pairing $((1, 2), (3, 4))$ has no inversions and its diagram has no nestings. Center: The pairing $((1,4), (3, 2))$ has one inversion and its diagram has one nesting. Right: The pairing $((1, 8), (3, 4), (5, 2), (7, 6))$ has four inversions and four nestings: $\{(1, 8), (3, 4 )\}$, $\{(1, 8), (5, 2 )\}$,  $\{(1, 8), (7, 6 )\}$, and $\{(3, 4), (5, 2)\}$. }
\label{fig:nesting}
\end{figure}
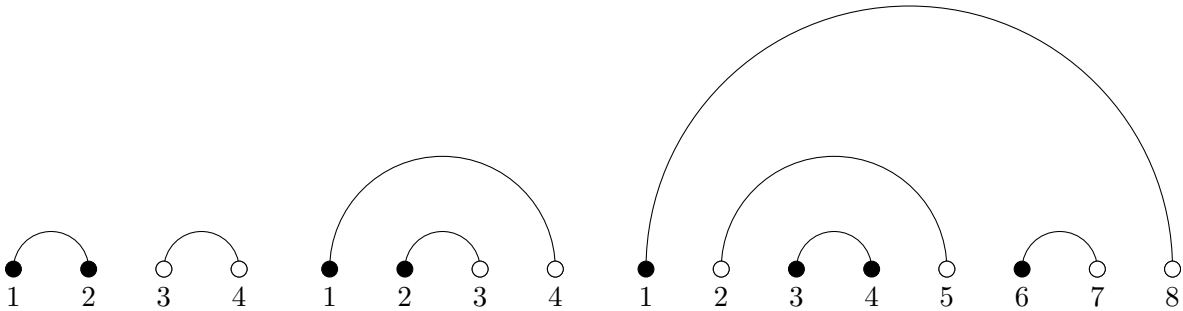

  \begin{rem}
\label{rem:inversionsarenestings}
 In the special case where an odd-even pairing $\pi$ is also planar, we remark that
an inversion of $\pi$ corresponds to a {\em nesting} in the  diagram constructed by placing the nodes in order on a line and linking pairs in the upper half-plane. More precisely, two arcs $(a_1, b_1)$, $(a_2, b_2)$ 
%with $a_1 < b_1$ and $a_2 < b_2$ 
are said to be {\em nesting} if $a_1 < a_2 < b_2 < b_1$  (see Figure~\ref{fig:nesting}). This correspondence between inversions and nestings follows immediately from the four node case, where the only planar pairings are $((1, 2), (3, 4))$ and $((1, 4), (3, 2))$. 
 %  If we have four nodes $\{b_1, b_2, w_1, w_2\}$ with $b_1 < b_2$ and $w_1 < w_2$, then the only black-white pairing with an inversion is $((b_1, w_2), (b_2, w_1))$. So we check that regardless of the ordering of these four nodes, this pairing is either nonplanar or a nesting. 
\end{rem}

\begin{proof}[Proof of Lemma \ref{firstlemma34}]

The proof of the lemma is technical, so we first identify a few easy cases. \\

\noindent {\bf Easy case 1.} If
\begin{itemize}
\item node 1 is black and $u_1 < s_1 < u_2 < s_2 < \cdots < u_k < s_k$, or
\item node 1 is white and $s_1 < u_1 < s_2 < u_2 < \cdots <s_k < u_k$,
\end{itemize} 
 then as in the proof of the first case of Lemma~\ref{lem:planarBWrho}, the pairing $((1, 2), (3, 4), \ldots, (2n-1, 2n))$ is a planar black-white pairing with $\sign_{BW}(\rho) = 1$, $\prod\limits_{(b, w) \in \rho} \sign(b, w) = 1$, and $\sign_{\cons}({\bf N}) = 1$. Thus the claim holds. \\

\noindent {\bf Easy case 2.} 
If
\begin{itemize}
\item node 1 is black and $s_1 < u_1 < s_2 < u_2 < \cdots <s_k < u_k$, or
\item node 1 is white and $u_1 < s_1 < u_2 < s_2 < \cdots < u_k < s_k$,
\end{itemize} 
then the pairing $\rho = ((2n, 1), (2, 3), \ldots, (2n-2, 2n-1))$ is black-white (as in the second case of Lemma~\ref{lem:planarBWrho}). 
In this case,
\begin{itemize}
\item $\sign_{BW}(\rho) =(-1)^{n-1}$, and 
\item $\prod\limits_{(b, w) \in \rho} \sign(b, w) = (-1)^{(2n-1 + a_{2n, 1} -1)/2} = (-1)^{n-1}(-1)^{a_{2n, 1}/2} = (-1)^{n-1}(-1)^{k}$,
\end{itemize}
so $\sign_{BW}(\rho) \prod\limits_{(b, w) \in \rho} \sign(b, w) = (-1)^{k} = \sign_{\cons} ({\bf N})$. \\

\noindent {\bf General case.} 
For the general case, we proceed by
induction on the number of couples of consecutive nodes of the same color. The base case is when there are two couples of consecutive nodes of the same color, which is Lemma \ref{lem:planarBWrho}. 
Assume the claim holds when we have a set of nodes that has $2(k-1)$ couples of consecutive nodes of the same color and let ${\bf N}$ be a set of nodes with $2k$ couples of consecutive nodes of the same color.
%\noindent {\Large{REWRITE \# 2}}
%\noindent {\bf General Case.} 

\noindent
\begin{minipage}{.75\textwidth}
\hspace{10pt} Using the notation from Remark~\ref{rem:notation},
%Let $(n_1, n_1 + 1), (n_2, n_2 + 1), \ldots, (n_{2k}, n_{2k} +1)$ be a list of the couples of consecutive nodes of the same color so that $n_1 < n_2 < \cdots < n_{2k}$.
 let $h$ be the smallest integer so that $n_{h-1}$ and $n_{h}$ are different colors. 
Then 
$\rho_1 = ((n_{h-1} + 1, n_{h-1} + 2), \ldots, (n_{h}-1, n_{h}))$
is a black-white pairing that contains at least one pair. 

\hspace{10pt} Throughout this proof, we will illustrate the main ideas with the example where ${\bf N}$ is a set of 12 nodes colored so that nodes 1, 3, 4, 5, 7, and 10 are black, as shown to the right. In this example, the couples of consecutive nodes of the same color are $(3, 4), (4, 5), (8, 9),$ and $(11, 12)$. Since $n_1 = 3$ and $n_2 = 4$ are black and $n_3 = 8$ is white, $h = 3$. So the pairing $\rho_1$ is $((5, 6), (7, 8))$.
\end{minipage} \hfill
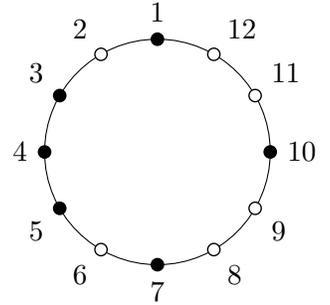
\begin{minipage}{.23\textwidth}
\begin{center}
 \begin{tikzpicture}[scale=.75]
  \draw (0,0) circle (2);
  \foreach \x in {1,2,...,12} {
   \node[shape=circle,fill=black, scale=0.5,label={{((\x-1)*360/12)+90}:\x}] (n\x) at ({((\x-1)*360/12)+90}:2) {}; };
     \foreach \x in {2, 6, 8, 9, 11, 12} {
     \node[shape=circle,fill=white, scale=0.4] (n\x) at ({((\x-1)*360/12)+90}:2) {};
  };
 \end{tikzpicture}
 \end{center}
\end{minipage}

Consider 
${\bf N'} = \{1, \ldots, |{\bf N}| - (n_h - n_{h-1}) \}$. Define
$\psi: {\bf N} - \{n_{h-1} + 1, \ldots, n_h \} \to {\bf N'}$ by
\begin{equation}
\label{eqn:relabeling}
\psi(\ell) =
\begin{cases}
 \ell &\mbox{ if }\ell \leq n_{h-1} \\
  \ell - (n_h - n_{h-1}) &\mbox{ if }\ell > n_{h} \\
 \end{cases}
\end{equation}
That is, $\psi$ defines a relabeling of the nodes of ${\bf N} - \{n_{h-1} + 1, \ldots, n_h \}$ so that node 1 is labeled $1,\ldots,$ node $n_{h-1}$ is labeled $n_{h-1}$, node $n_{h} + 1$ is labeled $n_{h-1} +1,\ldots,$ node $2n$ is labeled $2n - (n_{h} - n_{h-1})$. 
Since ${\bf N'}$ has $2k-2$ couples of consecutive nodes of the same color, by the induction hypothesis there is a black-white planar pairing $\rho_2$ of the nodes of ${\bf N'}$ such that
$$\sign_{BW}(\rho_2) \prod\limits_{(b, w) \in \rho_2} \sign(b, w) = \sign_{\cons}({\bf N'}).$$ 
Let $\psi^{-1}(\rho_2)$ denote the pairing that results from applying $\psi^{-1}$ to each node in $\rho_2$. 
That is, $\psi^{-1}(\rho_2)$ is the pairing obtained by returning the nodes of $\rho_2$ to their original labels in ${\bf N}$. Let $\rho = \rho_1 \cup \psi^{-1} (\rho_2)$. 
Observe that $\rho$ is a planar black-white pairing of ${\bf N}$.

In our example, the map $\psi$ defines a relabeling of ${\bf N} - \{5, 6, 7, 8\}$ so that node 9 is labeled $5,\ldots,$ node 12 is labeled 8. The node set ${\bf N'}$ has two couples of consecutive pairs of the same color. By Lemma \ref{lem:planarBWrho}, the pairing $\rho_2$ is $((1, 8), (3, 2), (5, 4), (7, 6))$, so the pairing $\psi^{-1}(\rho_2)$ is 
$((1, 12), (3, 2), (9,  4), (11, 10))$ and thus $\rho = ((1, 12), (3, 2), (5, 6), (7, 8), (9, 4), (11, 10))$, as shown in Figure~\ref{fig:lem5ex}.

We will next
\begin{enumerate}
\item[(1)] Compare $ \prod\limits_{(b, w) \in \rho_2} \sign(b, w)$
to
$ \prod\limits_{(b, w) \in \rho} \sign(b, w),$
\item[(2)] Compare $\sign_{BW}(\rho_2)$ to $\sign_{BW}(\rho)$, and
\item[(3)] Compare $\sign_{\cons}({\bf N'})$ to $\sign_{\cons}({\bf N})$.
\end{enumerate}

\begin{figure}
\centering
 \begin{tikzpicture}[scale=.75]
  \draw (0,0) circle (2);
  \foreach \x in {1,2,...,8} {
   \node[shape=circle,fill=black, scale=0.5,label={{((\x-1)*360/8)+90}:\x}] (n\x) at ({((\x-1)*360/8)+90}:2) {}; };
     \foreach \x in {2, 5, 7, 8} {
     \node[shape=circle,fill=white, scale=0.4] (n\x) at ({((\x-1)*360/8)+90}:2) {};
  };
    \foreach \x/\y in {1/8, 2/3, 4/5, 6/7} {
   \draw (n\x) -- (n\y);};
 \end{tikzpicture} \hspace{.5cm}
 \begin{tikzpicture}[scale=.75]
  \draw (0,0) circle (2);
  \foreach \x in {1,2,...,12} {
   \node[shape=circle,fill=black, scale=0.5,label={{((\x-1)*360/12)+90}:\x}] (n\x) at ({((\x-1)*360/12)+90}:2) {}; };
     \foreach \x in {2, 6, 8, 9, 11, 12} {
     \node[shape=circle,fill=white, scale=0.4] (n\x) at ({((\x-1)*360/12)+90}:2) {};
  };
      \foreach \x/\y in {5/6, 7/8} {
   \draw[color = red] (n\x) -- (n\y);};
         \foreach \x/\y in {1/12, 2/3, 4/9, 10/11} {
   \draw[densely dotted, thick, color = blue] (n\x) -- (n\y);};
 \end{tikzpicture}
 \caption{Left: The pairing $\rho_2$ of ${\bf N'}$ guaranteed by the induction hypothesis. Right: The pairing $\rho$ of ${\bf N}$. The pairing $\rho_1$ is shown with solid red lines and the pairing $\psi^{-1}(\rho_2)$ is shown with dotted blue lines.}
 \label{fig:lem5ex}
\end{figure}
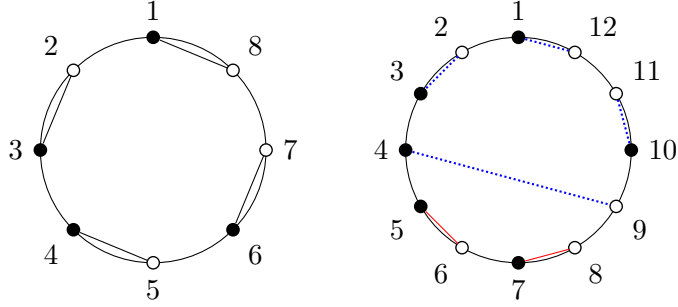

\paragraph{(1) Comparing $ \prod\limits_{(b, w) \in \rho_2} \sign(b, w)$
to
$\prod\limits_{(b, w) \in \rho} \sign(b, w)$.}

If $(b, w)$ is a pair in $\rho$ that is a pair of $\rho_1$, then $\sign(b, w) = 1$. If $(b, w)$ is a pair in $\rho$ that is a pair of $\psi^{-1} (\rho_2)$, then consider $(\psi(b), \psi(w))$ (the corresponding pair of $\rho_2$). If $b, w\leq n_{h-1}$ or $b, w\geq n_{h} + 1$, then $\sign(b, w) = \sign( \psi(b), \psi(w) )$ because $a_{b, w} = a_{\psi(b), \psi(w) }$. If $w\leq n_{h-1}$ and $b \geq n_{h} +1$ then
\begin{eqnarray*}
\sign( \psi(b), \psi(w) ) = (-1)^{ ( \psi(b) - \psi(w) + a_{\psi(b), \psi(w)} -1 )/2 }
& = &(-1)^{( b - (n_{h} - n_{h-1}) - w + a_{b, w} - 2 -1 )/2} \\
& = & (-1)^{ (-(n_{h} - n_{h-1}) - 2 ))/2} \sign(b, w) \\
& = & (-1)^{ (n_{h} - n_{h-1} + 2 )/2} \sign(b, w). 
\end{eqnarray*}
so
$$\prod\limits_{(b, w)\in \rho_2 } \sign(b, w)
= \prod\limits_{ \substack{ (b, w) \in \rho: \\  \min(b, w) \leq n_{h-1}  \text{ and} \\
 \max(b, w) \geq n_{h} + 1
 } } (-1)^{ (n_{h} - n_{h-1} + 2 )/2}
 \prod\limits_{(b, w) \in \rho} \sign(b, w).
 $$
 
 In the example, there are two pairs $(b, w)$ with $\min(b, w) \leq n_2$ and $\max(b, w) \geq n_3 + 1$: the pairs $(1, 12)$ and $(4, 9)$.
 
 \paragraph{(2) Comparing $\sign_{BW}(\rho_2)$ to $\sign_{BW}(\rho)$.}
 
Comparing $\sign_{BW}(\rho_2)$ to $\sign_{BW}(\rho)$ requires comparing the number of inversions of $\rho$ to the number of inversions of $\rho_2$ (see Definition \ref{def:invofrho}). Since $\rho_1$ contains only pairs of the form $(i, i+1)$, $\rho_1$ contains no inversions. Since the pairings under consideration are planar we can use the fact that inversions of a planar pairing correspond to nestings in the corresponding diagram (see Remark \ref{rem:inversionsarenestings}). 
Since there are $\frac{ n_h - n_{h-1} }{2}$ pairs in $\rho_1$, $\rho$ has 
$\frac{ n_h - n_{h-1} }{2}$ additional inversions compared to $\rho_2$ for each pair $(b, w)$ such that $\min(b, w) \leq n_{h-1}$ and $\max(b, w) \geq n_h + 1$. So, 
$$\sign_{BW}(\rho_2)  = \sign_{BW}(\rho) \prod\limits_{ \substack{ (b, w) \in \rho: \\  \min(b, w) \leq n_{h-1}  \text{ and} \\
 \max(b, w) \geq n_{h} + 1
 } } (-1)^{ (n_{h} - n_{h-1} )/2}$$
 
In the example, since there are two pairs $(b, w)$ with $\min\{b, w\} \leq n_2$ and $\max\{b, w\} \geq n_3 + 1$ and the pairing $\rho_1$ consists of two pairs, there are four more inversions in $\rho$ than in $\rho_2$. 
 
  \paragraph{(3) Comparing $\sign_{\cons}({\bf N'})$ to $\sign_{\cons}({\bf N})$.}
 %Finally we compare
 %$\sign_{\cons}({\bf N'})$ to $\sign_{\cons}({\bf N})$ (see Definition \ref{def:nodesign}). 
% Recall that $h$ is the smallest integer so that $n_{h-1}$ and $n_h$ are different colors.
 We will show that
   $$\sign_{\cons}({\bf N'}) = (-1)^{h-1} \sign_{\cons}({\bf N})$$
by comparing the number of inversions with respect to the node coloring of ${\bf N}$ to the number of inversions in with respect to the node coloring of ${\bf N'}$ (see Definition \ref{defn:inversioninnodecolors}). In particular, we show that there are $(h-2) + (h-1)$ inversions in ${\bf N}$ that do not have a counterpart in ${\bf N'}$ and that there are $h-2$ inversions in ${\bf N'}$ that do not have a counterpart in ${\bf N}$. 

Recall the notation from Remark~\ref{rem:notation}: $s_i$ is the first in a couple of consecutive black nodes in ${\bf N}$ and $u_i$ is the first in a couple of consecutive white nodes in ${\bf N}$. Define $s'_i$ and $u'_i$ analogously for ${\bf N'}$. 
%Recall also that we have the map $\psi: {\bf N} - \{n_{h-1} + 1, \ldots, n_h \} \to {\bf N'}$ which defines a relabeling of the nodes of ${\bf N} - \{n_{h-1} + 1, \ldots, n_h \}$ (see equation (\ref{eqn:relabeling})). 

First assume node 1 is black and that we have 
$$s_1 < \cdots < s_{h-1} < u_1 < \cdots.$$

\noindent {\em Inversions with respect to the node coloring of ${\bf N}$.}
By Remark \ref{rem:inversioninnodecolors}, there are two types of inversions with respect to the node coloring of {\bf N}. 
 \begin{enumerate}
 \item[(1)] Nodes $x$ and $y$ in ${\bf N}$ such that $x < y$, $x= s_a$, $y = u_b$, and $a \geq b$. 
 \item[(2)] Nodes $x$ and $y$ in ${\bf N}$ such that $x < y$, $x = u_a$, $y = s_b$, and $a > b$. 
 \end{enumerate}

Considering the first type of inversion, there are several cases:
\begin{enumerate}
\item[(a)] If $a \leq h-2$ and $b > 1$, then $\psi(x) = s'_a$ and $\psi(y) = u'_{b-1}$. Since $a \geq b$, $a \geq b-1$, so in this case there is a corresponding inversion in ${\bf N'}$. 
\item[(b)] If $a \leq h-2$ and $b = 1$, then $y \notin {\bf N} - \{n_{h-1} + 1, \ldots, n_h \}$, so in this case there is not a corresponding inversion in ${\bf N'}$. 
\item[(c)] If $a = h-1$ and $b \leq h-1$, then $x \notin {\bf N} - \{n_{h-1} + 1, \ldots, n_h \}$, so in this case there is not a corresponding inversion in ${\bf N'}$. 
\item[(d)] If $a > h-1$ and $b \leq a$, then $b > 1$ (since $u_1 < s_a$). In this case, $\psi(x) = s'_{a-1}$ and $\psi(y) = u'_{b-1}$, so there is a corresponding inversion in ${\bf N'}$. 
\end{enumerate}

Note that (b) gives $h-2$ inversions in ${\bf N}$ that are not in ${\bf N'}$ and (c) gives $h-1$ inversions in ${\bf N}$ that are not in ${\bf N'}$. 

Considering the second type of inversion, since $s_{h-1} < u_1$ we must have $a > h$. In this case, $\psi(x) = u'_{a-1}$ and $\psi(y) = s'_{b-1}$, so there is a corresponding inversion in ${\bf N'}$.

In the example, the pairs $(s_1, u_1), (s_2, u_1)$, and $(s_2, u_2)$ are inversions with respect to the node coloring of ${\bf N}$. 
%These are all inversions of the first type. 
Since $h = 3$, the inversion $(s_1, u_1)$ is in case (b) of the first type and the inversions $(s_2, u_1)$ and $(s_2, u_2)$ are in case (c) of the first type. So in this example, all of the inversions with respect to the node coloring of ${\bf N}$ do not have corresponding inversions in ${\bf N'}$. \\

\noindent {\em Inversions with respect to the node coloring of ${\bf N'}$.}
Similarly, there are two types of inversions in ${\bf N'}$. 
 \begin{enumerate}
 \item[(1)] Nodes $w$ and $z$ in ${\bf N'}$ such that $w < z$, $w = s'_a$, $z = u'_b$, and $a \geq b$. 
 \item[(2)] Nodes $w$ and $z$ in ${\bf N'}$ such that $w < z$, $w = u'_a$, $z = s'_b$, and $a > b$. 
 \end{enumerate}
 
 Considering the first type of inversion, there are two cases:
 \begin{enumerate}
 \item[(a)] If $a \leq h-2$, then $\psi^{-1}(w)= s_a$ and $\psi^{-1}(z) = u_{b+1}$.
 \begin{enumerate}
 \item[(i)] If $a \geq b+1$ then there is a corresponding inversion in {\bf N}.
 \item[(ii)] If $a = b$ there is not a corresponding inversion in {\bf N}. 
 \end{enumerate}
 \item[(b)] If $a \geq h-1$, then $\psi^{-1}(w) = s_{a+1}$ and $\psi^{-1}(z) = u_{b+1}$, so there is a corresponding inversion in ${\bf N}$. 
 \end{enumerate}
 
 We see that case (a)(ii) gives $h-2$ inversions in ${\bf N}'$ that are not in ${\bf N}$. 

 Considering the second type of inversion, since $s'_{h-2} < u'_1$ the only possibility is that $a > h-1$. In this case, $\psi^{-1}(w) = u_{a+1}$ and $\psi^{-1}(z) = s_{b+1}$, so there is a corresponding inversion in ${\bf N}$.

In the example, the only inversion with respect to the node coloring of ${\bf N'}$ is $(s'_1, u'_1)$, which is an example of case (a)(ii), so there is not a corresponding inversion in ${\bf N}$. 

We conclude that in the case where node 1 is black and we have 
$s_1 < \cdots < s_{h-1} < u_1 < \cdots,$  the equation $\sign_{\cons}({\bf N'}) = (-1)^{h-1} \sign_{\cons}({\bf N})$
 holds.

%We have shown that the equation
%  $$\sign_{\cons}({\bf N'}) = (-1)^{h-1} \sign_{\cons}({\bf N})$$
% holds. 
Combining this with
\begin{itemize}
\item $\sign_{BW}(\rho_2) \prod\limits_{(b, w) \in \rho_2} \sign(b, w) = \sign_{\cons}({\bf N'}),$ 
 \item $\prod\limits_{(b, w) \in \rho_2 } \sign(b, w)
= \prod\limits_{ \substack{ (b, w)\in \rho: \\  \min(b, w)\leq n_{h-1}  \text{ and} \\
 \max(b, w) \geq n_{h} + 1
 } } (-1)^{ (n_{h} - n_{h-1} + 2 )/2}
 \prod\limits_{(b, w) \in \rho} \sign(b, w)$, and
 \item $\sign_{BW}(\rho_2)  = \sign_{BW}(\rho) \prod\limits_{ \substack{ (b, w) \in \rho: \\  \min(b, w) \leq n_{h-1}  \text{ and} \\
 \max(b, w) \geq n_{h} + 1
 } } (-1)^{ (n_{h} - n_{h-1} )/2}$,
\end{itemize}
 
we have
 
 $$ \sign_{BW}(\rho) 
  \prod\limits_{(b, w) \in \rho} \sign(b, w)
  =    \sign_{\cons}({\bf N}) \cdot
  (-1)^{h-1}  \cdot 
  \prod\limits_{ \substack{ (b, w) \in \rho: \\  \min(b, w) \leq n_{h-1}  \text{ and} \\
 \max(b, w) \geq n_{h} + 1
 } } (-1).
 $$

So it remains to observe that the number of pairs $(b, w) \in \rho$ such that $\min(b, w) \leq s_{h-1}$ and $\max(b, w) \geq u_1 + 1$ has the same parity as $h-1$. There are exactly $h-1$ more black nodes than white nodes in the interval $[1, \ldots, s_{h-1}]$ because there are $h-1$ black nodes that are not followed by a white node in this interval. So there are $h-1$ black nodes that must all be paired with a white node with label $\geq u_1 + 1$. Therefore there are at least $h-1$ pairs $(b, w) \in \rho$ such that $\min(b, w) \leq s_{h-1}$ and $\max(b, w) \geq u_1 + 1$. There may be more than $h-1$ such pairs, but there must be $h-1+ 2m$ pairs for some $m \geq 0$. 

 There are three other cases:
when node 1 is white and we have $s_1 < \cdots < s_{h-1} < u_1 < \cdots$, when node 1 is black and we have $u_1 < \cdots < u_{h-1} < s_1 < \cdots$, and when node 1 is white and we have $u_1 < \cdots < u_{h-1} < s_1 < \cdots$.
These are omitted because the analyses are nearly identical to the case we just considered.

\end{proof}

\subsubsection{Proof of Lemma~\ref{lemma34}}

Recall that we want to show that
if $\rho$ is a black-white pairing on a graph $G$ with node set ${\bf N}$, 
\begin{equation}
\label{eqn0:lemma34}
\sign_{\cons}({\bf N}) \sign_{BW}(\rho) \prod\limits_{(b, w) \in \rho} \sign(b, w) =  (-1)^{\# \text{ crosses of } \rho}.
\end{equation}
By Lemma \ref{firstlemma34} there is a black-white planar pairing $\rho$ such that 
$$\sign_{BW}(\rho) \prod\limits_{(b, w) \in \rho} \sign(b, w) =\sign_{\cons} ({\bf N}).$$
Since $\rho$ is planar, $(-1)^{\# \text{crosses of } \rho} = 1$, so equation (\ref{eqn0:lemma34}) holds. \\

To prove equation (\ref{eqn0:lemma34}) holds for all black-white pairings we consider ways we can modify black-white pairings to obtain new black-white pairings and determine how these modifications affect equation (\ref{eqn0:lemma34}). 

\begin{defn}
Let $\sigma$ be a (not necessarily black-white) pairing on $\{1,\ldots,2n\}$, such that $x$ is not paired with $y$. When we {\em swap the locations of $x$ and $y$ in $\sigma$} we create a new pairing $\sigma'$ that is identical to $\sigma$ except that it contains the pairs $(x, \sigma(y))$ and $(y, \sigma(x))$ rather than $(x, \sigma(x))$ and $(y, \sigma(y))$. 
\end{defn}

\begin{example}
Suppose $\sigma$ is the pairing $((1, 3), (2, 4), (5, 6))$. If we swap the locations of $3$ and $4$ in $\sigma$ we obtain the pairing $\sigma' = ((1, 4), (2, 3), (5, 6))$. 
\end{example}

\begin{rem}
\label{rem:swappingandparity}
If $\rho$ is a black-white pairing and $\rho'$ is obtained from $\rho$ by swapping the locations of two nodes of the same color, $\sign_{BW}(\rho') = -\sign_{BW}(\rho)$. 
\end{rem}

Now we observe that we can obtain any black-white pairing on ${\bf N}$ from a given black-white pairing $\rho$ using the following types of swaps:
\begin{enumerate}
\item[(1)] Swapping the locations of $u$ and $u+1$ in $\rho$, where $(u, u+1)$ is a couple of consecutive white nodes.
\item[(2)] Swapping the locations of $x$ and $y$ in $\rho$, where $x < y$ are white nodes and all $\ell$ nodes appearing between $x$ and $y$ are black, where $\ell \geq 1$. 
\end{enumerate}

To see that these swaps are sufficient, let $w_1, \ldots, w_n$ be the white nodes in increasing order. The swaps described are the adjacent transpositions $(w_1, w_2), (w_2, w_3), \ldots, (w_{n-1}, w_{n})$.

We will show that equation (\ref{eqn0:lemma34}) holds after applying each type of swap.
%By Remark \ref{rem:swappingandparity}, each swap changes $\sign_{BW}(\rho)$. 
%Before proving Lemma \ref{lemma34}, 
This requires a few additional lemmas. Note that the proofs of Lemmas \ref{swapu} through \ref{lem:swapxy2} follow immediately from Definition \ref{def:signpair}.

\begin{lemma}
\label{swapu}
Let $b$ be a black node and let $(u, u+1)$ be a couple of consecutive white nodes. Then $\sign(b, u) = - \sign(b, u+1)$. 
\end{lemma}

\begin{proof}
If $b < u$, then $a_{b, u+1} = a_{b, u} + 1$. So 
\[
\sign(b, u) = (-1)^{(u - b + a_{b, u} - 1)/2}
=-  (-1)^{(u +1- b + a_{b, u+1} - 1)/2} = - \sign(b, u+1) .
\]
%and
%\[
%\sign(b, u+1) = (-1)^{(u +1- b + a_{b, u+1} - 1)/2} =  (-1)^{(u -b + a_{b, u} + 1)/2}.
%\]

If $b > u+1$, then $a_{b, u+1} = a_{b, u} - 1$. So 
\[
\sign(b, u) = (-1)^{(b-u + a_{b, u} - 1)/2}
=- (-1)^{(b-(u+1) + a_{b, u+1} - 1)/2} =- \sign(b, u+1) .
\]
%and 
%\[
%\sign(b, u+1) = (-1)^{(b-(u+1) + a_{b, u+1} - 1)/2} = (-1)^{(b-u + a_{b, u} - 3)/2}. 
%\]
%In both cases, $\sign(b, u) = - \sign(b, u+1)$. 
\end{proof}

\begin{lemma}
\label{lem:swapxy}
Assume the nodes $x$ and $y$ with $x < y$ are white and all $\ell$ nodes between $x$ and $y$ are black, where $\ell \geq 1$. If $b$ is a black node not in the interval $[x + 1, \ldots, y-1]$, then
$\sign(b, x) = (-1)^{\ell} \sign(b, y)$. 
\end{lemma}

\begin{proof}
If $b < x$, then $a_{b, y} = a_{b, x} + \ell - 1$. Then since $y = x + \ell + 1$, 
\begin{eqnarray*}
\sign(b, x) = (-1)^{ (x - b + a_{b, x} - 1)/2 } 
=  (-1)^{ (y - (\ell + 1) - b + a_{b, y} - \ell + 1 - 1)/2 } 
& = &  (-1)^{\ell}  (-1)^{ (y - b + a_{b, y} - 1)/2 }  \\
& = &  (-1)^{\ell}  \sign(b, y).
\end{eqnarray*}
If $b > y$, then $a_{b, y} = a_{b, x} - (\ell - 1)$. Then
\begin{eqnarray*}
\sign(b, x) = (-1)^{ (b-x + a_{b, x} - 1)/2 } 
 =   (-1)^{ (b - (y - (\ell + 1) ) + a_{b, y} + (\ell - 1) - 1)/2 }  
& = &  (-1)^{\ell}  (-1)^{ (b-y + a_{b, y} - 1)/2 }  \\
& = &  (-1)^{\ell}  \sign(b, y).
\end{eqnarray*}
\end{proof}

\begin{lemma}
\label{lem:swapxy2}
Assume the nodes $x$ and $y$ with $x < y$ are white and all $\ell$ nodes between $x$ and $y$ are black, where $\ell \geq 1$. If $b$ is a black node in the interval $[x + 1, \ldots, y-1]$, so $b = x+j$ for some $j \leq \ell$, then
$\sign(b, x) = (-1)^{\ell-1} \sign(b, y)$. 
\end{lemma}

\begin{proof}
Since $b = x+j$ and $a_{b, x} = j-1$, we see that
\[
\sign(b, x) =  (-1)^{ (b-x + a_{b, x} - 1)/2 } 
= (-1)^{ (j + j-1 - 1)/2 }  = (-1)^{j-1}. 
\]
Using the fact that $y - b = \ell + 1 - j$ and $a_{b, y} = \ell - j$, we have
\[
\sign(b, y) = 
(-1)^{ (y - b + a_{b, y} - 1)/2 }  = (-1)^{ (\ell + 1 - j + \ell - j - 1)/2} = (-1)^{\ell - j}. 
\]
So $\sign(b, x) = (-1)^{\ell-1} \sign(b, y)$. 
\end{proof}

%%%%%%%%%%%%%%%%%%%%%%%%%%%%%%%%
%%%%%%%%%%%%%%OLD PROOF%%%%%%%%%%%
%%%%%%%%%%%%%%%%%%%%%%%%%%%%%%%%

\begin{rem}
\label{crossingremark}
The symmetric group $S_{2n}$ acts on the set of pairings on $\{1, \ldots, 2n\}$: the transposition $(i,i+1)$ acts on a pairing $\rho$ by swapping the locations of $i$ and $i+1$.
If $i$ is paired with $i+1$, acting with $(i,i+1)$ leaves the pairing fixed; otherwise, $(i,i+1)$ acts nontrivially and changes the parity of the number of crossings.

Let $\rho$ be a (not necessarily black-white) pairing on $\{1,\ldots,2n\}$. Let $x$ and $y$ be two nodes such that $x <y$. Assume no node in the interval $[x, y]$ is paired with any other node in this interval. Then
$$(x, y) \rho = (x, x+1) \cdots (y-1, y) \cdots (x+1, x+2) (x, x+1) \rho$$
where each transposition of the form $(i, i+1)$ acts nontrivially. 
\end{rem}

\begin{lemma}
\label{crossinglemma3}
Let $\rho$ be a (not necessarily black-white) pairing on $\{1,\ldots,2n\}$. Let $x$ and $y$ be two nodes such that $x <y$ and $x$ is not paired with $y$. Assume that no node in the interval $[x+1, \ldots, y-1]$ is paired with any other node in this interval. Then when the locations of $x$ and $y$ in $\rho$ are swapped, 
\begin{itemize}
\item[(1)] if $x$ and $y$ were both paired with nodes in the interval $[x+1, \ldots, y-1]$,  the number of crossings of $\rho$ changes parity, 
\item[(2)] if exactly one of $x$ and $y$ was paired with a node in the interval $[x+1, \ldots, y-1]$,  then the number of crossings of $\rho$ does not change parity, and 
\item[(3)] if neither $x$ nor $y$ was paired with a node in the interval $[x+1, \ldots, y-1]$ then the number of crossings of $\rho$ changes parity. 
\end{itemize}
\end{lemma}

\begin{proof}
%The symmetric group $S_{2n}$ acts on the set of pairings on $\{1, \ldots, 2n\}$: the transposition $(i,i+1)$ acts on a pairing $\rho$ by swapping the locations of $i$ and $i+1$.
%If $i$ is paired with $i+1$, acting with $(i,i+1)$ leaves the pairing fixed; otherwise, $(i,i+1)$ acts nontrivially and changes the parity of the number of crossings.

Let $\rho$ be a pairing on $\{1,\ldots,2n\}$ and consider $(x, y) \rho$. There are several cases. The strategy is to factor $(x, y)$ into adjacent transpositions and determine which transpositions act nontrivially. \\

\noindent {\bf Case 1.} If the nodes $\rho(x)$ and $\rho(y)$ are both in the interval $[x+1, \ldots, y-1]$, then $(x, y) \rho = (\rho(x), \rho(y)) \rho$. Let $a = \min(\rho(x), \rho(y))$ and let $b = \max(\rho(x), \rho(y))$. Then
$$(\rho(x), \rho(y)) \rho =     (a, a+1) \cdots (b-1, b)\cdots (a+1, a+2) (a, a+1) \rho$$
We have written $(\rho(x), \rho(y))$ as a product of an odd number of transpositions of the form $(i, i+1)$. Since no node in the interval $[a, \ldots, b]$ is paired with any other node in this interval, all these transpositions act nontrivially by Remark \ref{crossingremark}. 
Thus the parity of the number of crossings changes. \\

\noindent {\bf Case 2.} If exactly one of the nodes $\rho(x)$ or $\rho(y)$ is in the interval $[x+1, \ldots, y-1]$, then 
$$(x, y) \rho =     (x, x+1) \cdots (y-1, y)\cdots (x+1, x+2) (x, x+1) \rho$$
and exactly one of these transpositions acts trivially. For if $x$ is paired with $x+k$, then after applying the transposition $(x, x+1)$ to $\rho$, $x+1$ and $x+k$ are paired. Similarly, after applying the transposition $(x+1, x+2)$ to $(x, x+1) \rho$, $x+2$ and $x+k$ are paired. It follows that the transposition $(x+k-1, x+k)$ acts trivially because when we reach this transposition, $x+k-1$ and $x+k$ are paired. Then, the transposition $(x+k-1, x+k-2)$ acts nontrivially and similarly we see that the remaining transpositions act nontrivially.   
Since an even number of transpositions of the form $(i, i+1)$ act nontrivially, the parity of the number of crossings does not change. \ \\

\noindent {\bf Case 3.} If neither of the nodes $\rho(x)$ and $\rho(y)$ are in the interval $[x+1, \ldots, y-1]$, then
$$(x, y) \rho =     (x, x+1) \cdots (y-1, y)\cdots (x+1, x+2) (x, x+1) \rho$$
so we have written $(\rho(x), \rho(y))$ as a product of an odd number of transpositions of the form $(i, i+1)$. Since no node in the interval $[x, y]$ is paired with any other node in this interval, all of these transpositions act nontrivially by Remark \ref{crossingremark}. Thus the parity of the number of crossings changes. \\
\end{proof}

Now that we have established  Lemmas \ref{swapu} through \ref{crossinglemma3}
we can show that equation (\ref{eqn0:lemma34}) holds after applying both types of swaps to $\rho$. 
By Remark \ref{rem:swappingandparity}, each swap changes $\sign_{BW}(\rho)$. \\

\noindent (1) {\bf Swapping the locations of $u$ and $u+1$.} \\

Let $b_1$ be the node paired with $u$ and let $b_2$ be the node paired with $u+1$. 
By Lemma \ref{swapu}, $\sign(b_1, u) = -\sign(b_1, u+1)$ and $\sign(b_2, u+1) = -\sign(b_2, u)$. 
So when we swap the locations of $u$ and $u+1$, $\prod\limits_{(b, w) \in \rho} \text{sign}(b, w)$ does not change. Since $\sign_{BW}(\rho)$ changes, the sign of the LHS of (\ref{eqn0:lemma34}) changes. Swapping $u$ and $u+1$ changes $(-1)^{\# \text{ crosses of } \rho}$,
% as swapping the locations of $w$ and $w+2$ in the case where neither $w$ nor $w+2$ is paired with $w+1$, so $(-1)^{\# \text{ crosses of } \rho}$ changes sign.
so swapping the locations of $u$ and $u+1$ does not affect equation (\ref{eqn0:lemma34}). \\

\noindent (2) {\bf Swapping the locations of $x$ and $y$, where $x < y$ are white nodes and all $\ell$ nodes between $x$ and $y$ are black.} \\

\noindent {\bf Case 1.} If $x$ and $y$ are both paired with black nodes in the interval $[x+1, x+2, \ldots, y-1]$, then $(-1)^{\# \text{ crosses of } \rho}$ changes sign by Lemma \ref{crossinglemma3}. 
By Lemma \ref{lem:swapxy2}, 
$$\sign(\rho(x), x) \sign(\rho(y), y) = ( (-1)^{\ell -1 })^2 \sign(\rho(x), y) \sign(\rho(y), x)$$
so $\prod\limits_{(b, w) \in \rho} \text{sign}(b, w)$ does not change.
Since $\sign_{BW}(\rho)$ changes, the sign of the LHS of (\ref{eqn0:lemma34}) changes.\\

\noindent {\bf Case 2.} If exactly one of $x$ and $y$ is paired with a black node in the interval $[x+1, x+2, \ldots, y-1]$, then 
$(-1)^{\# \text{ crosses of } \rho}$
%the number of crossings of $\rho$ 
does not change sign by Lemma \ref{crossinglemma3}. By Lemmas \ref{lem:swapxy} and \ref{lem:swapxy2}, 
$$\sign(\rho(x), x) \sign(\rho(y), y) = (-1)^{\ell - 1}(-1)^{\ell} \sign(\rho(x), y) \sign(\rho(y), x)$$
so
$\prod\limits_{(b, w) \in \rho} \text{sign}(b, w)$ changes. Since $\sign_{BW}(\rho)$ changes, the sign of the LHS of (\ref{eqn0:lemma34}) does not change sign. \\

\noindent {\bf Case 3.} If neither $x$ nor $y$ is paired with a black node in the interval $[x+1, x+2, \ldots, y-1]$, then $(-1)^{\# \text{ crosses of } \rho}$ changes sign. By Lemma \ref{lem:swapxy},
$$\sign(\rho(x), x) \sign(\rho(y), y) = ( (-1)^{\ell })^2 \sign(\rho(x), y) \sign(\rho(y), x)$$
so $\prod\limits_{(b, w) \in \rho} \text{sign}(b, w)$ does not change.
Since $\sign_{BW}(\rho)$ changes, the sign of the LHS of (\ref{eqn0:lemma34}) changes.\\

This completes the proof of Lemma~\ref{lemma34}. 
We conclude Section~\ref{sec:lem34} by proving Lemma~\ref{lem:OEandBWsigns}, which states that when a black-white pairing $\rho$ is also odd-even, $\sign_{OE}(\rho) = \sign_{BW}(\rho)$.

\subsubsection{Proof of Lemma~\ref{lem:OEandBWsigns}}
\label{sec:OEandBW}

Before we prove Lemma~\ref{lem:OEandBWsigns}, we prove the lemma in the case where $\rho$ is planar. 
\begin{lemma}
\label{lem:OEandBWsignsold}
%Let $G$ be a plane bipartite graph with a set of nodes {\bf N} with any node coloring that has an equal number of black and white nodes. 
When $\rho$ is a planar black-white pairing, 
$$\sign_{OE}(\rho) = \sign_{BW}(\rho)$$
\end{lemma}

%\begin{note*}
%Checked Lemma \ref{lem:OEandBWsigns} on all balanced node colorings on 6, 8, 10, 12 nodes
%\end{note*}

\begin{proof}
Let $\rho$ be a planar black-white pairing. 
Recall from Definition 
\ref{def:invofrho}
that all black-white pairings can be written
$\rho = ((b_1, \rho(b_1)), (b_2, \rho(b_2)), \ldots, (b_n, \rho(b_n)))$, where $b_1 < b_2 < \cdots < b_n$,
and we say that $(\rho(b_i), \rho(b_j))$ is an inversion of $\rho$ if $i < j$ and $\rho(b_i) > \rho(b_j)$.

All planar pairings are odd-even, and recall
from Definition~\ref{defn:invofpi} that
if\\
$\rho = ((1, \rho(1)), (3, \rho(3)), \ldots, (2n-1, \rho(2n-1)))$
is an odd-even pairing,
% that an inversion of an odd-even pairing is an inversion of the permutation
%$\begin{pmatrix}
% \frac{ \rho(1)}{2} & \frac{ \rho(3)}{2}  & \cdots & \frac{  \rho(2n-3) }{2} & \frac{ \rho(2n-1)}{2}
%\end{pmatrix}$. In this proof, 
we say $(\rho(i), \rho(j))$ is an inversion if $i < j$ and $\rho(i) > \rho(j)$.
%, even though the inversion is actually $(\frac{\rho(i)}{2}, \frac{\rho(j)}{2})$. 

We will show that there is a one-to-one correspondence between inversions of $\rho$ when it is considered as a black-white pairing (which we will call black-white inversions) and inversions of $\rho$ when it is considered as an odd-even pairing (which we will call odd-even inversions).

Consider a black-white inversion, that is, some $b_i < b_j$ such that $\rho(b_i) > \rho(b_j)$. There are several cases to consider: \\

\noindent {\bf Case 1.} $b_i, b_j$ are both odd. \\
In this case, $b_i = 2k-1$ and $b_j = 2 \ell -1$ for some $k < \ell$, so $(\rho(b_i), \rho(b_j))$ is an odd-even inversion. \\

\noindent {\bf Case 2.} $b_i, b_j$ are both even. \\
Since $b_i < b_j$ and $\rho(b_i) > \rho(b_j)$, $(b_j, b_i)$ is an odd-even inversion. \\

\noindent {\bf Case 3.} $b_i$ is odd and $b_j$ is even.\\
There are two subcases to consider. If $\rho(b_j) < b_i$, then it must be the case that $b_j > \rho(b_i)$. To see this, observe that if $b_j  < \rho(b_i)$, then 
$\rho(b_j) < b_i < b_j < \rho(b_i)$, but then we have a crossing, which contradicts the planarity of $\rho$. So $(b_j, \rho(b_i))$ is an odd-even inversion. 

If $\rho(b_j) > b_i$, then $\rho(b_i) > b_j$ (otherwise $b_i < \rho(b_j) < \rho(b_i) < b_j$, so $\rho$ has a crossing). So $(\rho(b_i), b_j)$ is an odd-even inversion. \\

\noindent {\bf Case 4.} $b_i$ is even and $b_j$ is odd.\\
If $\rho(b_i) > b_j$, then $\rho(b_j) > b_i$ (otherwise $\rho(b_j) < b_i < b_j < \rho(b_i)$ is a crossing), so $(\rho(b_j), b_i)$ is an odd-even inversion.
If $\rho(b_i) < b_j$, then $b_i > \rho(b_j)$ (otherwise $b_i < \rho(b_j) < \rho(b_i) < b_j$ is a crossing), so $(b_i, \rho(b_j))$ is an odd-even inversion. \\

A similar argument shows that for each odd-even inversion, there is a black-white inversion. Since there is a one-to-one correspondence between odd-even  inversions and black-white inversions, $\sign_{OE}(\rho) = \sign_{BW}(\rho)$. 
\end{proof}

\noindent {\bf Lemma~\ref{lem:OEandBWsigns}.}
%Let $G$ be a plane bipartite graph with a set of nodes {\bf N} with any node coloring that has an equal number of black and white nodes. 
When $\rho$ is a black-white pairing that is also odd-even, 
$$\sign_{OE}(\rho) = \sign_{BW}(\rho).$$

\begin{proof}
One can get from an odd-even black-white pairing $\rho_1$ to any other odd-even black white pairing $\rho_2$ by applying a series of moves, where each move swaps the locations of two nodes of the same color and parity. Since each of these moves changes $\sign_{OE}$ and $\sign_{BW}$, the claim follows from Lemma \ref{lem:OEandBWsignsold}. 
\end{proof}

\subsection{Lemmas 3.1 and 3.2 from Kenyon and Wilson}

\label{sec:lem31}
Throughout this section, $S$ denotes a balanced subset of nodes (a subset containing an equal number of black and white nodes). 
 In \cite{KW2006}, %\cite[Lemma 3.1]{KW2006}, 
Kenyon and Wilson show that  $Z^{D}(G \setminus S) Z^{D}(G \setminus S^{c})$ is a sum of double-dimer partition functions $Z^{DD}_{\pi}(G, {\bf N})$, where the sum is over all pairings $\pi$ that do not connect nodes in $S$ to nodes in $S^c$. 

\begin{lemma}\cite[Lemma 3.1]{KW2006}
\label{lem3.1}
If $S$ is a balanced subset of nodes 
%(a subset containing an equal number of black and white nodes) 
%and $Z^{D}(S) = Z^{D}( G \setminus S)$ denotes the weighted sum of dimer covers of $G \setminus S$, 
then 
$Z^{D}(G \setminus S) Z^{D}(G \setminus S^{c})$ is a sum of double-dimer configurations for all connection topologies $\pi$ for which $\pi$ connects no element of $S$ to an element of $S^{c}:={\bf N} \setminus S$. That is,
$$Z^{D}(G \setminus S) Z^{D}(G \setminus S^{c}) = Z^{DD} \sum\limits_{\pi} M_{S, \pi} \Pr (\pi),$$
where $M_{S, \pi}$ is 0 or 1 according to whether $\pi$ connects nodes in $S$ to $S^{c}$ or not.
\end{lemma}

This lemma relates the quantity $Z^{D}(G \setminus S) Z^{D}(G \setminus S^{c})$ to $\Pr(\pi)$. Next, Kenyon and Wilson show that $\dfrac{Z^{D}(G \setminus S) Z^{D}(G \setminus S^{c})}{(Z^{D}(G))^2}$ is a determinant in the quantities $X_{i, j}$. 

\begin{lemma}\cite[Lemma 3.2]{KW2006}
\label{lem:kwlem32}
Let $S$ be a balanced subset of $\{1, \ldots, 2n\}$. Then
$$\dfrac{Z^{D}(G \setminus S) Z^{D}(G \setminus S^{c})}{(Z^{D}(G))^{2}} = 
\det[(1_{i, j \in S} + 1_{i, j \notin S}) \times (-1)^{(|i - j| -1)/2} X_{i, j} ]^{i = 1, 3, \ldots, 2n-1}_{j = 2, 4, \ldots, 2n}.$$
\end{lemma}

The combination of these results shows that $\widehat{\Pr}(\pi)$ is a homogeneous polynomial in the $X_{i, j}$, since the matrix $M$ from \cite[Lemma 3.1]{KW2006} has full rank \cite[Lemma 3.3]{KW2006}. 
%Proving that the coefficients of the polynomial are integers requires more work. 
Our analogues of these lemmas have several differences (such as the additional global signs in our version of Lemma~\ref{lem:kwlem32}, see Lemma~\ref{lemma32gen}), but our proofs are quite similar to their proofs. 
 %concern the quantity $Z^{D}(G \setminus S) Z^{D}(G \setminus S^{c})$. 

 %concern the quantity $Z^{D}(G \setminus S) Z^{D}(G \setminus S^{c})$. 

We begin with Lemma~\ref{lem3.1}.
For a graph $G$ with node set ${\bf N}$ that does not necessarily have the property that all nodes are black and odd or white and even,
a statement very similar to Lemma \ref{lem3.1} holds. 
%Let $G$ be a graph with node set ${\bf N}$ that does not necessarily have the property that all odd nodes of $G$ are black and all even nodes are white. 
For the remainder of this section, we let $T \subseteq {\bf N}$ be the set of nodes that are odd and white or even and black. Since ${\bf N}$ is assumed to have an equal number of black and white nodes, $|T|$ is even. 

Let $\widetilde{G}$ be $G$ with an extra vertex and edge with weight 1 added to each node in $T$, so all of the nodes in $\widetilde{G}$ are black and odd or white and even. We note that
$Z^{D}(\widetilde{G} \setminus S) = Z^{D}(G \setminus (S \triangle T) ),$ where $S \triangle T$ denotes the symmetric difference of the sets $S$ and $T$. For example, if ${\bf N}$ is a set of 12 nodes colored so that nodes 1, 3, 4, 5, 7, and 10 are black (see the proof of Lemma~\ref{firstlemma34}) then $T = \{4, 9, 10, 11\}$. If $S =\{2, 3, 9, 10\}$, then $S \triangle T = \{2, 3, 4, 11\}$. 
Lemma \ref{lem3.1} implies the following.

\begin{cor}
Let $S$ be a balanced subset of nodes. 
$Z^{D}(G \setminus (S \triangle T)) Z^{D}(G \setminus (S \triangle T)^{c})$ is a sum of double-dimer configurations for all connection topologies $\pi$ for which $\pi$ connects no element of $S$ to an element of $S^{c}$. That is,
$$Z^{D}(G \setminus (S \triangle T)) Z^{D}(G \setminus  (S \triangle T)^{c}) = Z^{DD}(G) \sum\limits_{\pi} M_{S, \pi} \Pr (\pi),$$
where $M_{S, \pi}$ is 0 or 1 according to whether $\pi$ connects nodes in $S$ to $S^{c}$ or not.
\end{cor}

If $V = S \triangle T$, then $S = V\triangle T$, so we have:

\begin{cor}
\label{corlem31}
Let $V$ be a balanced subset of nodes.
$Z^{D}(G \setminus V) Z^{D}(G \setminus V^{c})$ is a sum of all connection topologies $\pi$ for which $\pi$ connects no elements of $V \triangle T$ to $(V \triangle T)^c$. That is,
$$Z^{D}(G \setminus V) Z^{D}(G \setminus V^{c}) = Z^{DD}(G) \sum\limits_{\pi} M_{V \triangle T, \pi} \Pr (\pi),$$
where $M_{V \triangle T, \pi}$ is 0 or 1 depending on whether $\pi$ connects nodes in $V \triangle T$ to $(V \triangle T)^{c}$.
\end{cor}

Corollary \ref{corlem31} is the version of Lemma \ref{lem3.1} that we will need to prove Theorem \ref{thm:thm1}.

Our version of  \cite[Lemma 3.2]{KW2006} is the following. 
\begin{lemma}
\label{lemma32gen}
Let $S$ be a balanced subset of ${\bf N} = \{1, \ldots, 2n\}$. Then 
\small
\begin{equation}
\label{lem32}
\dfrac{Z^{D}(G \setminus S) Z^{D}(G \setminus S^{c})}{(Z^{D}(G))^{2}} = \text{sign}_{\cons}({\bf N}) \text{sign}(S)
\det \left[(1_{i, j \in S} + 1_{i, j \notin S}) \times \text{sign}(i, j) Y_{i, j} \right]^{i = b_1, b_2, \ldots, b_n}_{j = w_1, w_2, \ldots, w_n}
\end{equation}
\normalsize
where $b_1, b_2, \ldots, b_n$ are the black nodes of $\{1, 2, \ldots, 2n\}$ listed in ascending order, $w_1, w_2, \ldots, w_n$ are the white nodes of $\{1, 2, \ldots, 2n\}$ listed in ascending order, $\text{sign}(i, j)$ is defined in Definition \ref{def:signpair}, $Y_{i, j} =\dfrac{Z^D(G_{i, j})}{Z^D(G)}$, and 
\[\text{sign}(S) = (-1)^{\# \text{ crosses of $\rho$} }, \]
where $\rho$ is a black-white pairing that does not connect\footnotemark~$S$ to $S^c$
%and $\rho |_{S}$ and $\rho |_{S^c}$ are planar.
and is planar when restricted to $S$ and planar when restricted to $S^c$. 
\end{lemma}
\footnotetext{The statement ``$\rho$ does not connect $S$ to $S^c$'' is an abbreviation for ``$\rho$ does not connect nodes in $S$ to nodes in $S^c$''.}

%\begin{rem}
%The fact that such a pairing $\rho$ always exists is a consequence of %Lemma \ref{firstlemma34}: 
%since in this lemma we constructed a black-white planar pairing $\rho$ given any node coloring. 
%since $S$ is balanced, Lemma \ref{firstlemma34} proves the existence of %a black-white planar pairing of $S$ and a black-white planar pairing of %$S^c$. 
%\end{rem}

\begin{rem}
The fact that such a pairing $\rho$ always exists is a consequence of Lemma \ref{firstlemma34}, which states that for any node coloring there is a planar black-white pairing $\rho$ satisfying \\ $\sign_{BW}(\rho) \prod\limits_{(b, w) \in \rho} \sign(b, w) = \sign_{c}({\bf N})$. 
%since in this lemma we constructed a black-white planar pairing $\rho$ given any node coloring. 
Since $S$ is balanced, the existence of a planar black-white pairing of $S$ and a planar black-white pairing of $S^c$ follows.  
\end{rem}

The proof of Lemma \ref{lemma32gen} requires some Kasteleyn theory. The reader familiar with basic facts about Kasteleyn matrices can skip the following section.

\subsubsection{Kasteleyn matrices}

Recall that  $G = (V_1, V_2, E)$ is a finite edge-weighted bipartite planar graph embedded in the plane. Let $\omega((i, j))$ denote the weight of an edge $(i, j) \in E$. 
%For the following definitions and lemmas, let $G = (V_1, V_2, E, \omega)$ be an edge-weighted bipartite planar graph with a fixed embedding in the plane. 

\begin{defn}
A Kasteleyn (or flat) weighting of $G$ is a choice of sign for each edge with the property that each face with 0 mod 4 edges has an odd number of $-$ signs and each face with 2 mod 4 edges has an even number of $-$ signs. 
\end{defn}

For the remainder of this section we will let $\sigma: E \to \pm 1$ denote the Kasteleyn weighting of $G$. 
A Kasteleyn matrix of $G$ is a weighted, signed bipartite adjacency matrix of $G$. More precisely, 
%given a Kasteleyn weighting, $\sigma: E \to \pm 1$, 
define a $|V_1| \times |V_2|$ matrix $K$ by 
$$
K_{i, j} = 
\begin{cases}
\sigma( (i, j) ) \omega( (i, j) ) & \mbox{ if } (i, j) \in E \\
0 & \mbox{ otherwise}
\end{cases}
$$
Kasteleyn showed that every bipartite planar graph with an even number of vertices has a 
Kasteleyn matrix \cite{Kas67}. Furthermore, if $|V_1| = |V_2|$ then $|\det K |$ is the weighted sum of all dimer configurations of $G$.

The proof of Lemma~\ref{lemma32gen} uses a few straightforward facts about Kasteleyn weightings. First, we will show that if $G = (V_1, V_2, E)$ has a Kasteleyn weighting $\sigma$, and we add edges to $G$ to obtain $G'$, we can choose weights for the added edges to obtain a Kasteleyn weighting $\sigma'$ of $G'$ with the property that $\sigma'(e) = \sigma(e)$ for all $e \in E$. 
%under certain conditions we can add edges to $G$ and choose their weights to obtain a Kasteleyn weighting 
\begin{lemma}
\label{kasteleynlemma1}
Let $b$ and $w$ be two vertices of opposite color on a face $F$ of $G = (V_1, V_2, E)$. 
Let $E' = E \cup \{ \tilde{e} \}$, where $\tilde{e} \notin E$ is an edge connecting $b$ and $w$ that separates $F$% the outer face of $G$
 into two faces and let $G' = (V_1, V_2, E')$. 
%Let $G = (V, E)$ be a planar bipartite graph. Fix a planar embedding of $G$ and let $\sigma: E \to \pm 1$ be a Kasteleyn weighting of $G$. Let $\tilde{e} \notin E$ be an edge connecting two vertices  
Define $\sigma': E \cup \{\tilde{e}\} \to \pm 1$ so that $\sigma'(e) = \sigma(e)$ for all $e \in E$ and choose $\sigma'(\tilde{e})$ so that one of the faces bounded by $\tilde{e}$ is flat (i.e., it has an odd number of $-$ signs if it has 0 mod 4 edges, and an even number of $-$ signs otherwise). Then $\sigma'$ is a Kasteleyn weighting of $G'$.
\end{lemma}

\begin{proof}
%Let $G = (V, E)$ be a planar bipartite graph. Fix a planar embedding of $G$ and let $\sigma: E \to \pm 1$ be a Kasteleyn weighting of $G$. Let $b$ and $w$ be two vertices of opposite color on the outer face of $G$. 
%Let $G = (V, E)$ be a planar bipartite graph. Fix a planar embedding of $G$ and let $\sigma: E \to \pm 1$ be a Kasteleyn weighting of $G$. Let $b$ and $w$ be two vertices of opposite color on the outer face of $G$. 
%There are two paths from $b$ to $w$ that use edges of $F$: the path $\mathcal{P}$ from $b$ to $w$ that goes clockwise, and the path $\mathcal{Q}$ that goes counterclockwise. 
By assumption, 
%We add 
the edge $\widetilde{e}$ 
separates $F$
%so that it separates a face of $G$
 into two faces: 
the face consisting of the edges of a path $\mathcal{Q}$ and the edge $\widetilde{e}$, and the face 
consisting of the edges of a path $\mathcal{P}$ and the edge $\widetilde{e}$.
% Call the new graph $G'$. 
 The path $\mathcal{P}$ consists of $1$ mod $4$ edges or $3$ mod $4$ edges. 
%There are two paths from $b$ to $w$ that use edges of the outer face of $G$: the path $\mathcal{P}$ from $b$ to $w$ that goes clockwise, and the path $\mathcal{Q}$ that goes counterclockwise. 
%We add the edge $\tilde{e}$ so that it separates the outer face of $G$ into two faces: 
%the face consisting of the edges of $\mathcal{Q}$ and the edge $\tilde{e}$, and the face 
%consisting of the edges of $\mathcal{P}$ and the edge $\tilde{e}$. Call the new graph $G'$. The path $\mathcal{P}$ consists of $1$ mod $4$ edges or $3$ mod $4$ edges. 
Define
\[
\sigma'(\tilde{e}) = 
\begin{cases}
\prod\limits_{e \in \mathcal{P}} \sigma(e) & \text{ if $\mathcal{P}$ has $1$ mod 4 edges } \\
-\prod\limits_{e \in \mathcal{P}} \sigma(e) & \text{ if $\mathcal{P}$ has $3$ mod 4 edges }
\end{cases}
\]
and define $\sigma'(e) = \sigma(e)$ for all $e \in E$. 
Now the face consisting of the path $\mathcal{P}$ and the edge $e$ is flat. It remains to check that 
the face $F'$ consisting of $\mathcal{Q}$ and $\widetilde{e}$
%the outer face of $G'$ 
is flat, 
which is done by breaking into cases based on whether the paths $\mathcal{P}$, $\mathcal{Q}$ have $1$ or $3$ edges mod 4.  
\end{proof}

\begin{lemma}
\label{kasteleynlemma2}
Let $W = \{v_1, \ldots, v_{2m}\}$ be a set of vertices on the outer face of $G = (V_1, V_2, E)$. Pair the vertices of $W$ so that we can add edges $e_1, \ldots, e_m$ connecting the pairs without introducing any edge crossings.
Let $E^{(m)} = E \cup  \{e_1, \ldots, e_m \}$ and let $G^{(m)} = (V_1, V_2, E^{(m)}, \omega)$.  
Define $\sigma_i: E \cup \{e_i\} \to \pm 1$ as in Lemma \ref{kasteleynlemma1}: $\sigma_i(e) = \sigma(e)$ for all $e \in E$ and $\sigma_i(e_i)$ is chosen so that one of the faces bounded by $e_i$ is flat. By Lemma \ref{kasteleynlemma1}, $\sigma_i$ is a Kasteleyn weighting for all $1 \leq i \leq m$. Then $\tau: E \cup \{e_1, \ldots, e_m \} \to \pm 1$ defined by $\tau(e) = \sigma(e)$ for all $e \in E$ and $\tau(e_i) = \sigma_i(e_i)$ for $1 \leq i \leq m$ is a Kasteleyn weighting of $G^{(m)}$.
\end{lemma}

\begin{proof}
We prove the claim by induction on $m$. When $m =1$, there is nothing to show. Assume the claim holds when we add $m-1$ edges to $G$. Now suppose we add $m$ edges $e_1, \ldots, e_m$. Choose an ``innermost" edge $e_j$, i.e. an edge with the property that one of its faces is bounded only by edges of $G$ and $e_j$. 
By the induction hypothesis, $\tau: E \cup \{e_1, \ldots,e_{j-1}, e_{j+1}, \ldots, e_{m} \} \to \pm 1$ defined by $\tau(e) = \sigma(e)$ for all $e \in E$ and $\tau(e_i) = \sigma_i(e_i)$ for $i = 1, 2, \ldots, j-1, j+1, \ldots, m$ is a Kasteleyn weighting of $G^{(m-1)}= (V_1, V_2, E \cup  \{ e_1, \ldots,e_{j-1}, e_{j+1}, \ldots e_{m} \} )$. 
Since $e_j$ is an innermost edge and $\sigma_j: E \cup \{e_j \} \to \pm 1$ was defined so that when $e_j$ is added to $G$, one of the faces bounded by $e_j$ is flat, we may apply Lemma \ref{kasteleynlemma1} to conclude that $\tau: E \cup \{e_1, \ldots, e_m \} \to \pm 1$ defined by $\tau(e) = \sigma(e)$ for all $e \in E$ and $\tau(e_i) = \sigma_i(e_i)$ for $1 \leq i \leq m$ is a Kasteleyn weighting of $G^{(m)} = (V_1, V_2, E \cup  \{e_1, \ldots, e_m \} )$. 
\end{proof}

We also need to show that if we delete rows and columns from a Kasteleyn matrix of a graph, the resulting submatrix is a Kasteleyn matrix of the corresponding graph.

\begin{lemma}
\label{kasteleynlemma3}
%Let $G = (V, E)$ be a planar bipartite graph. Fix a planar embedding of $G$ and 
Let $K$ be a Kasteleyn matrix of $G$. 
Let $S$ be a balanced subset of vertices on the outer face of $G$. 
Then $K_{\setminus S}$, the submatrix of $K$ formed by deleting the rows and columns from $S$,
is a Kasteleyn matrix of $G \setminus S$. 
\end{lemma}

To prove this, we need the following lemma and corollary, which are proven in \cite{Kuperberg}. 

\begin{lemma}\cite[Theorem 2.1]{Kuperberg}
If $G$ is a planar bipartite graph with an even number of vertices, there are an even number of faces with $4k$ sides. 
\end{lemma}

\begin{cor}\cite[Theorem 2.2]{Kuperberg}
\label{Kcor}
Every signed graph with an even number of vertices has an even number of non-flat faces. 
\end{cor}

\begin{proof}[Proof of Lemma \ref{kasteleynlemma3}]
$G \setminus S$ is flat at each internal face because $G$ is flat at each internal face, so it remains to show that it is flat on the outer face as well. 
Since $G \setminus S$ has an even number of vertices, it has an even number of non-flat faces by Corollary~\ref{Kcor}, so it must be flat on the outer face. 
\end{proof}

\subsubsection{Proof of Lemma~\ref{lemma32gen}}

\begin{proof}[Proof of Lemma~\ref{lemma32gen}]

Assume there are $2k$ couples of consecutive nodes of the same color. As in Remark~\ref{rem:notation}, we label the couples of consecutive white nodes $(u_i, u_{i} +1)$ and 
the couples of consecutive black nodes $(s_i, s_{i} +1)$ for $1 \leq i \leq k$. %Note that it is possible that $u_{i+1} = u_{i} + 1$ (for example, if there are three consecutive nodes of the same color). 

Following the proof of \cite[Lemma 3.2]{KW2006}, 
%Lemma 3.2 in Kenyon and Wilson, 
we adjoin to the graph $G$ $2n-2k$ edges connecting
all adjacent nodes except nodes $s_i$ and $s_i+1$ and nodes $u_i$ and $u_i+1$. 
The resulting graph is still bipartite by the assumption that the nodes alternate between black and white except for the nodes $s_i$ and $s_i+1$ 
and the nodes $u_i$ and $u_i+1$. 
Now add $4k$ more edges as follows. Since $G$ is bipartite, there is a white vertex $t_i$ on the outer face of $G$ between nodes $s_i$ and $s_i+1$ and a black vertex $v_i$ on the outer face of $G$ between nodes $u_i$ and $u_i+1$. 
Add edges connecting nodes $s_i$ and $t_i$ and $t_i$ and $s_i+1$, and edges connecting nodes $u_i$ and $v_i$ and $v_i$ and $u_i+1$. 
Give the $2n-2k + 4k = 2n + 2k$ edges we have added weight $\epsilon$ (and then take the limit $\epsilon \to 0$). Let $G'$ denote the resulting graph. 

Given a Kasteleyn weighting of a graph, the signs of edges incident to a vertex may be reversed, and each face will still have a correct number of minus signs. 
Fix a Kasteleyn weighting of the graph $G'$. List the vertices from the set ${\bf{N}} \cup \{ t_{i} \}_{i=1}^{k} \cup \{ v_{i} \}_{i=1}^{k}$ in counterclockwise order. 
 For each vertex in this list, 
 %of vertices $1, 2, \ldots, s, t, s+1, \ldots, u, v, u+1, \ldots, 2n$, 
% (or, if $u< s$, $1, 2, \ldots, u, v, u+1, \ldots, s, t, s+1, \ldots, 2n$)
 % (which are all nodes except $t$ and $v$), 
 if the edge from the vertex $i$ to the next vertex in the list $j$
has a minus sign, reverse the signs of all edges incident to vertex $j$. This ensures that the edges of weight $\epsilon$ we added to $G$ have positive sign, with the possible exception of the edge from node $2n$ to $1$, which must have sign $-(-1)^{n+k}$ for the outer face to have a correct number of minus signs (because if $n+k$ is even, the outer face has 0 mod 4 edges, and if $n+k$ is odd, the outer face has 2 mod 4 edges).

Let $S$ be a balanced subset of $\{1, \ldots, 2n\}$. 
Let $(w_1, b_1), \ldots, (w_j, b_j)$ be any noncrossing pairing of the nodes of $S$, where $w_1, \ldots, w_j$ are the white nodes of $S$ and $b_1, \ldots, b_j$ are the black nodes of $S$. Adjoin edges of weight $W$ connecting $w_i$ to $b_i$ for $1 \leq i \leq j$. Because of the edges of weight $\epsilon$ we adjoined to $G$, we let the sign of a new edge of weight $W$ connecting black node $b$ and white node $w$ be
\[
\text{sign}(b, w) = 
(-1)^{(|b-w|+ a_{b, w}-1)/2}, 
\]
where recall that $a_{b, w}$ is the number of couples of consecutive nodes of the same color in the interval $[\min\{b, w\}, \min\{b, w\} + 1, \ldots, \max\{b, w\}]$. 

Observe that with this choice of signs, when we add any one of the edges $e_i = (b_i, w_i)$ to $G'$ so that it separates the outer face of $G'$ into two faces, one of the faces bounded by $e_i$ is flat. So by Lemma \ref{kasteleynlemma2}, this is a Kasteleyn weighting. 

Let $K_W$ be the Kasteleyn matrix of the resulting graph, with rows and columns ordered so that $b_1, \ldots, b_j$ are the first $j$ rows and $w_1, \ldots, w_j$ are the first $j$ columns. 
Let $K = K_0$ be the corresponding Kasteleyn matrix when $W = 0$. 
Then $Z^{D}(G \setminus S) = \pm [W^j] \det(K_W)$
where $[W^j] \det(K_W)$ denotes the coefficient of $W^j$ in the polynomial $\det(K_W)$. 
(Because $[W^j] \det(K_W)$ is, up to a sign, the weighted sum of matchings that include all of the edges of weight $W,$ which is exactly
the weighted sum of matchings of $G \setminus S$.)
%Similarly, $Z^{D} = \det K_0$. 
Since each term of $\det(K_W)$ has the same sign, 
\[
\dfrac{ Z^{D}(G \setminus S)}{Z^{D}(G)} = \dfrac{ [W^j] \det(K_W)}{[W^0] \det(K_W)}.
\]
Next let $K_{\setminus S}$ denote the submatrix of $K$ formed by deleting the rows and columns from $S$. 
By Lemma \ref{kasteleynlemma3}, $K_{\setminus S}$ is a Kasteleyn matrix of $G \setminus S$.
The sign of $\det(K_{\setminus S})$ and the sign of 
$[W^j] \det(K_W)$ differ by the product of the signs of the edges of weight $W$. So, noting that $[W^0] \det(K_W)= \det(K)$, we have
\[
 \dfrac{ [W^j] \det(K_W)}{[W^0] \det(K_W)}
 = 
 \prod\limits_{\ell=1}^{j} \text{sign}(b_{\ell}, w_{\ell}) 
 \dfrac{ \det (K_{\setminus S} )}{ \det(K)}.
 \]
By Jacobi's determinant identity, 
\[
 \prod\limits_{\ell=1}^{j} \text{sign}(b_{\ell}, w_{\ell}) 
 \dfrac{ \det (K_{\setminus S}) }{ \det(K)}
 = 
  \prod\limits_{\ell=1}^{j} \text{sign}(b_{\ell}, w_{\ell}) 
  \det[ K_{b, w}^{-1}]_{w = w_1, \ldots, w_j}^{b = b_1, \ldots, b_j}.
 \]
 So we have
 \begin{equation}
 \label{proofeqn}
 \dfrac{ Z^{D}(G \setminus S)}{Z^{D}(G)} =  \prod\limits_{\ell=1}^{j} \text{sign}(b_{\ell}, w_{\ell}) 
  \det[ K_{b, w}^{-1}]_{w = w_1, \ldots, w_j}^{b = b_1, \ldots, b_j}.
 \end{equation}
Letting $S = \{b, w\}$ in equation (\ref{proofeqn}), we get
\[
Y_{b, w} = \dfrac{Z^{D}(G_{b, w})}{Z^{D}(G)}
 = \text{sign}(b, w) K_{b, w}^{-1}.
 \]
 From this and equation (\ref{proofeqn}) we find that when $\rho_1 = (w_1, b_1), \ldots, (w_j, b_j)$ is a noncrossing pairing of the nodes of $S$ and $\rho_2 = (w_{j+1}, b_{j+1}), \ldots, (w_n, b_n)$ is a noncrossing pairing of the nodes of $S^c$,
 \small
\begin{eqnarray*}
&& \dfrac{Z^{D}(G \setminus S) Z^{D}(G \setminus S^{c})}{(Z^{D}(G))^{2}}\\
 & =&  \prod\limits_{\ell=1}^{j} \text{sign}(b_{\ell}, w_{\ell}) \prod\limits_{\ell=j+1}^{n} \text{sign}(b_{\ell}, w_{\ell}) 
  \det \left[ \text{sign}(b, w) Y_{b, w}  \right]_{w = w_1, \ldots, w_j}^{b = b_1, \ldots, b_j}
 \det \left[ \text{sign}(b, w) Y_{b, w}  \right]_{w = w_{j+1}, \ldots, w_n}^{b = b_{j+1}, \ldots, b_n} \\
&=&  \prod\limits_{\ell=1}^{n} \text{sign}(b_{\ell}, w_{\ell}) 
%\prod\limits_{\ell=j+1}^{n} \text{sign}(b_{\ell}, w_{\ell}) 
  \det \left[ (1_{i, j \in S} + 1_{i, j \notin S})\text{sign}(b, w) Y_{b, w}  \right]_{w = w_1, \ldots, w_n}^{b = b_1, \ldots, b_n}
\end{eqnarray*}
 \normalsize
 which is equation (\ref{lem32}), except for the global sign and the order of the rows and columns (since $w_1, \ldots, w_n$ and $b_1, \ldots, b_n$
 are not necessarily in ascending order). 
 
 Let $\rho = \rho_1 \cup \rho_2$. 
 Reorder the rows so that the black nodes are in ascending order. For each row swap, make the corresponding column swap. Then $\rho$ pairs the node corresponding to row $i$ with the node corresponding to column $i$. Since the row swaps and column swaps we have made are in one-to-one correspondence, we have not changed the sign of the determinant. Finally, we need to put the columns in ascending order. The number of swaps required to do this is exactly $\sign_{BW}(\rho)$.
 
 So after reordering the rows and columns so that they are listed in ascending order, the global sign is:
 $$ \prod\limits_{\ell=1}^{n} \text{sign}(b_{\ell}, w_{\ell})  \sign_{BW}(\rho)$$
 which is equal to
 $\sign_{\cons}({\bf N}) (-1)^{\# \text{ crosses of $\rho$} } $
by Lemma \ref{lemma34}.

\end{proof}

%\subsection{Lemma 3.5 and Theorem 3.6 from Kenyon and Wilson}

\subsection{Defining $\mathcal{Q}^{(DD)}$.}
Let $Y'$ be the vector of monomials $Y'_{\rho}$ indexed by black-white pairings, where $Y'_{\rho} = (-1)^{\# \text{ crosses of $\rho$}} \prod\limits_{(i, j) \in \rho } Y_{i, j}$. 

In this section, we define $Q^{(DD)}$, which is the matrix satisfying
$P = Q^{(DD)}Y'$, where
$P$ is the vector indexed by planar pairings $\pi$ with entries
$\widetilde{\Pr}(\pi)$. Recall from Section~\ref{sec:organization} that $\widetilde{\Pr}(\pi) = \dfrac{ Z^{DD}_{\pi}(G, {\bf N}) }{ (Z^{D}(G))^{2} }$. 
%= \dfrac{\text{Pr}(\pi) Z^{DD}(G)}{(Z^D(G))^2 } = 
%\dfrac{\text{Pr}(\pi) Z^{DD}(G)}{(Z^D(G))^2 } $$
%$P_\pi = \text{Pr}(\pi) Z^{DD}/(Z^D)^2 = \text{Pr}(\pi) Z^{D}(\{1, 2, \ldots, 2n\})/Z^D$
% indexed by planar pairings $\pi$.
%In other words, the
% $\pi$th row of $Q^{(DD)}$ gives the polynomial $\widetilde{\Pr}(\pi)$. 

We begin with a few definitions. 

\begin{defn}
If $\sigma$ and $\tau$ are two pairings on a set of nodes $\{1, 2, \ldots, 2n\}$, construct the undirected multigraph $C$ with vertex set $\{1, 2, \ldots, 2n\}$ by adding an edge between vertices $i$ and $j$ for each pair $(i, j)$ of $\sigma$, and similarly for $\tau$. The {\em number of components in $\sigma \cup \tau$} is the number of connected components in $C$. Note that all connected components of $C$ are cycles. 
\end{defn}

\begin{example} If $\sigma = ((1, 2), (3, 4), (5, 6))$ and $\tau = ((1, 5), (2, 6), (3, 4))$ then there are two components in $\sigma \cup \tau$, as shown below.
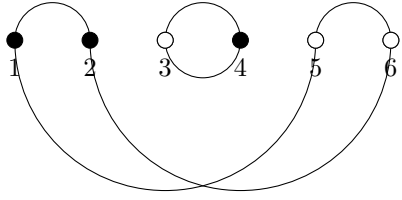
\[\begin{tikzpicture}[scale = .5]
	\vertex[fill] (n1) at (0, 0) [label=below:\small{1}] {};
	\vertex[fill] (n2) at (2,0) [label=below:\small{2}] {};
	\vertex (n3) at (4,0) [label=below:\small{3}]  {};
	\vertex[fill] (n4) at (6, 0) [label=below:\small{4}] {};
	\vertex (n5) at (8,0) [label=below:\small{5}] {};
	\vertex (n6) at (10,0) [label=below:\small{6}]  {};
	
	%arcs of pi
	\draw  (n2) arc (0:180:1cm);
	\draw  (n4) arc (0:180:1cm);
	\draw  (n6) arc (0:180:1cm);
	
	%arcs of rho
	\draw  (n1) arc (180:360:4cm);
	\draw  (n2) arc (180:360:4cm);
	\draw  (n3) arc (180:360:1cm);
	
		\vertex[fill=white] (n3) at (4,0) {};
	\vertex[fill=white] (n5) at (8,0) {};
	\vertex[fill=white] (n6) at (10,0)  {};

\end{tikzpicture}\]
\end{example}

\begin{defn}
\label{signrhopi}
If $\pi$ is an odd-even pairing and $\rho$ is a black-white pairing, define
\[ \text{sign}(\pi, \rho) = (-1)^{\# \text{nodes}/2} (-1)^{\# \text{ components in } \pi \cup \rho} \text{sign}_{OE}(\pi) \text{sign}_{BW}(\rho). \]
\end{defn}

\begin{defn}
\label{Bdefn}
Define the matrix $\mathcal{B}_{2}$ which has rows indexed by planar pairings and columns indexed by black-white pairings by
\[(\mathcal{B}_2)_{\pi, \rho} = \sign( \pi, \rho) 2^{\# \text{ components in } \pi \cup \rho }. \]
\end{defn}

Let $M$ be the matrix from Corollary~\ref{corlem31} and let $D$ be the vector indexed by balanced sets $S$ with entries $D_{S} = \dfrac{Z^{D}(G \setminus S)Z^{D}(G \setminus S^c)}{(Z^D(G))^2}$ (see Lemma~\ref{lemma32gen}). Following Kenyon and Wilson, we will show that 
$$M^{T}D = \mathcal{B}_{2} Y'$$ (Lemma~\ref{lem:mylem35}). 
This result is nontrivial, requiring several lemmas, but once it is established
 it is nearly immediate that
$$M^{T}M P = \mathcal{B}_{2} Y',$$
where $P$ is the vector indexed by planar pairings $\pi$ with entries
$\widetilde{\Pr}(\pi)$ (Theorem~\ref{thm36}). Kenyon and Wilson proved that $M^{T}M$ is invertible (\cite[Theorem 3.3]{KW2009}), so we conclude the section by defining $\mathcal{Q}^{(DD)}$ as $(M^{T}M)^{-1} \mathcal{B}_{2}$.

%In the following lemma,
%$M$ is the matrix from Corollary~\ref{corlem31}
%and $D$ is the vector indexed by balanced sets $S$ with entries $D_{S} = \dfrac{Z^{D}(G \setminus S)Z^{D}(G \setminus S^c)}{(Z^D(G))^2}$ (see Lemma~\ref{lemma32gen}). 

\begin{lemma}[analogue of Lemma 3.5 from \cite{KW2006}]
\label{lem:mylem35}
$M^{T}D = \mathcal{B}_2 Y'$.
\end{lemma}

In the proof of \cite[Lemma 3.5]{KW2006}, Kenyon and Wilson use the fact that if 
the nodes of $G$ are all either black and odd or white and even and 
$\pi$ and $\rho$ are odd-even pairings, then there are $2^{\text{\# components in } \pi \cup \rho}$ balanced sets $S$ such that
$\pi$ and $\rho$ do not connect $S$ to $S^c$ (for each component, either put all of its nodes in $S$ or all of its nodes in $S^c$). Recall from Section~\ref{sec:lem31} that $T \subseteq {\bf N}$ is the set of nodes that are odd and white or even and black; under Kenyon and Wilson's assumptions, $T = \emptyset$. 
It turns out that 
after removing the requirement that the nodes be black and odd or white and even, if
%when
 $\pi$ is an odd-even pairing and $\rho$ is a black-white pairing there are still $2^{\text{\# components in } \pi \cup \rho}$ sets $S$ such that
 $\rho$ does not connect $S$ to $S^c$ 
and $\pi$ does not connect $S \triangle T$ to $(S \triangle T)^c$.

 \begin{lemma}
 \label{alglemma}
 Let $\pi$ be an odd-even pairing and let $\rho$ be a black-white pairing. For each component of $\pi \cup \rho$ there are exactly two ways to put the nodes in this component into $S$ and $S^c$
so that $\rho$ does not connect $S$ to $S^c$ 
and $\pi$ does not connect $S \triangle T$ to $(S \triangle T)^c$. 
 \end{lemma}

\begin{proof}

We start by placing an initial node $a$ into $S$ or $S^c$, and then apply the algorithm below until all nodes in the
component have been placed into $S$ or $S^c$. \\

\noindent {\bf {\underline{Algorithm}}} \\

\noindent Step 1
\begin{itemize}
\item[(a)] If $a \in S \cap T^c$ or $a \in S^c \cap T$:

%If $a \in S$ and $a \notin T$ or $a \in S^c$ and $a \in T$:

\begin{itemize}
\item[(i)] If $\pi(a) \in T$:

\begin{itemize}
\item[] Put $\pi(a)$ in $S^c$.
\end{itemize}

\item[(ii)] Else if $\pi(a) \notin T$:
\begin{itemize}
\item[]  Put $\pi(a)$ in $S$.
\end{itemize}
\end{itemize}

\item[(b)] Else if $a \in S \cap T$ or $a \in S^c \cap T^c$:

%Else if $a \in S$ and $a \in T$ or $a \in S^c$ and $a \notin T$:

\begin{itemize}
\item[(i)]  If $\pi(a) \in T$:
\begin{itemize}
\item[] Put $\pi(a)$ in $S$.
\end{itemize}

\item[(ii)] Else if $\pi(a) \notin T$:
\begin{itemize}
\item[]  Put $\pi(a)$ in $S^c$.
\end{itemize}
\end{itemize}
\item[] Go to Step 2 with $a := \pi(a)$. 
\end{itemize}

\noindent { Step 2}
\begin{itemize}
\item[] If $a \in S$:
\begin{itemize}
\item[] Put $\rho(a)$ in $S$.
\end{itemize}
\item[] Else if $a \in S^c$:
\begin{itemize}
\item[] Put $\rho(a)$ in $S^{c}$.
\end{itemize}
\item[] Go to Step 1 with $a := \rho(a)$. 
\end{itemize}

\noindent {\bf Claim.} The set $S$ described in the algorithm is well-defined and balanced. 

\begin{proof} 
We will prove this claim by induction on the number of nodes in a component of $\pi \cup \rho$.

\noindent {\em {Base Cases.}}
First note that in the case where the nodes alternate between black and white,
% is done by Kenyon and Wilson in their proof of Lemma 3.3 \cite{KW2006}, 
%since when the nodes alternate between black and white 
$T = \emptyset$ or $T = {\bf N}$ so the algorithm reduces to putting all of the nodes in a component in $S$
or all of the nodes of a component in $S^{c}$, so $S$ is well-defined. Since in this case both pairings are black-white, $S$ is balanced as well. 

If there are two nodes in a component, since $\rho$ is a black-white pairing one of the nodes is black and the other is white, so by
the previous comment there is nothing to show. 
% We start by putting node 1 in $S$. Applying Step 1(a)(ii), we put node 2 in $S$. So $S = \{1, 2\}$ and $S^{c} = \emptyset$. 

If there are four nodes in a component, since $\rho$ is a black-white pairing two nodes must be black and two nodes must be white. 
By symmetry, it is enough to consider when nodes 1 and 2 are black and nodes 3 and 4 are white. 
There are two odd-even pairings: $((1, 2), (3, 4))$ and $((1, 4), (3, 2))$ and two black-white pairings: $((1, 4), (3, 2))$ and $((1, 3), (2, 4))$.

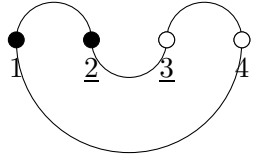
\begin{figure}[h!]

\[\begin{tikzpicture}[scale = .5]
	\vertex[fill] (n1) at (0, 0) [label=below:$1$] {};
	\vertex[fill]  (n2) at (2,0) [label=below:$\underline{2}$] {};
	\vertex (n3) at (4,0) [label=below:$\underline{3}$] {};
	\vertex (n4) at (6, 0) [label=below:$4$] {};
	
	%arcs of pi
	\draw  (n2) arc (0:180:1cm);
	\draw  (n4) arc (0:180:1cm);
	
	%arcs of rho
	\draw  (n1) arc (180:360:3cm);
	\draw  (n2) arc (180:360:1cm);

		\vertex[fill = white] (n3) at (4,0) {};
\vertex[fill = white]  (n4) at (6, 0) {};

\end{tikzpicture}\]

\caption{The diagram of $\pi \cup \rho$ when $\pi = ((1, 2), (3, 4))$ and $\rho = ((1, 4), (2, 3))$. The nodes in $T$ are underlined. If we start the algorithm by putting $1 \in S$, we get $S = \{1, 4\}$. }
\label{basecase}
\end{figure}

For example, when $\pi = ((1, 2), (3, 4))$ and $\rho = ((1, 4), (2, 3))$ (see Figure \ref{basecase}), we start by putting node 1 in $S$ . (We could also start by putting node 1 in $S^{c}$.) Then we run the algorithm: \\
Step 1.  Since $1 \notin T$ and $\pi(1) = 2 \in T$, we put $2 \in S^{c}$. \\
Step 2.  Since $2 \in S^{c}$ we put $\rho(2) = 3 \in S^{c}$. \\
Step 1.  Since $3 \in S^{c}$, $3 \in T$ and $4 \notin T$, we put $4 \in S$. \\
So we get $S = \{1, 4\}$, which is balanced. 
To check that $S$ is well-defined, it suffices to show that if we continue the algorithm for one more step, we do not get a contradiction.
% because this shows that it did not matter that we started at node 1. 
If we apply Step 2 starting at node 4, we get that we should put $\rho(4) = 1$ in $S$, as desired. 
%which is consistent. 
%$\rho$ does not bridge $S$ to $S^{c}$ and $\pi$ does not bridge $S \triangle T = \emptyset$ to to $(S \triangle T)^c = \{1, 2, 3, 4\}$. 
%Note that if we had started the algorithm by placing node 1 in $S^{c}$, then we would have
%$S = \{2, 3\}$.

In the table below are the results of applying the algorithm for each possible combination of odd-even pairings $\pi$ and black-white pairings $\rho$ that results in a component of size 4. When $\pi = \rho = ((1, 4), (3, 2))$, there are two components each of size 2, so this is omitted from the table. 

\begin{center}
\begin{tabular}{c |  c | c  | c | c | c}
$\pi$ & $\rho$   & $S$ & start & end & one more step \\ \hline
$((1, 2), (3, 4))$ & $((1, 4), (2, 3))$ & \{1, 4\} & $1 \in S$ & $4 \in S$ & $1 \in S$ \\
$((1, 2), (3, 4))$ & $((1, 3), (2, 4))$  & \{1, 3\} & $1 \in S$ & $3 \in S$ & $1 \in S$ \\
$((1, 4), (3, 2))$ & $((1, 3), (2, 4))$  & \{1, 2, 3, 4\} & $1 \in S$ & $3 \in S$ & $1 \in S$
\end{tabular}
\end{center}

In each case, $S$ is balanced, and continuing the algorithm for one more step does not create a contradiction.\\

\begin{figure}[h]
\centering
\begin{tikzpicture}[scale = .5]
	\vertex[fill] (n1) at (0, 0) [label=below:\small{$a$}] {};
	\vertex (n2) at (2,0) [label=below:\small{$\pi(a)$}] {};
	\vertex[fill] (n3) at (4,0) [label=above:$$] {};
	\vertex[fill] (n4) at (6, 0) [label=below:\small{$\rho(\pi(a))$}] {};
	\vertex (n5) at (8,0) [label=below:\small{$\rho(a)$}] {};
	\vertex[fill] (n6) at (10,0) [label=below:$$] {};
	
	%arcs of pi
	\draw  (n2) arc (0:180:1cm);
	\draw  (n4) arc (0:180:1cm);
	\draw  (n6) arc (0:180:1cm);
	
	%arcs of rho
	\draw  (n1) arc (180:360:4cm);
	\draw  (n2) arc (180:360:2cm);
	\draw  (n3) arc (180:360:3cm);
	
		\vertex[fill=white] (n2) at (2,0) {};
		\vertex[fill=white] (n5) at (8,0) {};

\end{tikzpicture} \hspace{1cm}
\begin{tikzpicture}[scale = .5]
	\vertex[fill] (n1) at (0, 0) [label=below:\small{$a$}] {};
	\vertex (n2) at (2,0) [label=below:\small{$\pi(a)$}] {};
	\vertex[fill] (n3) at (4,0) [label=above:$$] {};
	\vertex[fill] (n4) at (6, 0) [label=below:\small{$\rho(\pi(a))$}] {};
	\vertex (n5) at (8,0) [label=below:\small{$\rho(a)$}] {};
	\vertex[fill] (n6) at (10,0) [label=below:$$] {};
	
	%arcs of pi
	\draw  (n2) arc (0:180:1cm);
	\draw  (n4) arc (0:180:1cm);
	\draw  (n6) arc (0:180:1cm);
	
	%arcs of rho
	\draw  (n3) arc (180:360:3cm);
	\draw  (n4) arc (180:360:1cm);

	  \filldraw[color=white,draw = white] (0,-4) circle (0.05cm); 

	\vertex[fill=white] (n2) at (2,0) {};
		\vertex[fill=white] (n5) at (8,0) {};
\end{tikzpicture}

\caption{Illustration of case 1. Shown left is the odd-even pairing $\pi$ (top) and the black-white pairing $\rho$ (bottom). On the right we have replaced $\rho$ with $\tilde{\rho}$, the black white pairing on ${\bf N}$ $-$ $\{a, \pi(a) \}$ obtained by pairing $\rho(\pi(a))$ with $\rho(a)$.
}
\label{case1}
\end{figure}
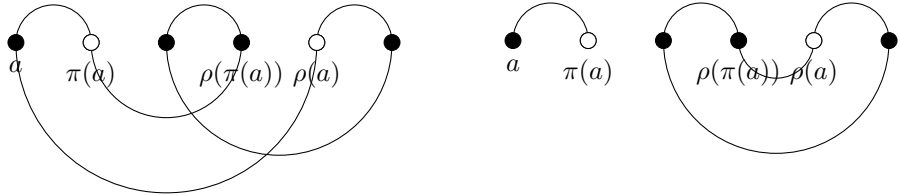

Now suppose that a component of $\pi \cup \rho$ has $2n$ nodes, where $2n > 4$. Assume that if a component has fewer than $2n$ nodes, the set $S$ is well-defined and balanced. Let ${\bf N}$ denote the set of nodes in this component. There are two cases to consider based on whether or not $\pi |_{\bf N}$ has a black-white pair.\\

\noindent {\bf Case 1.} (Illustrated in Figure \ref{case1}). 
Assume $\pi |_{\bf N}$ has at least one black-white pair $(a, \pi(a))$.
Since $\rho$ is a black-white pairing, $\rho(a)$ and $\rho(\pi(a))$ are opposite color.  
Consider the black-white pairing $\tilde{\rho}$ on ${\bf{N}} - \{a, \pi(a) \}$
obtained from $\rho$ by removing the pairs $(a, \rho(a))$ and $(\pi(a), \rho(\pi(a)))$ and adding the pair $(\rho(a), \rho(\pi(a)))$.
 %defined as follows: $\tilde{\rho}(i) = \rho(i)$ if $i \neq \rho(a), \rho(\pi(a))$ and $\tilde{\rho}(\rho(a)) = \rho(\pi(a))$. 
 Let $\tilde{\pi} = \pi |_{{\bf N} - \{a, \pi(a) \}}$. 
Now $\tilde{\pi} \cup \tilde{\rho}$ is a single component with $2n-2$ nodes. 
Start the algorithm by putting $\rho(\pi(a)) \in S$. By the induction hypothesis, the set $S$ produced by the algorithm is well-defined and balanced.
Note that the fact that $S$ is well-defined means that $\rho(a) \in S$. 

Considering the original component of $\pi \cup \rho$, when we start the algorithm at $\rho(\pi(a))$ it proceeds identically as it did with $\tilde{\pi} \cup \tilde{\rho}$ until we reach the node $\rho(a)$. Since $\rho(a) \in S$, applying Step 2 of the algorithm we add $a$ to $S$. 
(Note that we are guaranteed to be on Step 2 here by the fact that $\rho(a)$ is paired with $\rho(\pi(a))$ in $\tilde{\rho}$, and the algorithm starts with Step 1.)
Since $\pi$ is odd-even, black-white pairs of $\pi$ have the property that either both nodes are in $T$ or both are not in $T$. So after the next step of the algorithm (Step 1) we add $\pi(a)$ to $S$. Since we added $a$ and $\pi(a)$ to $S$, $S$ is still balanced.  Since $\pi(a) \in S$, continuing the algorithm for one more step would put $\rho(\pi(a)) \in S$, which is consistent. \\

\begin{figure}[h]
\centering
\begin{tikzpicture}[scale = .5]
	\vertex[fill] (n3) at (4, 0) [label=below:\small{$\pi(b)$}] {};
	\vertex (n4) at (6,0) [label=above:$$] {};
	\vertex[fill] (n5) at (8,0) [label=above:$$] {};
	\vertex[fill] (n6) at (10, 0) [label=below:\footnotesize{$\rho(\pi(a))$}] {};
	\vertex (n7) at (12,0) [label=below:\footnotesize{$\rho(\pi(b))$}] {};
	\vertex[fill] (n8) at (14,0) [label=below:\small{$b$}] {};
	\vertex (n11) at (16,0) [label=below:\small{$a$}] {};
	\vertex (n12) at (18, 0) [label=below:\small{$\pi(a)$}] {};
	
	%arcs of pi
	\draw  (n12) arc (0:180:1cm);
	\draw  (n8) arc (0:180:5cm);
	\draw  (n7) arc (0:180:3cm);
	\draw  (n6) arc (0:180:1cm);
	
	%arcs of rho
	\draw  (n3) arc (180:360:4cm);
	\draw  (n4) arc (180:360:1cm);
	\draw  (n6) arc (180:360:4cm);
	\draw  (n8) arc (180:360:1cm);
	
		\vertex[fill=white] (n4) at (6,0) {};
		\vertex[fill=white]  (n7) at (12,0) {};
			\vertex[fill=white] (n11) at (16,0)  {};
	\vertex[fill=white] (n12) at (18, 0) {};

\end{tikzpicture} \hspace{-1.75cm}
\begin{tikzpicture}[scale = .5]

	\vertex[fill] (n3) at (4, 0) [label=below:\small{$\pi(b)$}] {};
	\vertex (n4) at (6,0) [label=above:$$] {};
	\vertex[fill] (n5) at (8,0) [label=above:$$] {};
	\vertex[fill] (n6) at (10, 0) [label=below:\footnotesize{$\rho(\pi(a))$}] {};
	\vertex (n7) at (12,0) [label=below:\footnotesize{$\rho(\pi(b))$}] {};
	\vertex[fill] (n8) at (14,0) [label=below:\small{$b$}] {};
	\vertex (n11) at (16,0) [label=below:\small{$a$}] {};
	\vertex (n12) at (18, 0) [label=below:\small{$\pi(a)$}] {};
	
	%arcs of pi
	\draw  (n12) arc (0:180:1cm);
	\draw  (n8) arc (0:180:5cm);
	\draw  (n7) arc (0:180:3cm);
	\draw  (n6) arc (0:180:1cm);
	
	%arcs of rho
	\draw  (n4) arc (180:360:1cm);
	\draw  (n6) arc (180:360:1cm);

	 \filldraw[color=white,draw = white] (0,-4) circle (0.05cm); 

	\vertex[fill=white] (n4) at (6,0) {};
		\vertex[fill=white]  (n7) at (12,0) {};
					\vertex[fill=white] (n11) at (16,0)  {};
	\vertex[fill=white] (n12) at (18, 0) {};
	
\end{tikzpicture}

\caption{Illustration of case 2.
Shown left is the odd-even pairing $\pi$ (top) and the black-white pairing $\rho$ (bottom). On the right we have replaced $\rho$ with $\tilde{\rho}$, the black white pairing on ${\bf N}$ $-$ $\{a, \pi(a), b, \pi(b) \}$ obtained by pairing $\rho(\pi(a))$ with $\rho(\pi(b))$.
}
\label{case2a}
\end{figure}
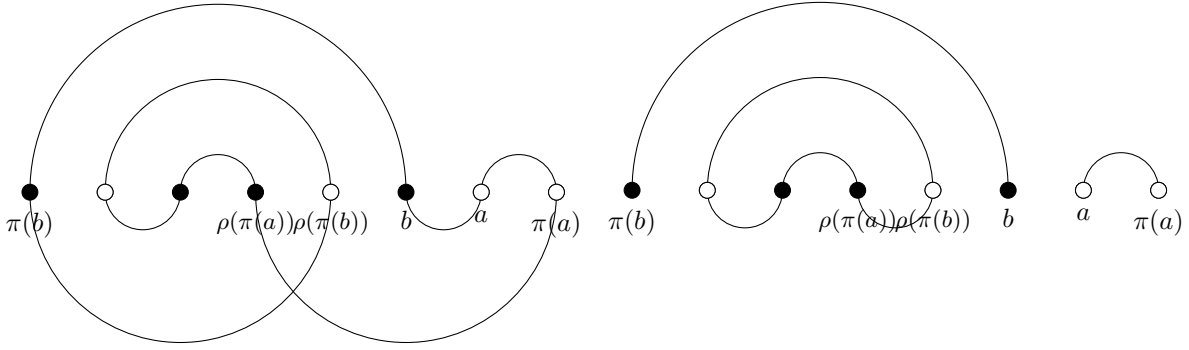

\noindent {\bf Case 2.} (Illustrated in Figure~\ref{case2a}). 
If $\pi |_{\bf N}$ does not have a black-white pair, then consider a white pair of $\pi |_{\bf N}$: $(a, \pi(a))$. Let $b = \rho(a)$. Since $a$ is white, $\rho(a)$ must be black, and $(b, \pi(b))$ is a black pair of $\pi |_{\bf N}$ by the assumption that $\pi |_{\bf N}$ does not have a black-white pair. 
Consider the black-white pairing $\tilde{\rho}$ on ${\bf{N}} - \{a, \pi(a), b, \pi(b)\}$ obtained from $\rho$ by removing the pairs $(a, b), (\pi(a), \rho(\pi(a))),$ and $(\pi(b), \rho(\pi(b)))$ and adding the pair $(\rho(\pi(a)), \rho(\pi(b)))$. 
%defined as follows: $\tilde{\rho}(i) = \rho(i)$ if $i \neq \rho(\pi(a)), \rho(\pi(b))$ and $\tilde{\rho}(\rho(\pi(a))) = \rho(\pi(b))$. 
Let $\tilde{\pi} = \pi |_{{\bf N} - \{a, \pi(a), b, \pi(b) \}}$. 
Now $\tilde{\pi} \cup \tilde{\rho}$ is a single component with $2n-4$ nodes. 
Start the algorithm by putting $\rho(\pi(a)) \in S$. By the induction hypothesis, the set $S$ produced by the algorithm is well-defined and balanced.
Note that the fact that $S$ is well defined means that $\rho(\pi(b)) \in S$. 

Considering the original component of $\pi \cup \rho$, when we start the algorithm by putting $\rho(\pi(a)) \in S$ it proceeds identically as it did with $\tilde{\pi} \cup \tilde{\rho}$ until we reach the node $\rho(\pi(b))$. Since $\rho(\pi(b)) \in S$, applying Step 2 of the algorithm we add $\pi(b)$ to $S$. Since $\pi$ is odd-even, exactly one of $\{b, \pi(b)\}$ is in $T$. This means that after applying Step 1 we put $b \in S^{c}$. Then we put $a \in S^{c}$ (since $\rho(a) = b$) and $\pi(a) \in S$ (since exactly one of $\{a, \pi(a)\}$ is in $T$). Since we added $\pi(b)$ and $\pi(a)$ to $S$, $S$ is still balanced. Since $\pi(a) \in S$, continuing the algorithm for one more step puts $\rho(\pi(a)) \in S$, which is consistent. 

\end{proof}

\noindent {\bf Claim.} After applying this algorithm, $\rho$ does not connect $S$ to $S^c$ and $\pi$ does not connect $S \triangle T$ to $(S \triangle T)^c$.

\begin{proof}
By Step 2, for each node $a$, $a$ and $\rho(a)$ will either both be in $S$ or $S^{c}$, so $\rho$ does not connect $S$ to $S^c$.
To show that $a$ and $\pi(a)$ are either both in $S \triangle T$ or both in $(S \triangle T)^c$, there are several cases to consider. 
\begin{itemize}
\item If $a$ and $\pi(a)$ are both not in $T$, then they are both placed into $S$ by Step 1(a)(ii) or both placed into $S^{c}$ by Step 1(b)(ii). 
In the first case, $a$ and $\pi(a)$ are both in $S \triangle T$, and in the second case $a$ and $\pi(a)$ are both in $(S \triangle T)^c$. 
%Then since $a$ and $\pi(a)$ are not in $T$, they will both be in $
\item If $a \in T$ and $\pi(a) \notin T$, then one of $a, \pi(a)$ is placed into $S$ and one is placed into $S^{c}$ by Step 1(a)(ii) or Step 1(b)(ii). If $a$ is placed into $S$ and $\pi(a)$ is placed into $S^{c}$, then $a$ and $\pi(a)$ are both in $(S \triangle T)^c$. The other case is similar. 
\item If $a \in T$ and $\pi(a) \in T$, then they are both placed into $S^{c}$ by Step 1(a)(i) or both placed into $S$ by Step 1(b)(i). 
\item If $a \notin T$ and $\pi(a) \in T$, then one is placed in $S$ and one is placed in $S^{c}$ by Step 1(a)(i) or Step 1(b)(i). 
\end{itemize}
\end{proof}

We have shown that the algorithm produces a well-defined balanced set $S$ with the desired properties.
We conclude that for each component of $\pi \cup \rho$ there are exactly two ways to put the nodes in this component into $S$ and $S^c$
so that $\rho$ does not connect $S$ to $S^c$ 
and $\pi$ does not connect $S \triangle T$ to $(S \triangle T)^c$. 

\end{proof}

We need two more facts to prove Lemma~\ref{lem:mylem35}.

\begin{lemma} 
\label{lem:signSdefn2}
Let $S$ be a balanced subset of nodes and let $\sign(S)$ be defined as in Lemma \ref{lemma32gen}. Then
\begin{equation*}
 \text{sign}(S) = (-1)^{\frac{\# \text{nodes}}{2}} (-1)^{\# \text{ comp in } \pi \cup \rho} \sign_{OE}(\pi) \sign_{BW}(\rho),
\end{equation*}
where $\pi$ is an odd-even pairing such that $\pi$ does not connect $S \triangle T$ to $(S \triangle T)^c$ and $\rho$ is a black-white pairing such that 
$\rho$ does not connect $S$ to $S^c$. 
\end{lemma}

The proof of Lemma~\ref{lem:signSdefn2} is lengthy and technical so we postpone it to Section~\ref{sec:signSformula} for ease of exposition. The following is an immediate consequence of this lemma.

\begin{cor}
\label{samesignlemma}
Let $\pi$ be an odd-even pairing and let $\rho$ be a black-white pairing. If $S_1$ and $S_2$ are balanced subsets of nodes such that
$\pi$ does not connect $S_i \triangle T$ to $(S_i \triangle T)^c$ and $\rho$ does not connect $S_i$ to $S_i^c$ for $i = 1, 2$, then
$\text{sign}(S_1) = \text{sign}(S_2)$.
\end{cor}

To see that Corollary~\ref{samesignlemma} follows from Lemma~\ref{lem:signSdefn2}, observe that if $\pi$ and $\rho$ satisfy the hypotheses of Lemma 2.3.7 for $S_1$ and $S_2$, then $S_1$ and $S_2$ must have the same sign, because all of the quantities on the right hand side of the equation in Lemma~\ref{lem:signSdefn2} depend only on $\pi$ and $\rho$.

\begin{proof}[Proof of Lemma \ref{lem:mylem35}]
Recall from Lemma \ref{lemma32gen} that
\begin{equation*}
D_{S} = \dfrac{Z^{D}(G \setminus S) Z^{D}(G \setminus S^{c})}{(Z^{D}(G))^{2}} = \sign_{\cons}({\bf N}) \text{sign}(S)
\det \left[(1_{i, j \in S} + 1_{i, j \notin S}) \times \text{sign}(i, j) Y_{i, j} \right]^{i = b_1,b_2, \ldots, b_n}_{j = w_1,w_2, \ldots, w_n},
\end{equation*}
where $b_1, b_2, \ldots, b_n$ are the black nodes listed in ascending order and $w_1, w_2, \ldots, w_n$ are the white nodes listed in ascending order.

When we expand the determinant, we get
\[D_{S} = \sign_{\cons}({\bf N}) \text{sign}(S) 
\sum_{ \substack{ \text{ BW pairings $\rho$: } \\ \text{ $\rho$ does not connect } \\ \text{ $S$ to $S^{c}$ } } }
\sign_{BW}(\rho) \prod\limits_{(i, j) \in \rho } \text{sign}(i, j)Y_{i, j}. \]
By Lemma \ref{lemma34}, 
\[D_{S} = \sign_{\cons}({\bf N})  \text{sign}(S) 
\sum_{ \substack{ \text{ BW pairings $\rho$: } \\ \text{ $\rho$ does not connect } \\ \text{ $S$ to $S^{c}$ } } }
\sign_{\cons}({\bf N}) (-1)^{\text{\# crosses of $\rho$} }
 \prod\limits_{(i, j) \in \rho } Y_{i, j}, \]
and thus by definition,
\begin{equation}
D_{S} = \text{sign}(S) 
\sum_{ \substack{ \text{ BW pairings $\rho$: } \\ \text{ $\rho$ does not connect } \\ \text{ $S$ to $S^{c}$ } } }  Y'_{\rho}.
\end{equation}

Let $\pi$ be a planar pairing and let $M$ be the matrix from Corollary \ref{corlem31}. The $\pi$th row of $M^{T}D$ is
\[ \sum_{ \substack{ S \subseteq \{1, 2, \ldots, 2n\} \\ \pi \text{ does not connect} \\ S \triangle T \text{ to } (S \triangle T)^c} } D_S. \]
%where $T$ is the set of toggled nodes (nodes that are odd and white or even and black). 
We see that
\begin{eqnarray*}
\sum_{ \substack{ S \subseteq \{1, 2, \ldots, 2n\}: \\ \pi \text{ does not connect} \\ S \triangle T \text{ to } (S \triangle T)^c} } D_S 
&= &  \sum_{ \substack{ S \subseteq \{1, 2, \ldots, 2n\}: \\ \pi \text{ does not connect} \\ S \triangle T \text{ to } (S \triangle T)^c} } \sign(S)
\sum_{ \substack{ \text{BW pairings } \rho: \\ \rho \text{ does not connect} \\ S \text{ to  }S^{c} }}  Y'_{\rho} \\
& = & \sum_{\text{BW pairings } \rho } \sum_{\substack{ S\text{: }\rho \text{ does not } \\ \text{ connect } S \text{ to  }S^{c} \text{ and }  \\ \pi \text{ does not connect} \\ S \triangle T \text{ to } (S \triangle T)^c} } \sign(S)Y'_{\rho}.
\end{eqnarray*}

By Lemma \ref{alglemma} and Corollary \ref{samesignlemma},
\[
 \sum_{\text{BW pairings } \rho } \sum_{\substack{ S\text{: }\rho \text{ does not } \\ \text{ connect } S \text{ to  }S^{c} \text{ and }  \\ \pi \text{ does not connect} \\ S \triangle T \text{ to } (S \triangle T)^c} } \sign(S)Y'_{\rho} 
 =  \sum_{\text{BW pairings } \rho } \sign( \pi, \rho ) 2^{ \# \text{ comp in } \pi \cup \rho } Y'_{\rho} .
 \]
Since this sum is the $\pi$th row of $\mathcal{B}_{2} Y'$, we have proven the claim. 

\end{proof}

\begin{thm}[analogue of Theorem 3.6 from \cite{KW2006}]
\label{thm36}
$ M^T M P = \mathcal{B}_{2} Y'$
\end{thm}

\begin{proof}
Noting that 
$\widetilde{\Pr}(\pi) = \Pr(\pi) Z^{DD}(G)/(Z^D(G))^2$, we see that
by Corollary~\ref{corlem31}, $MP = D$. 
Then, applying Lemma~\ref{lem:mylem35} we have 
%$\mathcal{M}_{2} P = 
$M^T M P = M^T D = \mathcal{B}_{2} Y'$.
\end{proof}

It remains to show that $M^T M$ is invertible. In fact, $M^{T}M$ is equal to the {\em meander matrix} $\mathcal{M}_{q}$ evaluated at $q = 2$.

% Kenyon and Wilson prove this for the matrix from Lemma \ref{lem3.1}, but as the only difference between the matrix from Corollary \ref{corlem31} and the matrix from Lemma \ref{lem3.1} is the ordering of the rows, Kenyon and Wilson's proof hold for the matrix $M$ from Corollary \ref{corlem31}  as well.

\begin{lemma}\cite[Lemma 3.3]{KW2006}
\label{lem3.3}
Let $M$ be the matrix from Lemma \ref{lem3.1}. Then
$M^{T}M = \mathcal{M}_{2}$, where $\mathcal{M}_{2}$
is a matrix with rows and columns indexed by planar pairings, with entries
$$(\mathcal{M}_{2})_{\sigma, \tau} = 2^{\# \text{ comp in } \sigma \cup \tau}$$
%the meander matrix $\mathcal{M}_{q}$ evaluated at $q = 2$. 
\end{lemma}

Since the only difference between the matrix from Lemma \ref{lem3.1} and the matrix from Corollary \ref{corlem31} is the ordering of the rows, Lemma~\ref{lem3.3} applies to the matrix $M$ from Corollary \ref{corlem31}  as well.

%\subsection{Proof of Theorem \ref{thm:thm1}}

%\label{sec:integervalued}

\begin{defn}
\label{matrixdefn}
Since $\mathcal{M}_{2}$ is invertible (see \cite{meanders}), define
$$\mathcal{Q}^{(DD)} = \mathcal{M}_{2}^{-1} \mathcal{B}_{2}$$
\end{defn}
Since $P = \mathcal{Q}^{(DD)} Y'$, $\mathcal{Q}^{(DD)}$ is the matrix of the $Y'$ polynomials: for a given planar pairing $\pi$, the $\pi$th row of $\mathcal{Q}^{(DD)}$ gives the polynomial $\widetilde{\Pr}(\pi)$.
% = \dfrac{\text{Pr}(\pi) Z^{DD}(G)}{(Z^D(G))^2} = \dfrac{Z^{DD}_{\pi}(G)}{(Z^D(G))^2} $. 

That is, 
\begin{equation*}
\widetilde{\Pr}(\sigma) 
=
 \sum_{\text{black-white pairings } \rho} \mathcal{Q}^{(DD)}_{\sigma, \rho} Y'_{\rho}.
 \end{equation*}

Our next aim is to prove that $\mathcal{Q}^{(DD)}$ is integer-valued.
%To prove that $\mathcal{Q}^{(DD)}$ is integer-valued, 
To this end, we will show that we can compute the columns combinatorially using a {\em transformation rule} from Kenyon and Wilson's study of {\em groves} \cite{KW2006}. 

\subsection{Groves}
\label{sec:groves}

\begin{defn}\cite{KW2006}
If $G$ is a finite edge-weighted planar graph embedded in the plane with a set of nodes, a {\em grove} is a spanning acyclic subgraph of $G$ such that each component tree contains at least one node. The weight of a grove is the product of the weights of the edges it contains. 
\end{defn}

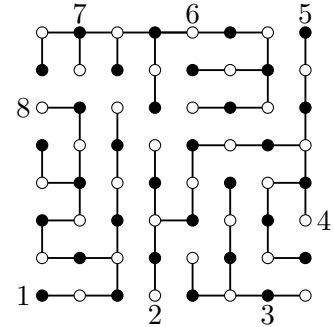
\begin{wrapfigure}{r}{.3\textwidth}
%%grove
\centering
\vspace{-.5cm}

\begin{tikzpicture}[scale=.5]
\def\maxX{7}
%\foreach \x [count=\n] in {0,...,\maxX}{
 %   \foreach \y in {0,...,7}{
 %       \draw[line width=.5pt] (\x,0) -- (\x,7);
 %       \draw[line width=.5pt] (0,\y) -- (\maxX,\y);
 %   }
%}

\node at  (-0.5,0) {$1$};
\node at  (3,-0.5) {$2$};
\node at  (6,-0.5) {$3$};
\node at  (7.5,2) {$4$};
\node at (7, 7.5) {$5$};
\node at (4, 7.5) {$6$};
\node at (1, 7.5) {$7$};

\node at  (-.5,5) {$8$};

%path from node 1 to node 8
   \draw[line width = .25mm]  (0, 0) -- (1, 0);
   \draw[line width = .25mm] (2, 0) -- (2, 1); 
   \draw[line width = .25mm]  (1, 1) -- (0, 1); 
   \draw[line width = .25mm] (1, 3) -- (1, 4); 
   \draw[line width = .25mm] (0, 2)-- (1, 2);
   \draw[line width = .25mm] (0, 5) -- (1, 5);
   
      \draw[line width = .25mm] (0, 3) -- (0, 4);
       \draw[line width = .25mm] (0, 3) -- (1, 3);
  
          \draw[line width = .25mm] (1, 6) -- (1, 7);
           \draw[line width = .25mm] (2, 6) -- (2, 7);
  
      \draw[line width = .25mm] (0, 6) -- (0, 7);
       \draw[line width = .25mm] (0, 7) -- (1, 7);
   
    \draw[line width = .25mm] (2, 2) -- (2, 1);
    \draw[line width = .25mm] (2, 2) -- (2, 3);
    \draw[line width = .25mm] (2, 3) -- (2, 4);
       \draw[line width = .25mm] (2, 4) -- (2, 5);
    
   \draw[line width = .25mm]  (3, 3) -- (3, 4);

%path from node 1 to node 8
   \draw[line width = .25mm] (1, 0) -- (2, 0);
 \draw[line width = .25mm]  (2, 1) -- (1, 1);
 \draw[line width = .25mm] (0, 1) -- (0, 2);
  \draw[line width = .25mm]  (1, 2) -- (1, 3);
   \draw[line width = .25mm] (1, 4) -- (1, 5);

    %node 7 to node 6 
   \draw[line width = .25mm]  (1, 7) -- (2, 7); 
       \draw[line width = .25mm] (3, 7) -- (4, 7); 
    \draw[line width = .25mm] (3, 7) -- (4, 7); 

\draw[line width = .25mm]   (3, 6) -- (3, 7);
\draw[line width = .25mm]   (6, 6) -- (6,  7);
   %\draw[line width = .25mm]   (3, 6) -- (4, 6);
     \draw[line width = .25mm]  (5, 6) -- (6, 6);
    %\draw[line width = .25mm] (3, 5) -- (4, 5);
    \draw[line width = .25mm]   (5, 5) -- (6, 5);
    
    \draw[line width = .25mm]   (3,2) -- (3, 3);

     \draw[line width = .25mm]   (5, 7) -- (6, 7);
          \draw[line width = .25mm]   (4, 7) -- (5, 7);

  \draw[line width = .25mm] (2, 7) -- (3, 7); 
\draw[line width = .25mm]  (4, 6) -- (5, 6);
   \draw[line width = .25mm]  (4, 5) -- (5, 5);
\draw[line width = .25mm]  (6, 6) -- (6, 5);
 \draw[line width = .25mm]  (3, 5) -- (3, 6);

     %node 2 to node 5
   \draw[line width = .25mm]  (3, 0) -- (3, 1);
   \draw[line width = .25mm]   (3, 2) -- (4, 2); 
   \draw[line width = .25mm]   (4,3)-- (4, 4);
   \draw[line width = .25mm]  (5, 4) -- (6, 4);
   \draw[line width = .25mm] (7, 4) -- (7,5);
   \draw[line width = .25mm]   (7, 6) -- (7,7);
      %node 2 to node 5
\draw[line width = .25mm] (3, 1) -- (3, 2);
\draw[line width = .25mm]  (4, 2) -- (4,3);
\draw[line width = .25mm]  (4, 4) -- (5, 4);
\draw[line width = .25mm]  (6, 4) -- (7, 4);
\draw[line width = .25mm]  (7, 5) -- (7,6);

   %node 3 to node 4
      \draw[line width = .25mm]  (6, 0) -- (7, 0);
   \draw[line width = .25mm]   (7, 1) -- (6, 1);
   \draw[line width = .25mm]   (6, 2) -- (6, 3);
   \draw[line width = .25mm]  (7, 3) -- (7, 2);

    \draw[line width = .25mm] (5, 2) -- (5, 3);
 \draw[line width = .25mm] (5, 2) -- (5, 1);
      
   \draw[line width = .25mm]  (5, 0) -- (5, 1);
   \draw[line width = .25mm] (4, 1) -- (4, 0);
   
\draw[line width = .25mm] (6, 1) -- (6, 2);
\draw[line width = .25mm]  (6, 3)-- (7, 3);
\draw[line width = .25mm] (4, 0) -- (5, 0);
\draw[line width = .25mm] (5, 0) -- (6, 0);

\draw[line width = .25mm] (7, 3) -- (7, 4);
   	
\foreach \x [count = \n] in {0, 2, 4, 6}{
\foreach \y in {0,  2, 4, 6}{
       \filldraw[fill=black, draw=black] (\x,\y) circle (0.15cm); 
         \filldraw[fill=white, draw=black] (\x+1,\y) circle (0.15cm); 
        }
        }
\foreach \x [count = \n] in {1, 3, 5, 7}{
\foreach \y in {1, 3, 5, 7}{
        \filldraw[fill=black, draw=black] (\x,\y) circle (0.15cm); 
          \filldraw[fill=white, draw=black] (\x-1,\y) circle (0.15cm); 
        }
        }
\end{tikzpicture}

\vspace{-.25cm}

\caption{A grove of a grid graph with $8$ nodes. The partition of the nodes is $\{ \{1, 8\}, \{2, 4, 5\}, \{3\}, \{6, 7\} \}$. }
\label{fig:grove}
\vspace{-1.5cm}
\end{wrapfigure}

The connected components of a grove partition the nodes into a planar partition. If $\sigma$ is a planar partition of $1, 2, \ldots, n$, let $\Pr(\sigma)$ be the probability that a random grove of $G$ partitions the nodes according to $\sigma$. Kenyon and Wilson showed that
 $\dddot{ \Pr}(\sigma) := \dfrac{ \Pr(\sigma) }{ \Pr( 1|2| \cdots |n)}$ is an integer-coefficient homogeneous polynomial in the variables $L_{i, j}$\footnotemark~\cite[Theorem 1.2]{KW2006}.
%The partition $\sigma = 1 | 2 | 3 | \cdots | n$ is called the {\underline {uncrossing}}. 

\footnotetext{When $G$ is viewed as a resistor network with conductances equal to the edge weights, $L_{i, j}$ is the current that would flow into node $j$ if node $i$ were set to one volt and all other nodes were set to zero volts \cite[Appendix A]{KW2006}.}

For example, the normalized probability  $\dddot{ \text{Pr}}(\sigma)$ that a random grove on four nodes partitions the nodes according to $1 | 234$ is
 $\dddot{ \text{Pr}}(\sigma) = L_{2, 3} L_{3, 4} + L_{2, 3} L_{2, 4} + L_{2, 4} L_{3, 4} + L_{1, 3} L_{2, 4}.$
 (See \cite[Section 1.2]{KW2006}). 
 
 Each monomial in the polynomial $\dddot{ \text{Pr}}(\sigma)$ is of the form $L_{\tau} = \sum\limits_{F} \prod\limits_{\{i, j\} \in F} L_{i, j}$. The sum is over spanning forests $F$ of the complete graph $K_n$ for which the trees of $F$ span the parts of $\tau$ and the product is over edges $\{i, j\}$ of the forest $F$.
 
 To compute these polynomials,
Kenyon and Wilson define a matrix $\mathcal{P}^{(t)}$ with rows indexed by planar partitions and columns indexed by all partitions and show how to compute the columns of this matrix combinatorially. The $\tau$th column of $\mathcal{P}^{(t)}$ is computed by writing the partition $\tau$ as a linear combination of planar partitions. So if $\tau$ is planar, then $\mathcal{P}^{(t)}_{\tau, \tau} = 1$ and $\mathcal{P}^{(t)}_{\sigma, \tau} = 0$ for all $\sigma \neq \tau$. If $\tau$ is nonplanar, the rule is a generalization of the rule for four nodes:
\begin{equation}
\label{eqn:4noderule}
13|24 \rightarrow
1 | 234 + 2 | 314 + 3 | 124 + 4 | 123- 12|34 - 14|23 
\end{equation}
This rule tells us, for example, that $P^{(t)}_{12|34, 13|24} = -1$.

%The generalization of the transformation rule is the following. 
In general, if a partition is nonplanar, then there will exist nodes $a < b < c < d$ such that $a$ and $c$ belong to one part, and $b$ and $d$ belong to another part. In Kenyon and Wilson's transformation rule, 1, 2, 3, and 4 in equation (\ref{eqn:4noderule})  are replaced with parts $A, B, C$ and $D$, which contain the nodes $a, b, c$ and $d$, respectively. 
%Kenyon and Wilson define the following rule for transforming a nonplanar partition into a linear combination of planar partitions. 

%\begin{defn}\cite{KW2006}
%\label{ltau}
%Let $L_{i, j}$ denote the current that would flow into node $j$ if node $i$ were sent to one volt and the remaining nodes set to zero volts. Note that $L_{i, j} = L_{j, i}$. 
%For a partition $\tau$ on $1, \ldots, n$, define $L_{\tau} = \sum\limits_{F} \prod\limits_{\{i, j\} \in F} L_{i, j}$ where the sum is over spanning forests $F$ of the complete graph $K_n$ for which the trees of $F$ span the parts of $\tau$ and the product is over edges $\{i, j\}$ of the forest $F$.
%\end{defn}

%If a partition is nonplanar, then there will exist nodes $a < b < c < d$ such that $a$ and $c$ belong to one part, and $b$ and $d$ belong to another part.
%Kenyon and Wilson define the following rule for transforming a nonplanar partition into a linear combination of planar partitions. 

\begin{trule}\cite[Rule 1]{KW2006} 
\label{kwrule1}
Arbitrarily subdivide the part containing $a$ and $c$ into two sets $A$ and $C$ such that $a \in A$ and $c \in C$, and similarly subdivide the part containing $b$ and $d$ into $B \ni b$ and $D \ni d$. Let the remaining parts of the partition be denoted by ``rest." Then the transformation rule is
$$AC|BD|\text{rest} \to A|BCD|\text{rest}  +  B|ACD|\text{rest} +   C|ABD|\text{rest} +   D|ABC|\text{rest} -   AB|CD|\text{rest} - AD|BC|\text{rest}$$
\end{trule}

\begin{rem}
\label{rem:rule1error}
If we arbitrarily subdivide the part containing $a$ and $c$ into two sets $A$ and $C$ such that $a \in A$ and $c \in C$, and similarly for the part containing $b$ and $d$, it is possible to repeat Rule~\ref{kwrule1} indefinitely without ever obtaining a linear combination of planar partitions. 

For example, consider the partition $1235|46$. One crossing is $a = 1$, $b = 4$, $c= 5$, $d= 6$. If we choose $A = \{1, 2, 3\}$, $B = \{4\}, C = \{5\},$ and $D = \{6\}$, then after applying Rule~\ref{kwrule1}, all of the partitions are planar. But if we choose $A = \{1, 2\}$, $B = \{4\}$, $C = \{3, 5\}$, and $D = \{6\}$ then after applying Rule~\ref{kwrule1} we get
$$12|3456 + 4|12356 + 35|1246 + 6|12345 - 124|356 - 126|345$$
which includes nonplanar partitions. For example, the partition $124|356$ has crossing $a = 1$, $b = 3$, $c= 4$, $d= 6$. If we choose $A = \{1, 2\}$, $B = \{3, 5\}$, $C = \{4\}$, and $D = \{6\}$ then after applying Rule~\ref{kwrule1} to $124|356$ we get
 $$12|3456 + 35|1246 + 4|12356 + 6|12345 - 1235|46 - 126|345.$$
So after applying Rule~\ref{kwrule1} twice, all partitions cancel except for the partition $1235|46$, which is the partition we started with. We could continue this process indefinitely.

\end{rem}

Remark~\ref{rem:rule1error} motivates the following modification of Rule~\ref{kwrule1}. 

\begin{trule}
\label{hrule}
% If a partition $\tau$ is nonplanar, then there will exist items $a < b < c < d$ such that $a$ and $c$ belong to one part, and $b$ and $d$ belong to another part. 
Subdivide the part containing $a$ and $c$ into two sets $A$ and $C$ such that $A$ contains all the items in this part less than $b$, and $C$ contains all other items. Similarly, subdivide the part containing $b$ and $d$ into two sets $B$ and $D$ so that $B$ contains all items in this part less than $c$, and $D$ contains all other items.
Then the transformation rule is
$$AC|BD|\text{rest} \to A|BCD|\text{rest}  +  B|ACD|\text{rest} +   C|ABD|\text{rest} +   D|ABC|\text{rest} -   AB|CD|\text{rest} - AD|BC|\text{rest}$$
\end{trule}

Applying Rule~\ref{hrule} repeatedly will result in a linear combination of planar partitions. 

%Kenyon and Wilson's main result concerning groves is the following theorem. 
%As we mentioned above, Kenyon and Wilson's main result about groves concerns the computation of the polynomials $\dddot{ \Pr}(\sigma) := \dfrac{ \Pr(\sigma) }{ \Pr( 1|2| \cdots |n)}$.
We have now presented all the definitions needed to state
Kenyon and Wilson's main result for groves.

\begin{thm}\cite[Theorem 1.2]{KW2006} 
\label{thm:KWgrovethm}
Any partition $\tau$ may be transformed into a formal linear combination of planar partitions by repeated application of Rule~\ref{hrule}\footnotemark, and the resulting linear combination does not depend on the choices made when applying Rule~\ref{hrule}, so that we may write
$$\tau \to \sum_{\text{ planar partitions } \sigma}  \mathcal{P}_{\sigma, \tau}^{(t)} \sigma.$$
For any planar partition $\sigma$, the same coefficients  $\mathcal{P}_{\sigma, \tau}^{(t)}$ satisfy the equation
$$\dddot{ \Pr}(\sigma)  = \dfrac{ \Pr(\sigma) }{ \Pr( 1|2| 3| \cdots |n)} = \sum_{\text{ partitions } \tau} \mathcal{P}_{\sigma, \tau}^{(t)} L_{\tau}$$
for bipartite edge-weighted planar graphs. 
%More generally, for any graph these coefficients satisfy
% $$\sum_{\text{ partitions } \tau}  \mathcal{P}_{\sigma, \tau}^{(t)} \dfrac{ \text{Pr}(\tau) }{ \text{Pr}( 1|2| 3| \cdots |n)} = \sum_{\text{ partitions } \tau} \mathcal{P}_{\sigma, \tau}^{(t)} L_{\tau}$$ 
%The matrix $\mathcal{P}^{(t)}$ is called the {\underline{projection matrix from partitions to planar partitions}}.
\end{thm}
\footnotetext{In \cite[Theorem 1.2]{KW2006}, Rule~\ref{kwrule1} is used in the theorem statement, but for the reasons stated in Remark \ref{rem:rule1error}, we have changed it to Rule~\ref{hrule}.}

%\begin{rem}
%\label{rem:proj=1}
%When $\sigma$ is planar, $\mathcal{P}^{(t)}_{\sigma, \sigma} = 1$. 
%\end{rem}

\subsection{Proof that $\mathcal{Q}^{DD}$ is integer-valued}

\label{sec:firstmajorproof}

We will complete the proof of Theorem~\ref{thm:thm1} by showing that 
we can use the transformation rule introduced in the previous section to compute the columns of the matrix $\mathcal{Q}^{(DD)}$.

\begin{rem}
\label{defn:myrule1}
For pairings, both Rule~\ref{kwrule1} and Rule~\ref{hrule} become the following:
If a pairing $\rho$ is nonplanar, then there will exist items $a < b < c < d$ such that $a$ and $c$ are paired, and $b$ and $d$ are paired. Then the transformation rule is
\begin{equation}
\label{eqn:myrule1}
ac|bd|\text{rest} \to -  ab|cd|\text{rest} - ad|bc|\text{rest}.
\end{equation}
\end{rem}

\begin{trule}
\label{myrule2}
For a black-white pairing $\rho$, repeatedly apply (\ref{eqn:myrule1}) until we have written $\rho$ as a linear combination of planar pairings. Then multiply each planar pairing in this sum by $\sign_{OE}(\sigma)\sign_{BW}(\rho)$.
\end{trule}

The fact that Rule~\ref{myrule2} is well-defined follows from Theorem \ref{thm:KWgrovethm}.
%Let $\widetilde{Q}$ be the matrix obtained by the procedure from Rule~\ref{myrule2}, so the $(\sigma, \rho)$th entry of $\widetilde{\mathcal{Q}}$  is the product of $\sign_{OE}(\sigma)\sign_{BW}(\rho)$ with the coefficient of $\sigma$ when $\rho$ is written as a linear combination of planar pairings using (\ref{eqn:myrule1}). That is,
%\begin{equation}
%\label{eqn:Qtilde}
%\widetilde{\mathcal{Q}}_{\sigma, \rho} = \sign_{OE}(\sigma) \sign_{BW}(\rho) \mathcal{P}_{\sigma, \rho}^{(t)}.
% \end{equation}’’
Proving that Rule~\ref{myrule2} computes the columns of $\mathcal{Q}^{(DD)}$ will prove that the matrix $\mathcal{Q}^{(DD)}$ is integer-valued and gives us the desired theorem, which is stated in full below.

\begin{customthm}{\ref{thm:thm1}}
Any black-white pairing $\rho$ can be transformed into a formal linear combination of planar pairings by repeated application of Rule~\ref{myrule2}, and the resulting linear combination does not depend on the choices we made when applying Rule~\ref{myrule2}, so that we may write
$$\rho \to \sum\limits_{\text{planar pairings } \sigma} \mathcal{Q}^{(DD)}_{\sigma, \rho} \sigma.$$
For any planar pairing $\sigma$, these same coefficients $\mathcal{Q}^{(DD)}_{\sigma, \rho}$ satisfy the equation
$$\widetilde{\Pr}(\sigma) :=
\dfrac{Z^{DD}_{\sigma}(G, {\bf N})}{(Z^D(G))^2}
=
 \sum_{\text{black-white pairings } \rho} \mathcal{Q}^{(DD)}_{\sigma, \rho} Y'_{\rho}.$$
 \end{customthm}

\begin{rem}
The fact that the resulting linear combination does not depend on the choices we made when applying Rule~\ref{myrule2} is an immediate consequence of Theorem~\ref{thm:KWgrovethm}.
\end{rem}

The proof of Theorem~\ref{thm:thm1} requires two additional lemmas.

\begin{lemma}
\label{lem:componentlemma}
Let $\pi$ be a pairing and 
let $\rho$ be a pairing with nodes $a < b < c < d$ that form a crossing in $\rho$. Let $\rho_1$ be the pairing obtained from $\rho$ by replacing the pairs $(a, c)$ and $(b, d)$ with $(a, b)$ and $(c, d)$ and let $\rho_2$ be the pairing obtained from $\rho$ by replacing the pairs
$(a, c)$ and $(b, d)$ with $(a, d)$ and $(b, c)$. 
Then either
\begin{enumerate}
\item[(1)]  $\pi \cup \rho$ has one more component than both $\pi \cup \rho_1$ and $\pi \cup \rho_2$,
%$ \# \text{ comp in }  \pi \cup \rho_1 - \# \text{ comp in } \pi \cup \rho = -1$ and 
%$ \# \text{ comp in }  \pi \cup \rho_2 - \# \text{ comp in } \pi \cup \rho = -1$
\item[(2)] $\pi \cup \rho_1$ has one more component than $\pi \cup \rho$, and $\pi \cup \rho_2$ and $\pi \cup \rho$ have the same number of components, or
%$ \# \text{ comp in }  \pi \cup \rho_1 - \# \text{ comp in } \pi \cup \rho = 1$ and 
%$ \# \text{ comp in }  \pi \cup \rho_2 - \# \text{ comp in } \pi \cup \rho = 0$ 
\item[(3)] $\pi \cup \rho_2$ has one more component than $\pi \cup \rho$, and $\pi \cup \rho_1$ and $\pi \cup \rho$ have the same number of components.
%$ \# \text{ comp in }  \pi \cup \rho_1 - \# \text{ comp in } \pi \cup \rho = 0$ and 
%$ \# \text{ comp in }  \pi \cup \rho_2 - \# \text{ comp in } \pi \cup \rho = 1$
\end{enumerate}
\end{lemma}

\begin{proof}
Observe that either $a, b, c,$ and $d$ are all in the same component of $\pi \cup \rho$ or $a$ and $c$ are in the same component and $b$ and $d$ are in a different component. If $a$ and $c$ are in the same component and $b$ and $d$ are in a different component, then pairing $a$ with $b$ and $c$ with $d$ merges these two components.
% so $ \# \text{ comp in }  \pi \cup \rho_1 - \# \text{ comp in } \pi \cup \rho = -1$. 
Similarly, pairing $a$ with $d$ and $b$ with $c$ merges these two components.
% so $ \# \text{ comp in }  \pi \cup \rho_2 - \# \text{ comp in } \pi \cup \rho = -1$. 

 If $a, b, c$ and $d$ are in the same component, then we consider the following path in $\pi \cup \rho$:
\begin{equation}
\label{eqn:path}
c - a - \pi(a) - \rho(\pi(a)) - \cdots
\end{equation}
This path reaches $b$ or $d$ before it reaches $c$ since by assumption $a, b, c, d$ are all in the same component. If it reaches $b$ before $d$, then in $\rho_1$, $a$ and $b$ are in a different component than $c$ and $d$. This is because path (\ref{eqn:path}) is replaced with
$$b - a  - \pi(a) - \rho(\pi(a)) - \cdots - b,$$
so $\pi \cup \rho_1$ has one more component than $\pi \cup \rho$.
%So $ \# \text{ comp in }  \pi \cup \rho_1 - \# \text{ comp in } \pi \cup \rho = 1$.
 In $\rho_2$, $a, b, c,$ and $d$ are all in the same component, because path (\ref{eqn:path}) is replaced with
$$d - a  - \pi(a) - \rho(\pi(a)) - \cdots - b - c,$$
so $\pi \cup \rho_2$ and $\pi \cup \rho$ have the same number of components.
% $ \# \text{ comp in }  \pi \cup \rho_2 - \# \text{ comp in } \pi \cup \rho = 0$. 
If the path reaches $d$ before $b$, then in $\rho_2$, $a$ and $d$ are in a different component than $b$ and $c$, so 
$\pi \cup \rho_2$ has one more component than $\pi \cup \rho$. 
%$ \# \text{ comp in }  \pi \cup \rho_2 - \# \text{ comp in } \pi \cup \rho = 1$. 
In $\rho_1$, $a, b, c,$ and $d$ are in the same component, so 
 $\pi \cup \rho_1$ and $\pi \cup \rho$ have the same number of components.
%$ \# \text{ comp in }  \pi \cup \rho_1 - \# \text{ comp in } \pi \cup \rho = 0$. 
\end{proof}

\begin{lemma} 
\label{lem:decompintoplanar}
Let $\rho$ be a pairing (not necessarily black-white). Then for any planar pairing $\pi$, 
\begin{equation}
\label{eqn:decompintoplanar}
\sign_{OE}(\pi) (-1)^{ C_{\rho}} (-1)^{\# \text{nodes}/2} 2^{C_{\rho}} =
 \sum_{\text{ planar pairings }  \sigma} \mathcal{P}_{\sigma, \rho}^{(t)} \sign_{OE}( \sigma)2^{C_{\sigma} }.
 \end{equation}
Here, $C_{\rho}$ denotes the number of components in $ \pi \cup \rho$ and  $C_{\sigma}$ denotes the number of components in $ \pi \cup \sigma$.
\end{lemma}

%Lemma~\ref{lem:decompintoplanar} requires one additional lemma.

\begin{proof}
We will prove the claim by induction on the number of crossings in $\rho$.

\noindent {\em Base Case.} When $\rho$ has 0 crossings, equation (\hyperref[eqn:decompintoplanar]{\ref{eqn:decompintoplanar}}) becomes
\[ \sign_{OE}(\pi) (-1)^{ C_{\rho}} (-1)^{\# \text{nodes}/2} 2^{  C_{\rho}} =   \mathcal{P}_{\rho, \rho}^{(t)} \sign_{OE}( \rho)2^{ C_{\rho}},\]
%\[ \sign_{OE}(\pi) (-1)^{\# \text{ comp in } \pi \cup \rho} (-1)^{\# \text{nodes}/2} 2^{\# \text{ comp in } \pi \cup \rho} =   \mathcal{P}_{\rho, \rho}^{(t)} \sign_{OE}( \rho)2^{\# \text{ comp in } \pi \cup \rho}, \]
which %by Remark~\ref{rem:proj=1} 
is equivalent to
\begin{equation}
\label{eqn:decompintoplanarbasecase}
 \sign_{OE}(\pi) (-1)^{C_{\rho}} (-1)^{\# \text{nodes}/2}  \sign_{OE}( \rho) = 1 .
\end{equation}
%\begin{equation}
%\label{eqn:decompintoplanarbasecase}
% \sign_{OE}(\pi) (-1)^{\# \text{ comp in } \pi \cup \rho} (-1)^{\# %\text{nodes}/2}  \sign_{OE}( \rho) = 1 .
%\end{equation}
First suppose $\rho = \pi$. Since $(-1)^{\# \text{ comp in } \pi \cup \pi} =  (-1)^{\# \text{nodes}/2}$, equation~(\ref{eqn:decompintoplanarbasecase}) holds.
We can obtain any planar pairing from any other planar pairing by a sequence of moves, where each move consists of swapping the locations of two nodes of the same parity. So we will show that when $\rho$ is a planar pairing, $x$ and $y$ are two nodes of the same parity, and $\rho'$ is the pairing obtained from $\rho$ by swapping the locations of $x$ and $y$, 
replacing $\rho$ with $\rho'$ does not change
the left hand side of equation~(\ref{eqn:decompintoplanarbasecase}).
%the left hand side of equation~(\ref{eqn:decompintoplanarbasecase}) does not change.
Since $\sign_{OE}(\rho) = - \sign_{OE}(\rho')$, we must show that $(-1)^{\# \text{ comp in } \pi \cup \rho} = - (-1)^{\# \text{ comp in } \pi \cup \rho'}$.

If $x$ and $\rho(x)$ are in a different component than $y$ and $\rho(y)$ in $\pi \cup \rho$, then $\pi \cup \rho'$ has one fewer component than $\pi \cup \rho$. If $x, \rho(x), y,$ and $\rho(y)$ are all in the same component in $\pi \cup \rho$, then without loss of generality assume that $x$ and $y$ are both even, so $\rho(x)$ and $\rho(y)$ are both odd, and 
consider the following path in $\pi \cup \rho$:
$$\rho(x) - x- \pi(x) - \rho(\pi(x)) - \cdots.$$
Since $\rho$ and $\pi$ are both odd-even, segments $n \!\frown\! \rho(n)$ go from an odd node to an even node. Since $\rho(y)$ is odd and $y$ is even, this means that we must reach the node $\rho(y)$ before the node $y$. Therefore we have the path
$$\rho(x) - x - \pi(x) - \rho(\pi(x)) - \cdots -  \rho(y) -  y - \cdots$$
When we replace the pairs $(x, \rho(x)), (y, \rho(y))$ with $(x, \rho(y))$ and $(y, \rho(x))$, this path is replaced with
$$\rho(y) - x - \pi(x) - \rho(\pi(x)) - \cdots - \rho(y)$$
so $(x, \rho(y))$ and $(y, \rho(x))$ are in different components of $\pi \cup \rho'$. 
We conclude that equation (\ref{eqn:decompintoplanarbasecase}) holds for all planar pairings $\rho$.

Now assume that equation (\ref{eqn:decompintoplanar}) holds for pairings $\rho$ with $\leq k$ crossings. 

Let $\rho$ be a pairing with $k+1$ crossings. Let $a < b < c < d$ be nodes that form a crossing in $\rho$. Let $\rho_1$ be the pairing obtained by replacing the pairs $(a, c)$ and $(b, d)$ with $(a, b)$ and $(c, d)$ and let $\rho_2$ be the pairing obtained by replacing the pairs
$(a, c)$ and $(b, d)$ with $(a, d)$ and $(b, c)$. We claim that both $\rho_1$ and $\rho_2$ have fewer than $k+1$ crossings. Observe that if a chord connecting two nodes $n_1$ and $n_2$ crosses the chord connecting $a$ and $b$ in $\rho_1$, it also crosses the chord connecting $a$ and $c$ or the chord connecting $b$ and $d$ in $\rho$. 
Similarly, if a chord connecting two nodes crosses the chord connecting $c$ and $d$ in $\rho_1$, it also crosses the chord connecting $a$ and $c$ or the chord connecting  $b$ and $d$ in $\rho$. It follows that $\rho_1$ has at least one fewer crossing than $\rho$. A similar argument shows that $\rho_2$ has at least one fewer crossing than $\rho$.
%: if a chord connecting two nodes crosses the chord connecting $a$ and $d$ in $\rho_2$, it crosses the chord connecting $b$ and $d$ and/or the chord connecting $a$ and $c$ in $\rho$. If a chord crosses the chord connecting $b$ and $c$ in $\rho_2$, it crosses the chord connecting $b$ and $d$ and/or the chord connecting $a$ and $c$ in $\rho$. 
By the induction hypothesis, 
\[
 \sign_{OE}(\pi) (-1)^{C_{\rho_{1}}} (-1)^{\# \text{nodes}/2} 2^{C_{\rho_{1}}} =
 \sum_{\text{ planar pairings } \sigma} \mathcal{P}_{\sigma, \rho_1}^{(t)} \sign_{OE}( \sigma)2^{C_{\sigma}}
\]
and
\[ \sign_{OE}(\pi) (-1)^{ C_{\rho_{2}}  } (-1)^{\# \text{nodes}/2} 2^{ C_{\rho_{2}}   } =
 \sum_{\text{ planar pairings } \sigma}  \mathcal{P}_{\sigma, \rho_2}^{(t)} \sign_{OE}( \sigma)2^{  C_{\sigma}   }.\]
%By the induction hypothesis, 
%\[ \sign_{OE}(\pi) (-1)^{\# \text{ comp in } \pi \cup \rho_1} (-1)^{\# \text{nodes}/2} 2^{\# \text {comp in } \pi \cup \rho_1} =
% \sum_{\text{ planar pairings } \sigma} \mathcal{P}_{\sigma, \rho_1}^{(t)} \sign_{OE}( \sigma)2^{\# \text{ comp in } \pi \cup \sigma}\]
%and
%\[ \sign_{OE}(\pi) (-1)^{\# \text{ comp in } \pi \cup \rho_2} (-1)^{\# \text{nodes}/2} 2^{\# \text{comp in } \pi \cup \rho_2} =
 %\sum_{\text{ planar pairings } \sigma} \mathcal{P}_{\sigma, \rho_2}^{(t)} \sign_{OE}( \sigma)2^{\# \text{ comp in } \pi \cup \sigma}.\]
By the transformation rule (\ref{eqn:myrule1}),
%By Remark~\ref{defn:myrule1},
$$\mathcal{P}_{\sigma, \rho_1}^{(t)} + \mathcal{P}_{\sigma, \rho_2}^{(t)} = -\mathcal{P}_{\sigma, \rho}^{(t)}$$
so we have
\[
\sum_{\substack{ \text{ planar} \\ \text{pairings } \sigma }} \mathcal{P}_{\sigma, \rho}^{(t)} \sign_{OE}( \sigma)2^{C_{\sigma}  } = - \sign_{OE}(\pi)(-1)^{\# \text{nodes}/2} \left( (-1)^{C_{\rho_{1}}} 2^{C_{\rho_{1}} } + 
(-1)^{ C_{\rho_{2}}    } 2^{ C_{\rho_{2}}  } \right).
\]
%\begin{eqnarray*}
% & & \sum_{\text{ planar pairings } \sigma} \mathcal{P}_{\sigma, \rho}^{(t)} \sign_{OE}( \sigma)2^{\# \text{ comp in } \pi \cup \sigma}\\
%& & = - \sign_{OE}(\pi)(-1)^{\# \text{nodes}/2} \left( (-1)^{\# \text{ comp in } \pi \cup \rho_1} 2^{\# \text{ comp in } \pi \cup \rho_1} + 
%(-1)^{\# \text{ comp in } \pi \cup \rho_2} 2^{\# \text{ comp in } \pi \cup \rho_2} \right).
%\end{eqnarray*}
By Lemma \ref{lem:componentlemma} there are three cases to consider:
\begin{enumerate}
\item[(1)]  $\pi \cup \rho$ has one more component than both $\pi \cup \rho_1$ and $\pi \cup \rho_2$,
%$ \# \text{ comp in }  \pi \cup \rho_1 - \# \text{ comp in } \pi \cup \rho = -1$ and 
%$ \# \text{ comp in }  \pi \cup \rho_2 - \# \text{ comp in } \pi \cup \rho = -1$
\item[(2)] $\pi \cup \rho_1$ has one more component than $\pi \cup \rho$, and $\pi \cup \rho_2$ and $\pi \cup \rho$ have the same number of components, and
%$ \# \text{ comp in }  \pi \cup \rho_1 - \# \text{ comp in } \pi \cup \rho = 1$ and 
%$ \# \text{ comp in }  \pi \cup \rho_2 - \# \text{ comp in } \pi \cup \rho = 0$ 
\item[(3)] $\pi \cup \rho_2$ has one more component than $\pi \cup \rho$, and $\pi \cup \rho_1$ and $\pi \cup \rho$ have the same number of components.
%$ \# \text{ comp in }  \pi \cup \rho_1 - \# \text{ comp in } \pi \cup \rho = 0$ and 
%$ \# \text{ comp in }  \pi \cup \rho_2 - \# \text{ comp in } \pi \cup \rho = 1$
\end{enumerate}

\noindent{\bf Case (1).}
Since $C_{\rho_{i}} - C_{\rho} = -1$ for $i = 1, 2$, 
%$ \# \text{ comp in }  \pi \cup \rho_i - \# \text{ comp in } \pi \cup \rho = -1$ for $i = 1, 2$,
\begin{eqnarray*}
 (-1)^{C_{\rho_{1}} } 2^{ C_{\rho_{1}} } + 
(-1)^{C_{\rho_{2}} } 2^{C_{\rho_{2}} } &= & 
 -(-1)^{  C_{\rho} }\cdot  \frac{1}{2} \cdot 2^{ C_{\rho} } + -(-1)^{ C_{\rho} }\cdot  \frac{1}{2} \cdot 2^{C_{\rho}} \\
 & = & 
 (-1)^{C_{\rho} } 2^{C_{\rho}} \left(-\frac{1}{2} - \frac{1}{2} \right)\\
 & = & -(-1)^{C_{\rho}} 2^{C_{\rho}}.
 \end{eqnarray*}
 %So
%\begin{eqnarray*}
%& & - \sign_{OE}(\pi)(-1)^{\# \text{nodes}/2} \left( (-1)^{\# \text{comp in } \pi \cup %\rho_1} 2^{\# \text{comp in } \pi \cup \rho_1} + 
%(-1)^{\# \text{comp in } \pi \cup \rho_2} 2^{\# \text{comp in } \pi \cup \rho_2} \right) \\
% &=& \sign_{OE}(\pi)(-1)^{\# \text{nodes}/2}(-1)^{\# \text{comp in } \pi \cup \rho} 2^{\# \text{comp in } \pi \cup \rho} 
% \end{eqnarray*}
 
\noindent{\bf Cases (2) and (3).}
We will only include the proof for case (2), since case (3) is completely analogous. 
Since $C_{\rho_{1}} - C_{\rho} = 1$ and
%Since $ \# \text{ comp in }  \pi \cup \rho_1 - \# \text{ comp in } \pi \cup \rho = 1$ and 
 $C_{\rho_{2}} - C_{\rho} = 0$, 
%$ \# \text{ comp in }  \pi \cup \rho_2 - \# \text{ comp in } \pi \cup \rho = 0$, 
\begin{eqnarray*}
  (-1)^{ C_{\rho_{1}} } 2^{ C_{\rho_{1}} } + 
(-1)^{ C_{\rho_{2}}} 2^{ C_{\rho_{2}} }
 & = & 
 -(-1)^{ C_{\rho} } 2 \cdot 2^{ C_{\rho}} + (-1)^{ C_{ \rho}} 2^{ C_{ \rho}} \\
 & = & 
 (-1)^{C_{ \rho} } 2^{ C_{\rho}} (-2 + 1)\\
 & = & -(-1)^{ C_{ \rho}} 2^{ C_{ \rho}} .
 \end{eqnarray*}
So in all cases, 
\begin{eqnarray*}
 - \sign_{OE}(\pi)(-1)^{\# \text{nodes}/2} 
\left( (-1)^{ C_{ \rho_{1}}   }    2^{ C_{ \rho_{1}}  } + 
(-1)^{ C_{ \rho_{2}}    } 2^{C_{ \rho_{2}}  } \right) 
 = \sign_{OE}(\pi)(-1)^{\# \text{nodes}/2}
 (-1)^{C_{ \rho}} 
 2^{C_{ \rho}},
 \end{eqnarray*}
and thus
$$\sign_{OE}(\pi) (-1)^{ C_{\rho}   } (-1)^{\# \text{nodes}/2} 2^{ C_{\rho} } =
 \sum_{\text{ planar pairings } \sigma} \mathcal{P}_{\sigma, \rho}^{(t)} \sign_{OE}( \sigma)2^{ C_{\sigma} }.$$

\end{proof}

\begin{proof}[Proof of Theorem \ref{thm:thm1}]
Let $\widetilde{\mathcal{Q}}$ be the matrix obtained by the procedure from Rule~\ref{myrule2}, so the $(\sigma, \rho)$th entry of $\widetilde{\mathcal{Q}}$  is the product of $\sign_{OE}(\sigma)\sign_{BW}(\rho)$ with the coefficient of $\sigma$ when $\rho$ is written as a linear combination of planar pairings using (\ref{eqn:myrule1}). That is,
\begin{equation}
\label{eqn:Qtilde}
\widetilde{\mathcal{Q}}_{\sigma, \rho} = \sign_{OE}(\sigma) \sign_{BW}(\rho) \mathcal{P}_{\sigma, \rho}^{(t)}.
 \end{equation}

We will show that 
$$\mathcal{M}_{2} \widetilde{\mathcal{Q}} {\bf e}_{i} =  \mathcal{B}_{2} {\bf e}_{i}$$
for all $i$. This will show that $\mathcal{M}_{2} \widetilde{\mathcal{Q}} = \mathcal{M}_{2} \mathcal{Q}^{(DD)}$, which proves the theorem since $\mathcal{M}_{2}$ is invertible.

Let $\rho$ be a black-white pairing. Recall from Definition
\hyperref[Bdefn]{\ref{Bdefn}}
that $(\mathcal{B}_2)_{\pi, \rho} = \sign(\pi, \rho) 2^{C_{ \rho}  }$.
% By assumption, 
% $(\widetilde{\mathcal{Q}})_{\pi, \rho} = \sign_{OE}(\pi) \sign_{BW}(\rho) \mathcal{P}_{\pi, \rho}^{(t)}$. 
Then by equation (\ref{eqn:Qtilde}), to show that $\mathcal{M}_{2} \widetilde{\mathcal{Q}} {\bf e}_{i} =  \mathcal{B}_{2} {\bf e}_{i}$,
%by equation (\ref{eqn:Qtilde})
we need to show that for each planar pairing $\pi$, 
$$\sign(\pi, \rho) 2^{C_{\rho}} =
 \sum_{\text{ planar pairings } \sigma} \mathcal{P}_{\sigma, \rho}^{(t)} \sign_{OE}( \sigma) \sign_{BW}( \rho) 2^{C_{\sigma}}. $$
By Definition \ref{signrhopi},
 $$\sign(\pi, \rho) = (-1)^{\# \text{nodes}/2} (-1)^{C_{\rho}} \sign_{OE}(\pi) \sign_{BW} (\rho).$$
Applying Lemma~\ref{lem:decompintoplanar} completes the proof. 
\end{proof}
% it suffices to prove the following. 

%This proves that
%\begin{equation}
%\label{eqn:intval}
%(\mathcal{Q}^{(DD)})_{\pi, \rho} = \sign_{OE}(\pi) \sign_{BW}(\rho) \mathcal{P}^{(t)}_{\pi, \rho},
%\end{equation}
%which completes the proof of
%Theorem \ref{thm:thm1}. Its full statement is below. \\

%\begin{thm}(analogue of Theorem 1.4 from \cite{KW2006}). 

%\end{thm}

\subsection{Another characterization of $\sign(S)$}

\label{sec:signSformula}

In this section, we prove Lemma \ref{lem:signSdefn2}, which was key in establishing Lemma~\ref{lem:mylem35}.\\

%Let $T \subseteq {\bf N}$ be the set of nodes that are odd and white or even and black.\\

\noindent {\bf Lemma \ref{lem:signSdefn2}.} Let $S$ be a balanced subset of nodes and let $\sign(S)$ be defined as in Lemma \ref{lemma32gen}. Then
\begin{equation}
\label{eqn:signSdefn2}
 \text{sign}(S) = (-1)^{\# \text{nodes}/2} (-1)^{\# \text{ comp in } \pi \cup \rho} \sign_{OE}(\pi) \sign_{BW}(\rho)
\end{equation}
where $\pi$ is an odd-even pairing such that $\pi$ does not connect $S \triangle T$ to $(S \triangle T)^c$ and $\rho$ is a black-white pairing such that 
$\rho$ does not connect $S$ to $S^c$. \\

%\begin{note*}
%I have checked that this is equal to the definition of sign(S) in Lemma \ref{lemma32gen} for all node colorings on 6, 8, and 10 nodes where node 1 is black and node $2n$ is white. 
%\end{note*}

Proving Lemma~\ref{lem:signSdefn2} requires
\begin{itemize}
\item[(1)] proving that such pairings $\pi$ and $\rho$ always exist,
\item[(2)] proving that equation (\ref{eqn:signSdefn2}) is well-defined, and
% if there are two pairs of pairings $(\pi_1, \rho_1)$ and $(\pi_2, \rho_2)$ such that $\pi_i$ does not bridge $S \triangle T$ to $(S \triangle T)^c$ and $\rho_i$ does not bridge $S$ to $S^c$, that the right hand side of equation (\ref{eqn:signSdefn2}) is the same, and
\item[(3)] proving that equation (\ref{eqn:signSdefn2}) holds. 
\end{itemize}

We will postpone the proof of (1) because the fact that such pairings $\pi$ and $\rho$ always exist will follow quickly from the proofs of (2) and (3). 

\subsubsection{Proof that equation (\ref{eqn:signSdefn2}) is well-defined}

The strategy of the proof is to define local moves that allow us to get from a pair $(\pi_1, \rho_1)$ such that $\pi_1$ does not connect $S \triangle T$ to $(S \triangle T)^c$ and $\rho_1$ does not connect $S$ to $S^c$
to any other pair $(\pi_2, \rho_2)$ with this property, and to show that these moves do not change the right hand side of equation~(\ref{eqn:signSdefn2}). 

Specifically, we will define two types of local moves. First, we define moves that modify $\pi$ by swapping the locations of two nodes of the same parity under certain conditions but leave $\rho$ fixed, called moves of type $A_{OE}$. Next, we define moves that modify $\rho$ by swapping the locations of two nodes of the same color under similar conditions but leave $\pi$ fixed, called moves of type $A_{BW}$.

In order to describe the conditions under which we can swap the locations of two nodes, we need the following definition.

\begin{defn}
We call a pair of nodes $(a, \eta(a))$ a {\em transition pair} if exactly one of the nodes $a, \eta(a)$ is in $T$. 
\end{defn}

%\begin{rem}
%If $S$ is a set such that $\pi$ does not bridge $S \triangle T$ to $(S \triangle T)^c$
%and $\rho$ does not bridge $S$ to $S^c$ and $(a, \pi(a))$ is a transition pair, then exactly one of the nodes $a, \pi(a)$ is in $S$. 
%\end{rem}

\begin{rem}
If $a$ and $b$ are two nodes in the same component of $\pi \cup \rho$, there are two paths from $a$ to $b$. Since the algorithm in Lemma~\ref{alglemma} is well-defined, the parity of the number of transition pairs is independent of the path. 
% number of transitions pairs in each path has the same parity.
%the parity of the number of transition pairs in the different paths is the same. 
\end{rem}

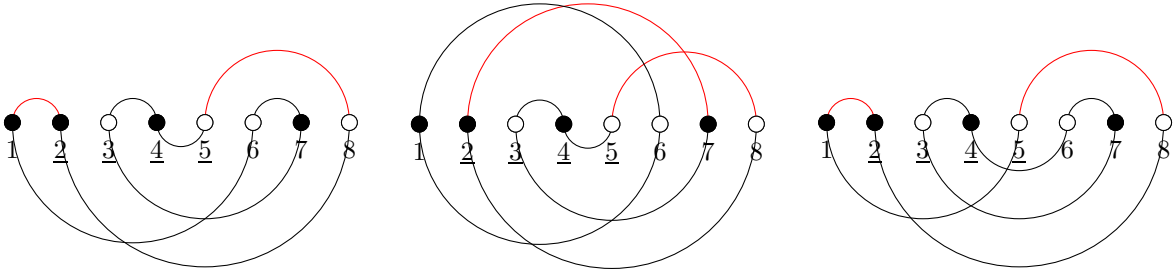
\begin{figure}[h]
\centering
\begin{minipage}{.32\textwidth}
\vspace{.6cm}
\[\begin{tikzpicture}[scale = .32]

	\vertex[fill] (n1) at (4, 0) [label=below:\small{$1$}] {};
	\vertex[fill] (n2) at (6,0) [label=below:\small{${\underline{ 2}}$}] {};
	\vertex (n3) at (8,0) [label=below:\small{$\underline{3}$}] {};
	\vertex[fill] (n4) at (10, 0) [label=below:\small{$\underline{4}$}] {};
	\vertex (n5) at (12,0) [label=below:\small{$\underline{5}$}] {};
	\vertex (n6) at (14,0) [label=below:\small{$6$}] {};
	\vertex[fill] (n7) at (16,0) [label=below:\small{$7$}] {};
	\vertex (n8) at (18, 0) [label=below:\small{$8$}] {};

	%arcs of pi
	\draw[color = red]  (n8) arc (0:180:3cm);
	\draw  (n7) arc (0:180:1cm);
	\draw[color = red]  (n2) arc (0:180:1cm);
	\draw  (n4) arc (0:180:1cm);
	
	%arcs of rho
	\draw  (n1) arc (180:360:5cm);
	\draw  (n2) arc (180:360:6cm);
	\draw  (n4) arc (180:360:1cm);
	\draw  (n3) arc (180:360:4cm);

	%So arcs are below vertices
	\vertex[fill]  at (4, 0) {};
	\vertex[fill] at (6,0) {};
	\vertex[fill = white] at (8,0) {};
	\vertex[fill] at (10, 0) {};
	\vertex[fill = white] at (12,0) {};
	\vertex[fill = white] at (14,0)  {};
	\vertex[fill]  at (16,0)  {};
	\vertex[fill = white] at (18, 0) {};

\end{tikzpicture}\]
\end{minipage}
\begin{minipage}{.32\textwidth}
\[\begin{tikzpicture}[scale = .32]
	\vertex[fill] (n1) at (4, 0) [label=below:\small{$1$}] {};
	\vertex[fill] (n2) at (6,0) [label=below:\small{$\underline{2}$}] {};
	\vertex (n3) at (8,0) [label=below:\small{$\underline{3}$}] {};
	\vertex[fill] (n4) at (10, 0) [label=below:\small{$\underline{4}$}] {};
	\vertex (n5) at (12,0) [label=below:\small{$\underline{5}$}] {};
	\vertex (n6) at (14,0) [label=below:\small{$6$}] {};
	\vertex[fill] (n7) at (16,0) [label=below:\small{$7$}] {};
	\vertex (n8) at (18, 0) [label=below:\small{$8$}] {};
	
	%arcs of pi
	\draw[color = red]  (n8) arc (0:180:3cm);
	\draw[color = red]  (n7) arc (0:180:5cm);
	\draw  (n6) arc (0:180:5cm);
	\draw  (n4) arc (0:180:1cm);
	
	%arcs of rho
	\draw  (n1) arc (180:360:5cm);
	\draw  (n2) arc (180:360:6cm);
	\draw  (n4) arc (180:360:1cm);
	\draw  (n3) arc (180:360:4cm);
	
	%So arcs are below vertices
	\vertex[fill]  at (4, 0) {};
	\vertex[fill] at (6,0) {};
	\vertex[fill = white] at (8,0) {};
	\vertex[fill] at (10, 0) {};
	\vertex[fill = white] at (12,0) {};
	\vertex[fill = white] at (14,0)  {};
	\vertex[fill]  at (16,0)  {};
	\vertex[fill = white] at (18, 0) {};

\end{tikzpicture}\]
\end{minipage}
\begin{minipage}{.32\textwidth}
\vspace{.6cm}
\[\begin{tikzpicture}[scale = .32]
	
	\vertex[fill] (n1) at (4, 0) [label=below:\small{$1$}] {};
	\vertex[fill] (n2) at (6,0) [label=below:\small{${\underline{ 2}}$}] {};
	\vertex (n3) at (8,0) [label=below:\small{$\underline{3}$}] {};
	\vertex[fill] (n4) at (10, 0) [label=below:\small{$\underline{4}$}] {};
	\vertex (n5) at (12,0) [label=below:\small{$\underline{5}$}] {};
	\vertex (n6) at (14,0) [label=below:\small{$6$}] {};
	\vertex[fill] (n7) at (16,0) [label=below:\small{$7$}] {};
	\vertex (n8) at (18, 0) [label=below:\small{$8$}] {};

	%arcs of pi
	\draw[color = red]  (n8) arc (0:180:3cm);
	\draw  (n7) arc (0:180:1cm);
	\draw[color = red]  (n2) arc (0:180:1cm);
	\draw  (n4) arc (0:180:1cm);
	
	%arcs of rho
	\draw  (n1) arc (180:360:4cm);
	\draw  (n2) arc (180:360:6cm);
	\draw  (n4) arc (180:360:2cm);
	\draw  (n3) arc (180:360:4cm);

	%So arcs are below vertices
	\vertex[fill]  at (4, 0) {};
	\vertex[fill] at (6,0) {};
	\vertex[fill = white] at (8,0) {};
	\vertex[fill] at (10, 0) {};
	\vertex[fill = white] at (12,0) {};
	\vertex[fill = white] at (14,0)  {};
	\vertex[fill]  at (16,0)  {};
	\vertex[fill = white] at (18, 0) {};

\end{tikzpicture}\]
\end{minipage}
\caption{{\em Left}: The diagram of $\pi \cup \rho$, where $\pi = ((1, 2), (3, 4), (5, 8), (7, 6))$ and $\rho = ((1, 6), (2, 8), (4, 5), (7,3))$. Nodes that are in $T$ are underlined and arcs between transition pairs are red. 
{\em Center}: Since $1$ and $7$ are two nodes of the same parity and color and a path from $1$ to $7$ contains an even number of transition pairs, if we let $\pi' = ((1, 6), (3, 4), (5, 8), (7, 2))$  then $(\pi, \rho)$ and $(\pi', \rho)$ differ by a move of type $A_{OE}$. {\em Right}: Since $1$ and $4$ are two nodes of the same color and a path from $1$ to $4$ contains an even number of transition pairs, if we let $\rho' = ((1, 5), (2, 8), (4, 6), (7, 3))$  then $(\pi, \rho)$ and $(\pi', \rho)$ differ by a move of type $A_{BW}$.}
\label{fig:typeAOE}
\end{figure}

\begin{defn}
\label{defn:AOE}
Suppose $\pi$ is an odd-even pairing and $\rho$ is a black-white pairing. Let $a$ and $b$ be two nodes of the same parity. If
\begin{itemize}
\item $a$ and $b$ are in different components,
\item $a$ and $b$ are the same color and a path from $a$ to $b$ contains an even number of transition pairs, or
\item $a$ and $b$ are different colors and a path from $a$ to $b$ contains
 an odd number of transition pairs, 
\end{itemize}
let $\pi'$ be the pairing obtained from $\pi$ by swapping the locations of $a$ and $b$ in $\pi$. We say that $(\pi', \rho)$ and $(\pi, \rho)$ differ by a move of type $A_{OE}$. See Figure \ref{fig:typeAOE} for an example. 
\end{defn}

\begin{defn}
\label{defn:ABW}
Let $\pi$ be an odd-even pairing and let $\rho$ be a black-white pairing. 
Suppose $a$ and $b$ are the same color 
and either $a$ and $b$ are in different components, or a path in $\pi \cup \rho$ from $a$ to $b$ contains an even number of transition pairs.

Suppose we swap the locations of $a$ and $b$ in $\rho$ to obtain the pairing $\rho'$. 
Then
we say that $(\pi, \rho')$ and $(\pi, \rho)$ differ by a {\em move of type $A_{BW}$}.  See Figure \ref{fig:typeAOE} for an example. 
\end{defn}

\begin{lemma}
\label{AOElemma}
Let $\pi, \pi'$ be odd-even pairings and let $\rho$ be a black-white pairing such that $( \pi, \rho)$ and $(\pi', \rho)$ differ by a move of type $A_{OE}$. Then the number of components in $\pi \cup \rho$ and the number of components in $\pi' \cup \rho$ differ by one.  
%Suppose $\pi$ is an odd-even pairing and $\rho$ is a black-white pairing. 
%Let $a$ and $b$ be two nodes of the same parity. Let $\pi'$ be the pairing obtained from $\pi$ by swapping the locations of $a$ and $b$ in $\pi$. If
%\begin{itemize}
%\item[(1)] $a$ and $b$ are in different components, 
%\item[(2)] $a$ and $b$ are the same color and a path from $a$ to $b$ contains an even %number of transition pairs, or
%\item[(3)] $a$ and $b$ are different colors and a path from $a$ to $b$ contains an n odd number of transition pairs, 
%\end{itemize}
%then the number of components in $\pi \cup \rho$ and the number of components in $\pi' \cup \rho$ differ by one.  
\end{lemma}

\begin{proof}
%This proof is similar to the proof of Lemma \ref{ABWlemma}. 
%Consider a move of type $A_{OE}$ that swaps the locations of two nodes $a$ and $b$ of the same parity in $\pi$. 
If $a$ and $b$ are in different components of $\pi \cup \rho$, swapping the locations of $a$ and $b$ in $\pi$ merges these two components, 
so the number of components decreases by one. 

If $a$ and $b$ are in the same component, 
without loss of generality assume that node $a$ is white. 
%If $a$ and $b$ are in the same component of $\pi \cup \rho$, as in the proof of Lemma \ref{ABWlemma} we again rely on the following observation. 
Consider the following path from $a$ to $b$, which starts by traversing the edge connecting $a$ to $\pi(a)$:
$$a - \pi(a) - \cdots - b.$$
 We claim that
we always reach $b$ before $\pi(b)$. 
This follows from the observation that because $\rho$ is black-white and $\pi$ is odd-even,
a path in $\pi \cup \rho$ alternates between black
and white nodes unless a pair $(d, \pi(d))$ in the path is a transition pair. 
So since our path starts at a white node by traversing the edge in $\pi$,
if we consider an edge $d \!\frown\!  \pi(d)$ of the path,  $d$ is white and $\pi(d)$ is black if and only if we traverse this edge
after passing through an even number 
of transition pairs.
So, if we were to reach $\pi(b)$ before $b$, $b$ is black if and only if there are an even number
of transition pairs between $a$ and $b$, a contradiction since $a$ is white.
It follows that we must reach $b$ before $\pi(b)$. 

Thus we have the following path in $\pi \cup \rho$: 
$$a - \pi(a) - \cdots - b - \pi(b)$$
When we replace the pairs $(a, \pi(a))$ and $(b, \pi(b))$ in $\pi$ with $(a, \pi(b))$ and $(b, \pi(a))$ to obtain $\pi'$
%a path in $\pi \cup \rho$ goes from $a$ to $\rho(b)$, and 
the middle portion of the path above $\pi(a) - \cdots - b$ becomes a new component, so the number of components increases by one. 
\end{proof}

\begin{cor}
\label{AOEcor}
A move of type $A_{OE}$ does not change the right hand side of equation (\ref{eqn:signSdefn2}). 
\end{cor}

\begin{proof}
If $(\pi, \rho)$ and $(\pi', \rho)$ differ by a move of type $A_{OE}$, then
 $(-1)^{\# \text{ comp in } \pi \cup \rho}  = 
 -(-1)^{\# \text{ comp in } \pi' \cup \rho}$ by Lemma~\ref{AOElemma}
 and $\sign_{OE}(\pi) = -\sign_{OE}(\pi')$, so replacing $\pi$ with $\pi'$ does not change the right hand side of equation (\ref{eqn:signSdefn2}). 
%By Lemma \ref{AOElemma}, a move of type $A_{OE}$ changes the number of components in $\pi \cup \rho$ by 1. This move changes $\sign_{OE}(\pi)$, but does not change $\sign_{BW}(\rho)$.
\end{proof}

\begin{cor}
\label{ABWcor}
A move of type $A_{BW}$ does not change the right hand side of equation (\ref{eqn:signSdefn2}). 
\end{cor}

\begin{proof}
The proof that a move of type $A_{BW}$ changes the number of components in $\pi \cup \rho$ by one is analogous to the proof of Lemma \ref{AOElemma}. The claim follows as it did in the proof of Corollary~\ref{AOEcor}.
%An analogous proof to the proof of Lemma \ref{AOElemma} shows that a move of type $A_{BW}$ changes the number of components in $\pi \cup \rho$ by one. 
%This move changes $\sign_{OE}(\pi)$, and does not change $\sign_{BW}(\rho)$. 
\end{proof}

\begin{proof}[Proof that equation (\ref{eqn:signSdefn2}) is well-defined]
By Corollaries \ref{AOEcor} and \ref{ABWcor}, moves of type $A_{OE}$ and type $A_{BW}$ do not change the right hand side of equation (\ref{eqn:signSdefn2}). So to prove that the formula for $\sign(S)$ is well-defined, it suffices to show that these two types of moves are enough to get from a pair $(\pi_1, \rho_1)$ such that $\pi_1$ does not connect $S \triangle T$ to $(S \triangle T)^c$ and $\rho_1$ does not connect $S$ to $S^c$ to any other pair $(\pi_2, \rho_2)$ with this property.

We can get from any pairing of nodes in $S$ to any other pairing of nodes in $S$ using moves of type $A_{BW}$ because type $A_{BW}$ moves allow us to exchange any nodes of the same color in $S$. By the same reasoning, we can get from any pairing of nodes in $S^c$ to any other pairing of nodes in $S^c$. 
So, if $\rho$ and $\rho'$ are two pairings that both do not connect $S$ to $S^c$, then we can get from $\rho$ to $\rho'$ using a sequence of moves of type $A_{BW}$.

Similarly, we can get from any odd-even pairing of nodes in $S \triangle T$ to any other odd-even pairing of nodes in $S \triangle T$ by swapping nodes of the same parity in $S \triangle T$. We can also get from any odd-even pairing of nodes in $(S \triangle T)^c$ to any other odd-even pairing of nodes in $(S \triangle T)^c$. So, if $\pi$ and $\pi'$ are two odd-even pairings that both do not connect $S \triangle T$ to $(S \triangle T)^c$, then we can get from $\pi$ to $\pi'$ using a sequence of moves of type $A_{OE}$.

We have thus shown if we have two pairs
of pairings $(\pi_1, \rho_1)$ and $(\pi_2, \rho_2)$ such that $\pi_i$ is odd-even and does not connect $S \triangle T$ to $(S \triangle T)^c$ and $\rho_i$ is black-white and does not connect
$S$ to $S^c$, that the right hand side of equation (\ref{eqn:signSdefn2}) is unchanged when we replace $(\pi_1, \rho_1)$ with $(\pi_2, \rho_2)$. 

\end{proof}

\subsubsection{Proof that equation (\ref{eqn:signSdefn2}) holds}

\label{sec:eqnholds}

First assume that $S$ is a balanced set of size $2j$ such that there is a planar black-white pairing $\rho$ that does not connect $S$ to $S^c$. Although it may not be obvious that such a set always exists, recall from Lemma \ref{firstlemma34} that regardless of the node coloring of {\bf N}, there exists a planar black-white pairing $\rho$ of ${\bf N}$. So we choose $S$ to be $2j$ of the arcs of $\rho$. 

Then by definition,
$$\sign(S) = (-1)^{\# \text{ crosses of } \rho} = 1.$$
Let
 $\pi = \rho$. Since $\pi$ is odd-even and black-white, for all pairs in $\pi$, either both nodes of the pair are in $T$ or both are not in $T$, so $\pi$ does not connect $S \triangle T$ to $(S \triangle T)^c$. 
Since $\pi = \rho$, $(-1)^{\# \text{nodes}/2}= (-1)^{\# \text{ comp in } \pi \cup \rho}$. Also,
$\sign_{OE}(\pi) = \sign_{BW}(\rho)$ by Lemma~\ref{lem:OEandBWsignsold}, so equation (\ref{eqn:signSdefn2}) holds.

We can obtain any balanced set of size $2j$ from $S$ by making a sequence of the following types of replacements: 
\begin{itemize}
\item[(1)] Replace $x \in S$ with $x+1 \in S^c$, where $(x, x+1)$ is a couple of consecutive nodes of the same color. (Or replace $x+1 \in S$ with $x \in S^c$). 
%\item[(2)] Replacing $s$ with $s+1$, where $(s, s+1)$ is a couple of consecutive black nodes.
\item[(2)] Replace $x \in S$ with $y \in S^c$, where $x < y$ are the same color and all $\ell$ nodes in the interval $[x+1, x+2, \ldots, y-1]$ are the opposite color of $x$ and $y$ ($\ell \geq 1$). (Or replace $y \in S$ with $x \in S^c$). 
\end{itemize}

Therefore it suffices to show the following. Assume we're given a balanced set $S$, an odd-even pairing $\pi$ that does not connect $S \triangle T$ to $(S \triangle T)^c$, and 
a black-white pairing $\rho$ that does not connect $S$ to $S^c$
such that $\rho |_{S}$ and $\rho |_{S^c}$ are planar.
% and is planar when restricted to $S$ and when restricted to $S^c$.
After making either of the above two types of replacements to obtain $S'$, we can construct
an odd-even pairing $\pi'$ that does not connect $S' \triangle T$ to $(S' \triangle T)^c$ and a black-white pairing
 $\rho'$ that does not connect $S'$ to $S'^c$
 such that $\rho' |_{S'}$ and $\rho' |_{S'^c}$ are planar.
 % is planar when restricted to $S'$ and when restricted to $S'^c$. 
  After replacing $S$, $\pi$, $\rho$ in equation (\ref{eqn:signSdefn2}) with $S'$, $\pi'$ and $\rho'$, equation (\ref{eqn:signSdefn2}) still holds.

%We will show that equation (\ref{eqn:signSdefn2}) holds after applying each type of replacement to obtain $S'$. More precisely, given $S$, $\rho$ that does not bridge $S$ to $S^c$ and $\pi$ that does not bridge $S \triangle T$ to $(S \triangle T)^c$, and $S'$ meaning that we will construct $\pi'$ and $\rho'$ 
This requires several lemmas.

\begin{lemma}
\label{lem:pix=y}
Let $S$ be a balanced subset of nodes. 
Let $x$ and $y$ be two nodes of the same color and opposite parity with $x < y$ such that $x \in S$ and $y \notin S$. 
Let
$\rho$ be a black-white pairing such that $\rho$ does not connect $S$ to $S^c$ and let $\pi$ be an odd-even pairing such that $\pi$ does not connect $S \triangle T$ to $(S \triangle T)^c$. 
Let $S' = S \setminus \{x\} \cup \{y\}$ and let $\rho'$ be the pairing obtained by swapping the locations of $x$ and $y$ in $\rho$. Then
%Assume Hypotheses~\ref{hyp} and that $\pi(x) = y$. 
\begin{itemize}
\item[(a)] if $\pi(x) = y$,
%$\rho'$ does not bridge $S'$ to $S'^c$ and 
\begin{itemize}
\item[(i)]$\pi$ does not connect $S' \triangle T$ to $(S' \triangle T)^c$, and 
\item[(ii)] when $\rho$ is replaced with $\rho'$, the right hand side of equation (\ref{eqn:signSdefn2}) changes sign. 
\end{itemize}
\item[(b)] if $\pi(x) \neq y$, let $\pi'$ be the pairing obtained from $\pi$ by pairing $x$ with $y$, $\pi(x)$ with $\pi(y)$, and leaving the remaining pairs the same. Then
\begin{itemize}
\item[(i)] $\pi'$ does not connect $S' \triangle T$ to $(S' \triangle T)^c$. 
\item[(ii)] when $\rho$ is replaced with $\rho'$ and $\pi$ is replaced with $\pi'$, the right hand side of equation (\ref{eqn:signSdefn2}) changes sign. 
\end{itemize}
\end{itemize}
\end{lemma}

\begin{proof}
We will first prove part (a).
%By Lemma \ref{lem:swapab}, $\rho'$ does not bridge $S'$ to $S'^c$. 
The fact that $\pi$ does not connect $S' \triangle T$ to $(S' \triangle T)^c$ follows from the observation that since $\pi(x) = y$, both $x$ and $y$ are in $S \triangle T$ or both are in $(S \triangle T)^c$. If both $x, y$ are in $S \triangle T$ then since we assumed $x \in S$ and $y \notin S$, $y$ must be in $T$, so both $x, y$ are in $(S' \triangle T)^c$. 
So $\pi$ does not connect $S' \triangle T$ to $(S' \triangle T)^c$. 
%the difference between $S \triangle T$ and $S' \triangle T$ is that one of $S \triangle T$, $S' \triangle T$ contains $u$ and $u+1$ and the other does not. 

Since we obtained $\rho'$ from $\rho$ by swapping the locations of $x$ and $y$, 
 $\sign_{BW}(\rho') = - \sign_{BW}(\rho)$.
The number of components in $\pi \cup \rho$ is the same as the number of components in $\pi \cup \rho'$ because when we replace $\rho$ with $\rho'$ the path $\pi(\rho(x)) - \rho(x) - x - y- \rho(y)$ is replaced with 
$\pi(\rho(x)) - \rho(x) - y - x- \rho(y)$.
%$\pi(\rho(u)) - \rho(u+1) - u - u+1 - \rho(u)$.
% (Note that $\pi(\rho(x))$ could be equal to $\rho(y)$). 
 So the right hand side of equation (\ref{eqn:signSdefn2}) changes sign.

Next, we prove part (b). 
The proof of (i) relies on the observation that since $x$ and $y$ are the same color but opposite parity, exactly one of the nodes $x, y$ is in $T$. 
This implies that 
$x$ and $y$ are both in $S \triangle T$ or both in $(S \triangle T)^c$ and that
$x$ and $y$ are both in $S' \triangle T$ or both in $(S' \triangle T)^c$.

Since $x$ and $y$ are both in $S \triangle T$ or both in $(S \triangle T)^c$, $\pi(x)$ and $\pi(y)$ are both in $S \triangle T$ or both in $(S \triangle T)^c$. Since neither $\pi(x)$ nor $\pi(y)$ is $x$ or $y$, 
$\pi(x)$ and $\pi(y)$ are both in $S' \triangle T$ or both in $(S' \triangle T)^c$. We conclude that $\pi'$ does not connect $S' \triangle T$ to $(S' \triangle T)^c$.

For the proof of (ii), first note that pairing $x$ and $y$ and $\pi(x)$ with $\pi(y)$ is the same as swapping the locations of $y$ and $\pi(x)$. It follows that $\sign_{OE}(\pi') = - \sign_{OE}(\pi)$, and since $\sign_{BW}(\rho') = - \sign_{BW}(\rho)$, it remains to show that the number of components in $\pi' \cup \rho'$ and the number of components in $\pi \cup \rho$ differ by 1.

By letting $a = \pi(x)$ and $b = y$ in Definition~\ref{defn:AOE}, 
we see that $(\pi', \rho)$ and $(\pi, \rho)$ differ by a move of type $A_{OE}$. 
%We use Lemma \ref{AOElemma} with $a = \pi(x)$ and $b = y$ to show that the number of components in $\pi' \cup \rho$ and the number of components in $\pi \cup \rho$ differ by one.
If $\pi(x)$ and $y$ are in different components, this is clear, since
%we can clearly apply the lemma, 
%since 
$\pi(x)$ and $y$ have the same parity.
If $\pi(x)$ and $y$ are in the same component,
% it is not immediately obvious that we can apply the lemma, so we show
% It's not immediately obvious that we can apply this lemma.
%so we first note that $\pi(x)$ and $y$ have the same parity. 
 %So it remains to show 
% that in this case, 
we must show that
 they are the same color if and only if there are an even number of transition pairs between them.
%If $\pi(x)$ and $y$ are in different components, we can clearly apply the lemma, since they have the same parity.
%If $\pi(x)$ and $y$ are in the same component, it is not immediately obvious that we can apply the lemma, so we show
% It's not immediately obvious that we can apply this lemma.
%so we first note that $\pi(x)$ and $y$ have the same parity. 
 %So it remains to show 
 %that in this case, they are the same color if and only if there is an even number of transition pairs between them. 
 This is because
\begin{itemize}
\item $y$ and $x$ are the same color
\item a path from $y$ to $x$ contains an odd number of transition pairs (since $x \in S$ and $y \notin S$)
\item $x$ and $\pi(x)$ are the same color if and only if $(x, \pi(x))$ is a transition pair
\end{itemize}
So, 
by Lemma \ref{AOElemma},
% to show that 
the number of components in $\pi' \cup \rho$ and the number of components in $\pi \cup \rho$ differ by one. Then, since $\pi'(x) = y$, by the proof of part (a), the number of components in $\pi' \cup \rho'$ is the same as the number of components in $\pi' \cup \rho$.

%To prove (b), note that $\sign_{BW}(\rho') = - \sign_{BW}(\rho)$. 
%and $\sign_{OE}(\pi') = - \sign_{OE}(\pi)$ because the two even nodes from the nodes $\{x, y, \pi(x), %\pi(y)\}$ have been swapped. 
%It remains to show that the number of components in $\pi' \cup \rho'$ and the number of components in $\pi \cup \rho$ differ by 1. 

We conclude that when $\rho$ is replaced with $\rho'$ and $\pi$ is replaced with $\pi'$ the right hand side of equation (\ref{eqn:signSdefn2}) changes sign. 
\end{proof}

\begin{lemma} 
\label{lem:planarwhenrestricted}
Let $S \subseteq {\bf N}$ be a balanced set. Let $x, y$ be nodes of the same color such that $x \in S$, $y \in S^c$, $x < y$ and all 
$\ell$ nodes
in the interval $[x+1, x+2, \ldots, y-1]$
 %appearing between $x$ and $y$ 
 are the opposite color of $x$ and $y$ ($\ell \geq 1$). 

Let $\rho$ be a black-white pairing such that $\rho$ does not connect $S$ to $S^c$ and $\rho |_{S}$ and $\rho |_{S^c}$ are planar.
%$\rho$ is planar when restricted to $S$, and planar when restricted to $S^c$.

\begin{itemize}
\item[(1)] If $\rho(x)$ is not in the interval $[x+1, \ldots, y-1]$ and there is a node in this interval that is in $S$, let $k$ be the smallest integer such that $x+k$ is in $S$ and let $\rho'$ be the pairing obtained from $\rho$ by replacing the pairs $(x, \rho(x))$ and $(x+k, \rho(x+k))$ with the pairs $(x, x+k)$ and $(\rho(x), \rho(x+k))$. Then
$\rho' |_{S}$ and $\rho' |_{S^c}$ are planar.
%$\rho'$ is planar when restricted to $S$, and when restricted to $S^c$. 
Also, replacing $\rho$ with $\rho'$ does not change the right hand side of equation (\ref{eqn:signSdefn2}). 
\item[(2)] If $\rho(y)$ is not in the interval $[x+1, \ldots, y-1]$ and there is a node in this interval that is in $S^c$, let $k$ be the smallest integer such that $y -  k$ is in $S^c$ and let $\rho'$ be the pairing obtained from $\rho$ by replacing the pairs $(y, \rho(y))$ and $(y-k, \rho(y-k))$ with the pairs $(y,y -  k)$ and $(\rho(y), \rho(y-k))$. Then
$\rho' |_{S}$ and $\rho' |_{S^c}$ are planar.
%$\rho'$ is planar when restricted to $S$, and when restricted to $S^c$. 
Also, replacing $\rho$ with $\rho'$ does not change the right hand side of equation (\ref{eqn:signSdefn2}). 
\end{itemize}

\end{lemma}

\begin{proof}
Since the proofs of (1) and (2) are completely analogous, we only prove (1). 

We first show that $\rho'|_{S}$ is planar. 
Since we chose the smallest integer $k$ such that $x+k$ is in $S$, there are no chords connecting two nodes in $S$ that cross the chord $x \!\frown\! (x+k)$.
%pairing $x$ with $x+k$ to obtain $\rho'$ did not create crossings in $\rho' |_{S}$
We need to check that there are no chords connecting two nodes in $S$ that cross the chord $\rho(x) \!\frown\! \rho(x+k)$. 
%that pairing $\rho(x)$ with $\rho(x+k)$ did not create any crossings in $\rho' |_{S}$. 
If there was such a crossing,
%If pairing $\rho(x)$ with $\rho(x+k)$ resulted in a crossing in $\rho' |_{S}$, 
that means that there is a node $a \in S$
such that one of the following holds:
\begin{itemize}
\item[(1)] $a < \rho(x+k) < \rho(a) < \rho(x)$,
\item[(2)] $\rho(x+k) < a <  \rho(x) < \rho(a)$,
\item[(3)] $a < \rho(x) < \rho(a) < \rho(x+k)$, or
\item[(4)] $\rho(x) < a <  \rho(x+k) < \rho(a)$.
\end{itemize}

We use the facts that if $a > x$ then $a > x+k$ (since otherwise $a \in S^c$, a contradiction) or, similarly, if $\rho(a) > x$ then $\rho(a) > x+k$, to show that if the inequalities in (1), (2), (3), or (4) hold, then $\rho |_{S}$ is not planar.

For example, in case (1), if $a > x$ then $a > x+k$. So we have
$$x+k < a < \rho(x+k) < \rho(a),$$
which contradicts that $\rho |_{S}$ is planar.
If $a < x$ then there are two cases. If $\rho(a) < x$, we have
$a < \rho(x+k) < \rho(a) < x+k.$
If instead $\rho(a) > x$, we have
$a < x < \rho(a) < \rho(x).$
In both cases, we have a contradiction.

In case (2), if $a > x$, then we have
$x < a <  \rho(x) < \rho(a).$
If $a < x$ and $\rho(a) < x$, then 
$a < \rho(x) < \rho(a) <x.$
If $a < x$ and $\rho(a) > x$, then 
%$\rho(a) > x+k$. In this case,
$\rho(x+k) < a < x+k < \rho(a).$
In all cases, we have a contradiction. 

Case (3) is similar to case (2), and case (4) is similar to case (1).

We conclude that $\rho' |_{S}$ is planar. Since $\rho |_{S^c}$ was planar and the nodes $x, x+k, \rho(x), \rho(x+k)$ are all in $S$, $\rho' |_{S^c}$ is also planar. %when restricted to $S^c$. 
% did not have an impact on its planarity
%We conclude that $\rho'$ is planar when restricted to $S$ and $S^c$. 

Next, we observe that the number of components in $\pi \cup \rho$ and the number of components in $\pi \cup \rho'$ differ by 1. 
This is because to obtain the pairing $\rho'$ from $\rho$, we swapped the locations of $x$ and $\rho(x+k)$. Since $x$ and $\rho(x+k)$ are both in $S$ and both the same color, $(\pi, \rho')$ and $(\pi, \rho)$ differ by a move of type $A_{BW}$. So by Corollary \ref{ABWcor}, the number of components in $\pi \cup \rho$ and the number of components in $\pi \cup \rho'$ differ by 1. Since $\sign(\rho') = - \sign(\rho)$, replacing $\rho$ with $\rho'$ does not change the right hand side of equation (\ref{eqn:signSdefn2}). 
\end{proof}

The following useful observation is immediate from the definitions.

 \begin{rem}
\label{lem:swapab}
Let $\sigma$ be a pairing such that
$x$ and $y$ are two nodes that are not paired in $\sigma$, and let $\sigma'$ be the pairing obtained by swapping the locations of $x$ and $y$ in $\sigma$. 
%(so $\sigma'$ contains the pairs $(x, \sigma(y))$ and $(y, \sigma(x))$.)
Suppose $S$ is a balanced subset of nodes such that $x \in S$ and $y \in S^c$. Let $S' = (S \setminus \{x \} ) \cup \{y \}$. If $\sigma$ does not connect $S$ to $S^c$, then $\sigma'$ does not connect $S'$ to $S'^c$. 
\end{rem}

\noindent {\em{Proof that equation (\ref{eqn:signSdefn2}) holds.}}
%First assume that $S$ is a balanced set of size $2j$ such that there exists a planar black-white pairing $\rho$ that does not bridge $S$ to $S^c$. Although it may not seem immediately obvious that such a set always exists, recall from Lemma \ref{firstlemma34} that regardless of the node coloring of {\bf N}, there exists a planar black-white pairing $\rho$ of ${\bf N}$. So we can choose $S$ to be $2j$ of the arcs of $\rho$. 
%Then by definition,
%$$\sign(S) = (-1)^{\# \text{crossings of } \rho} = 1$$
%Let
% $\pi = \rho$. Since $\pi$ is odd-even and black-white, for all pairs in $\pi$, either both %nodes of the pair are in $T$ or both are not in $T$, so $\pi$ does not bridge $S \triangle T$ to $(S \triangle T)^c$. 
%Since $\pi = \rho$, $(-1)^{\frac{\# \text{nodes}}{2}} = (-1)^{\# \text{comp in } \pi \cup \rho}$ and
%$\sign_{OE}(\pi) = \sign_{BW}(\rho)$, so equation (\ref{eqn:signSdefn2}) holds. 
%\noindent {\bf Assumptions.} 
Throughout this proof, we assume that we are given a balanced set $S$, an odd-even pairing $\pi$ that does not connect $S \triangle T$ to $(S \triangle T)^c$, and 
a black-white pairing $\rho$ that does not connect $S$ to $S^c$ and is planar when restricted to $S$ and when restricted to $S^c$.

Recall from the beginning of Section~\ref{sec:eqnholds} that we are considering two types of replacements that we can make to $S$ to obtain $S'$: (1) replacing $x \in S$ with $x+1 \in S^c$, where $(x, x+1)$ is a couple of consecutive nodes of the same color, and (2) replacing $x \in S$ with $y \in S^c$, where $x < y$ are the same color and all $\ell$ nodes appearing between $x$ and $y$ are the opposite color of $x$ and $y$ for some $\ell \geq 1$.
%We will make
%one of the previously described types of replacements to obtain $S'$. %Then, 
For both types of replacements,
we will construct a black-white pairing
 $\rho'$ that does not connect $S'$ to $S'^c$ such that
$\rho |_{S'}$ and $\rho |_{S'^c}$ are planar
% is planar when restricted to $S'$ and when restricted to $S'^c$ 
 and
an odd-even pairing $\pi'$ that does not connect $S' \triangle T$ to $(S' \triangle T)^c$. We will show that after replacing $S$, $\pi$, $\rho$ in equation (\ref{eqn:signSdefn2}) with $S'$, $\pi'$ and $\rho'$, equation (\ref{eqn:signSdefn2}) still holds. \\

\noindent {\bf (1) Replace $x \in S$ with $x+1$ $\in S^c$. } \\

%Without loss of generality we assume that $x$ and $x+1$ are both white. 

%Let $\pi$ be an odd-even pairing such that $\pi$ does not bridge $S \triangle T$ to $(S \triangle T)^c$ and let
%$\rho$ be a black-white pairing such that $\rho$ does not bridge $S$ to $S^c$ and $\rho$ is planar when restricted to $S$ and when restricted to $S^c$.
Suppose we replace $x \in S$ with $x +1 \in S^c$ to obtain $S'$. 
There are two cases to consider based on whether or not $\pi(x) = x+1$. In both cases, we let $\rho'$ be the pairing obtained by swapping the locations of $x$ and $x+1$ in $\rho$. By Remark \ref{lem:swapab}, $\rho'$ is a black-white pairing that does not connect $S'$ to $S'^c$. Also note that since $\rho |_{S}$
and $\rho |_{S^c}$ are planar, 
%when restricted to $S$ and restricted to $S^c$, 
$\rho' |_{S'}$ and $\rho' |_{S'^c}$ are planar.\\
% is planar when restricted to $S'$ and when restricted to $S'^c$. \\

\noindent {\bf Case 1.} If $\pi(x) = x+1$, 
%let $\rho'$ be the pairing obtained by swapping the locations of $x$ and $x+1$ in $\rho$. 
%By Lemma \ref{lem:pix=y}, 
%$\rho'$ is a black-white pairing that does not bridge $S'$ to $S'^c$ and 
$\pi$ does not connect $S' \triangle T$ to $(S' \triangle T)^c$ and when we replace $\rho$ with $\rho'$, the right hand side of equation (\ref{eqn:signSdefn2}) changes sign by Lemma \ref{lem:pix=y}.
%We will show that 
%\begin{itemize}
%\item[(a)] $\rho'$ is a black-white pairing that does not bridge $S'$ to $S'^c$ and $\pi$ does not bridge $S' \triangle T$ to $(S' \triangle T)^c$. 
%\item[(b)] Equation (\ref{eqn:signSdefn2}) holds. 
%\end{itemize}
%We will first address (a). 
%By Lemma \ref{lem:swapab}, $\rho'$ does not bridge $S'$ to $S'^c$. The fact that $\pi$ does not bridge $S' \triangle T$ to $(S' \triangle T)^c$ follows from the observation that since $\pi(x) = x+1$, both $x$ and $x+1$ are in $S \triangle T$ or both are in $(S \triangle T)^c$. If both $x, x+1$ are in $S \triangle T$ then since we assumed $x \in S$, $x+1 \in T$, so both $x, x+1$ are in $(S' \triangle T)^c$. 
%So $\pi$ does not bridge $S' \triangle T$ to $(S' \triangle T)^c$. 
%the difference between $S \triangle T$ and $S' \triangle T$ is that one of $S \triangle T$, $S' \triangle T$ contains $u$ and $u+1$ and the other does not. 
%Since we obtained $\rho'$ from $\rho$ by swapping the locations of $x$ and $x+1$, 
% $\sign_{BW}(\rho') = - \sign_{BW}(\rho)$.
%The number of components in $\pi \cup \rho$ is the same as the number of components in $\pi \cup \rho'$ because the path $\pi(\rho(x)) - \rho(x) - x - x+1 - \rho(x+1)$ is replaced with 
%$\pi(\rho(x)) - \rho(x) - x+1 - x- \rho(x+1)$.
%$\pi(\rho(u)) - \rho(u+1) - u - u+1 - \rho(u)$.
% (Note that $\pi(\rho(x))$ coxld be equal to $\rho(x+1)$). So the right hand side of equation (\ref{eqn:signSdefn2}) changes sign. 
Since we swapped the locations of $x$ and $x+1$ in $\rho$ to obtain $\rho'$, $(-1)^{\text{\# of crosses of }  \rho} =-(-1)^{\text{\# of crosses of }  \rho'}$. So equation (\ref{eqn:signSdefn2}) holds. \\

\noindent {\bf Case 2.} If $\pi(x) \neq x+1$, 
%observe that since $x \in S$ and $x+1 \notin S$, either $x$ and $x+1$ are in different components or 
%$x$ and $x+1$ are in the same component and a path from $x$ to $x+1$ contains an odd number of transition pairs.
let $\pi'$ be the pairing obtained from $\pi$ by pairing $x$ with $x+1$, $\pi(x)$ with $\pi(x+1)$, and leaving the remaining pairs the same. 
%Let $\rho'$ be the pairing obtained from $\rho$ by swapping the locations of $x$ and $x+1$. 
By Lemma \ref{lem:pix=y}, $\pi'$ does not connect $S' \triangle T$ to $(S' \triangle T)^c$ and when we replace $\pi$ with $\pi'$ and $\rho$ with $\rho'$, the right hand side of equation (\ref{eqn:signSdefn2}) changes sign. 
As in Case 1, $(-1)^{\text{\# of crosses of }  \rho} =-(-1)^{\text{\# of crosses of }  \rho'}$, so equation (\ref{eqn:signSdefn2}) holds. \\
%Since we swapped the locations of $x$ and $x+1$ in $\rho$, $(-1)^{\text{\# of crosses of }  \rho}$ changes as well. \\

\noindent {\bf (2) Replace $x$ with $y$, where $x < y$ are the same color and all $\ell$ nodes in the interval $[x+1, x+2, \ldots, y-1]$ are the opposite color of $x$ and $y$ ($\ell \geq 1$). }\\

%Without loss of generality, we assume $x$ and $y$ are white nodes and all $\ell$ nodes between $x$ and $y$ are black.
% Let $\pi$ be an odd-even pairing such that $\pi$ does not bridge $S \triangle T$ to $(S \triangle T)^c$ and let
%$\rho$ be a black-white pairing such that $\rho$ does not bridge $S$ to $S^c$ and $\rho$ is planar when restricted to $S$ and when restricted to $S^c$.

Suppose we replace $x \in S$ with $y \in S^c$ to obtain $S'$. There are several cases to consider based on whether $x$ and $y$ are paired with nodes in the interval $[x+1, x+2, \ldots, y-1]$. \\

\noindent {\bf Case 1.} We first consider the case when both $x$ and $y$ are paired with a node in the interval $[x+1, x+2, \ldots, y-1]$. \\

\noindent {\em Construction of $\rho'$}.
Let $\rho^{(1)}$ be the pairing obtained by swapping the locations of $x$ and $y$. By Remark \ref{lem:swapab}, $\rho^{(1)}$ does not connect $S'$ to $S'^c$.

\begin{wrapfigure}{r}{.3\textwidth}
\vspace{.2cm}
  \begin{tikzpicture}[decoration={brace, mirror, raise=4pt}]
        \draw (1*60:2.4) node {$\rho(x)$};
                \draw (0*60:2.4) node {$x$};
          \draw (-4*60:2.4) node {$\rho(y)$};
              \draw (-30-3*60:2.4) node {$y$};
               \draw (9*15:2.3) node {\small{$S$}};
     \node[shape=circle,fill=black,  scale=0.5] (1) at (0*30:2) {}; %y
        \draw (1) arc (0:360:2); %CIRCLE
    
      \node[shape=circle,fill=white,  scale=0.4] at (0*15:2) {}; %y
        \node[shape=circle,fill=black, scale=0.5] (2) at (1*15:2) {};
              \node[shape=circle,fill=black, scale=0.5] (3) at (2*15:2) {}; %rho(y)
                    \node[shape=circle,fill=black, scale=0.5] (4) at (3*15:2) {};
                          \node[shape=circle,fill=black, scale=0.5] (5) at (4*15:2) {}; %rho(x)
                                \node[shape=circle,fill=black, scale=0.5] (6) at (5*15:2) {}; %x
                                 \node[shape=circle,fill=black, scale=0.5] (7) at (6*15:2) {}; 
                                     \node[shape=circle,fill=black, scale=0.5] (0) at (7*15:2) {}; 
                                    \node[shape=circle,fill=black, scale=0.5] (8) at (8*15:2) {}; %rho(i)
                                      \node[shape=circle,fill=black, scale=0.5] (9) at (9*15:2) {}; 
                                        \node[shape=circle,fill=black, scale=0.5] (10) at (10*15:2) {}; 
                                        % \node[shape=circle,fill=black, scale=0.5] (11) at (11*15:2) {};
                                     \node[shape=circle,fill=white, scale=0.4]  at (10*15:2) {};

        \draw (1) -- (5);
          \draw  (10) -- (8);
          %={brace,amplitude=10pt,mirror,raise=4pt},
          \draw[decorate] (2) -- (4) node [black,midway,xshift=0.5cm,yshift = 0.25cm] {\small $\in S^c$};

\end{tikzpicture}
\caption{A possible configuration of the nodes in Case 1.}
\label{fig:bothpaired}
\vspace{-.75cm}
 \end{wrapfigure}
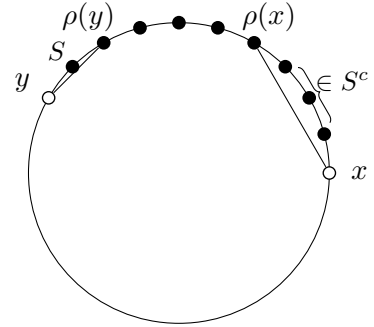

We observe that if $\ell > 2$, at least one of $\rho^{(1)} |_{S'}$,
 $\rho^{(1)} |_{S'^c}$ is not planar.
%when restricted to $S'$, not planar when restricted to $S'^c$, or both. 
To see this, observe that since $\rho |_{S}$ and $\rho |_{S^c}$ are planar,
%when restricted to $S$ and when restricted to $S^c$, 
the nodes 
in the interval $[x+1, \ldots, \rho(x)-1]$
%between 
%$x$ and $\rho(x)$ 
are in $S^c$ and the nodes in the interval
$[\rho(y)+1, \ldots, y-1]$
%between $\rho(y)$ and $y$
 are in $S$ (see Figure~\ref{fig:bothpaired}). 

Suppose towards a contradiction that $\rho^{(1)} |_{S'}$ and $\rho^{(1)} |_{S'^c}$ are planar. Since $\rho^{(1)} |_{S'^c}$ is planar, all nodes 
in the interval $[x+1, \ldots, \rho(y)-1]$
%between $x$ and $\rho(y)$ 
are in $S'$. This means that either 
\begin{itemize}
\item[(1)] $\rho(y) = x+1$, or
\item[(2)] $\rho(x) = x+1$ and $\rho(y) = x+2$. 
\end{itemize}
If (1) holds, there is at least one node 
%between $\rho(y)$ and $y$ 
in the interval $[\rho(y) + 1, \ldots, y-1]$
other than $\rho(x)$. By the observation in the previous paragraph, this node is in $S'$. If it is 
in the interval $[\rho(x) + 1, \ldots, y-1]$
%between $\rho(x)$ and $y$
 its chord crosses the $\rho(x) \!\frown\! y$ chord, contradicting the assumption that $\rho^{(1)}|_{S'}$ is planar. If it is 
% between $\rho(x)$ and $x$ 
 in the interval $[x+1, \ldots, \rho(x)-1]$
 it crossed the $\rho(x) \!\frown\! x$ chord, contradicting the planarity of $\rho|_{S}$. If (2) holds, there is at least one node 
 in the interval $[\rho(y)+1, \ldots, y-1]$,
% between $\rho(y)$ and $y$, 
 this node is in $S'$, and its chord crosses the $\rho(x) \!\frown\! y$ chord, contradicting the assumption that $\rho^{(1)} |_{S'}$ is planar.
% when restricted to $S'$. 

Observe that since $\rho$ pairs $x$ and $y$ with nodes in the interval $[x+1, x+2, \ldots, y-1]$, any crossings in $\rho^{(1)} |_{S'}$ must involve nodes in the interval $[x+1, x+2, \ldots, y-1]$.

We claim that we can undo the crossings in $\rho^{(1)} |_{S'}$ one at a time without changing the right hand side of equation (\ref{eqn:signSdefn2}). To prove the claim, we will describe a procedure for constructing $\rho^{(m+1)}$ from $\rho^{(m)}$ so that $\rho^{(m+1)} |_{S'}$ has one fewer crossing than
$\rho^{(m)} |_{S'}$. \\

\begin{figure}[h]
\centering
 %\begin{minipage}{.4\linewidth}
 \begin{tikzpicture}[scale = 0.9]
        \draw (1*60:2.4) node {$i_1 \in S'$};
         %\draw (3*30:2.4) node {$\notin S'$};
                \draw (1*30:2.6) node {$\rho(x)$};
                \draw (0*60:2.4) node {$x$};
          \draw (-4*60:2.4) node {$i_2 \in S'$};
              \draw (-30-3*60:2.4) node {$y$};
               \draw (8*30:2.4) node {$\rho(i_2)$};
                \draw (10*30:2.4) node {$\rho(i_1)$};

     \node[shape=circle,fill=black,  scale=0.5] (1) at (0*30:2) {}; %x

        \draw (1) arc (0:360:2); %CIRCLE
    
      \node[shape=circle,fill=white,  scale=0.4] at (0*30:2) {}; %x
        \node[shape=circle,fill=black, scale=0.5] (2) at (1*30:2) {};
              \node[shape=circle,fill=black, scale=0.5] (3) at (2*30:2) {}; %i
                    \node[shape=circle,fill=black, scale=0.5] (4) at (3*30:2) {};
      \node[shape=circle,fill=black, scale=0.5] (5) at (4*30:2) {}; %i_2
    \node[shape=circle,fill=black, scale=0.5] (6) at (5*30:2) {}; %y
    \node[shape=circle,fill=white, scale=0.4] (7) at (5*30:2) {}; %y
   \node[shape=circle,fill=black, scale=0.5] (8) at (8*30:2) {}; %rho(i_2)
     \node[shape=circle,fill=black, scale=0.5] (9) at (10*30:2) {}; %rho(i_1)
     \node[shape=circle,fill=white, scale=0.4]  at (10*30:2) {}; %rho(i_1)
                                     \node[shape=circle,fill=white, scale=0.4]  at (8*30:2) {};

        \draw (6) -- (2); %y to rho(x)
          \draw (5) -- (8); %i_2 to rho(i_2)
            \draw (3) -- (9); %i_1 to rho(i_1)
\end{tikzpicture}
%\end{minipage}
%\begin{minipage}{.4\linewidth}
 \begin{tikzpicture}[scale = 0.9]
        \draw (1*60:2.4) node {$i_1 \in S'$};
         %\draw (3*30:2.4) node {$\notin S'$};
                \draw (1*30:2.6) node {$\rho(x)$};
                \draw (0*60:2.4) node {$x$};
          \draw (-4*60:2.4) node {$i_2 \in S'$};
              \draw (-30-3*60:2.4) node {$y$};
               \draw (8*30:2.4) node {$\rho(i_2)$};
                \draw (10*30:2.4) node {$\rho(i_1)$};

     \node[shape=circle,fill=black,  scale=0.5] (1) at (0*30:2) {}; %x

        \draw (1) arc (0:360:2); %CIRCLE
    
      \node[shape=circle,fill=white,  scale=0.4] at (0*30:2) {}; %y
        \node[shape=circle,fill=black, scale=0.5] (2) at (1*30:2) {};
              \node[shape=circle,fill=black, scale=0.5] (3) at (2*30:2) {}; %i
                    \node[shape=circle,fill=black, scale=0.5] (4) at (3*30:2) {};
                          \node[shape=circle,fill=black, scale=0.5] (5) at (4*30:2) {}; %rho(x)
\node[shape=circle,fill=black, scale=0.5] (6) at (5*30:2) {}; %x
               \node[shape=circle,fill=white, scale=0.4] (7) at (5*30:2) {}; %x
               \node[shape=circle,fill=black, scale=0.5] (8) at (8*30:2) {}; %rho(i)
            \node[shape=circle,fill=white, scale=0.4]  at (8*30:2) {};
                                          \node[shape=circle,fill=black, scale=0.5] (9) at (10*30:2) {}; %rho(i_2)
     \node[shape=circle,fill=white, scale=0.4]  at (10*30:2) {}; %rho(i_2)

        \draw (6) -- (3); %y to i_1
          \draw (5) -- (8); %i_2 to rho(i_2)
            \draw (2) -- (9); %rho(x) to rho(i_1)

      %  \draw (1) -- (3);
       %   \draw (5) -- (8);
       %          \draw (2) -- (9);
\end{tikzpicture}
%\end{minipage}
%\end{center}
 \begin{tikzpicture}[scale = 0.9]
        \draw (1*60:2.4) node {$i_1 \in S'$};
         %\draw (3*30:2.4) node {$\notin S'$};
                \draw (1*30:2.6) node {$\rho(x)$};
                \draw (0*60:2.4) node {$x$};
          \draw (-4*60:2.4) node {$i_2 \in S'$};
              \draw (-30-3*60:2.4) node {$y$};
               \draw (8*30:2.4) node {$\rho(i_2)$};
                \draw (10*30:2.4) node {$\rho(i_1)$};

     \node[shape=circle,fill=black,  scale=0.5] (1) at (0*30:2) {}; %x
      
        \draw (1) arc (0:360:2); %CIRCLE
    
      \node[shape=circle,fill=white,  scale=0.4] at (0*30:2) {}; %y
        \node[shape=circle,fill=black, scale=0.5] (2) at (1*30:2) {};
              \node[shape=circle,fill=black, scale=0.5] (3) at (2*30:2) {}; %i
                    \node[shape=circle,fill=black, scale=0.5] (4) at (3*30:2) {};
                          \node[shape=circle,fill=black, scale=0.5] (5) at (4*30:2) {}; %rho(x)
\node[shape=circle,fill=black, scale=0.5] (6) at (5*30:2) {}; %x
               \node[shape=circle,fill=white, scale=0.4] (7) at (5*30:2) {}; %x
               \node[shape=circle,fill=black, scale=0.5] (8) at (8*30:2) {}; %rho(i)
            \node[shape=circle,fill=white, scale=0.4]  at (8*30:2) {};
                                          \node[shape=circle,fill=black, scale=0.5] (9) at (10*30:2) {}; %rho(i_2)
     \node[shape=circle,fill=white, scale=0.4]  at (10*30:2) {}; %rho(i_2)

        \draw (6) -- (5); %y to i_1
          \draw (3) -- (8); %i_2 to rho(i_2)
           
            \draw (2) -- (9); %rho(x) to rho(i_1)
     %   \draw (1) -- (2);
     %     \draw (5) -- (8);
      %           \draw (3) -- (9);
\end{tikzpicture}
\caption{Illustration of the procedure for undoing the crossings in $\rho^{(1)} |_{S'}$. Left: Choose the smallest node $i_1$ in $S'$ greater than $\rho^{(1)}(y$) $=$ $\rho(x)$ whose chord crosses the chord $y$ $\!\frown\!$ $\rho(x)$. Center: Swap the locations of $i_1$ and $\rho(x)$ in $\rho^{(1)}$ to obtain $\rho^{(2)}$. Right: Repeat this procedure to obtain $\rho^{(2)}$.}
\label{fig:iterative}
\end{figure}
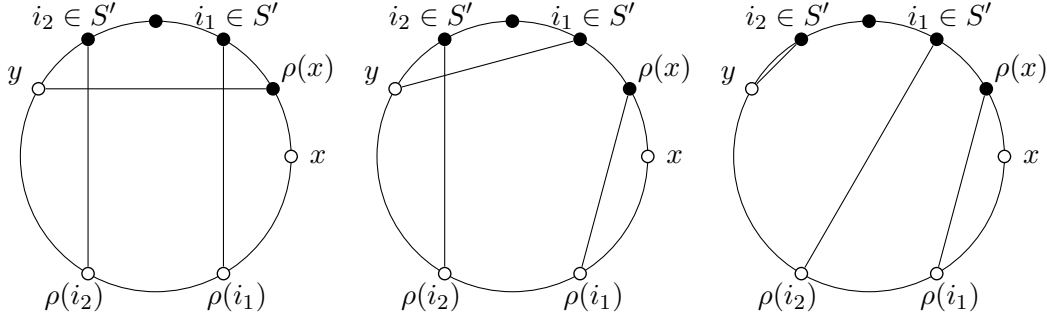

\begin{procedure}
\label{procedure}
(Illustrated in Figure~\ref{fig:iterative}). Choose the smallest node $i_m \in S'$ greater than $\rho^{(m)}(y)$ such that $i_m \!\frown\! \rho^{(m)}(i_m)$ crosses the chord $y \!\frown\! \rho^{(m)}(y)$. Note that $\rho^{(m)}(i_m) = \rho(i_m)$ for all $m$ and $\rho^{(1)}(y) = \rho(x)$.
 Since $i_m$ and $\rho^{(m)}(y)$ are the same color and both in $S'$, we can swap the locations of $i_m$ and $\rho^{(m)}(y)$ in $\rho^{(m)}$ to obtain $\rho^{(m+1)}$, and this is a move of type $A_{BW}$.
By Corollary \ref{ABWcor}, replacing $\rho^{(m)}$ with $\rho^{(m+1)}$ does not change the right hand side of equation (\ref{eqn:signSdefn2}). We claim that $\rho^{(m+1)} |_{S'}$ has one fewer crossing than $\rho^{(m)} |_{S'}$. First observe that since all nodes between $x$ and $y$ are the same color, any chord that crosses the chord  $i_m \!\frown\! y$ must have also crossed the chord $y\!\frown\!\rho^{(m)}(y)$. 
So we just need to check that pairing $\rho^{(m)}(y)$ with $\rho(i_m)$ did not create any crossings. 
 If a black-white chord $a \!\frown\! \rho(a)$  with $a, \rho(a) \in S'$ crosses $\rho^{(m)}(y) \!\frown\! \rho(i_m)$, then if
 % either one of $a, \rho(a)$ is in the interval $[i+1, \ldots, y-1]$, it would have also crossed $i - \rho(i)$. If
  either one of $a, \rho(a)$ is in the interval $[\rho^{(m)}(y) + 1, \ldots, i_m-1]$, it would have crossed $y \!\frown\!\rho^{(m)}(y)$, contradicting the assumption that $i_m$ is the node in $S'$ closest to $\rho^{(m)}(y)$ that crossed $y\!\frown\!\rho^{(m)}(y)$. 
%  Note that neither $a$ nor $\rho(a)$ can be in the interval $[x+1, \ldots, \rho(x) -1]$ since every node in this interval is in $S'^c$. Neither $a$ nor $\rho(a)$ can be in the interval $[\rho(x), \ldots, \rho^{(m)}(y)]$ since nodes in this interval are either in $S'^c$ or were chosen as $i_j$ for $j  < m$. 
So both $a, \rho(a)$ are outside the interval $[\rho^{(m)}(y) + 1, \ldots, i_m-1]$, meaning $a \!\frown\! \rho(a)$ crosses $\rho^{(m)}(y) \!\frown\! \rho(i_m)$ if and only if it crosses $i_m \!\frown\! \rho(i_m)$. \\
%We next check that this move did not create any additional crossings. 
\end{procedure}
%Since we chose $i$ to be the node in $S'$ closest to $\rho(x)$ that crossed $y\!\frown\!\rho(x)$, every node between $i$ and $y$ whose chord crossed the chord $y\!\frown\!\rho(x)$ crosses the chord $i \!\frown\! y$. 

%We conclude that $\rho^{(2)} |_{S'}$ has one fewer crossing than $\rho^{(1)} |_{S'}$ and replacing $\rho^{(1)}$ with $\rho^{(2)}$ did not change the right hand side of equation (\ref{eqn:signSdefn2}). 
Note that if $\rho^{(m)}$ does not connect $S'$ to $S'^c$, then $\rho^{(m+1)}$ does not connect $S'$ to $S'^c$. 
We repeat this procedure until we have a pairing $\rho^{(n)}$ such that $\rho^{(n)}|_{S'}$ is planar.
%when restricted to $S'$. 

% the claim follows from the induction hypothesis. \\

%Now suppose the claim holds when there are $k-1$ crossings when $\rho'$ is restricted to $S'$. Assume there are $k$ crossings when $\rho'$ is restricted to $S'$. 

Similarly, we can undo the crossings in $\rho^{(n)} |_{S'^c}$ one at a time without changing the right hand side of equation (\ref{eqn:signSdefn2}). The resulting pairing is $\rho'$. \\

\noindent {\em Construction of $\pi'$ and analysis of equation (\ref{eqn:signSdefn2})}.

We break into subcases based on the parity of $\ell$.

\noindent {\bf Case 1a.} $\ell$ is odd

\noindent {\em Analysis of LHS of (\ref{eqn:signSdefn2})}.
Since $x$ and $y$ are both paired with black nodes in the interval $[x+1, x+2, \ldots, y-1]$,
$(-1)^{\# \text{ crosses of } \rho^{(1)}} = - (-1)^{\# \text{ crosses of } \rho}$ by Lemma \ref{crossinglemma3}. 
We will show that when we undo crossings to obtain $\rho'$ as described, we apply Procedure~\ref{procedure} an odd number of times. Recall that every node between $x$ and $\rho(x)$ is in $S^c$ and every node between $\rho(y)$ and $y$ is in $S$. It follows that $\rho(x) < \rho(y)$ or $\rho(x) = \rho(y) +1$. Putting these facts together, we see that every node in $S'^c \cap \{x+1, \ldots, y-1\}$ crosses the $x \!\frown\! \rho(y)$ chord, and every node in $S' \cap \{x+1, \ldots, y-1\}$ crosses the $y \!\frown\! \rho(x)$ chord. Since there are an odd number of nodes in $\{x+1, \ldots, y-1\} \setminus \{ \rho(x), \rho(y) \}$, we must apply Procedure~\ref{procedure} an odd number of times. 
We conclude that $(-1)^{\# \text{ crosses of } \rho'} = (-1)^{\# \text{ crosses of } \rho}$.

\noindent {\em Construction of $\pi'$ and analysis of RHS of (\ref{eqn:signSdefn2})}.
Since $\ell$ is odd, $x$ and $y$ are the same parity, so we let $\pi'$ be the pairing obtained by swapping the locations of $x$ and $y$. 
%We can swap the locations of $x$ and $y$ because $\ell$ is odd, so $x$ and $y$ are the same parity.
We claim that $\pi'$ does not connect $S' \triangle T$ to $(S' \triangle T)^c$. 
Since $x$ and $y$ are the same parity and the same color, either both of $x, y$ are in $T$ or neither $x$ nor $y$ are in $T$.
Since $x, y$ are either both in $T$ or both not in $T$, exactly one of $x, y$ is in $S \triangle T$. 
So by Remark \ref{lem:swapab}, $\pi'$ does not connect $S' \triangle T$ to $(S' \triangle T)^c$. 

Also, $\pi' \cup \rho^{(1)}$ has the same number of components as $\pi \cup \rho$ because when we replace $\pi$ with $\pi'$ and $\rho$ with $\rho^{(1)}$,
%only difference between $\pi \cup \rho$ and $\pi' \cup \rho^{(1)}$ is 
the path
$\cdots - \pi(x) - x - \rho(x) - \cdots$  in
$\pi \cup \rho$
is replaced with
$\cdots -\pi(y) - x - \rho(y) - \cdots$ in $\pi' \cup \rho^{(1)}$ and the path $\cdots - \pi(y) - y - \rho(y) - \cdots$ in $\pi \cup \rho$
is replaced with
$\cdots -\pi(x) - y - \rho(x) - \cdots$.

Since we applied Procedure~\ref{procedure} an odd number of times and each application of Procedure~\ref{procedure} is a move of type $A_{BW}$, by Lemma~\ref{AOElemma}, 
$$ (-1)^{\# \text{ comp in } \pi' \cup \rho'}
= -(-1)^{\# \text{ comp in } \pi' \cup \rho^{(1)}}
= - (-1)^{\# \text{ comp in } \pi \cup \rho}.$$
Since $\sign_{BW}(\rho^{(1)}) = -\sign_{BW}(\rho)$ and 
$\sign_{BW}(\rho^{(m+1)}) =- \sign_{BW}(\rho^{(m)})$, 
$\sign_{BW}(\rho') = \sign_{BW}(\rho)$. 
Finally, since $\sign_{OE}(\pi') = - \sign_{OE}(\pi)$, 
we conclude that equation (\ref{eqn:signSdefn2}) holds when $\pi$ is replaced with $\pi'$ and $\rho$ is replaced with $\rho'$.

%We have thus shown that when we replace $\pi$ with $\pi'$ and $\rho$ with $\rho'$ the right hand side of equation (\ref{eqn:signSdefn2}) does not change. 

\noindent {\bf Case 1b.} $\ell$ is even

\noindent {\em Analysis of LHS of (\ref{eqn:signSdefn2})}. As in Case 1a, $(-1)^{\# \text{ crosses of } \rho^{(1)}} = - (-1)^{\# \text{ crosses of } \rho}$. 
%$(-1)^{\# \text{ crossings of } \rho}$ changes.
 If $\ell = 2$ then we let $\rho' = \rho^{(1)}$ and both $\rho'  |_{S'}$ and
 $\rho'  |_{S'^c}$ are planar.
 If $\ell > 2$, then we will show that
% , so equation (\ref{eqn:signSdefn2}) holds. If $\ell > 2$, then we need to show that 
 when we undo crossings in $\rho^{(1)}$ to obtain $\rho'$ we apply Procedure~\ref{procedure} an even number of times. The reasoning is analogous to the $\ell$ is odd case: the claim follows from the fact that there are an even number of nodes in $\{x + 1, \ldots, y-1\} \setminus \{ \rho(x), \rho(y) \}$. We conclude that
 $(-1)^{\# \text{ crosses of } \rho'} = - (-1)^{\# \text{ crosses of } \rho}$. 
 
 \noindent {\em Construction of $\pi'$ and analysis of RHS of (\ref{eqn:signSdefn2})}.
We break into cases based on whether $\pi(x) = y$ or $\pi(x) \neq y$. 
If $\pi(x) = y$, we let $\pi' = \pi$. 
If $\pi(x) \neq y$, we let $\pi'$ be the pairing obtained from $\pi$ by pairing $x$ with $y$, $\pi(x)$ with $\pi(y)$, and leaving the remaining pairs the same. In both cases
$\pi'$ does not connect $S' \triangle T$ to $(S' \triangle T)^c$ 
and 
$\sign_{OE}(\pi) (-1)^{\# \text{ comp in } \pi \cup \rho}
=  \sign_{OE}(\pi') (-1)^{\# \text{ comp in } \pi' \cup \rho^{(1)}}$ by Lemma \ref{lem:pix=y}. 
Since we applied Procedure~\ref{procedure} an even number of times and each application of Procedure~\ref{procedure} is a move of type $A_{BW}$, by Lemma~\ref{AOElemma}, $ (-1)^{\# \text{ comp in } \pi' \cup \rho^{(1)}}
= (-1)^{\# \text{ comp in } \pi' \cup \rho'}$. 
Finally, since $\sign_{BW}(\rho^{(1)}) = -\sign_{BW}(\rho)$ and 
$\sign_{BW}(\rho^{(m+1)}) =- \sign_{BW}(\rho^{(m)})$, 
$\sign_{BW}(\rho') = -\sign_{BW}(\rho)$. We conclude that
 when $\rho$ is replaced with $\rho'$ and $\pi$ is replaced with $\pi'$, the right hand side of equation (\ref{eqn:signSdefn2}) changes sign. Thus equation (\ref{eqn:signSdefn2}) holds. \\

\noindent {\bf Case 2.} We next consider the case where exactly one of $x$ or $y$ is paired with a black node in the interval $[x+1, x+2, \ldots, y-1]$. 

Without loss of generality, suppose that $x$ is the node that is paired with a black node in the interval $[x+1, \ldots, y-1]$. There are two subcases to consider. 

\noindent {\bf Case 2a.} If one of the $\ell$ nodes between $x$ and $y$ is in $S^c$, then let $k$ be the smallest integer such that $y-k$ is in $S^c$ and let  $\rho'$ be the pairing obtained by pairing $y$ with $y-k$ and $\rho(y)$ with $\rho(y-k)$. 
By Lemma
\ref{lem:planarwhenrestricted},  $\rho' |_{S}$ and $\rho' |_{S^c}$ are planar, and replacing $\rho$ with $\rho'$ does not change the right hand side of equation (\ref{eqn:signSdefn2}). 

To show that replacing $\rho$ with $\rho'$ does not change the left hand side of equation (\ref{eqn:signSdefn2}), we must show that $(-1)^{\# \text{ crosses of } \rho'} = (-1)^{\# \text{ crosses of } \rho}$. 
This follows from the observations that:
\begin{itemize}
\item since $\rho |_{S^c}$ is planar, the chords $(y-k) \!\frown\! \rho(y-k)$ and $y \!\frown\! \rho(y)$ do not cross, and 
\item a chord $a \!\frown\! \rho(a)$ crosses exactly one of
$(y-k) \!\frown\! \rho(y-k)$, $y \!\frown\! \rho(y)$ if and only if it crosses exactly one of $\rho(y-k) \!\frown\! \rho(y)$, $(y-k) \!\frown\! y$. 
%\item If a chord $a \!\frown\! \rho(a)$ crosses both of $(y-k) \!\frown\! \rho(y-k)$, $y \!\frown\! \rho(y)$, it crosses neither of $\rho(y-k) \!\frown\! \rho(y)$, $(y-k) \!\frown\! y$. 
%\item If a chord $a \!\frown\! \rho(a)$ crosses both of $\rho(y-k) \!\frown\! \rho(y)$, $(y-k) \!\frown\! y$, it crosses neither of $(y-k) \!\frown\! \rho(y-k)$, $y \!\frown\! \rho(y)$.
\end{itemize}

 Thus we have reduced Case 2a to Case 1, where both $x$ and $y$ are paired with nodes in the interval $[x+1, x+2, \ldots, y-1]$. 

\noindent {\bf Case 2b.} If all of the $\ell$ nodes between $x$ and $y$ are in $S$ (this includes the case where the only node between $x$ and $y$ is $x+1$), then since $\rho |_{S}$ is planar, $x$ is paired with $x+1$.  When we swap the locations of $x$ and $y$ to obtain $\rho^{(1)}$, $\rho^{(1)}|_{S'^c}$ is planar but $\rho^{(1)}|_{S'}$ is not planar. In fact, every node between $x+1$ and $y$ is in $S$ (and therefore in $S'$) and crosses the $y\!\frown\! (x+1)$ chord. As in Case 1, we obtain $\rho'$ by applying Procedure~\ref{procedure} to undo the crossings in $\rho^{(1)}|_{S'}$, and this does not change the right hand side of equation (\ref{eqn:signSdefn2}). We break into cases based on whether $\ell$ is odd or $\ell$ is even before constructing $\pi'$.

\noindent {\bf Case 2bi.} $\ell$ is odd

Since exactly one of $x$ and $y$ is paired with a node in the interval $[x+1, x+2, \ldots, y-1]$,
$(-1)^{\# \text{ crosses of } \rho^{(1)}} = (-1)^{\# \text{ crosses of } \rho}$ by Lemma \ref{crossinglemma3}. 
We claim that when we undo crossings to obtain $\rho'$, there are an even number of crossings to undo. This is because every node between $x+1$ and $y$ crosses the $(x+1)\!\frown\! y$ chord, and since $\ell$ is odd there are an even number of such nodes. So
$(-1)^{\# \text{ crosses of } \rho'} = (-1)^{\# \text{ crosses of } \rho}$.
 %equation (\ref{eqn:signSdefn2}) holds. 

We let $\pi'$ be the pairing obtained by swapping the locations of $x$ and $y$. By the type of arguments used in Case 1, $\pi'$ does not connect $S' \triangle T$ to $(S' \triangle T)^c$, $(-1)^{\# \text{ comp in }\pi' \cup \rho'} = (-1)^{\# \text{ comp in } \pi \cup \rho }$, and $\sign_{BW}(\rho') = -\sign_{BW}(\rho)$. 
We conclude that equation (\ref{eqn:signSdefn2}) holds when $\pi$ is replaced with $\pi'$ and $\rho$ is replaced with $\rho'$.

 %and replacing $\pi$ with $\pi'$ and $\rho$ with $\rho'$ does not change the

\noindent {\bf Case 2bii.} $\ell$ is even

As in Case 2bi, $(-1)^{\# \text{ crosses of } \rho^{(1)}} =(-1)^{\# \text{ crosses of } \rho }$. When we undo crossings to obtain $\rho'$, there are an odd number of crossings to undo, so $(-1)^{\# \text{ crosses of } \rho'} = - (-1)^{\# \text{ crosses of } \rho }$. 

We break into cases based on whether $\pi(x) = y$ or $\pi(x) \neq y$. If $\pi(x) = y$, we let $\pi'  = \pi$.
If $\pi(x) \neq y$, we let $\pi'$ be the pairing obtained from $\pi$ by pairing $x$ with $y$, $\pi(x)$ with $\pi(y)$, and leaving the remaining pairs the same.
In both cases
$\pi'$ does not connect $S' \triangle T$ to $(S' \triangle T)^c$, 
and 
$\sign_{OE}(\pi) (-1)^{\# \text{ comp in } \pi \cup \rho}
=  \sign_{OE}(\pi') (-1)^{\# \text{ comp in } \pi' \cup \rho^{(1)}}$.
 %by Lemma \ref{lem:pix=y}. 
By the type of arguments used in Case 1,
%Since we applied Procedure~\ref{procedure} an odd number of times and each application of Procedure~\ref{procedure} is a move of type $A_{BW}$, by Lemma~\ref{AOElemma}, 
$ (-1)^{\# \text{ comp in } \pi' \cup \rho^{(1)}}
= -(-1)^{\# \text{ comp in } \pi' \cup \rho'}$ and
%Finally, since $\sign_{BW}(\rho^{(1)}) = -\sign_{BW}(\rho)$ and 
%$\sign_{BW}(\rho^{(m+1)}) =- \sign_{BW}(\rho^{(m)})$, 
$\sign_{BW}(\rho') = \sign_{BW}(\rho)$. We conclude that when 
 when $\rho$ is replaced with $\rho'$ and $\pi$ is replaced with $\pi'$, the right hand side of equation (\ref{eqn:signSdefn2}) changes sign. Thus equation (\ref{eqn:signSdefn2}) holds. \\

%the right hand side of equation (\ref{eqn:signSdefn2}) changes. Since exactly one of $x$ and $y$ is paired with black nodes in the interval $[x+1, x+2, \ldots, y-1]$, $(-1)^{\# \text{ crossings of } \rho}$ does not change.

%There are an odd number of crossings to undo because there are an odd number of nodes in the interval $[x + 2, \ldots, y-1]$. \\

\noindent {\bf Case 3.} Finally, we observe that we can reduce the case where neither $x$ nor $y$ is paired with a black node in the interval $[x+1, x+2, \ldots, y-1]$ to the case where exactly one of $x$ or $y$ is paired with a black node in the interval $[x+1, x+2, \ldots, y-1]$. 

First assume that at least one of the $\ell$ nodes between $x$ and $y$ is in $S$.
Choose the smallest integer $k$ such that $x+k$ is in $S$.
Let $\rho'$ be the pairing that pairs $x$ with $x+k$ and $\rho(x)$ with $\rho(x+k)$. 
%We claim that since $\rho$ was planar when restricted to $S$ and $S^c$, 
By Lemma \ref{lem:planarwhenrestricted}, 
$\rho' |_{S}$ is planar and $\rho|_{S^c}$ are planar 
and replacing $\rho$ to $\rho'$ does not change the right hand side of equation (\ref{eqn:signSdefn2}). 
The argument that $(-1)^{\# \text{ crosses of } \rho'} = (-1)^{\# \text{ crosses of } \rho}$ is the same as the argument in Case 2a.

Finally, if all of the $\ell$ nodes between $x$ and $y$ are in $S^c$, pair $y$ with $x+ \ell$. The argument then proceeds identically.

\subsubsection{Proof that $(\pi, \rho)$ exists}

We conclude by proving the
%Finally, we address the
 existence of an odd-even pairing $\pi$ and a black-white pairing $\rho$ such that $\pi$ does not connect $S \triangle T$ to $(S \triangle T)^c$ and $\rho$ does not connect $S$ to $S^c$.

Recall that at the beginning of Section~\ref{sec:eqnholds} we showed that for all $j$ there is a balanced set $S$ of size $2j$ with a planar black-white pairing $\rho$ that does not connect $S$ to $S^c$, and by choosing $\pi = \rho$ we also have an odd-even pairing $\pi$ that does not connect $S \triangle T$ to $(S \triangle T)^c$. 

 %Recall that for $S = \{2n, 1, 2, \ldots, |S| - 1 \}$ we gave an odd-even pairing $\pi$ and a black-white pairing $\rho$ with this property. 
 
We also showed that any balanced set of size $2j$ can be obtained from $S$ by making a sequence of replacements of types (1) and (2) discussed in the beginning of Section~\ref{sec:eqnholds}. 
Furthermore, we showed that given an odd-even pairing $\pi$ and a black-white pairing $\rho$ such that $\pi$ does not connect $S \triangle T$ to $(S \triangle T)^c$
and $\rho$ does not connect $S$ to $S^c$,
and a set $S'$ obtained from $S$ by making a replacement of the form (1) or (2),
we can modify $\pi$ and $\rho$ to obtain $\pi'$ and $\rho'$ so that $\pi'$ does not connect $S' \triangle T$ to $(S' \triangle T)^c$
and $\rho'$ does not connect $S'$ to $S'^c$. 

We conclude that for each balanced subset $S$, there is an odd-even pairing $\pi$ and a black-white pairing $\rho$ with the desired properties.

\section{A recurrence for tripartite double-dimer configurations}

\subsection{Kenyon and Wilson's determinant formula}

In this section we prove our analogue of Kenyon and Wilson's determinant formula for tripartite pairings. Recall the statement of their theorem from Section \ref{sec:KWwork}:

\begin{customthm}{\ref{thm:kw61}}\cite[Theorem 6.1]{KW2009}
Suppose that the nodes are contiguously colored red, green, and blue (a color may occur zero times), and that $\sigma$ is the (unique) planar pairing in which like colors are not paired together. We have
$$
\widehat{\Pr}(\sigma)
= \sign_{OE}(\sigma) \det [1_{i, j \text{ RGB-colored differently } } X_{i, j} ]^{i = 1, 3, \ldots, 2n-1}_{j =2, 4 \ldots, 2n}.$$
\end{customthm}

%We recall from Section~\ref{sec:KWwork} that $\widehat{ \Pr }(\sigma) = \dfrac{Z^{DD}_{\sigma}(G, {\bf N}) }{(Z^D(G^{BW}))^2}$.

Kenyon and Wilson proved Theorem~\ref{thm:kw61} by combining two key results. The first is from their study of groves (see Section~\ref{sec:groves}). Recall that Kenyon and Wilson showed in Theorem~\ref{thm:KWgrovethm} that $\dddot{ \Pr}(\sigma)$ is an integer-coefficient homogeneous polynomial in the variables $L_{i, j}$. %More precisely, they showed, 
%$\dddot{ \Pr}(\sigma) = \sum\limits_{\text{ partitions } \tau} \mathcal{P}_{\sigma, \tau}^{(t)} L_{\tau}$,
%where
%$L_{\tau} = \sum\limits_{F} \prod\limits_{\{i, j\} \in F} L_{i, j}$. Here the sum is over spanning forests $F$ of the complete graph $K_n$ for which the trees of $F$ span the parts of $\tau$ and the product is over edges $\{i, j\}$ of the forest $F$.
Furthermore, they showed that when $\sigma$ is a partition that is a tripartite pairing, the grove  polynomial $\dddot{ \Pr}(\sigma)$ can be expressed as a Pfaffian whose entries are $L_{i, j}$ or 0. 

\begin{thm}\cite[Theorem 3.1]{KW2009}
\label{KWthm31}
Let $\sigma$ be the tripartite pairing partition defined by circularly contiguous sets of nodes $R, G,$ and $B$, where $|R|, |G|,$ and $|B|$ satisfy the triangle inequality. Then 
$$\dddot{ \Pr }(\sigma) 
= \text{Pf} 
\begin{pmatrix}
0 & L_{R, G} & L_{R, B} \\
-L_{G, R} & 0 & L_{G, B} \\
-L_{B, R} & -L_{B, G} & 0 
\end{pmatrix}$$
where $L$ is the matrix with entries $L_{i, j}$ whose rows and columns are indexed by the nodes, and $L_{R,G}$ is the submatrix of $L$ whose rows are the red nodes and columns are the green nodes.
\end{thm}

%\footnotetext{The $L$ matrix is the matrix with entries $L_{i, j}$, with rows and columns indexed by the nodes. It is the negative of the Dirichlet-to-Neumann matrix from }

The second result they needed is
a theorem which allows one to compute the double-dimer polynomials 
$\widehat{\Pr}(\sigma)$ using the grove polynomials. 

%model using the polynomials for groves. 
\begin{thm}\cite[Theorem 4.2]{KW2006}
\label{thm42}
 If a planar partition $\sigma$ only contains pairs and we make the following substitutions to the grove partition polynomial $\dddot{ \Pr}(\sigma)$:
$$L_{i, j} \to 
\begin{cases} 
0, & \text{ if } i \text{ and } j \text{ have the same parity, } \\
(-1)^{(|i-j| -1)/2} X_{i, j}, & \text{ otherwise,} 
\end{cases}$$
then the result is $\sign_{OE}(\sigma)$ times the double-dimer pairing polynomial $\widehat{ \Pr }(\sigma)$, when we interpret $\sigma$ as a pairing.
%, and $(-1)^{\sigma}$ is the signature of the permutation $\sigma_1, \sigma_3, \ldots, \sigma_{2n-1}$. \\
\end{thm}

We prove Theorem~\ref{thm61} (our version of Theorem~\ref{thm:kw61}) similarly. We can use Theorem~\ref{KWthm31} as stated, but we need the following analogue of Theorem~\ref{thm42}:

%Our version of this theorem is very similar, except for the additional global sign $\sign_{\cons}({\bf N})$. 
\begin{thm}
\label{mythm42}
%Suppose $G$ is a graph with node set ${\bf N}$.  
If a planar partition $\sigma$ only contains pairs and we make the following substitutions to the grove partition polynomial
$\dddot{ \Pr}(\sigma)$:
$$L_{i, j} \to 
\begin{cases} 
0, & \text{ if  i and  j are the same color,} \\
\sign(i, j) Y_{i, j}, & \text{ otherwise,} 
\end{cases}$$
 then the result is $\sign_{c}({\bf N}) \sign_{OE}(\sigma) \widetilde{\Pr }(\sigma).$
%$$L_{\rho} \to 
%\begin{cases} 
%0, & \text{ if } i \text{ and } j \text{ are both black or both white for some $(i, j) \in \rho$} \\
%\sign_{\cons} ({\bf N}) \prod\limits_{(i, j) \in \rho} \sign(i, j) Y_{i, j}, & \text{ otherwise,} 
%\end{cases}$$
% then the result is $\sign_{OE}(\sigma)$ times the double-dimer pairing polynomial $\widehat{\text{Pr}}(\sigma)$ when we interpret $\sigma$ as a pairing. 
\end{thm}

\begin{proof}
%By Theorem~\ref{thm:KWgrovethm}, $\dddot{ \Pr}(\sigma) = \sum\limits_{\text{ partitions } \tau} \mathcal{P}_{\sigma, \tau}^{(t)} L_{\tau}$.
In Theorem~\ref{thm:thm1}, we established that
$$\widetilde{\Pr}(\sigma) :=
\dfrac{Z^{DD}_{\sigma}(G, {\bf N})}{(Z^D(G))^2}
=
 \sum_{\text{black-white pairings } \rho} \mathcal{Q}^{(DD)}_{\sigma, \rho} Y'_{\rho}.$$
In the proof of Theorem~\ref{thm:thm1}, we showed
%We proved in Section~\ref{sec:firstmajorproof}
% that 
 $\mathcal{Q}^{(DD)}_{\sigma, \rho} = \sign_{OE}(\sigma) \sign_{BW}(\rho) \mathcal{P}_{\sigma, \rho}^{(t)}$ (see equation~(\ref{eqn:Qtilde})). 
 This connects the polynomials $\widetilde{\Pr}(\sigma)$ to the grove polynomials $\dddot{ \Pr}(\sigma)$, since 
 $\dddot{ \Pr}(\sigma) = \sum\limits_{\text{ partitions } \tau} \mathcal{P}_{\sigma, \tau}^{(t)} L_{\tau}$ by Theorem~\ref{thm:KWgrovethm}.
 
 Specifically, in the case where $\sigma$ is a pair, we have
  $$\dddot{ \Pr}(\sigma) = \sign_{OE}(\sigma) \sum\limits_{\text{ pairs } \rho} \sign_{BW}(\rho) \mathcal{Q}^{(DD)}_{\sigma, \rho}  L_{\rho}.$$

Observe the sum is over all pairs $\rho$ rather than all partitions. This is because by Rule~\ref{hrule}, 
when we express a partition as a linear combination of planar partitions, any singleton parts of that partition show up in each planar partition with nonzero coefficient. 
%Therefore if $\sigma$ is a planar partition that is a pairing and $\tau$ contains at least one singleton part, $\mathcal{P}_{\sigma, \tau}^{(t)} = 0$. 
Also observe that when we apply Rule~\ref{hrule} to a partition, each of the resulting partitions contains the same number of parts as the original partition.
%write $\tau$ as a linear combination of planar parittions, eachtransform partitions using Rule~\cite{kwrule1}, the number of parts is conserved. 
It follows that if $\sigma$ is a pairing, and $\mathcal{P}_{\sigma, \rho}^{(t)} \neq 0$ for some partition $\rho$, then $\rho$ is also pairing.

% Recall that by Definition~\ref{matrixdefn} the matrix $\mathcal{Q}^{(DD)}$ gives coefficients for the monomials $Y_{\rho} = \prod\limits_{(i, j) \in \rho} Y_{i, j}$ weighted by $(-1)^{\# \text{ crosses of }\rho}$.  
Finally, we recall that $Y'_{\rho} = (-1)^{\text{\# crosses of }\rho} \prod\limits_{i \text{ black} } Y_{i,\rho(i)}$
 and by Lemma \ref{lemma34},
\begin{equation*}
\sign_{\cons}({\bf N}) \sign_{BW}(\rho) \prod\limits_{(i, j) \in \rho} \sign(i, j) =  (-1)^{\# \text{ crosses of } \rho}.
\end{equation*}
The theorem follows. 
\end{proof}

%The theorem follows because in order to get the polynomial from the matrix you have to multiply each term $Y_{\rho}$ by $(-1)^{# crosses}$. This is equivalent to muultiplying by 
%$\sign_{\cons}({\bf N}) \sign_{BW}(\rho) \prod\limits_{(i, j) \in \rho} \sign(i, j)$, and the latter two terms cancel. 

%\subsection{Theorem 6.1 from Kenyon and Wilson}

%The following theorems are from \cite{KW2009}. 
%To prove their determinant formula, 
%Kenyon and Wilson combine Theorem~\ref{mythm42} with the following theorem that expresses the grove polynomial

%Remarkably, the analogous theorem in our more general setting has no additional global sign. 

The remainder of this section will be devoted to proving the following theorem.

\begin{customthm}{\ref{thm61}}
Suppose that the nodes are contiguously colored red, green, and blue (a color may occur zero times), and that $\sigma$ is the (unique) planar pairing in which like colors are not paired together.  We have
$$\widetilde{\Pr}(\sigma)= \sign_{OE}(\sigma) \det [1_{i, j \text{ RGB-colored differently } } Y_{i, j} ]^{i = b_1, b_2, \ldots, b_{n}}_{j = w_1, w_2, \ldots, w_{n} }.$$
where $b_1 < b_2 < \ldots < b_n$ are the black nodes listed in increasing order and $w_1 < w_2 < \ldots < w_n$ are the white nodes listed in increasing order. 
\end{customthm}
While our proof of Theorem~\ref{thm61} is very similar to Kenyon and Wilson's proof of Theorem~\ref{thm:kw61}, we do require the following technical lemma. 

\begin{lemma}
\label{lem:lifesaver}
 Let ${\bf N}$ be a set of $2n$ nodes and
let $(n_1, n_1 +1), \ldots, (n_{2k}, n_{2k} + 1)$ be a complete list of couples of consecutive nodes of the same color. Define $(-1)^{i > j}$ to be $-1$ if $i > j$, and $1$ otherwise, and let
\[M =  [(-1)^{i > j} \sign(i, j) Y_{i, j} ]^{i = b_1, b_2, \ldots, b_{n}}_{j = w_1, w_2, \ldots, w_{n} }, \]
where $b_1 < b_2 < \cdots < b_n$ are the black nodes listed in increasing order and $w_1 < w_2 < \cdots < w_n$ are the white nodes listed in increasing order. 
Then $M$ is a block matrix where within each block, the signs of the entries are staggered in a checkerboard pattern. 

Furthermore,
let $t$ be the total number of
rows and columns of $M$ that we need to multiply by $-1$ to obtain a matrix with entries whose signs are staggered in a checkerboard pattern where the upper left entry is positive. If node 1 is black, 
$$(-1)^{t} =  \sign_{\cons}({\bf N}) 
(-1)^{\sum\limits_{i=1}^{2k} \lfloor \frac{n_i}{2} \rfloor} $$
and if node 1 is white, 
$$(-1)^{t} =  (-1)^{n} \sign_{\cons}({\bf N}) 
(-1)^{\sum\limits_{i=1}^{2k} \lfloor \frac{n_i}{2} \rfloor}.$$
\end{lemma}

\begin{proof}
We will first prove the claim that $M$ is a block matrix where within each block, the signs of the entries are staggered in a checkerboard pattern. 

\noindent 
\begin{minipage}{.75\textwidth}
\hspace{10pt} We begin with an example. 
Suppose we have 20 nodes colored as shown right. Then there are four couples of consecutive nodes of the same color: $(4, 5), (8, 9), (13, 14)$, and $(17, 18)$ and $M$ is the matrix shown below. We see that the blocks of $M$ correspond to consecutive nodes of the same color. More precisely, the last column in a block corresponds to a white node that precedes at least two consecutive black nodes. The first column in the next block corresponds to the first white node after these consecutive black nodes. Similarly, the nodes corresponding to the last row in a block and the first row in the next block are separated by at least two consecutive white nodes.
\end{minipage} \hspace{.1cm}
\begin{minipage}{.2\textwidth}
\begin{center}
 \begin{tikzpicture}[scale=.65]
  \draw (0,0) circle (2);
  \foreach \x in {1,2,...,20} {
   \node[shape=circle,fill=black, scale=0.5,label={{((\x-1)*360/20)+90}:\small{\x}}] (n\x) at ({((\x-1)*360/20)+90}:2) {}; };
     \foreach \x in {2, 4, 5, 7, 10, 12, 15, 17, 18, 20} {
     \node[shape=circle,fill=white, scale=0.4] (n\x) at ({((\x-1)*360/20)+90}:2) {};
  };
 \end{tikzpicture}
 \end{center}
\end{minipage}

$$ \left(\begin{array}{cc cc  | cc |  cc cc}
Y_{1, 2} & -Y_{1, 4} & Y_{1, 5} & -Y_{1, 7} & -Y_{1, 10} & Y_{1, 12} & Y_{1, 15} & -Y_{1, 17} & Y_{1, 18} & -Y_{1, 20} \\
-Y_{3, 2} & Y_{3, 4} & -Y_{3, 5} & Y_{3, 7} & Y_{3, 10} & -Y_{3, 12} & -Y_{3, 15} & Y_{3, 17} & -Y_{3, 18} & Y_{3, 20} \\
\hline
-Y_{6, 2} & Y_{6, 4} & -Y_{6, 5} & Y_{6, 7} & Y_{6, 10} & -Y_{6, 12} & -Y_{6, 15} & Y_{6, 17} & -Y_{6, 18} & Y_{6, 20} \\
Y_{8, 2} & -Y_{8, 4} & Y_{8, 5} & -Y_{8, 7} & -Y_{8, 10} & Y_{8, 12} & Y_{8, 15} & -Y_{8, 17} & Y_{8, 18} & -Y_{8, 20} \\
-Y_{9, 2} & Y_{9, 4} & -Y_{9, 5} & Y_{9, 7} & Y_{9, 10} & -Y_{9, 12} & -Y_{9, 15} & Y_{9, 17} & -Y_{9, 18} & Y_{9, 20} \\
Y_{11, 2} & -Y_{11, 4} & Y_{11, 5} & -Y_{11, 7} & -Y_{11, 10} & Y_{11, 12} & Y_{11, 15} & -Y_{11, 17} & Y_{11, 18} & -Y_{11, 20} \\
-Y_{13, 2} & Y_{13, 4} & -Y_{13, 5} & Y_{13, 7} & Y_{13, 10} & -Y_{13, 12} & -Y_{13, 15} & Y_{13, 17} & -Y_{13, 18} & Y_{13, 20} \\
Y_{14, 2} & -Y_{14, 4} & Y_{14, 5} & -Y_{14, 7} & -Y_{14, 10} & Y_{14, 12} & Y_{14, 15} & -Y_{14, 17} & Y_{14, 18} & -Y_{14, 20} \\
-Y_{16, 2} & Y_{16, 4} & -Y_{16, 5} & Y_{16, 7} & Y_{16, 10} & -Y_{16, 12} & -Y_{16, 15} & Y_{16, 17} & -Y_{16, 18} & Y_{16, 20} \\
\hline
-Y_{19, 2} & Y_{19, 4} & -Y_{19, 5} & Y_{19, 7} & Y_{19, 10} & -Y_{19, 12} & -Y_{19, 15} & Y_{19, 17} & -Y_{19, 18} & Y_{19, 20}
\end{array}\right).$$

Since in the matrix above, row $i$ does not correspond to node $i$, we introduce the following notation. 
We define the map $B: \{1, 2, \ldots, n\} \to \{b_1, \ldots, b_n\}$ by letting $B(i)$ be the node corresponding to row $i$. 
%In the example above, $R(4) = 8$. 
 Similarly, we define $W: \{1, 2, \ldots, n\} \to \{w_1, \ldots, w_n\}$ by letting $W(j)$ be the node corresponding to column $j$.
In the example above, $B(4) = 8$ and $W(8) = 17$. 

We will show that $M$ has the form
$$
\kbordermatrix{      & W(j) < s_1    & s_1 < W(j) < s_2 & \cdots & s_{k-1} < W(j) < s_k \\
B(i) < u_1 & A_{1, 1} &A _{1, 2} & \cdots & A_{1, k}   \\
u_1 < B(i) < u_2&   A_{2, 1} & A_{2, 2} & \cdots & A_{2, k}   \\
\vdots &        \vdots   &    \vdots     &     \ddots      &    \vdots         \\
u_{k-1} < B(i) < u_k &   A_{k, 1} & A_{k, 2} & \cdots & A_{k, k}   \\
},$$
where in each block $A_{i, j}$, the signs of the entries are staggered in a checkerboard pattern. Note that a block could be empty. 

We first show that within a block, rows $i$ and $i+1$ have opposite sign. There are two cases to consider: 
\begin{itemize}
\item[(1)] $B(i+1) - B(i) =2$, and 
\item[(2)] $B(i+1) - B(i) = 1$. 
\end{itemize} 
These are the only cases because if $B(i+1) - B(i) > 2$, then there is at least one couple of consecutive white nodes between $B(i)$ and $B(i+1)$, so rows $i$ and $i+1$ are in different blocks.

In case (1), there is not a couple of consecutive nodes of the same color between $B(i)$ and $B(i+1)$, so $a_{B(i), w} = a_{B(i+1), w}$ for all $w$. It follows immediately from the definition
$\sign(b,w) = (-1)^\frac{|b-w| + a_{b, w} -1}{2}$
that $\sign(B(i+1), w) = -\sign(B(i), w)$ unless $B(i) < w < B(i+1)$. But if $B(i) < w < B(i+1)$, the sign $(-1)^{b > w}$ flips. So in case (1), rows $i$ and $i+1$ have opposite sign.

In case (2), $(B(i), B(i+1))$ is a couple of consecutive black nodes, so $|a_{B(i+1), w} - a_{B(i), w}| = 1$. 
If $B(i+1) > w$, 
$$\sign(B(i+1), w) = 
(-1)^{\frac{B(i+1) - w+ a_{B(i+1), w} -1}{2}}
= (-1)^{\frac{B(i) +1 - w+ a_{B(i), w} +1  -1}{2}}
= - \sign(B(i), w).$$
The case where $B(i+1) < w$ is completely analogous. 

We conclude that within a block, rows $i$ and $i+1$ have opposite sign. The proof that within a block columns $j$ and $j+1$ have opposite sign is identical. So, within each block, the signs of the entries are staggered in a checkerboard pattern.\\

Since $M$ is a block matrix where the signs of each block are staggered in a checkerboard pattern, we can always choose rows and columns to multiply by $-1$ so that the signs of the matrix entries are staggered in a checkerboard pattern and the upper left entry is positive. 
Let $t$ be the total number of
rows and columns we need to multiply by $-1$. We claim that if node 1 is black, 
$(-1)^{t} =  \sign_{\cons}({\bf N}) 
(-1)^{\sum\limits_{i=1}^{2k} \lfloor \frac{n_i}{2} \rfloor} $
and if node 1 is white, 
$(-1)^{t} =  (-1)^{n} \sign_{\cons}({\bf N}) 
(-1)^{\sum\limits_{i=1}^{2k} \lfloor \frac{n_i}{2} \rfloor}.$ 
%where recall that $2n$ is the total number of nodes. 

We will prove the claim by induction on $n$, where $2n$ is the total number of nodes. 
The base case is when there are 4 nodes. In this case, $(-1)^{n} = 1$. We check all possible node colorings in the table below. 

\begin{center}
\begin{tabular}{ |c | c |  c |   c | c |}
\hline
& & & & \\[-1em]
black nodes & $M$ & $t$ & $\sign_{\cons}({\bf N})$ & $(-1)^{\sum  \lfloor \frac{n_i}{2} \rfloor}$ \\
\hline
1, 2 &$ \left(\begin{array}{ c  c  }
-Y_{1, 3} & Y_{1, 4} \\
Y_{2, 3} &  -Y_{2, 4} 
\end{array}
\right)$&  2 & $-1$ & $-1$ \\
\hline
1, 3 &$ \left(\begin{array}{ c  c  }
Y_{1, 2} & -Y_{1, 4} \\
-Y_{3, 2} &  Y_{3, 4} 
\end{array}
\right)$&  0 & 1 & 1 \\
\hline
1, 4 &$ \left(\begin{array}{ c  c  }
Y_{1, 2} & -Y_{1, 3} \\
Y_{4, 2} &  -Y_{4, 3} 
\end{array}
\right)$&  1 & 1 & $-1$ \\
\hline
3, 4 &$ \left(\begin{array}{ c  c  }
Y_{3, 1} & -Y_{3,2} \\
-Y_{4, 1} &  Y_{4, 2} 
\end{array}
\right)$&  0& $-1$ & $-1$ \\
\hline
2, 4 &$ \left(\begin{array}{ c  c  }
-Y_{2, 1} & Y_{2, 3} \\
Y_{4, 1} &  -Y_{4, 3} 
\end{array}
\right)$&  2& 1 & 1 \\
\hline
2, 3 &$ \left(\begin{array}{ c  c  }
-Y_{2, 1} & -Y_{2,4} \\
Y_{3, 1} &  Y_{3, 4} 
\end{array}
\right)$&  1& 1 & $-1$ \\
\hline
\end{tabular}
\end{center}

Now assume the claim holds when there are $2n-2$ nodes and suppose that $|{\bf N}| = 2n$. 
Choose the largest nodes $x, x+1$ such that $x, x+1$ are different colors. 

%As in the proof of Lemma~\ref{firstlemma34},
Let
${\bf N'} = \{1, \ldots, 2n-2 \}$. Define
$\psi: {\bf N} - \{x, x+1 \} \to {\bf N'}$ by
\begin{equation*}
\psi(\ell) =
\begin{cases}
 \ell &\mbox{ if }\ell  < x \\
  \ell - 2 &\mbox{ if }\ell > x+1 \
 \end{cases}.
\end{equation*}
That is, $\psi$ defines a relabeling of the nodes of ${\bf N} - \{x, x+1 \}$ so that node 1 is labeled 1, $\ldots,$ node $x-1$ is labeled $x-1$, node $x+2$ is labeled $x,\ldots,$ node $2n$ is labeled $2n - 2$. 

Recall that $(n_1, n_1 +1), \ldots, (n_{2k}, n_{2k} + 1)$ is a complete list of couples of consecutive nodes of the same color in ${\bf N}$. Let 
$(n'_1, n'_1 +1), \ldots, (n'_{2j}, n'_{2j} + 1)$ be a complete list of couples of consecutive nodes of the same color in ${\bf N'}$.

Let $M'$ denote the matrix corresponding to ${\bf N'}$. 
Let $t'$ denote the total number of rows and columns we need to multiply by $-1$ to get a matrix $M'_{(1)}$ with entries whose signs are staggered in a checkerboard pattern so that the upper left entry is positive. 
By the induction hypothesis, 
$$(-1)^{t'} = \sign_{\cons}({\bf N'}) (-1)^{\sum\limits_{i=1}^{2j} \lfloor \frac{n'_i}{2} \rfloor}$$
if node 1 is black and
$$(-1)^{t'} = (-1)^{n-1} \sign_{\cons}({\bf N'}) (-1)^{\sum\limits_{i=1}^{2j} \lfloor \frac{n'_i}{2} \rfloor}$$
if node 1 is white.

There are several cases to consider based on whether or not ${\bf N}$ and ${\bf N'}$ have the same number of couples of consecutive nodes of the same color.
In each case, we will assume that node $2n$ is white. 
When $2n$ is black, the argument is completely analogous. 

 Each case will involve two steps:
\begin{enumerate}
\item[(i)] comparing $\sign_{\cons}({\bf N'})$ to $\sign_{\cons}({\bf N})$ and $(-1)^{\sum \lfloor \frac{n'_i}{2} \rfloor}$ to
$(-1)^{\sum \lfloor \frac{n_i}{2} \rfloor}$, and
\item[(ii)] comparing $t$ to $t'$.
%\item[(ii)] comparing the signs of the entries of $M$ to the signs of the entries of $M'$, 
%\item[(iii)] returning the nodes of $M'_{(1)}$ to their original labels to obtain $M'_{(2)}$, 
%\item[(iv)] adding the column and row corresponding to nodes $x+1$ and $x$ to $M'_{(2)}$ to obtain $M'_{(3)}$, 
%\item[(v)] multiplying additional rows and columns of $M'_{(3)}$ (if necessary) to obtain a matrix $M'_{(4)}$ with entries whose signs are staggered in a checkerboard pattern, with the top left entry positive, and 
%\item[(vi)] letting $\tilde{M}$ be the matrix obtained from $M$ by doing all the multiplications we did to $M'$, and considering the relationship between $\tilde{M}$ and $M'_{(4)}$. 
\end{enumerate}

%Since the proof is lengthy, we illustrate the main ideas with examples.

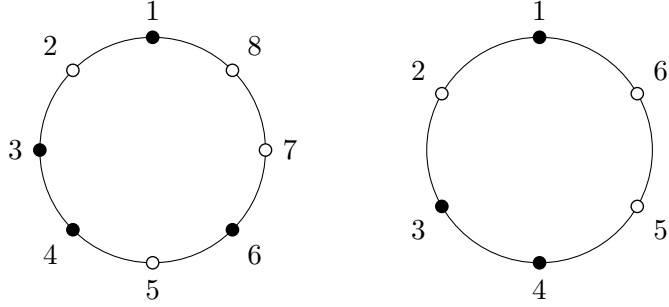
\begin{figure}[htb]
\centering
 \begin{tikzpicture}[scale=.75]
  \draw (0,0) circle (2);
  \foreach \x in {1,2,...,8} {
   \node[shape=circle,fill=black, scale=0.5,label={{((\x-1)*360/8)+90}:\x}] (n\x) at ({((\x-1)*360/8)+90}:2) {}; };
     \foreach \x in {2, 5, 7, 8} {
     \node[shape=circle,fill=white, scale=0.4] (n\x) at ({((\x-1)*360/8)+90}:2) {};
  };
 \end{tikzpicture} \hspace{1cm}
  \begin{tikzpicture}[scale=.75]
  \draw (0,0) circle (2);
  \foreach \x in {1,2,...,6} {
   \node[shape=circle,fill=black, scale=0.5,label={{((\x-1)*360/6)+90}:\x}] (n\x) at ({((\x-1)*360/6)+90}:2) {}; };
     \foreach \x in {2, 5, 6} {
     \node[shape=circle,fill=white, scale=0.4] (n\x) at ({((\x-1)*360/6)+90}:2) {};
  };
 \end{tikzpicture}
 \caption{Shown left is an example of a possible node coloring ${\bf N}$ that could occur in Case 1(a). Nodes 6 and 7 are deleted from ${\bf N}$ and relabeled to obtain ${\bf N'}$, which is shown right.}
 \label{fig:signlemmaex1}
\end{figure}

\noindent {\bf Case 1.}
In the first case, ${\bf N'}$ has the same number of couples of consecutive nodes of the same color as ${\bf N}$. There are two ways this can occur: $x+1 < 2n$, or $x+1 = 2n$. \\

\noindent {\bf Case 1(a).} $x+1 < 2n$

We first assume that node 1 is black. 

\noindent{\bf (i) Comparing $\sign_{\cons}({\bf N})$ to  $\sign_{\cons}({\bf N'})$ and $(-1)^{\sum \lfloor \frac{n'_i}{2} \rfloor}$ to
$(-1)^{\sum \lfloor \frac{n_i}{2} \rfloor}$.}

Since ${\bf N'}$ has the same number of couples of consecutive nodes as ${\bf N}$, $\sign_{\cons}({\bf N}) = \sign_{\cons}({\bf N'})$. 
Since we assumed that $x$ and $x+1$ are the largest nodes such that $x$ and $x+1$ are different colors and ${\bf N'}$ has the same number of couples of consecutive nodes as ${\bf N}$, node $x-1$ and all nodes in the interval $[x+1, \ldots, 2n]$ are white. 
%First assume that node 1 is black. 
Since each node in the interval $[x+1, \ldots, 2n]$ of ${\bf N}$ is white, each node in the interval $[x-1, \ldots, 2n-2]$ of ${\bf N}$ is white, and node 1 is black in both ${\bf N}$ and ${\bf N'}$ we have
% and the nodes  there are $2n-x-1$ nodes $n_i$ that are first in a couple of consecutive nodes such that $\psi(n_i) = n_i - 2$.
%It follows that
\[
(-1)^{\sum\limits_{i=1}^{2k} \lfloor \frac{n'_i}{2} \rfloor}
= 
(-1)^{2n-x-1} (-1)^{\sum\limits_{i=1}^{2k} \lfloor \frac{n_i}{2} \rfloor}
.\]
%= 
%(-1)^{\sum\limits_{i=1}^{2k} \lfloor \frac{\psi(n_i)}{2} \rfloor}
We conclude that
\begin{equation}
\label{eqn:comparison}
\sign_{\cons}({\bf N'}) (-1)^{\sum\limits_{i=1}^{2j} \lfloor \frac{n'_i}{2} \rfloor}
=
(-1)^{2n-x-1} \sign_{\cons}({\bf N}) (-1)^{\sum\limits_{i=1}^{2j} \lfloor \frac{n_i}{2} \rfloor}.
\end{equation}

\noindent {\bf (ii) Comparing $t$ to $t'$.}

%Let $t$ be the total number of row and column multiplications needed to obtain $\tilde{M}$ from $M$ where $\tilde{M}$ has entries whose signs are staggered in a checkerboard pattern. 
Comparing the parity of $t$ and $t'$ is a multi-step process. 
Recall that we obtained $M'_{(1)}$ from $M'$ by multiplying $t'$ rows and columns of $M'$.
We start by returning the nodes of $M'_{(1)}$ to their original labels to obtain $M'_{(2)}$.
Then, we add the row and column corresponding to nodes $x$ and $x+1$ to $M'_{(2)}$ to get $M'_{(3)}$. Finally, we let $\widetilde{M}$ be the matrix obtained from $M$ by doing all the row and column multiplications we did to $M'$ to get $M'_{(1)}$, and consider the relationship between $\widetilde{M}$ and $M'_{(3)}$. 

%For ease of exposition, as above we define the map $R: \{1, 2, \ldots, n\} \to \{b_1, \ldots, b_n\}$ by letting $R(i)$ be the node corresponding to row $i$ of $M$ and we define $C: \{1, 2, \ldots, n\} \to \{w_1, \ldots, w_n\}$ by letting $C(j)$ be the node corresponding to column $j$ of $M$. We define $R': \{1, 2, \ldots, n-1\} \to \{b_1, \ldots, b_{n-1}\}$ and $C': \{1, 2, \ldots, n-1\} \to \{w_1, \ldots, w_{n-1}\}$ analogously for $M'$. 

Previously we defined the map $B: \{1, 2, \ldots, n\} \to \{b_1, \ldots, b_n\}$ by letting $B(i)$ be the node corresponding to row $i$ of $M$ and we defined $W: \{1, 2, \ldots, n\} \to \{w_1, \ldots, w_n\}$ by letting $W(j)$ be the node corresponding to column $j$ of $M$. 
It will be convenient to let $R:= B^{-1}$ and $C:= W^{-1}$, so for example $R(6)$ is the row corresponding to the black node $6$, and $C(7)$ is the column corresponding to the white node 7. 

We define $B'$, $W'$, $R'$ and $C'$
%: \{1, 2, \ldots, n-1\} \to \{b_1, \ldots, b_{n-1}\}$ and $C': \{1, 2, \ldots, n-1\} \to \{w_1, \ldots, w_{n-1}\}$ 
analogously for $M'$.

Because this portion of the proof is long,
we will illustrate the main ideas with an example. 
Let $G$ be a graph with 8 nodes where nodes $1, 3, 4$ and $6$ are colored black (see Figure~\ref{fig:signlemmaex1}). In this example, $x = 6$. So ${\bf N'} = \{1, 2, 3, 4, 5, 6\}$ where nodes $1, 3,$ and $4$ are black. This means that
$$M' = 
\left(
\begin{array}{ c  c c }
Y_{1, 2} & Y_{1, 5} & - Y_{1, 6}  \\
-Y_{3, 2} & - Y_{3, 5} &  Y_{3, 6}  \\
  Y_{4, 2} & Y_{4, 5} & - Y_{4, 6} \\
\end{array}
\right).$$
To obtain
$$M'_{(1)} = 
\left(
\begin{array}{ c  c c }
Y_{1, 2} & -Y_{1, 5} &  Y_{1, 6}  \\
-Y_{3, 2} &  Y_{3, 5} &  -Y_{3, 6}  \\
  Y_{4, 2} & -Y_{4, 5} &  Y_{4, 6} \\
\end{array}
\right)$$
we multiply the second and third columns of $M$ by $-1$, so $t' = 2$. 

In general, to get from $M'$ to $M'_{(1)}$,
either we multiply all of the columns in a block or none of the columns in a block, 
because within each block, the signs of the entries are staggered in a checkerboard pattern. The same is true for the rows. 

\noindent {\bf Return the nodes of $M'_{(1)}$ to their original labels.}
Next, we return the nodes to their original labels (equivalently, we apply the map $\psi^{-1}$) to get $M'_{(2)}$. Note that the only entries that are affected are the entries in the columns corresponding to nodes $\psi(x+2), \ldots, \psi(2n)$. 

In the example, we return node 6 to its original label of 8, resulting in the matrix
$$M'_{(2)} = 
\left(
\begin{array}{ c  c c }
Y_{1, 2} & -Y_{1, 5} &  Y_{1, 8}  \\
-Y_{3, 2} &  Y_{3, 5} &  -Y_{3, 8}  \\
  Y_{4, 2} & -Y_{4, 5} &  Y_{4, 8} \\
\end{array}
\right).$$

\noindent {\bf Add the row and column corresponding to nodes $x$ and $x+1$ to $M'_{(2)}$.}
Now, add to $M'_{(2)}$ the column corresponding to node $x+1$ (i.e. the column with entries $(-1)^{i > x+1} \sign(i, x+1) Y_{i, x+1}$) in between the columns corresponding to nodes $x-1$ and $x+2$. 
%to $M'_{(2)}$. 
Also add the row corresponding to node $x$ as the last row. 
%consider the matrix obtained from $M'_{(2)}$ by 
%Change the sign of the entries in the new column in rows $R(a)$ if $R'(\psi(a))$ was a row we multiplied by $-1$. 
%Similarly, change the sign of the entries in the new row in columns $C(b)$ if $C'(\psi(b))$ was a column we multiplied by $-1$. 
Change the sign of the entries in the new column 
in the rows of $M'$ that we multiplied by $-1$.
Similarly, change the sign of the entries of the new row in the columns that we multiplied by $-1$. 
Call the resulting matrix $M'_{(3)}$.

In the example, we get
$$M'_{(3)} = 
\left(
\begin{array}{ c  c c c}
Y_{1, 2} & -Y_{1, 5} &  -Y_{1, 7}  &Y_{1, 8}  \\
-Y_{3, 2} &  Y_{3, 5} & Y_{3, 7}  & -Y_{3, 8}  \\
  Y_{4, 2} & -Y_{4, 5} &  -Y_{4, 7}  & Y_{4, 8} \\
   -Y_{6, 2} & Y_{6, 5} &  Y_{6, 7}  & -Y_{6, 8} \\ 
\end{array}
\right),$$
%where note that we changed the sign of entries $Y_{6, \psi^{-1}(5)} = Y_{6, 5}$ and $(R(6), C(8))$ since we multiplied columns  $C'(\psi(5)) = C'(5)$ and  $C'(\psi(8)) = C'(6)$ of $M'$ by $-1$. 
where note that we changed the sign of entries $Y_{6, \psi^{-1}(5)} = Y_{6, 5}$ and $Y_{6, \psi^{-1}(6)} = Y_{6, 8} $ because we multiplied the columns of $M'$ corresponding to nodes 5 and 6 
%($C'(5)$ and $C'(6)$) 
by $-1$. 
%corresponding to nodes $5$ and $6$ by $-1$. 
% and the column of $M'$ corresponding to the node $\psi(8) = 6$ by $-1$. 
%$(4, 2)$ and $(4, 4)$ since we multiplied columns $2$ and $4$ of $M'$ by $-1$
%where note that we changed the sign of entries $(R^{-1}(6), C^{-1}(5))$ and $(R^{-1}(6), C^{-1}(8))$ since we multiplied columns  $C'^{-1}(\psi(5)) = C'^{-1}(5)$ and  $C'(\psi(8)) = C'(6)$ of $M'$ by $-1$. 

%It is important to note that 
Since 
%$M$ is a block matrix with checkerboard blocks and 
we changed the signs of entries in the row $R(x)$ and the column $C(x+1)$ as described above,
$M'_{(3)}$ is a block matrix with checkerboard blocks with the following additional properties: 
\begin{itemize}
\item[(1)] All columns strictly to the left of column $C(x+1)$ and all rows strictly above row $R(x)$ are in the same block. 
\item[(2)] The $j$th entry of $C(x-1)$ and $C(x+2)$ have opposite sign because they were adjacent in $M'$, which is checkerboard.
\item[(3)] All columns strictly to the right of $C(x+1)$ are in the same block(s). 
\item[(4)]  $C(x+1)$ is either in same block as $C(x+2)$ or in the same block as $C(x-1)$.
\item[(5)] $R(x)$ is either in the same block as all other rows, or in its own block.
\end{itemize}
%  Therefore it is possible to obtain a except for the last $2n-x+2$ columns and possibly row $R(x)$ and column $C(x+1)$...\\

%Now we consider the multiplications required to get $M'_{(3)}$ to be a checkerboard matrix. 

%Note that we must multiply the column corresponding to $x+1$ if and only if we multiplied every column in its block. Call the matrix after doing this multiplication (or omitting this multiplication, if not necessary) $M'_{(4)}$. In our example, 
%$$M'_{(4)} = 
%\left(
%\begin{array}{ c  c c c}
%Y_{1, 2} & -Y_{1, 5} &  Y_{1, 7}  & Y_{1, 8}  \\
%-Y_{3, 2} &  Y_{3, 5} & -Y_{3, 7}  & -Y_{3, 8}  \\
%  Y_{4, 2} & -Y_{4, 5} &  Y_{4, 7}  & Y_{4, 8} \\
%   -Y_{6, 2} & Y_{6, 5} &  -Y_{6, 7}  & -Y_{6, 8} \\ 
%\end{array}
%\right)$$
%Observe that $M'_{(4)}$ is checkerboard except for the last $2n-x+2$ columns and possibly the last row. 

\noindent {\bf Compare $\widetilde{M}$ to  $M'_{(3)}$ and conclusion.}
Observe that if $i< x$ and $j > x + 1$, then
\begin{equation}
\label{eqn:signsoff}
\sign( \psi(i), \psi(j) ) = (-1)^{ ( \psi(j) - \psi(i) + a_{\psi(i), \psi(j)} -1 )/2 }
 = (-1)^{( j - i - 2 + a_{i, j} -1 )/2} \\
 =  -\sign(i, j),
\end{equation}
so the entries in the columns $C(x+2), \ldots, C(2n)$ are opposite in sign in $M$ compared to the entries in columns $C'(\psi(x+2)), \ldots, C'(\psi(2n))$ in $M'$. 

Returning to our example, we see that
\[ M = 
\left(
\begin{array}{ c  c c c }
Y_{1, 2} & Y_{1, 5} & - Y_{1, 7} &    Y_{1, 8} \\
-Y_{3, 2} & - Y_{3, 5} &  Y_{3, 7} & -  Y_{3, 8} \\
  Y_{4, 2} & Y_{4, 5} & - Y_{4, 7} &   Y_{4, 8} \\
 - Y_{6, 2} & -  Y_{6, 5} & Y_{6, 7} & -Y_{6, 8} \\
\end{array}
\right) \text{ and  }
M' = 
\left(
\begin{array}{ c  c c }
Y_{1, 2} & Y_{1, 5} & - Y_{1, 6}  \\
-Y_{3, 2} & - Y_{3, 5} &  Y_{3, 6}  \\
  Y_{4, 2} & Y_{4, 5} & - Y_{4, 6} \\
\end{array}
\right), \]
%so indeed the entries in column $C(8) = 4$ are opposite signs of the entries in column $C(6) = 3$. 
so indeed each entry in column $C(8) = 4$ has sign opposite of the corresponding
entry of column $C(6) = 3$.

%Next we consider the relationship between the matrices $M'$ and $M$.

%Recall that in the process of getting from $M'$ to $M'_{(1)}$, we multiplied $t'$ rows and columns by $-1$. Let $a_1, \ldots, a_r$ be the nodes in ${\bf N'}$ corresponding to the rows we multiplied by $-1$. Let $b_1, \ldots, b_c$ be the nodes in ${\bf N'}$ corresponding to the columns we multiplied by $-1$. 

%Now we consider the multiplications required to get $M'_{(3)}$ to be a checkerboard matrix. 
%Call the matrix after doing this multiplication (or omitting this multiplication, if not necessary) $M'_{(4)}$. In our example, 
%$$M'_{(4)} = 
%\left(
%\begin{array}{ c  c c c}
%Y_{1, 2} & -Y_{1, 5} &  Y_{1, 7}  & Y_{1, 8}  \\
%-Y_{3, 2} &  Y_{3, 5} & -Y_{3, 7}  & -Y_{3, 8}  \\
%  Y_{4, 2} & -Y_{4, 5} &  Y_{4, 7}  & Y_{4, 8} \\
 %  -Y_{6, 2} & Y_{6, 5} &  -Y_{6, 7}  & -Y_{6, 8} \\ 
%\end{array}
%\right)$$
%Observe that $M'_{(4)}$ is checkerboard except for the last $2n-x+2$ columns and possibly the last row. 

Now let $\widetilde{M}$ be the matrix $M$ obtained by doing all of the $t'$ row and column multiplications we did to $M'$ to obtain $M'_{(1)}$. 
%Then the entries columns $C(x+2), \ldots, C(2n)$ of $\tilde{M}$ are opposite sign in $M$ compared to the entries in columns $C'(\psi(x+2)), \ldots, C'(\psi(2n))$ in $M'_{(3)}$. 
In our example, 
\[  \widetilde{M} = 
\left(
\begin{array}{ c  c c c }
Y_{1, 2} & -Y_{1, 5} & - Y_{1, 7} &   - Y_{1, 8} \\
-Y_{3, 2} &  Y_{3, 5} &  Y_{3, 7} &  Y_{3, 8} \\
  Y_{4, 2} & -Y_{4, 5} & - Y_{4, 7} &  - Y_{4, 8} \\
 - Y_{6, 2} &   Y_{6, 5} & Y_{6, 7} & Y_{6, 8} \\
\end{array}
\right).\]

Since we changed the signs of entries in the row $R(x)$ and the column $C(x+1)$ as described in the previous step, by equation~(\ref{eqn:signsoff}) $\widetilde{M}$ is identical to $M'_{(3)}$ except for the columns $C(x+2), \ldots, C(2n)$. 
%Since $M'_{(3)}$ is checkerboard except for the last $2n-x+2$ columns and possibly the row $R(x)$ and column $C(x+1)$, 
%(The number of multiplications is $t'+1$ if we multiplied the column corresponding to $x+1$, and $t'$ otherwise). 
Combining this fact with observations (1), (2), and (3) about $M'_{(3)}$ above, we conclude that
that $\widetilde{M}$ is checkerboard except possibly for the row $R(x)$ and column $C(x+1)$. Both the entries in $R(x)$ and the entries in $C(x+1)$ alternate in signs, so it remains to determine whether or not we need to multiply $R(x)$ and/or $C(x+1)$ by $-1$. 

%After we are finished multiplying rows and columns and have obtained a checkerboard matrix, the entry $(R(x), C(x+1))$ will have positive sign if and only if $x$ is odd. This follows from the observations that in a matrix with checkerboard entries, the entry $(R(x), C(2n))$ has positive sign, and all nodes $x+1, \ldots, 2n$ are white. Since $(-1)^{x > x+1} = 1$ and $\sign(x+1, x) = 1$, this means that $x$ is odd if and only if we must multiply both $R(x)$ and $C(x+1)$ by $-1$ or neither by $-1$ to achieve a checkerboard matrix. 

%Since $C(x+1)$ is in the same block as $C(x-2)$, we must multiply $C(x+1)$ by $-1$ if and only if we multiplied every column in its block. 

%We conclude that $x$ is odd if and only if the parity of $t$ is the same as the parity of $t'$. Since $x$ is odd if and only if $(-1)^{2n-x-1} = 1$, by equation (\ref{eqn:comparison}), $t$ has the same parity as 
%$$\sign_{\cons}({\bf N})
%(-1)^{\sum\limits_{i=1}^{2j} \lfloor \frac{n_i}{2} \rfloor},$$
%as desired.

Since $(-1)^{x > x+1} = 1$ and $\sign(x+1, x) = 1$, the entry $(R(x), C(x+1))$ of $\widetilde{M}$ is positive. 
Also, since in a matrix with checkerboard entries, the entry $(R(x), C(2n))$ has positive sign, and all nodes $x+1, \ldots, 2n$ are white, 
the entry $(R(x), C(x+1))$ of the final checkerboard matrix we get after multiplying $R(x)$ and/or $C(x+1)$ by $-1$ has positive sign if and only if $x$ is odd. 
%After we are finished multiplying $R(x)$ and/or $C(x+1)$ by $-1$ to obtain a checkerboard matrix, the entry $(R(x), C(x+1))$ 
%will have positive sign if and only if $x$ is odd. This is immediate from the observations that 
%Then, since $(-1)^{x > x+1} = 1$ and $\sign(x+1, x) = 1$,  
This means that $x$ is odd if and only if we must multiply both $R(x)$ and $C(x+1)$ by $-1$ or neither by $-1$ to achieve a checkerboard matrix. 

%Since $C(x+1)$ is in the same block as $C(x-2)$, we must multiply $C(x+1)$ by $-1$ if and only if we multiplied every column in its block. 
We conclude that $x$ is odd if and only if the parity of $t$ is the same as the parity of $t'$. Since $x$ is odd if and only if $(-1)^{2n-x-1} = 1$, by equation (\ref{eqn:comparison}), $t$ has the same parity as 
$$\sign_{\cons}({\bf N})
(-1)^{\sum\limits_{i=1}^{2j} \lfloor \frac{n_i}{2} \rfloor},$$
as desired.

When node 1 is white, the argument is very similar, but we have
\[
(-1)^{\sum\limits_{i=1}^{2j} \lfloor \frac{n'_i}{2} \rfloor}
= 
(-1)^{2n-x} (-1)^{\sum\limits_{i=1}^{2j} \lfloor \frac{n_i}{2} \rfloor}
\]
since $(8, 1)$ is a couple of consecutive nodes of the same color. 
It follows that 
\begin{equation*}
(-1)^{n-1} \sign_{\cons}({\bf N'}) (-1)^{\sum\limits_{i=1}^{2j} \lfloor \frac{n'_i}{2} \rfloor}
=
(-1)^{n} (-1)^{2n-x-1} \sign_{\cons}({\bf N}) (-1)^{\sum\limits_{i=1}^{2j} \lfloor \frac{n_i}{2} \rfloor}.
\end{equation*}
The rest of the argument is identical. \\

\noindent {\bf Case 1(b).} $x+1 = 2n$. 

If ${\bf N'}$ has the same number of couples of consecutive nodes as ${\bf N}$ and $x+1 = 2n$, there are two possibilities: either nodes $2n-2$ and $1$ are black, or nodes $2n-2$ and $1$ are white. 
%There are two ways this can occur, since we assumed node $2n$ is white. We could have that node $2n-1$ and nodes $2n-2$ are black and node $1$ is black, or we could have that nodes 1 and $2n-2$ are white. 

We first assume that node $1$ is black. 

\noindent{\bf (i) Comparing $\sign_{\cons}({\bf N})$ to  $\sign_{\cons}({\bf N'})$ and $(-1)^{\sum \lfloor \frac{n'_i}{2} \rfloor}$ to
$(-1)^{\sum \lfloor \frac{n_i}{2} \rfloor}$.}

Since ${\bf N'}$ has the same number of couples of consecutive nodes as ${\bf N}$ and $x+1 = 2n$, 
$\sign_{\cons}({\bf N}) = \sign_{\cons}({\bf N'})$ and 
$$(-1)^{\sum\limits_{i=1}^{2k} \lfloor \frac{n'_i}{2} \rfloor} = 
(-1)^{\sum\limits_{i=1}^{2k} \lfloor \frac{n_i}{2} \rfloor},$$
so
$$(-1)^{\sum\limits_{i=1}^{2k} \lfloor \frac{n_i}{2} \rfloor} \sign_{\cons}({\bf N}) =(-1)^{\sum\limits_{i=1}^{2k} \lfloor \frac{n'_i}{2} \rfloor} \sign_{\cons}({\bf N'}).$$

\noindent {\bf (ii) Comparing $t$ to $t'$.}

In this case $\psi$ is the identity map, so $M'_{(1)} = M'_{(2)}$. We obtain $M'_{(3)}$ as described in Case 1(a). By the same reasoning as in Case 1(a), all columns to the left of column $C(2n)$ and all rows above row $R(2n-1)$ are in the same block. $C(2n)$ is either in the same block as the other columns or in its own block, and similarly for row $R(2n-1)$. 

Let $\widetilde{M}$ be the matrix $M$ obtained by doing all of the $t'$ multiplications we did to $M'$ to obtain $M'_{(3)}$. 
Since $\psi$ is the identity map, $\widetilde{M} = M'_{(3)}$. 
It remains to determine whether or not we need to multiply $R(2n-1)$ and/or $C(2n)$ by $-1$. 

Since $(-1)^{2n-1 > 2n} = 1$ and $\sign(2n-1, 2n) = 1$, we need to multiply both $C(2n)$ and $R(2n-1)$ or neither in order for $\widetilde{M}$ to be checkerboard. So $t$ has the same parity as $t'$ and therefore the same parity as
$$(-1)^{\sum\limits_{i=1}^{2k} \lfloor \frac{n_i}{2} \rfloor} \sign_{\cons}({\bf N}).$$
This proves the claim when node 1 is black. 
If node 1 is white, the only difference is that
$$(-1)^{\sum\limits_{i=1}^{2k} \lfloor \frac{n'_i}{2} \rfloor}
= -(-1)^{\sum\limits_{i=1}^{2k} \lfloor \frac{n_i}{2} \rfloor}
$$
because node $2n$ was the first in a couple of consecutive white nodes in ${\bf N}$, and in ${\bf N'}$, node $2n-2$ is the first in a couple of consecutive white nodes. 
It follows that
$$(-1)^{n-1} (-1)^{\sum\limits_{i=1}^{2k} \lfloor \frac{n'_i}{2} \rfloor} \sign_{\cons}({\bf N'}) =(-1)^{n} (-1)^{\sum\limits_{i=1}^{2k} \lfloor \frac{n_i}{2} \rfloor} \sign_{\cons}({\bf N}).$$
The rest of the proof is the same. \\

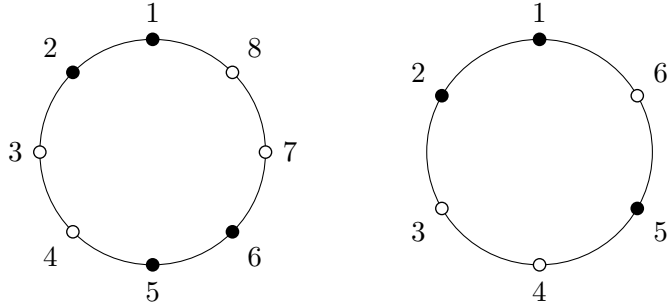
\begin{figure}[htb]
\centering
  \begin{tikzpicture}[scale=.75]
  \draw (0,0) circle (2);
  \foreach \x in {1,2,...,8} {
   \node[shape=circle,fill=black, scale=0.5,label={{((\x-1)*360/8)+90}:\x}] (n\x) at ({((\x-1)*360/8)+90}:2) {}; };
     \foreach \x in {3, 4, 7, 8} {
     \node[shape=circle,fill=white, scale=0.4] (n\x) at ({((\x-1)*360/8)+90}:2) {};
  };
 \end{tikzpicture} \hspace{1cm}
   \begin{tikzpicture}[scale=.75]
  \draw (0,0) circle (2);
  \foreach \x in {1,2,...,6} {
   \node[shape=circle,fill=black, scale=0.5,label={{((\x-1)*360/6)+90}:\x}] (n\x) at ({((\x-1)*360/6)+90}:2) {}; };
     \foreach \x in {3, 4, 6} {
     \node[shape=circle,fill=white, scale=0.4] (n\x) at ({((\x-1)*360/6)+90}:2) {};
  };
 \end{tikzpicture}
 \caption{Shown left is an example of a possible node coloring ${\bf N}$ that could occur in Case 2(a). Nodes 6 and 7 are deleted from ${\bf N}$ and relabeled to obtain ${\bf N'}$, which is the image on the right.}
 \label{fig:signlemmaex2}
\end{figure}

\noindent {\bf Case 2.}
In the second case, ${\bf N'}$ has two fewer couples of consecutive nodes of the same color compared to ${\bf N}$. Again, there are two ways this can occur: $x+1 < 2n$, or $x+1 = 2n$. \\

\noindent {\bf Case 2(a).} $x+1 < 2n$

We first assume that node 1 is black. 
As in Case 1(a), 
%Since this case is lengthy, 
we illustrate the main ideas with an example. 
Let $G$ be a graph with $8$ nodes where nodes $1, 2, 5$ and $6$ are colored black (see Figure~\ref{fig:signlemmaex2}). 
In this example, $x = 6$, so ${\bf N'} = \{1, 2, 3, 4, 5, 6\}$ where nodes $1, 2,$ and $5$ are black. 

\noindent{\bf (i) Comparing $\sign_{\cons}({\bf N})$ to  $\sign_{\cons}({\bf N'})$ and $(-1)^{\sum \lfloor \frac{n'_i}{2} \rfloor}$ to
$(-1)^{\sum \lfloor \frac{n_i}{2} \rfloor}$.}

%Since we assumed node $2n$ is white, node $x$ is black. 
Since we assumed ${\bf N'}$ has two fewer couples of consecutive nodes of the same color compared to ${\bf N}$, nodes $x-1$  and $x$ are both black. 

%As in the proof of Proposition~\ref{firstlemma34}, we 
%define $L: \{n_1, n_2, \ldots, n_{2k} \} \to \{s_1, \ldots, s_k,  u_1, \ldots, u_k \}$ to be the function that gives each of the nodes in ${\bf N}$ that is the first in a couple of consecutive nodes of the same color the label $s_j$ or $u_j$ appropriately. Similarly, we define $L'$ to be the function that gives
%each of the nodes in ${\bf N'}$ that is the first in a couple of consecutive nodes of the same color the label $s_i'$ or $u_i'$ appropriately. 

Recall that $(s_i, s_i + 1)$ denotes a couple of consecutive black nodes of the same color and $(u_i, u_i + 1)$ denotes a couple of consecutive white nodes of the same color. 
By our assumptions, we have
$$\cdots < s_k < u_{k - (2n - x - 2)} < \cdots < u_{k-1} < u_k.$$
When we remove nodes $x$ and $x+1$, there is a one-to-one correspondence between the inversions with respect to the node coloring of ${\bf N}$ and the inversions with respect to the node coloring of ${\bf N'}$ except for the inversions in {\bf N} of the form $(s_k, u_i)$ for $k - (2n-x - 2) \leq i \leq k$. Thus we have that
$$\sign_{\cons}({\bf N}) = (-1)^{2n - x -1} \sign_{\cons}({\bf N'}).$$

In our example, in ${\bf N}$ we have $s_1 < u_1 < s_2 < u_2$ and in ${\bf N'}$ we have $s'_1 < u'_1$. There is one fewer inversion in ${\bf N'}$ compared to ${\bf N}$ and so  $\sign_{\cons}({\bf N}) = - \sign_{\cons}({\bf N'})$.

Next we compare $(-1)^{\sum\limits_{i=1}^{2k-2} \lfloor \frac{n'_i}{2} \rfloor}$ to
$(-1)^{\sum\limits_{i=1}^{2k} \lfloor \frac{n_i}{2} \rfloor}$. In our example, nodes $5$ and $7$ (nodes $x-1$ and $x+1$, respectively) are in $\{n_1, \ldots, n_{2k}\}$. In ${\bf N'}$, the black node $\psi(5) =5$ is adjacent to the white node $\psi(8)=6$ and the couple of consecutive white nodes $(7, 8)$ is not replaced by a new couple of consecutive nodes, so 
$$(-1)^{\sum\limits_{i=1}^{2k} \lfloor \frac{n_i}{2} \rfloor}
 = (-1)^{\lfloor \frac{1}{2} \rfloor}
 (-1)^{\lfloor \frac{3}{2} \rfloor}
  (-1)^{\lfloor \frac{5}{2} \rfloor}
    (-1)^{\lfloor \frac{7}{2} \rfloor} = 1$$ while
 $(-1)^{\sum\limits_{i=1}^{2k-2} \lfloor \frac{n'_i}{2} \rfloor} 
  = (-1)^{\lfloor \frac{1}{2} \rfloor}
 (-1)^{\lfloor \frac{3}{2} \rfloor} = -1$.

In general, since the nodes $x+1, \ldots, 2n$ are white and $x-1$ is black, the node $x-1$ and all nodes in the interval
$[x+1, \ldots, 2n-1]$ are equal to $n_i$ for some $i$. 
%$L(x-1) = s_k$ and  $L(x+1) = u_{k-(2n-x+2)}, \ldots, 
%L(2n-1) =  u_k$. 
%, \ldots, 2n-1$ are equal to $ \ldots,$ in ${\bf N}$ and $x-1$ is $s_k$ in ${\bf N}$. 
In ${\bf N'}$, all nodes in the interval $[\psi(x+2), \ldots, \psi(2n-1)]$ are equal to $n'_i$ for some $i$. 

% $L'(\psi(x+2)) = u'_{k-(2n-x+2)}, \ldots, L'(\psi(2n-1)) =  u'_{k-1}$. 
From the observations that
\begin{itemize}
\item to obtain ${\bf N'}$ we deleted nodes $x$ and $x+1$ from ${\bf N}$, 
\item $\psi(x-1)$ is adjacent to the white node $\psi(x+2)$ in ${\bf N'}$, and
\item  $\psi(y) = y - 2$ for $y > x+1$,
\end{itemize}
we get
$$(-1)^{\sum\limits_{i=1}^{2k-2} \lfloor \frac{n'_i}{2} \rfloor}
=
(-1)^{2n - x - 2} (-1)^{ \lfloor \frac{x-1}{2} \rfloor} (-1)^{ \lfloor \frac{x+1}{2} \rfloor} (-1)^{ \sum\limits_{i=1}^{2k} \lfloor \frac{n_i}{2} \rfloor}.$$
It follows that
$$(-1)^{\sum\limits_{i=1}^{2k-2} \lfloor \frac{n'_i}{2} \rfloor}
= (-1)^{x+1} (-1)^{ \sum\limits_{i=1}^{2k} \lfloor \frac{n_i}{2} \rfloor},$$
so we conclude that
$$\sign_{\cons}({\bf N'})  
(-1)^{\sum\limits_{i=1}^{2k-2} \lfloor \frac{n'_i}{2} \rfloor}
= \sign_{\cons}({\bf N})  
(-1)^{ \sum\limits_{i=1}^{2k} \lfloor \frac{n_i}{2} \rfloor}.$$

\noindent {\bf (ii) Comparing $t$ to $t'$.}
This portion of the proof is similar in structure to part (ii) of Case 1(a). 
%The structure of this portion of the proof will be similar to the structure of part (ii) of Case 1(a). 
In our example, 
$$M' = 
\left(
\begin{array}{ c  c c }
-Y_{1, 3} & Y_{1, 4} &  -Y_{1, 6}  \\
Y_{2, 3} & - Y_{2, 4} &  Y_{2, 6}  \\
  Y_{5, 3} & -Y_{5, 4} &  Y_{5, 6}  \\
\end{array}
\right).$$
We multiply all three columns and the last row of $M'$ by $-1$ to obtain $M'_{(1)}$ and return node 6 to its original label of 8 to obtain $M'_{(2)}$, so 
\[ M'_{(1)} = 
\left(
\begin{array}{ c  c c }
Y_{1, 3} & -Y_{1, 4} &  Y_{1, 6}  \\
-Y_{2, 3} & Y_{2, 4} &  -Y_{2, 6}  \\
  Y_{5, 3} & -Y_{5, 4} &  Y_{5, 6}  \\
\end{array}
\right) 
\text{ and  }
M'_{(2)} = 
\left(
\begin{array}{ c  c c }
Y_{1, 3} & -Y_{1, 4} &  Y_{1, 8}  \\
-Y_{2, 3} & Y_{2, 4} &  -Y_{2, 8}  \\
  Y_{5, 3} & -Y_{5, 4} &  Y_{5, 8}  \\
\end{array}
\right).
\]

%Recall that in the process of getting from $M'$ to $M'_{(1)}$ we need to multiply all of the columns in a block or none of the columns in a block, and likewise for the rows. Note that $C'(\psi(x+2)), \ldots, C'(\psi(2n))$ are in the same blocks of $M'$. Unlike in Case 1(a), the first column to the left of $C'(\psi(x+2))$ may or may not be in the same block of $M'$ as $C'(\psi(x+2))$. 

\noindent {\bf Add the column and row corresponding to nodes $x+1$ and $x$ to $M'_{(2)}$.}
Now, add the column corresponding to node $x+1$ immediately to the left of the column corresponding to the node $x+2$ in $M'_{(2)}$. 
%to $M'_{(2)}$. 
Also add the row corresponding to node $x$ as the last row. 
%consider the matrix obtained from $M'_{(2)}$ by 
Change the sign of the entries in the new column in rows $R(a)$ if $R'(\psi(a))$ was a row we multiplied by $-1$. 
Similarly, change the sign of the entries in the new row in columns $C(b)$ if $C'(\psi(b))$ was a column we multiplied by $-1$, and call the resulting matrix $M'_{(3)}$, which is a block matrix with checkerboard blocks with properties (1) and (3)-(5) from Case 1(a). Property (2) has to be slightly modified:
\begin{itemize}
\item[(2)] The $j$th entry of the first column to the left of $C(x+1)$ and the $j$th entry of $C(x+2)$ have opposite sign because they were adjacent in $M'$, which is checkerboard.
\end{itemize}
%Now, as we did in Case 1(a), consider the matrix $M'_{(3)}$ obtained from $M'_{(2)}$ by adding the column corresponding to node $x+1$ and the row corresponding to node $x$.
%Change the sign of this entries in this new column in rows $R'(\psi(a)) =R(a)$ if $R'(\psi(a))$ was a row we multiplied by $-1$. 
%Similarly, change the sign of the entries in this new row in columns $C(b))$ if $C'(\psi(b))$ was a column we multiplied by $-1$. 
In our example, 
\[
M'_{(3)} = 
\left(
\begin{array}{ c  c c c }
Y_{1, 3} & -Y_{1, 4} &  Y_{1, 7} & Y_{1, 8}  \\
-Y_{2, 3} & Y_{2, 4} & - Y_{2, 7} &   -Y_{2, 8}  \\
  Y_{5, 3} & -Y_{5, 4} & Y_{5, 7} &   Y_{5, 8}  \\
    Y_{6, 3} & -Y_{6, 4} & Y_{6, 7} &   Y_{6, 8}  \\
\end{array}
\right) .
\]

%As in Case 1(a), $M'_{(3)}$ is a block matrix with checkerboard blocks with the following additional properties: 
%\begin{itemize}
%\item[(1)] All columns to the left of column $C(x+1)$ and all rows above row $R(x)$ are in the same block. 
%\item[(2)] The $j$th entry of the first column to the left of $C(x+1)$ and $C(x+2)$ have opposite sign because they were adjacent in $M'$, which is checkerboard.
%\item[(3)] All columns to the right of $C(x+1)$ are in the same block(s). 
%\item[(4)]  $C(x+1)$ is either in same block as $C(x+2)$ or in the same block as the first column to its left. 
%\item[(5)] $R(x)$ is either in the same block as all other rows, or in its own block.
%\end{itemize}

%Note that $M'_{(3)}$ is a block matrix with checkerboard blocks. All columns to the left of column $C(x+1)$ and all rows above row $R(x)$ are in the same block. 
%All columns to the right of $C(x+1)$ are in the same block(s). 
%$C(x+1)$ is either in the same block as $C(x+2)$ or in the same block as the first column to its left. $R(x)$ is either in the same block as all other rows, or in its own block. 

\noindent{\bf Compare $\widetilde{M}$ to the entries of $M'_{(3)}$ and conclusion.}
Observe that if $i< x$ and $j > x + 1$ then
\begin{equation*}
\sign( \psi(i), \psi(j) ) = (-1)^{ ( \psi(j) - \psi(i) + a_{\psi(i), \psi(j)} -1 )/2 }
 = (-1)^{( j-2 - i  + a_{i, j} - 2 -1 )/2} \\
 =  \sign(i, j),
\end{equation*}
so unlike in Case 1(a), the entries in the columns $C(x+2), \ldots, C(2n)$ are the same sign in $M$ as the entries in columns $C'(\psi(x+2)), \ldots, C'(\psi(2n))$ in $M'$. 

Returning to our example, we see that
the entries in column $C(8) = 4$ have the same signs as the entries in column $C(6) = 3$, as

\[M = 
\left(
\begin{array}{ c  c c c }
-Y_{1, 3} & Y_{1, 4} &  Y_{1, 7} &    -Y_{1, 8} \\
Y_{2, 3} & - Y_{2, 4} &  -Y_{2, 7} &  Y_{2, 8} \\
  Y_{5, 3} & -Y_{5, 4} & - Y_{5, 7} &   Y_{5, 8} \\
 - Y_{6, 3} & Y_{6, 4} & Y_{6, 7} & -Y_{6, 8} \\
\end{array}
\right) \text{ and }
M' = 
\left(
\begin{array}{ c  c c }
-Y_{1, 3} & Y_{1, 4} &  -Y_{1, 6}  \\
Y_{2, 3} & - Y_{2, 4} &  Y_{2, 6}  \\
  Y_{5, 3} & -Y_{5, 4} &  Y_{5, 6}  \\
\end{array}
\right).
\]

Let $\widetilde{M}$ be the matrix $M$ obtained by doing all of the $t'$ multiplications we did to $M'$ to obtain $M'_{(3)}$. We see that $\widetilde{M} = M'_{(3)}$, so $\widetilde{M}$ is checkerboard except for the columns $C(x+2), \ldots, C(2n)$ and also possibly the row $R(x)$ and/or the column $C(x+1)$. 

%Now we consider the multiplications required to get $M'_{(3)}$ to be a checkerboard matrix, which we will call $M'_{(4)}$. 

There are two cases to consider. In the first case, $C(x+1)$ is not in the same block as the first column to its left, so we need to multiply $C(x+1)$ by $-1$. Then, since $C(x+1)$ is in the same block as $C(x+2), \ldots, C(2n)$, we need to multiply the remaining $2n-x-1$ columns by $-1$ as well. So we have done $t' + 2n - x$ total multiplications. It remains to consider whether or not we need to multiply row $R(x)$ by $-1$. Recall from Case 1(a) that after we are finished multiplying rows and columns and have obtained a checkerboard matrix, the entry $(R(x), C(x+1))$ must have positive sign if and only if $x$ is odd. 
Since $(-1)^{x > x+1} = 1$, $\sign(x+1, x) = 1$, and we multiplied $C(x+1)$ by $-1$, we multiply $R(x)$ by $-1$ if and only if $x$ is odd. Therefore if $x$ is odd, we have done $t' + 2n - x + 1$ multiplications, and if $x$ is even, we have done $t' + 2n - x$ multiplications. We have thus shown that $t$ has the same parity as $t'$.

If $C(x+1)$ is in the same block as the first column to its left, we do not need to multiply $C(x+1)$ by $-1$ but we still need to multiply the remaining $2n-x-1$ columns by $-1$. So we have done  $t' + 2n - x -1$ total multiplications. Since we did not multiply $C(x+1)$ by $-1$, we multiply $R(x)$ by $-1$ if and only if $x$ is even. Therefore, if $x$ is even, we have done $t' + 2n - x$ total multiplications and if $x$ is odd we have done $t' + 2n - x + 1$ multiplications. Again, $t$ has the same parity as $t'$. 

In both cases, $t$ has the same parity as
$\sign_{\cons}({\bf N})  
(-1)^{ \sum\limits_{i=1}^{2k} \lfloor \frac{n_i}{2} \rfloor},$ which completes the proof when node 1 is black.

When node 1 is white, we have
$$\cdots < s_k < u_{k - (2n - x - 1)} < \cdots < u_{k-1} < u_k$$
but $(s_k, u_k)$ is not an inversion since node 1 is white. 
So we still have
$\sign_{\cons}({\bf N}) = (-1)^{2n - x-1} \sign_{\cons}({\bf N'})$. 
Since $(2n, 1)$ is a couple of consecutive white nodes,
$$(-1)^{\sum\limits_{i=1}^{2k-2} \lfloor \frac{n'_i}{2} \rfloor}
=
(-1)^{2n - x - 1} (-1)^{ \lfloor \frac{x-1}{2} \rfloor} (-1)^{ \lfloor \frac{x+1}{2} \rfloor} (-1)^{ \sum\limits_{i=1}^{2k} \lfloor \frac{n_i}{2} \rfloor}.$$
It follows that
$$(-1)^{\sum\limits_{i=1}^{2k-2} \lfloor \frac{n'_i}{2} \rfloor}
= (-1)^{x} (-1)^{ \sum\limits_{i=1}^{2k} \lfloor \frac{n_i}{2} \rfloor},$$
so we conclude that
$$(-1)^{n-1} \sign_{\cons}({\bf N'})  
(-1)^{\sum\limits_{i=1}^{2k-2} \lfloor \frac{n'_i}{2} \rfloor}
=-(-1)^{n} (-1)^{2n-1} \sign_{\cons}({\bf N})  
(-1)^{ \sum\limits_{i=1}^{2k} \lfloor \frac{n_i}{2} \rfloor}
= (-1)^{n}\sign_{\cons}({\bf N})  
(-1)^{ \sum\limits_{i=1}^{2k} \lfloor \frac{n_i}{2} \rfloor}. 
$$
The rest of the argument is the same. \\

\noindent {\bf Case 2(b).} $x+1 = 2n$

If ${\bf N'}$ has two fewer couples of consecutive nodes of the same color and compared to ${\bf N}$ and $x+1 = 2n$, it must be the case that nodes $2n-1$ and $2n-2$ are both black and node 1 is white.

\noindent{\bf (i) Comparing $\sign_{\cons}({\bf N})$ to  $\sign_{\cons}({\bf N'})$ and $(-1)^{\sum \lfloor \frac{n'_i}{2} \rfloor}$ to
$(-1)^{\sum \lfloor \frac{n_i}{2} \rfloor}$.}

Removing nodes $2n$ and $2n-1$ does not remove any inversions with respect to the node coloring of ${\bf N}$ (recall that $(s_k, u_k)$ is not an inversion when node $1$ is white). Thus $\sign_{\cons}({\bf N}) = \sign_{\cons}({\bf N'})$. 

Next observe that
$$(-1)^{\sum\limits_{i=1}^{2k-2} \lfloor \frac{n'_i}{2} \rfloor} 
= (-1)^{\lfloor \frac{2n}{2} \rfloor} (-1)^{\lfloor \frac{2n-2}{2} \rfloor} (-1)^{\sum\limits_{i=1}^{2k} \lfloor \frac{n_i}{2} \rfloor} 
= - (-1)^{\sum\limits_{i=1}^{2k} \lfloor \frac{n_i}{2} \rfloor} .
$$
Since node 1 is white, we have 
\begin{equation*}
(-1)^{n-1} \sign_{\cons}({\bf N'}) (-1)^{\sum\limits_{i=1}^{2j} \lfloor \frac{n'_i}{2} \rfloor}
=
(-1)^{n}  \sign_{\cons}({\bf N}) (-1)^{\sum\limits_{i=1}^{2j} \lfloor \frac{n_i}{2} \rfloor}.
\end{equation*}

\noindent {\bf (ii) Comparing $t$ to $t'$.}

This argument is identical to (ii) in Case 1(b), and we conclude that $t$ has the same parity as $t'$, and therefore the same parity as $(-1)^{n}  \sign_{\cons}({\bf N}) (-1)^{\sum\limits_{i=1}^{2j} \lfloor \frac{n_i}{2} \rfloor}$. 

%As in Case 1(b), $\psi$ is the identity map, so $M_{(1)}' = M_{(2)}'$.

%Clearly the entries of $M'$ are the same as the entries of $M$. We multiply $t'$ rows and columns by $-1$ to obtain $M'_{(1)}$. No relabeling is necessary. We add $R(x)$ and $C(x+1)$, which is the final row and column, respectively. Since $(-1)^{2n-1 > 2n} = 1$ and $\sign(2n-1, 2n) = 1$, we must multiply either both the last row and column by $-1$ or neither to obtain $M'_{(4)}$, which is equal to $\tilde{M}$. It follows that $t$ has the same parity as $t'$. 
\end{proof}

Now that we have established Lemma~\ref{lem:lifesaver}, the proof of Theorem~\ref{thm61} is very similar to Kenyon and Wilson's proof of Theorem~\ref{thm:kw61}.

\begin{proof}[Proof of Theorem~\ref{thm61}]{\color{white} a}
%(Similar to the proof of Theorem 6.1 in \cite{KW2009}).
%By Remark~\ref{rem:cyclicreordering}, 

%\vspace{1pt}

\noindent
\begin{minipage}{.65\textwidth}
\hspace{10pt}  Without loss of generality\footnotemark, we may assume that when we list the nodes in counterclockwise order starting with the red ones, they are in the order $1, 2, \ldots, 2n$. 
Combining Theorems \ref{mythm42} and \ref{KWthm31} immediately gives a Pfaffian formula for the double-dimer model. For example, let $G$ be a graph with eight nodes where nodes 1, 3, 4, and 6 are black. Assume the nodes are colored red, green and blue as shown to the right, so $\sigma = ((1, 8), (3, 4), (5, 6), (7, 2))$. Then by Theorem~\ref{KWthm31}, 
\end{minipage}
\begin{minipage}{.3\textwidth}
\begin{center}
 \begin{tikzpicture}[scale=.65]
  \draw (0,0) circle (2);
  \foreach \x in {1,2,...,8} {
   \node[shape=circle,fill=black, scale=0.5,label={{((\x-1)*360/8)+90}:\color{red}{\small{\x}}}] (n\x) at ({((\x-1)*360/8)+90}:2) {}; };
     \foreach \x in {1,2,3} {
   \node[shape=circle,fill=black, scale=0.5,label={{((\x-1)*360/8)+90}:\color{red}{\small{\x}}}] (n\x) at ({((\x-1)*360/8)+90}:2) {}; };
     \foreach \x in {4, 5} {
   \node[shape=circle,fill=black, scale=0.5,label={{((\x-1)*360/8)+90}:\color{green}{\small{\x}}}] (n\x) at ({((\x-1)*360/8)+90}:2) {}; };
        \foreach \x in {6, 7, 8} {
   \node[shape=circle,fill=black, scale=0.5,label={{((\x-1)*360/8)+90}:\color{blue}{\small{\x}}}] (n\x) at ({((\x-1)*360/8)+90}:2) {}; };
     \foreach \x in {2, 5, 7, 8} {
     \node[shape=circle,fill=white, scale=0.4] (n\x) at ({((\x-1)*360/8)+90}:2) {};
  };
  
        \foreach \x/\y in {1/8, 3/4, 5/6, 7/2} {
   \draw (n\x) -- (n\y);};
 \end{tikzpicture}
 \end{center}
 \end{minipage}

\begin{equation}
\label{eqn:Pfex}
\dddot{\Pr }(18 | 34 | 56 | 72) = 
\begin{pmatrix}
0 & 0 & 0 & L_{1, 4} & L_{1, 5} & L_{1, 6} & L_{1, 7} & L_{1, 8} \\
0 & 0 & 0 & L_{2, 4} & L_{2, 5} & L_{2, 6} & L_{2, 7} & L_{2, 8} \\
0 & 0 & 0 & L_{3, 4} & L_{3, 5} & L_{3, 6} & L_{3, 7} & L_{3, 8} \\
-L_{4, 1} & -L_{4, 2} & -L_{4, 3} & 0 & 0 & L_{4, 6} & L_{4, 7} & L_{4, 8} \\
-L_{5, 1} & -L_{5, 2} & -L_{5, 3} & 0 & 0 & L_{5, 6} & L_{5, 7} & L_{5, 8} \\
-L_{6, 1} & -L_{6, 2} & -L_{6, 3} & -L_{6, 4} &  -L_{6, 5} & 0 & 0& 0\\
-L_{7, 1} & -L_{7, 2} & -L_{7, 3} & -L_{7, 4} &  -L_{7, 5} & 0 & 0& 0\\
-L_{8, 1} & -L_{8, 2} & -L_{8, 3} & -L_{8, 4} &  -L_{8, 5} & 0 & 0& 0\\
\end{pmatrix}.
\end{equation}
%So by making the substitution in Theorem~\ref{mythm42}, we can get a Pfaffian formula for the polynomial $\widetilde{\text{Pr} }(18 | 34 | 56 | 72)$. 
So making the substitution in Theorem~\ref{mythm42} expresses 
$\widetilde{\Pr }(18 | 34 | 56 | 72)$ as a Pfaffian, up to a global sign.

Presently, we explain how we can obtain a determinant formula from this Pfaffian formula. We make the substitution $L_{i, j} \to 0$ when $i$ and $j$ are both black or both white and we reorder the rows and columns so the black nodes are listed first. In the example, the above matrix becomes 
\[
\begin{pmatrix}
0 & 0 & 0 & 0 & 0 & L_{1, 5} & L_{1, 7} & L_{1, 8} \\
0 & 0 & 0 & 0 & 0 & L_{3, 5} & L_{3, 7} & L_{3, 8} \\
0 & 0 & 0 & 0 & -L_{4, 2} & 0 & L_{4, 7} & L_{4, 8} \\
0 & 0 & 0 & 0 & -L_{6, 2} & -L_{6, 5} & 0 & 0 \\
0 & 0& L_{2, 4} &L_{2, 6} &  0& 0 & 0& 0\\
-L_{5, 1} & -L_{5, 3} & 0 & L_{5, 6} &  0& 0 & 0& 0\\
-L_{7, 1} & -L_{7, 3} & -L_{7, 4} & 0 &  0& 0 & 0& 0\\
-L_{8, 1} & -L_{8, 3} & -L_{8, 4} & 0 &  0& 0 & 0& 0\\
\end{pmatrix}.
\]

%old idea for global sign
%Let $(n_1, n_1 + 1), (n_2, n_2 + 1), \ldots, (n_{2k}, n_{2k} + 1)$ be a list of the couples of consecutive nodes of the same color. Let $h_1, \ldots, h_j$ be a complete list of all of the indices so that $n_{h_{i} - 1}$ and $n_{h_{i}}$ are different colors. Then the global sign is
%$$\prod_{i=1}^{j} (-1)^{ (n_{h_{i}} - n_{h_{i} - 1} )/2}$$
%because the idea is that we get the nodes to alternate by swapping $s_i$ with the closest $u_k$?
%Need to check this and see if there's a better way to write the global sign. 

\footnotetext{We can renumber the nodes while preserving their cyclic order without changing the global sign of the Pfaffian in 
Theorem~\ref{KWthm31}. 
This is because if we move the last row and column to be the first row and column, the sign of the Pfaffian changes. But since the entries above the diagonal must be non-negative, we negate the new first row and column and the Pfaffian changes sign again.}

Simultaneous swaps of two different rows and corresponding columns changes the sign of the Pfaffian. 
Assuming the graph has $2k$ couples of consecutive nodes of the same color, 
we claim that the number of swaps needed so that the black nodes are listed first has the same parity as
$$ \dfrac{n(n-1)}{2} + \sum\limits_{i=1}^{2k} \left\lfloor \frac{n_i}{2} \right\rfloor,$$
if node 1 is black. If node 1 is white, the number of swaps needed has the same parity as
$$ \dfrac{n(n+1)}{2} + \sum\limits_{i=1}^{2k} \left\lfloor \frac{n_i}{2} \right\rfloor.$$

%To prove this, we will first show that when we have $2k$ couples of consecutive nodes of the same color, the number of node swaps needed to get to a node coloring that alternates black and white has the same parity as
%$$\sum\limits_{i=1}^{2k} \left\lfloor \frac{n_i}{2} \right\rfloor.$$
To prove this, we will first show that 
%when we have $2k$ couples of consecutive nodes of the same color, 
the number of node swaps needed to get
from a node coloring with $2k$ couples of consecutive nodes of the same color
 to a node coloring that alternates black and white has the same parity as
$\sum\limits_{i=1}^{2k} \left\lfloor \frac{n_i}{2} \right\rfloor.$

We will prove this by induction on $k$. When $k = 0$, $0$ swaps are needed, so the claim holds trivially.
Assume the claim holds when ${\bf N}$ has $2(k-1)$ couples of consecutive nodes of the same color and suppose we have a set of nodes that has $2k$ couples of consecutive nodes of the same color. 
Let $h$ be the smallest integer so that $n_{h-1}$ and $n_h$ are different colors. Then $n_{h-1}$ and $n_h$ are the same parity and there are an even number of nodes in the interval $[n_{h-1} +1, \ldots, n_{h}]$, which alternate in color. If we swap $n_h$ with $n_{h} -1$, $n_{h} - 2$ with $n_{h} -3$ $,\ldots,$ $n_{h-1} + 2$ with $n_{h-1} +1$, we will have done 
$\frac{n_{h} - n_{h-1}}{2}$ swaps and we will have a node coloring with $2(k-1)$ couples of consecutive nodes of the same color.
If $n_h$ and $n_{h-1}$ are both even then $\frac{n_{h} - n_{h-1}}{2}$ clearly has the same parity as $\left\lfloor \frac{ n_{h} }{2} \right\rfloor + \left\lfloor \frac{ n_{h-1} }{2} \right\rfloor$. If $n_h$ and $n_{h-1}$ are both odd then by writing $\frac{n_{h} - n_{h-1}}{2} = 
 \frac{n_{h} -1 - (n_{h-1} -1)}{2}$ we see that $\frac{n_{h} - n_{h-1}}{2}$ and $\left\lfloor \frac{ n_{h} }{2} \right\rfloor + \left\lfloor \frac{ n_{h-1} }{2} \right\rfloor$ have the same parity. By the induction hypothesis, the number of swaps needed to get to a node coloring that alternates black and white has the same parity as
 $$\sum\limits_{\substack{ 1 \leq i \leq 2k \\ i \neq h, h-1} } \left\lfloor \frac{n_i}{2} \right\rfloor.$$
The claim follows.

Assume node 1 is black. 
If there are no couples of consecutive nodes of the same color, 
the number of swaps needed to put the black nodes first is
\[ 1 + 2 +  3 + \cdots + (n-1) = \dfrac{n(n-1)}{2} \]
because the third node requires 1 swap, the fifth node requires 2 swaps, the seventh node requires 3 swaps$,\ldots,$ and the $(2n-1)$st node requires $n-1$ swaps. 
So if there are $2k$ couples of consecutive nodes of the same color, since 
$\sum\limits_{i=1}^{2k} \left\lfloor \frac{n_i}{2} \right\rfloor$
node swaps are needed to get to a node coloring that alternates black and white,
the number of swaps needed so that the black nodes are listed first has the same parity as
$$ \dfrac{n(n-1)}{2} + \sum\limits_{i=1}^{2k} \left\lfloor \frac{n_i}{2} \right\rfloor.$$

If node 1 is white, 
the number of swaps needed to put the black nodes first is
\[ 1 + 2 +  3 + \cdots + n = \dfrac{n(n+1)}{2} \]
because the second node requires 1 swap, the fourth node requires 2 swaps$,\ldots,$ and the $(2n)$th node requires $n$ swaps. So the number of swaps needed so that the black nodes are listed first has the same parity as
$$ \dfrac{n(n+1)}{2} + \sum\limits_{i=1}^{2k} \left\lfloor \frac{n_i}{2} \right\rfloor.$$

Next, observe that after the rows and columns have been sorted, the matrix has the form
\[
\begin{pmatrix}
0  & \pm L_{B, W} \\
\mp L_{W, B} & 0 
\end{pmatrix}
\]
where $B$ represents the black nodes, $W$ the white nodes, and the
 signs of the entries in $\pm L_{B, W}$
 % individual signs 
  are $+$ if the black node has a smaller label than the white node and $-$ otherwise. 
The Pfaffian of this matrix is the determinant of the upper right submatrix times $(-1)^{\frac{n(n-1)}{2}}$. To summarize, after making the substitution $L_{i, j} \to 0$ when $i$ and $j$ are both black or both white and sorting the rows and columns so the black nodes are listed first,
\[
\text{Pf} 
\begin{pmatrix}
0 & L_{R, G} & L_{R, B} \\
-L_{G, R} & 0 & L_{G, B} \\
-L_{B, R} & -L_{B, G} & 0 
\end{pmatrix}
=
(-1)^{\sum\limits_{i=1}^{2k} \lfloor \frac{n_i}{2} \rfloor} 
\det
\begin{pmatrix}
L_{B, W}
\end{pmatrix}, 
\]
when node 1 is black. When node 1 is white, 
\[
\text{Pf} 
\begin{pmatrix}
0 & L_{R, G} & L_{R, B} \\
-L_{G, R} & 0 & L_{G, B} \\
-L_{B, R} & -L_{B, G} & 0 
\end{pmatrix}
=
(-1)^{n} (-1)^{\sum\limits_{i=1}^{2k} \lfloor \frac{n_i}{2} \rfloor} 
\det
\begin{pmatrix}
L_{B, W}
\end{pmatrix}.
\]

In the example, after this substitution and reordering, the Pfaffian of matrix (\ref{eqn:Pfex}) is equal to
\[
(-1)^{1 + 3}
\det
\begin{pmatrix}
0 & L_{1, 5} & L_{1, 7} & L_{1, 8} \\
0 & L_{3, 5} & L_{3, 7} & L_{3, 8} \\
 -L_{4, 2} & 0 & L_{4, 7} & L_{4, 8} \\
 -L_{6, 2} & -L_{6, 5} & 0 & 0 \\
\end{pmatrix}
\]
because $n_1 = 3$ and $n_2 = 7$. 

Next we do the substitution $L_{i, j} \to  \sign(i, j) Y_{i, j}$. 
The result is the matrix 
$$M =  [1_{i, j  \text{ RGB-colored differently}} (-1)^{i > j} \sign(i, j) Y_{i, j} ]^{i = b_1, b_2, \ldots, b_{n}}_{j = w_1, w_2, \ldots, w_{n} }$$
where $b_1 < b_2 < \cdots < b_n$ are the black nodes listed in increasing order and $w_1 < w_2 < \cdots < w_n$ are the white nodes listed in increasing order.
By Theorem~\ref{mythm42},
$$\widetilde{\Pr}(\sigma) = \sign_{OE}(\sigma) \sign_{\cons}({\bf N}) (-1)^{\sum\limits_{i=1}^{2k} \lfloor \frac{n_i}{2} \rfloor} \det(M)$$
if node 1 is black and 
$$\widetilde{\Pr}(\sigma) = \sign_{OE}(\sigma) \sign_{\cons}({\bf N}) (-1)^{n} (-1)^{\sum\limits_{i=1}^{2k} \lfloor \frac{n_i}{2} \rfloor} \det(M)$$
is node 1 is white. 
By Lemma~\ref{lem:lifesaver},
$M$ is a block matrix where within each block, the signs of the entries are staggered in a checkerboard pattern. 
%In this example, we get 
%\begin{eqnarray*}
%& & 
%\det
%\begin{pmatrix}
%0 & (-1)^{\frac{5-1}{2}} Y_{1, 5} &  (-1)^{\frac{7-1}{2}} Y_{1, 7} &  %(-1)^{\frac{(8-1) + 1}{2}} Y_{1, 8} \\
%0 & (-1)^{\frac{5-3}{2}} Y_{3, 5} &  (-1)^{\frac{7-3}{2}} Y_{3, 7} & %(-1)^{\frac{(8-3) + 1}{2}}  Y_{3, 8} \\
% -  (-1)^{\frac{4-2}{2}} Y_{4, 2} & 0 & (-1)^{\frac{(7-4) - 1}{2}}  Y_{4, 7} %&  (-1)^{\frac{8-4}{2}} Y_{4, 8} \\
% -  (-1)^{\frac{6-2}{2}}Y_{6, 2} & -  (-1)^{\frac{(6-5) - 1}{2}} Y_{6, 5} & %0 & 0 \\
%\end{pmatrix} \\
%&=&
%\[
% \det
%\left(
%\begin{array}{ c | c c c }
%0 & Y_{1, 5} & - Y_{1, 7} &    Y_{1, 8} \\
%0 & - Y_{3, 5} &  Y_{3, 7} & -  Y_{3, 8} \\
%  Y_{4, 2} & 0 & - Y_{4, 7} &   Y_{4, 8} \\
%  \hline
% - Y_{6, 2} & -  Y_{6, 5} & 0 & 0 \\
%\end{array}
%\right)
%\]
%\end{eqnarray*}
%Observe that this is a block matrix where the signs of the entries in each block are staggered in a checkerboard pattern. 
Next, we multiply rows and columns of $M$ by $-1$ so that the signs of the matrix entries are staggered in a checkerboard pattern and the upper left entry is positive. Call the resulting matrix $\widetilde{M}$. 
%We obtain the matrix $\widetilde{M}$ from $M$ by multiplying rows and columns by $-1$
%We can always choose rows and columns to multiply by $-1$ 
%Call the resulting matrix $\til{M}$. 
By Lemma \ref{lem:lifesaver}, 
$$\widetilde{\Pr}(\sigma) = \sign_{OE}(\sigma) \det(\widetilde{M})$$
regardless of whether node 1 is black or white.
Then, if we multiply every other row by $-1$ and every other column by $-1$, the signs of all matrix entries are positive and the determinant is unchanged. 
We conclude that
$$\widetilde{\Pr}(\sigma)= \sign_{OE}(\sigma) \det [1_{i, j \text{ RGB-colored differently } } Y_{i, j} ]^{i = b_1, b_2, \ldots, b_{n}}_{j = w_1, w_2, \ldots, w_{n} }.$$

%Indeed, by Lemma,  this is always the case. 
%Also by Lemma, 
%if $t$ be the total number of
%rows and columns we need to multiply by $-1$ to get a matrix with entries whose signs are staggered in a checkerboard pattern where the upper left entry is positive. If node 1 is black, 
%$$(-1)^{t} =  \sign_{\cons}({\bf N}) 
%(-1)^{\sum\limits_{i=1}^{2k} \lfloor \frac{n_i}{2} \rfloor} $$
%and if node 1 is white, 
%$$(-1)^{t} =  (-1)^{n} \sign_{\cons}({\bf N}) 
%(-1)^{\sum\limits_{i=1}^{2k} \lfloor \frac{n_i}{2} \rfloor}.$$
%We will prove that this is always the case. 

%Once we have a matrix with signs staggered in a checkerboard pattern so that the upper left entry is positive, we can multiply every other row and every other column by $-1$ to get a matrix with the property that the sign of each entry is positive.
%From Theorem \ref{mythm42}, we have a global sign of $\sign_{\cons}({\bf N} )\sign_{OE}(\sigma)$. The sign of $\sign_{OE}(\sigma)$ comes from the fact that after making this substitution described in Theorem \ref{mythm42} what we have is $\sign_{OE}(\sigma)$ times the double-dimer pairing polynomial. 

%All of the global signs cancel except for $\sign_{OE}(\sigma)$. 

Returning to our example, we find that 
\[\widetilde{ \Pr }(18 | 34 | 56 |72  ) = 
\sign_{OE}(18|34|56|72)
\det
\begin{pmatrix}
0 & Y_{1, 5} & Y_{1, 7} &    Y_{1, 8} \\
0 &  Y_{3, 5} &  Y_{3, 7} &   Y_{3, 8} \\
  Y_{4, 2} & 0 &  Y_{4, 7} &   Y_{4, 8} \\
  Y_{6, 2} &   Y_{6, 5} & 0 & 0 \\
\end{pmatrix}
=
-\det
\begin{pmatrix}
0 & Y_{1, 5} & Y_{1, 7} &    Y_{1, 8} \\
0 &  Y_{3, 5} &  Y_{3, 7} &   Y_{3, 8} \\
  Y_{4, 2} & 0 &  Y_{4, 7} &   Y_{4, 8} \\
  Y_{6, 2} &   Y_{6, 5} & 0 & 0 \\
\end{pmatrix}.
\]

\end{proof}

\subsection{Proof of Theorem~\ref{thm:cond}}

%\begin{example}
%\label{ex:cond}
%Let $G$ be a bipartite graph with a set {\bf N} of $2n$ nodes. 
%Assume the nodes of $G$ are black and odd or white and even. Assume also that the nodes are contiguously colored red, green, and blue (a color may occur zero times), and when we list the nodes in counterclockwise order starting with the red nodes, they are in the order $1, 2, \ldots, 2n$. Let $\sigma$ be the (unique) planar pairing in which like colors are not paired together. 

Now that we have established Theorem~\ref{thm61}, Theorem~\ref{thm:cond} follows from the proof method described in Section~\ref{sec:proofsketch}.

\begin{customthm}{\ref{thm:cond}}
Let $G = (V_1, V_2, E)$ be a finite edge-weighted planar bipartite graph with a set of nodes {\bf N}.
Divide the nodes into three circularly contiguous sets $R$, $G$, and $B$ such that $|R|, |G|,$ and $|B|$ satisfy the triangle inequality and let $\sigma$ be the corresponding tripartite pairing.
If $x, w \in V_1$ and $y, v \in V_2$ then
\begin{eqnarray*}
& & 
 \sign_{OE}(\sigma) \sign_{OE}(\sigma'_{xywv})Z^{DD}_{\sigma}(G, {\bf N}) Z^{DD}_{\sigma_{xywv}}(G, {\bf N} - \{x, y, w, v\}) \hspace{.4cm}\\ &=& 
\sign_{OE}(\sigma'_{xy}) \sign_{OE}(\sigma'_{wv})
Z^{DD}_{\sigma_{xy}}(G, {\bf N} - \{x, y\})  Z^{DD}_{\sigma_{wv}}(G, {\bf N} - \{w, v\}) \\ 
&& -  \sign_{OE}(\sigma'_{xv}) \sign_{OE}(\sigma'_{wy})
 Z^{DD}_{\sigma_{xv}}(G, {\bf N} - \{x, v\})  Z^{DD}_{\sigma_{wy}}(G, {\bf N} - \{w, y\}) 
 \end{eqnarray*}
  where for $i, j \in \{x, y, w, v\}$, $\sigma_{ij}$ is the unique planar pairing on ${\bf N} - \{i, j\}$ %corresponding node set 
 in which like colors are not paired together, and $\sigma_{ij}'$ is the pairing after the the node set ${\bf N} - \{i, j\}$ has been relabeled so that the nodes are numbered consecutively. 

\end{customthm}

\begin{proof}
First we assume that all pairings in the theorem statement exist. 
%Without loss of generality, assume $x < w$ and $y < v$.
Let 
$$M = [1_{i ,j \text{ RGB-colored differently } } Y_{i, j} ]^{i = b_1, b_2, \ldots, b_n}_{j = w_1, w_2, \ldots, w_n}$$ 
and
let $r_x$ and $r_w$ denote the rows corresponding to nodes $x$ and $w$, respectively (so we are assuming that $x$ and $w$ are colored black).

We first move the columns corresponding to $y$ and $v$ (i.e. the columns with entries $Y_{i, y}$ and $Y_{i, v}$, respectively) to the columns $r_x$ and $r_w$. We observe that we can do this without exchanging the column with entries $Y_{i, y}$ with the column with entries $Y_{i,v}$. 
For example, if $r_x < c_y < r_w < c_v$, we swap column $c_y$ with column $c_y - 1$, then column $c_y -1$ with column $c_y-2, \ldots,$ and column $r_{x} + 1$ with column $r_x$. 
%Let $s_{y}$ denote the number of column swaps we make in this process. 
Next, we swap column $c_v$ with column $c_v - 1, \ldots,$ and column $r_{w} + 1$ with column $r_w$. 
%Let $s_{v}$ denote the number of column swaps we make in this process. 
If instead $c_y < c_v < r_x < r_w$, we swap column $c_v$ with column $c_v +1, \ldots,$ and column $r_{w} - 1$ with column $r_w$ before swapping $c_y$ with column $c_y +1, \ldots,$ and column $r_{x} - 1$ with column $r_x$. 

Without loss of generality we assume that we move the column with entries $Y_{i, y}$ to column $r_x$ and the column with entries $Y_{i, v}$ to column $r_w$ to obtain the matrix $\widetilde{M}$. 
Let $s_{y}$ denote the number of column swaps we make in the process of moving the column with entries $Y_{i, y}$.
% to column $r_x$.
Let $s_{v}$ denote the number of column swaps we make in the process of moving the column with entries $Y_{i, v}$.
% to column $r_w$. 
 Note that $s_{y}$ and $s_v$ are well-defined up to parity. 
Note also that after making these swaps, the columns are still in ascending order, aside from the columns with entries $Y_{i, y}$ and $Y_{i, v}$.

%Note that the parity 
%of column swaps required to get $y$ in column $r_x$ is $\sign(y, x)$ and 
%the parity of 
%column swaps we need to get $v$ to be in column $r_w$ is $\sign(v, w)$. Let $\widetilde{M}$ be the matrix after making these swaps. 

By the Desnanot-Jacobi identity, 
\begin{equation}
\label{eqn:condpf}
\det(\widetilde{M}) \det(\widetilde{M}^{r_{x}, r_{w}}_{r_{x}, r_{w}}) = \det (\widetilde{M}_{r_{x}}^{r_{x}}) \det(\widetilde{M}_{r_{w}}^{r_{w}}) - \det(\widetilde{M}^{r_{x}}_{r_{w}}) \det(\widetilde{M}^{r_{w}}_{r_{x}}),
\end{equation}
where recall that $M_{s}^{t}$ is the matrix $M$ with row $s$ and column $t$ removed.

We apply Theorem \ref{thm61} to each term in equation (\ref{eqn:condpf}). 
First consider $\det (\widetilde{M}_{r_{x}}^{r_{x}})$. In order to apply Theorem~\ref{thm61} we must reorder the columns. Since we have removed the column $r_x$ which had entries $Y_{i, y}$, $s_{v}$ column swaps will put the columns in the correct (ascending) order. 
%The fact that $s_v$ column swaps will put the columns in the correct order follows 
This follows from the previous observation that we moved the columns corresponding to $y$ and $v$ without exchanging the column with entries $Y_{i,y}$ with the column with entries $Y_{i,v}$. 
%Thus the removal of the column $r_x$ does not affect the number of column swaps required to put the columns in ascending order. 

We must also relabel the nodes ${\bf N} - \{x, y \}$
so that they are numbered consecutively.
Recall that $\sigma_{xy}$ denotes the unique planar pairing of ${\bf N} - \{x, y \}$ in which like colors are not paired together. When we relabel ${\bf N} - \{x, y \}$ we relabel $\sigma_{xy}$ as well. Call the resulting node set ${\bf N'}$ and the resulting pairing $\sigma'_{xy}$. Then by Theorem~\ref{thm61}, 
$$\det (\widetilde{M}_{r_{x}}^{r_{x}}) 
= (-1)^{s_{v}} \sign_{OE}(\sigma'_{xy}) \dfrac{ Z^{DD}_{\sigma'_{xy}}(G, {\bf N'} )}{ (Z^{D}(G))^2 }, $$
and thus
$$\det (\widetilde{M}_{r_{x}}^{r_{x}}) 
=(-1)^{s_{v}} \sign_{OE}(\sigma'_{xy}) \dfrac{ Z^{DD}_{\sigma_{xy}}(G, {\bf N} - \{x, y \})}{ (Z^{D}(G))^2 }.  $$

%We apply Theorem \ref{thm61} to each term in equation (\ref{eqn:condpf}). 
%First consider $\det (\widetilde{M}_{r_{x}}^{r_{x}})$. In order to apply Theorem~\ref{thm61}, we must reorder the columns, which introduces a global sign of $\sign(v, w)$, and relabel the nodes ${\bf N} - \{x, y \}$
%so that they are numbered consecutively.
%Let $\sigma_1$ be the unique planar pairing of ${\bf N} - \{x, y \}$ in which like colors are not paired together. When we relabel ${\bf N} - \{x, y \}$ we relabel $\sigma_1$ as well. Call the resulting node set ${\bf N'}$ and the resulting pairing $\sigma'_1$. Then by Theorem~\ref{thm61}, 
%$$\det (\widetilde{M}_{r_{x}}^{r_{x}}) 
%= \sign(v, w) \sign_{OE}(\sigma'_1) \dfrac{ Z^{DD}_{\sigma'_1}(G, {\bf N'} )}{ (Z^{D}(G))^2 }.  $$
 %It is immediate that
%$$\det (\widetilde{M}_{r_{x}}^{r_{x}}) 
%= \sign(v, w) \sign_{OE}(\sigma'_1) \dfrac{ Z^{DD}_{\sigma_1}(G, {\bf N} - \{x, y \})}{ (Z^{D}(G))^2 }.  $$

%and comes from reordering the columns of
%$\widetilde{M}_{r_{x}}^{r_{x}}$
% so that they are in the order required by Theorem \ref{thm61}. 

%Note that since Theorem \ref{thm61} requires that the nodes are labeled consecutively, we are using the fact that the weighted sum of double dimer configurations on $(G, {\bf N} - \{x, y \})$ is equal to the weighted sum of double dimer configurations on $(G, {\bf N'})$,  where ${\bf N'} = \{1, 2, \ldots, 2n-2 \}$. 

Similarly, we have
%\begin{equation}
\begin{equation}
\label{eqn:RHS}
\begin{split}
\det (\widetilde{M}_{r_{w}}^{r_{w}}) 
&= (-1)^{s_{y}}
 \sign_{OE}(\sigma'_{wv}) \dfrac{ Z^{DD}_{\sigma_{wv}}(G, {\bf N} - \{w, v \})}{ (Z^{D}(G))^2 }, \\
\det (\widetilde{M}_{r_{x}}^{r_{w}}) 
&=  (-1)^{s_{y}}
 \sign_{OE}(\sigma'_{xv}) \dfrac{ Z^{DD}_{\sigma_{xv}}(G, {\bf N} - \{x, v \})}{ (Z^{D}(G))^2 }, \text{ and} \\
\det (\widetilde{M}_{r_{w}}^{r_{x}}) 
&=  (-1)^{s_{v}}
 \sign_{OE}(\sigma'_{yw}) \dfrac{ Z^{DD}_{\sigma_{yw}}(G, {\bf N} - \{y, w \})}{ (Z^{D}(G))^2 }. 
\end{split}
\end{equation}
It follows that the right hand side of equation (\ref{eqn:condpf}) is 
\begin{eqnarray*}
& &
% \sign(y, x)
%\sign(v, w)
% \cdot 
  (-1)^{s_{y}}
    (-1)^{s_{v}}
 \bigg{(}
\sign_{OE}(\sigma'_{xy}) \sign_{OE}(\sigma'_{wv})\dfrac{ Z^{DD}_{\sigma_{xy}}(G, {\bf N} - \{x, y \})}{ (Z^{D}(G))^2 }
\dfrac{ Z^{DD}_{\sigma_{wv}}(G, {\bf N} - \{w, v \})}{ (Z^{D}(G))^2 }\\
 & & - \sign_{OE}(\sigma'_{xv}) \sign_{OE}(\sigma'_{yw})\dfrac{Z^{DD}_{\sigma_{xv}}(G, {\bf N} - \{x, v \})}{ (Z^{D}(G))^2 } 
\dfrac{ Z^{DD}_{\sigma_{yw}}(G, {\bf N} - \{y, w \})}{ (Z^{D}(G))^2 } \bigg{). }
\end{eqnarray*}
Applying Theorem~\ref{thm61} to the left hand side of equation (\ref{eqn:condpf}), we have
\begin{equation}
\label{eqn:LHS}
\begin{split}
  \det(\widetilde{M})   &=  
  (-1)^{s_{y}}
    (-1)^{s_{v}}
 \sign_{OE}(\sigma) \dfrac{Z^{DD}_{\sigma}(G, {\bf N})}{ (Z^{D}(G))^2 }, \text{ and} \\
\det(\widetilde{M}^{r_{x}, r_{w}}_{r_{x}, r_{w}}) 
&=
 \sign_{OE}(\sigma'_{xywv})
\dfrac{ Z^{DD}_{\sigma_{xywv}}(G, {\bf N} - \{x, y, w, v\})}{ (Z^{D}(G))^2 }.
\end{split}
\end{equation}

We conclude that
%\begin{eqnarray*}
%& &   \sign_{OE}(\sigma) \sign_{OE}(\sigma'_5)Z^{DD}_{\sigma}(G, {\bf N}) Z^{DD}_{\sigma_5}(G, {\bf N} - \{x, y, w, v\}) \\ &=& 
%\sign_{OE}(\sigma'_1) \sign_{OE}(\sigma'_2)
%Z^{DD}_{\sigma_1}(G, {\bf N} - \{x, y\})  Z^{DD}_{\sigma_2}(G, {\bf N} - \{w, v\}) \\
%& & -  \sign_{OE}(\sigma'_3) \sign_{OE}(\sigma'_4)
% Z^{DD}_{\sigma_3}(G, {\bf N} - \{x, v\})  Z^{DD}_{\sigma_4}(G, {\bf N} - \{w, y\}) 
% \end{eqnarray*}
\begin{eqnarray*}
& & 
 \sign_{OE}(\sigma) \sign_{OE}(\sigma'_{xywv})Z^{DD}_{\sigma}(G, {\bf N}) Z^{DD}_{\sigma_{xywv}}(G, {\bf N} - \{x, y, w, v\}) \hspace{.4cm}\\ &=& 
\sign_{OE}(\sigma'_{xy}) \sign_{OE}(\sigma'_{wv})
Z^{DD}_{\sigma_{xy}}(G, {\bf N} - \{x, y\})  Z^{DD}_{\sigma_{wv}}(G, {\bf N} - \{w, v\}) \\ 
&& -  \sign_{OE}(\sigma'_{xv}) \sign_{OE}(\sigma'_{wy})
 Z^{DD}_{\sigma_{xv}}(G, {\bf N} - \{x, v\})  Z^{DD}_{\sigma_{wy}}(G, {\bf N} - \{w, y\}).
 \end{eqnarray*}

It is not necessarily the case that the pairings $\sigma_{xy},\sigma_{wv}, \sigma_{xv}, \sigma_{wy}$, and $\sigma_{xywv}$ all exist. First consider the case where
one of the pairings
$\sigma_{xy},\sigma_{wv}, \sigma_{xv}, \sigma_{wy}$ does not exist.
% $\sigma_i$ does not exist for $i \in \{1, 2, 3, 4\}$. 
 Without loss of generality, assume that $\sigma_{xy}$ does not exist. This means that the number of nodes of different colors in ${\bf N} - \{x, y\}$ do not satisfy the triangle inequality. Then $\det(\widetilde{M}_{r_{x}}^{r_{x}}) = 0$ since every black-white pairing contains an $RGB$-monochromatic pair. There are two possibilities in this case: either the theorem statement holds trivially, or
\begin{equation}
\label{eqn:reallyspecialcase}
Z^{DD}_{\sigma}(G, {\bf N}) Z^{DD}_{\sigma_{xywv}}(G, {\bf N} - \{x, y, w, v\}) =
 Z^{DD}_{\sigma_{xv}}(G, {\bf N} - \{x, v\})  Z^{DD}_{\sigma_{wy}}(G, {\bf N} - \{w, y\}).
 \end{equation}

Since the numbers of nodes of different colors in ${\bf N} - \{x, y\}$ do not satisfy the triangle inequality, 
without loss of generality, we may assume there are more red nodes than the combined number of blue and green nodes in ${\bf N} - \{x, y\}$.
Since we assumed that in ${\bf N}$, $|R|, |G|,$ and $|B|$ satisfy the triangle inequality, it follows that $|R| = |G| + |B|$ in ${\bf N}$. 
%this is only possible if in ${\bf N}$ the number of red nodes is equal to the total number of green and blue nodes. 
Assuming without loss of generality that when we list the nodes in counterclockwise order starting with the red ones, they are in the order $1, 2, \ldots, 2n$, this means $\sigma$ is the pairing $((1, 2n), (2, 2n-1), \ldots, (n, n+1))$. It must be the case that $x$ and $y$ are both green or blue, so if either $w$ or $v$ is green, then $\sigma_{xywv}$ does not exist and either $\sigma_{xv}$ or $\sigma_{wv}$ does not exist, so the equality holds trivially. If both $w$ and $v$ are red, then $\sigma_{xv}$, $\sigma_{wv}$ and $\sigma_{xywv}$ all exist.
 In this case, the rest of the proof proceeds as above and we have
 \begin{eqnarray*}
  & & \sign_{OE}(\sigma) \sign_{OE}(\sigma'_{xywv})Z^{DD}_{\sigma}(G, {\bf N}) Z^{DD}_{\sigma_{xywv}}(G, {\bf N} - \{x, y, w, v\}) \\
  &=& -  \sign_{OE}(\sigma'_{xv}) \sign_{OE}(\sigma'_{wy})
 Z^{DD}_{\sigma_{xv}}(G, {\bf N} - \{x, v\})  Z^{DD}_{\sigma_{wy}}(G, {\bf N} - \{w, y\}).
 \end{eqnarray*}
Recall that inversions in a planar pairing correspond to nestings (see Remark~\ref{rem:inversionsarenestings}). 
Because $\sigma$ is the pairing $((1, 2n), (2, 2n-1), \ldots, (n, n+1))$, $\sign(\sigma'_{xv}) = \sign(\sigma'_{wy})$ and 
$\sign(\sigma'_{xywv}) = \sign(\sigma) \cdot (-1)^{n-1} \cdot (-1)^{n-2}$. Equation (\ref{eqn:reallyspecialcase}) follows.

 If $\sigma_{xywv}$ does not exist, then this means that the numbers of nodes of different colors in ${\bf N} - \{x, y, w, v\}$ do not satisfy the triangle inequality. 
Without loss of generality, we may assume there are more red nodes than the combined number of blue and green nodes in ${\bf N} - \{x, y, w, v\}$. By the same reasoning as above, 
$\det(\widetilde{M}_{r_{x}, r_{w}}^{r_{x}, r_{w}}) = 0$.
%The set $\{x, y, w, v\}$ must contain more green and blue nodes than red nodes, otherwise the original set of nodes ${\bf N}$ would not satisfy the triangle inequality. 
%There are two possibilities. Either $\{x, y, w, v\}$ consists only of green and blue nodes, in which case
%Or $\{x, y, w, v\}$ contains one red node. 
%There are two cases.
%There are either four more red nodes than the combined %number of blue and green nodes in ${\bf N} - \{x, y, w, v\}$, or there are two more red nodes than the combined number of blue and green nodes in ${\bf N} - \{x, y, w, v\}$.
% If there are four more red nodes than the combined number of blue and green nodes, 
%  the equation holds trivially. 
%If there are two more red nodes than the combined number %of blue and green nodes, 
%there are two possibilities. 
 If any one of $x, y, w,$ or $v$ is red, then the equation holds trivially.
If all of $x, y, w$, and $v$ are green or blue, then in the original node set ${\bf N}$, $|R| +2= |G| + |B|$.
% If all of $x, y, w,$ and $v$ are green or blue, then 
 So each of the pairings $\sigma'_{xy}, \sigma'_{wv},  \sigma'_{xv}$, and $\sigma'_{wy}$ is $((1, 2n-2), (2, 2n-1), \ldots)$. Then we have
\begin{eqnarray*}
& & \sign_{OE}(\sigma'_{xy}) \sign_{OE}(\sigma'_{wv}) Z^{DD}_{\sigma_{xy}}(G, {\bf N} - \{x, y\})  Z^{DD}_{\sigma_{wv}}(G, {\bf N} - \{w, v\}) \\
&=  &
\sign_{OE}(\sigma'_{xv}) \sign_{OE}(\sigma'_{wy}) Z^{DD}_{\sigma_{xv}}(G, {\bf N} - \{x, v\})  Z^{DD}_{\sigma_{wy}}(G, {\bf N} - \{w, y\}).
\end{eqnarray*}
Since inversions in a planar pairing correspond to nestings, all pairings have the same sign. So in this case, 
\[
Z^{DD}_{\sigma_{xy}}(G, {\bf N} - \{x, y\})  Z^{DD}_{\sigma_{wv}}(G, {\bf N} - \{w, v\}) 
= 
 Z^{DD}_{\sigma_{xv}}(G, {\bf N} - \{x, v\})  Z^{DD}_{\sigma_{wy}}(G, {\bf N} - \{w, y\}) .
\]

\end{proof}

\begin{rem}
\label{rem:caseanalysis}
To simplify the expression in Theorem~\ref{thm:cond}, it suffices to know the RGB-coloring of the nodes $x, y, w$, $v$. 

Without loss of generality, assume that when we list the nodes in counterclockwise order starting with the red ones, they are in the order $1, 2, \ldots, 2n$. 
Let $|RG(\sigma)|$ be the number of red-green pairs in $\sigma$. Define $|GB(\sigma)|$ and $|RB(\sigma)|$ similarly. 
Assume that $|RG(\sigma)|, |GB(\sigma)|,$ and $|RB(\sigma)|$ are all nonzero. 

If the set of nodes $\{x, y\}$ contains one red node and one blue node, then
 $\sigma_{xy}$ has one fewer red-blue pair than $\sigma$, but the number of red-green and green-blue pairs is the same (see Figure~\ref{fig:specialcase}). 
 By Remark~\ref{rem:inversionsarenestings}, to determine the relationship between $\sign_{OE}(\sigma)$ and $\sign_{OE}(\sigma'_{xy})$ it suffices to count the number of nestings in the diagram of $\sigma$ that involve a red-blue pair $(n_r, n_b)$.
 There is one nesting for each red-blue pair other than $(n_r, n_b)$, one nesting for each red-green pair, and one nesting for each green-blue pair (see Figure~\ref{fig:specialcase1}). Therefore, 
$$\sign(\sigma'_{xy}) = \sign_{OE}(\sigma) \cdot (-1)^{|RG(\sigma)|} \cdot (-1)^{|GB(\sigma)|} \cdot  (-1)^{|RB(\sigma)|-1}. $$

\begin{figure}[htb]
\centering
 \begin{tikzpicture}[scale=.75]
  \draw (0,0) circle (2);
  \foreach \x in {1,2,...,16} {
   \node[shape=circle,fill=black, scale=0.5,label={{((\x-1)*360/16)+90}:\x}] (n\x) at ({((\x-1)*360/16)+90}:2) {}; };
     \foreach \x in {1, 2, 3, 4, 5} {
     \node[shape=circle,fill=red, scale=0.4] (n\x) at ({((\x-1)*360/16)+90}:2) {};
     };
        \foreach \x in {6, 7, 8, 9, 10, 11} {
     \node[shape=circle,fill=green, scale=0.4] (n\x) at ({((\x-1)*360/16)+90}:2) {};
     };
               \foreach \x in {12, 13, 14, 15, 16} {
     \node[shape=circle,fill=blue, scale=0.4] (n\x) at ({((\x-1)*360/16)+90}:2) {};
     };
       \foreach \x/\y in {1/16, 3/8, 5/6, 7/4, 9/14, 11/12, 13/10, 15/2} {
   \draw (n\x) -- (n\y);
  };
 \end{tikzpicture}
  \begin{tikzpicture}[scale=.75]
  \draw (0,0) circle (2);
  \foreach \x in {1,2,3} {
   \node[shape=circle,fill=black, scale=0.5,label={{((\x-1)*360/14)+90}:\x}] (n\x) at ({((\x-1)*360/14)+90}:2) {}; };

   \node[shape=circle,fill=black, scale=0.5,label={{((4-1)*360/14)+90}:5}] (n4) at ({((4-1)*360/14)+90}:2) {}; 
   \node[shape=circle,fill=black, scale=0.5,label={{((5-1)*360/14)+90}:6}] (n5) at ({((5-1)*360/14)+90}:2) {}; 
      \node[shape=circle,fill=black, scale=0.5,label={{((6-1)*360/14)+90}:7}] (n6) at ({((6-1)*360/14)+90}:2) {}; 
      \node[shape=circle,fill=black, scale=0.5,label={{((7-1)*360/14)+90}:8}] (n7) at ({((7-1)*360/14)+90}:2) {}; 
            \node[shape=circle,fill=black, scale=0.5,label={{((8-1)*360/14)+90}:9}] (n8) at ({((8-1)*360/14)+90}:2) {}; 
           \node[shape=circle,fill=black, scale=0.5,label={{((9-1)*360/14)+90}:10}] (n9) at ({((9-1)*360/14)+90}:2) {}; 
     \node[shape=circle,fill=black, scale=0.5,label={{((10-1)*360/14)+90}:11}] (n10) at ({((10-1)*360/14)+90}:2) {}; 
          \node[shape=circle,fill=black, scale=0.5,label={{((11-1)*360/14)+90}:13}] (n11) at ({((11-1)*360/14)+90}:2) {}; 
      \node[shape=circle,fill=black, scale=0.5,label={{((12-1)*360/14)+90}:14}] (n12) at ({((12-1)*360/14)+90}:2) {}; 
  \node[shape=circle,fill=black, scale=0.5,label={{((13-1)*360/14)+90}:15}] (n13) at ({((13-1)*360/14)+90}:2) {}; 
   \node[shape=circle,fill=black, scale=0.5,label={{((14-1)*360/14)+90}:16}] (n14) at ({((14-1)*360/14)+90}:2) {}; 
   
     \foreach \x in {1, 2, 3, 4} {
     \node[shape=circle,fill=red, scale=0.4] (n\x) at ({((\x-1)*360/14)+90}:2) {};
     };
        \foreach \x in {5, 6, 7, 8, 9, 10} {
     \node[shape=circle,fill=green, scale=0.4] (n\x) at ({((\x-1)*360/14)+90}:2) {};
     };
               \foreach \x in {11, 12, 13, 14} {
     \node[shape=circle,fill=blue, scale=0.4] (n\x) at ({((\x-1)*360/14)+90}:2) {};
     };
       \foreach \x/\y in {1/14, 3/6, 5/4, 7/2, 9/12, 11/10, 13/8} {
   \draw (n\x) -- (n\y);
  };
 \end{tikzpicture}
   \begin{tikzpicture}[scale=.75]
  \draw (0,0) circle (2);
  \foreach \x in {1,2,...,6} {
   \node[shape=circle,fill=black, scale=0.5,label={{((\x-1)*360/14)+90}:\x}] (n\x) at ({((\x-1)*360/14)+90}:2) {}; };

    \node[shape=circle,fill=black, scale=0.5,label={{((7-1)*360/14)+90}:8}] (n7) at ({((7-1)*360/14)+90}:2) {}; 
            \node[shape=circle,fill=black, scale=0.5,label={{((8-1)*360/14)+90}:9}] (n8) at ({((8-1)*360/14)+90}:2) {}; 
           \node[shape=circle,fill=black, scale=0.5,label={{((9-1)*360/14)+90}:10}] (n9) at ({((9-1)*360/14)+90}:2) {}; 
     \node[shape=circle,fill=black, scale=0.5,label={{((10-1)*360/14)+90}:12}] (n10) at ({((10-1)*360/14)+90}:2) {}; 
          \node[shape=circle,fill=black, scale=0.5,label={{((11-1)*360/14)+90}:13}] (n11) at ({((11-1)*360/14)+90}:2) {}; 
      \node[shape=circle,fill=black, scale=0.5,label={{((12-1)*360/14)+90}:14}] (n12) at ({((12-1)*360/14)+90}:2) {}; 
  \node[shape=circle,fill=black, scale=0.5,label={{((13-1)*360/14)+90}:15}] (n13) at ({((13-1)*360/14)+90}:2) {}; 
   \node[shape=circle,fill=black, scale=0.5,label={{((14-1)*360/14)+90}:16}] (n14) at ({((14-1)*360/14)+90}:2) {}; 
   
     \foreach \x in {1, 2, 3, 4, 5} {
     \node[shape=circle,fill=red, scale=0.4] (n\x) at ({((\x-1)*360/14)+90}:2) {};
     };
        \foreach \x in {6, 7, 8, 9} {
     \node[shape=circle,fill=green, scale=0.4] (n\x) at ({((\x-1)*360/14)+90}:2) {};
     };
               \foreach \x in {10, 11, 12, 13, 14} {
     \node[shape=circle,fill=blue, scale=0.4] (n\x) at ({((\x-1)*360/14)+90}:2) {};
     };
       \foreach \x/\y in {1/14, 3/12, 5/6, 7/4, 9/10, 11/8, 13/2} {
   \draw (n\x) -- (n\y);
  };
 \end{tikzpicture}
 \caption[The tripartite pairing on ${\bf N}$ $-$ $\{4, 12\}$ compared to the pairing on ${\bf N}$.]
 {The tripartite pairing on ${\bf N}$ $-$ $\{4, 12\}$ (shown center) has one fewer red-blue pair than  the tripartite pairing on ${\bf N}$ (shown left). The tripartite pairing on ${\bf N}$ $-$ $\{7, 11\}$ (shown right) has one fewer red-green pair, one fewer green-blue pair, and one more red-blue pair compared to the pairing on ${\bf N}$. }
 
% Left: The pairing $\sigma$ of the nodes ${\bf N}$. Center: The pairing $\sigma_6$, which is the unique planar pairing on ${\bf N} - \{4, 12\}$ in which like colors are not paired together. Right: The pairing $\sigma_7$, which is the unique planar pairing on ${\bf N} - \{7, 11\}$ in which like colors are not paired together.} 
 \label{fig:specialcase}
 \end{figure}
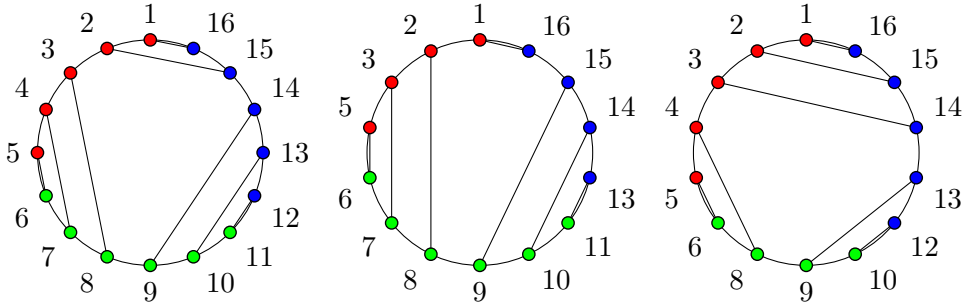

If the set of nodes $\{x, y\}$ contains one red node and one green node, then
$\sigma_{xy}$ has one fewer red-green pair than $\sigma$, but the number of 
red-blue and green-blue pairs is the same.
% By Remark~\ref{rem:inversionsarenestings}, to determine the relationship between $\sign_{OE}(\sigma)$ and $\sign_{OE}(\sigma'_{xy})$ it suffices to count the number of nestings in the diagram of $\sigma$ that involve a red-green pair $(n_r, n_g)$. 
  So we count the number of nestings in the diagram of $\sigma$ that involve a red-green pair $(n_r, n_g)$. 
 There is one nesting for each red-green pair other than $(n_r, n_g)$, and one nesting for each red-blue pair.
Therefore, 
$$\sign_{OE}(\sigma'_{xy}) = \sign_{OE}(\sigma) \cdot (-1)^{|RB(\sigma)|} \cdot (-1)^{|RG(\sigma)| - 1}.$$

Similarly, if the set of nodes $\{x, y\}$ contains one green node and one blue node, then
$$\sign_{OE}(\sigma'_{xy}) = \sign_{OE}(\sigma) \cdot (-1)^{|RB(\sigma)|} \cdot (-1)^{|GB(\sigma)| - 1}.$$

\begin{figure}[htb]
 \centering
 \begin{tikzpicture}[scale = .75]
	\vertex[fill = red] (n1) at (1, 0) [label=below:1] {};
	\vertex[fill = red] (n2) at (2,0) [label=below:2] {};
	\vertex[fill = red] (n3) at (3,0) [label=below:3] {};
	\vertex[fill = red] (n4) at (4, 0) [label=below:4] {};
	\vertex[fill = red] (n5) at (5,0) [label=below:5] {};
	\vertex[fill = green] (n6) at (6,0) [label=below:6] {};
	\vertex[fill = green] (n7) at (7,0) [label=below:7] {};
	\vertex[fill = green] (n8) at (8, 0) [label=below:8] {};
	\vertex[fill = green] (n9) at (9, 0) [label=below:9] {};
	\vertex[fill = green] (n10) at (10, 0) [label=below:10] {};
	\vertex[fill = green] (n11) at (11, 0) [label=below:11] {};
	\vertex[fill = blue] (n12) at (12, 0) [label=below:12] {};
	\vertex[fill = blue] (n13) at (13, 0) [label=below:13] {};
	\vertex[fill = blue] (n14) at (14, 0) [label=below:14] {};
	\vertex[fill = blue] (n15) at (15, 0) [label=below:15] {};
	\vertex[fill = blue] (n16) at (16, 0) [label=below:16] {};
	
	%arcs of pi
	\draw  (n16) arc (0:180:7.5cm);
	\draw[thick,dotted]  (n15) arc (0:180:6.5cm);
	\draw  (n14) arc (0:180:2.5cm);
	\draw  (n13) arc (0:180:1.5cm);
	\draw  (n12) arc (0:180:0.5cm);
	\draw  (n8) arc (0:180:2.5cm);
	\draw  (n7) arc (0:180:1.5cm);
	\draw  (n6) arc (0:180:0.5cm);
	
	\vertex[fill = red] (n1) at (1, 0) {};
	\vertex[fill = red] (n2) at (2,0)  {};
	\vertex[fill = red] (n3) at (3,0) {};
	\vertex[fill = red] (n4) at (4, 0) {};
	\vertex[fill = red] (n5) at (5,0) {};
	\vertex[fill = green] (n6) at (6,0) {};
	\vertex[fill = green] (n7) at (7,0) {};
	\vertex[fill = green] (n8) at (8, 0) {};
	\vertex[fill = green] (n9) at (9, 0) {};
	\vertex[fill = green] (n10) at (10, 0) {};
	\vertex[fill = green] (n11) at (11, 0) {};
	\vertex[fill = blue] (n12) at (12, 0) {};
	\vertex[fill = blue] (n13) at (13, 0) {};
	\vertex[fill = blue] (n14) at (14, 0) {};
	\vertex[fill = blue] (n15) at (15, 0) {};
	\vertex[fill = blue] (n16) at (16, 0) {};

\end{tikzpicture}
\caption[A red-blue pair is removed from $\sigma$.]
{When a red-blue pair is removed from $\sigma$, the number of nestings in the diagram of $\sigma$ decreases by $|RB(\sigma)| -1 + |RG(\sigma)| + |GB(\sigma)|$.}
%\caption{The pairing $\sigma =((1, 12), (3, 6), (5,4), (7, 10), (9, 8), (11, 2))$ is shown above. Let $\sigma'$ be the pairing obtained by removing the pair $(4, 7)$ and relabeling. To determine the relationship between $\sign_{OE}(\sigma)$ and $\sign_{OE}(\sigma')$ we count the number of nestings that include the pair $(4, 7)$.}
\label{fig:specialcase1}
\end{figure}
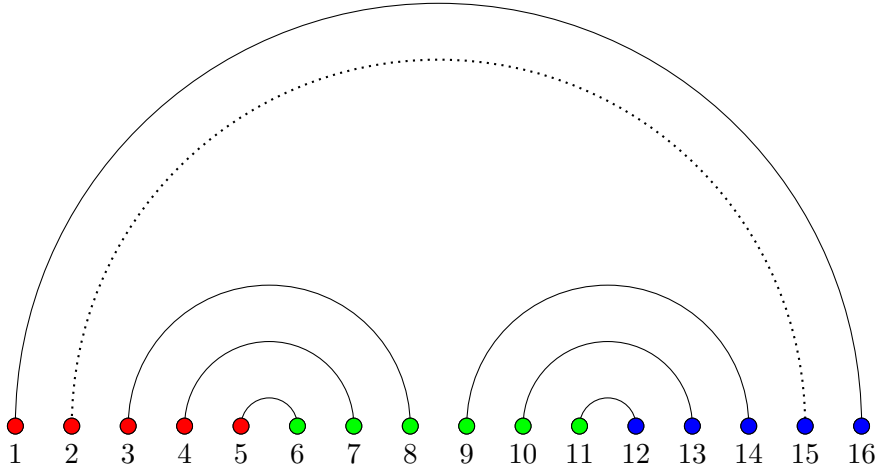

If both $x$ and $y$ are green nodes, then $\sigma_{xy}$
has one fewer red-green pair, one fewer green-blue pair, and one more red-blue pair, as shown in Figure \ref{fig:specialcase}.
%which like colors are not paired. Note that ${\bf N} - \{t_2, t_3\}$ has two fewer green nodes. It follows that $\sigma_{t_{2} t_{3}}$ has one more red-blue pair than $\sigma$, one fewer red-green pair, and one fewer green-blue pair, as shown in Figure \ref{fig:specialcase}.
 Removing a red-green pair from $\sigma$ removes $|RB(\sigma)| + |RG(\sigma)| -1$ nestings. Then, removing a green-blue pair removes $|RB(\sigma)| + |GB(\sigma)| -1$ nestings. After these pairs have been removed, adding a red-blue pair results in $|RB(\sigma)| + |GB(\sigma)|-1 + |RG(\sigma)|-1$ additional nestings. 
Therefore, 
$$ \sign_{OE}(\sigma'_{xy}) = \sign_{OE}(\sigma) \cdot (-1)^{|RB(\sigma)|}.$$

If both $x$ and $y$ are red nodes, $\sigma_{xy}$
has one fewer red-blue pair, one fewer red-green pair, and one more green-blue pair. Removing a red-blue pair from $\sigma$ removes $|RG(\sigma)| + |GB(\sigma)| + |RB(\sigma)| -1$ nestings. Then, removing a red-green pair removes $|RB(\sigma)| -1 + |RG(\sigma)| -1$ nestings. After these pairs have been removed, adding a green-blue pair results in $|RB(\sigma)| -1  + |GB(\sigma)|$ additional nestings. Thus
$$ \sign_{OE}(\sigma'_{xy}) = \sign_{OE}(\sigma) \cdot (-1)^{|RB(\sigma)|}.$$

Similarly, if both $x$ and $y$ are blue, 
$$ \sign_{OE}(\sigma'_{xy}) = \sign_{OE}(\sigma) \cdot (-1)^{|RB(\sigma)|}.$$
\end{rem}

If we assume that the nodes $x, y, w, v$ alternate black and white and the set $\{x,y,w,v\}$ contains at least one node of each RGB color, we can use Remark~\ref{rem:caseanalysis} to show that the all the signs in Theorem~\ref{thm:cond} are positive.
% signs simplify nicely in Theorem~\ref{thm:cond}. 
%Let $x, y, w, v$ be nodes appearing in a cyclic order such that at least one of these four nodes is red, at least one is green, and at least one is blue.

\begin{customthm}{\ref{cor:cond}}
%Let $G = (V_1, V_2, E)$ be a planar bipartite graph with a set of nodes {\bf N}. 
Divide the nodes into three circularly contiguous sets $R$, $G$, and $B$ such that $|R|, |G|$ and $|B|$ satisfy the triangle inequality and let $\sigma$ be the corresponding tripartite pairing. 
Let $x, y, w, v$ be nodes appearing in a cyclic order such that the set
$\{x,y,w,v\}$ contains at least one node of each RGB color.
% such that $\{x, y\}$ and $\{v, w\}$ are pairs of $\sigma$.
If $x$ and $w$ are both black and $y$ and $v$ are both white, then
% If $x, w \in V_1$ and $y, v \in V_2$ then
%Contiguously color the nodes red, green and blue using the minimum number of colors so that no pair contains nodes of the same color, and assume that the minimum number of colors required is 3. 
%Assume that the nodes are contiguously colored red, green, and blue (a color may occur zero times). 
%Assume the nodes are are numbered consecutively in counterclockwise order starting with the red nodes so that node 1 is the first red node.
%Let the nodes $x,y, w, v$ appear in a cyclic order. If $x, w$ are black and $y, v$ are white, and if in addition $\sigma(x) = y$ and $\sigma(w) = v$ then
\begin{eqnarray*}
Z^{DD}_{\sigma}(G, {\bf N}) Z^{DD}_{\sigma_{xywv}}(G, {\bf N} - \{x, y, w, v\}) &=& 
Z^{DD}_{\sigma_{xy}}(G, {\bf N} - \{x, y\})  Z^{DD}_{\sigma_{wv}}(G, {\bf N} - \{w, v\}) \\
& &  +
 Z^{DD}_{\sigma_{xv}}(G, {\bf N} - \{x, v\})  Z^{DD}_{\sigma_{wy}}(G, {\bf N} - \{w, y\}).
 \end{eqnarray*}
 \normalsize
 %where $\sigma_{xy}, $ is the unique planar pairing on the corresponding node set in which like colors are not paired together. 
 \end{customthm}
%Let $G$ be a bipartite graph with a set {\bf N} of $2n$ nodes. Fix a tripartite pairing $\sigma$. Contiguously color the nodes red, green and blue using the minimum number of colors so that no pair contains nodes of the same color, and assume that the minimum number of colors required is 3. 
%Let the nodes $x,y, w, v$ appear in a cyclic order. If $x, w$ are black and $y, v$ are white, and if in addition $\sigma(x) = y$ and $\sigma(w) = v$ then
%\begin{eqnarray*}
%& & 
%Z^{DD}_{\sigma}(G, {\bf N}) Z^{DD}_{\sigma_5}(G, {\bf N} - \{x, y, w, v\}) \hspace{.4cm}\\ &=& 
%Z^{DD}_{\sigma_1}(G, {\bf N} - \{x, y\})  Z^{DD}_{\sigma_2}(G, {\bf N} - \{w, v\})  +
 %Z^{DD}_{\sigma_3}(G, {\bf N} - \{x, v\})  Z^{DD}_{\sigma_4}(G, {\bf N} - \{w, y\}) 
 %\end{eqnarray*}
% where $\sigma_i$ is the unique planar pairing on the corresponding node set in which like colors are not paired together. 

\begin{proof}
Without loss of generality, assume that when we list the nodes in counterclockwise order starting with the red ones, they are in the order $1, 2, \ldots, 2n$. Assume also that one of the nodes $x, y, w, v$ is red, two are green, and one is blue. The other cases are very similar. 
 By the assumption that the nodes are in cyclic order,
 % this means that either
%\begin{itemize}
%\item $x$ is red, $y$ and $w$ are green, and $v$ is blue
%\item $x$ is blue, $y$ is red, and $w$ and $v$ are green
%\item $x$ is green, $y$ is blue, $w$ is red, and $v$ is green
%\item $x$ and $y$ are green, $w$ is blue, and $v$ is red
%\end{itemize}
%
%Since we only need to compare the products of signs $\sign_{OE}(\sigma'_{xy})\sign_{OE}(\sigma'_{wv})$, 
%$\sign_{OE}(\sigma'_{xv})\sign_{OE}(\sigma'_{wy})$, and $\sign_{OE}(\sigma)\sign_{OE}(\sigma'_{xywv})$
%Thus 
there are two possibilities\footnotemark : 
\footnotetext{The assumption that the nodes $x, y, w, v$ are in cyclic order is required. Otherwise, it would be possible for $x$ to be red, $y$ to be green, $w$ to be blue, and $v$ to be green. In this case, the sets $\{x, y\}$ and $\{x, v\}$ consist of one red node and one green node, and the sets $\{w, v\}$ and $\{y, w\}$ consists of one green node and one blue node.}
\begin{itemize}
\item[(i)] One of the sets $\{x, y\}, \{w, v\}$ consists of a red node and a green node and the other consists of a green node and a blue node. Also, one of the sets $\{x, v\},\{y, w\}$ consists consists of one red node and one blue node, and the other consists of two green nodes. 
\item[(ii)] One of the sets $\{x, v\},\{y, w\}$ consists of a red node and a green node and the other consists of a green node and a blue node. Also, one of the sets $\{x, y\}, \{w, v\}$ consists consists of one red node and one blue node, and the other consists of two green nodes. 
\end{itemize}

% and the set $\{x,v\}$ would consist of one red node and one green node and the set $\{y, w\}$ would consist of a green node and a blue node). 

We only prove case (i), as case (ii) is essentially the same.  
%Since these two cases are so similar, we only consider case (i). 
By Remark~\ref{rem:caseanalysis},
\begin{eqnarray*}
\sign_{OE}(\sigma'_{xy}) \sign_{OE}(\sigma'_{wv}) 
&=& (-1)^{|RB(\sigma)|} (-1)^{|RG(\sigma)| -1} (-1)^{|RB(\sigma)|} (-1)^{|GB(\sigma)| -1} \\
&=&(-1)^{|RG(\sigma)|}(-1)^{|GB(\sigma)|}, \text{ and} \\
\sign_{OE}(\sigma'_{xv}) \sign_{OE}(\sigma'_{wy}) 
&=&  (-1)^{|RG(\sigma)|} (-1)^{|GB(\sigma)|}  (-1)^{|RB(\sigma)|-1} (-1)^{|RB(\sigma)|} \\
&=&- (-1)^{|RG(\sigma)|}(-1)^{|GB(\sigma)|}.
\end{eqnarray*}
Since we can obtain $\sigma_{xywv}$ by first removing the nodes $x, y$, and then removing the nodes $w, v$, 
$\sign_{OE}(\sigma'_{xywv}) = \sign_{OE}(\sigma'_{xy}) \sign_{OE}(\sigma'_{wv}) $.
%\begin{eqnarray*}
%\sign_{OE}(\sigma'_{xy}) \sign_{OE}(\sigma'_{wv}) 
%&=& (-1)^{|RB(\sigma)|} (-1)^{|RG(\sigma)| -1} (-1)^{|RB(\sigma)|} (-1)^{|GB(\sigma)| -1}
%\\
%&=&(-1)^{|RG(\sigma)|}(-1)^{|GB(\sigma)|} .
%\end{eqnarray*}
Thus by Theorem~\ref{thm:cond}, 
\begin{eqnarray*}
Z^{DD}_{\sigma}(G, {\bf N}) Z^{DD}_{\sigma_{xywv}}(G, {\bf N} - \{x, y, w, v\})  &=& 
Z^{DD}_{\sigma_{xy}}(G, {\bf N} - \{x, y\})  Z^{DD}_{\sigma_{wv}}(G, {\bf N} - \{w, v\}) \\
&&+
 Z^{DD}_{\sigma_{xv}}(G, {\bf N} - \{x, v\})  Z^{DD}_{\sigma_{wy}}(G, {\bf N} - \{w, y\}) .
 \end{eqnarray*}
\normalsize
\end{proof}

\section{Acknowledgements}
 
I thank my advisor, Benjamin Young, for posing this problem and his invaluable guidance throughout this project. Thanks also to Richard Kenyon, who provided thoughtful comments on the introduction and to Gregg Musiker for a helpful discussion. Finally, thank you to the anonymous referee who provided numerous helpful comments.

\bibliographystyle{alpha}
\bibliography{sample}

\end{document}